\documentclass[10pt,reqno]{article}
\usepackage{amsfonts,amsrefs,latexsym,amsmath, amssymb, mathrsfs, verbatim}
\usepackage{url,color}
\usepackage{pifont}
\usepackage{upgreek}
\usepackage{fancyhdr}
\usepackage[linktocpage]{hyperref}
\usepackage{calligra}
\usepackage{marvosym}
\usepackage[percent]{overpic}
\usepackage{pict2e}
\usepackage{graphicx}
\usepackage{tensor}
\usepackage{amsthm}
\usepackage{subfigure}
\usepackage[hypcap=false]{caption}
\usepackage{xcolor}
\usepackage{wasysym}
\usepackage{accents}
\usepackage{enumitem}
\usepackage{tikz}
\usetikzlibrary{shapes.geometric, arrows}
\usetikzlibrary{mindmap}
\setlist[enumerate]{label*=\arabic*.}
\numberwithin{equation}{section}

\topmargin - .23 in

\newtheorem{theorem}{Theorem}[section]
\newtheorem{proposition}{Proposition}[section]
\newtheorem{lemma}[proposition]{Lemma}
\newtheorem{corollary}[proposition]{Corollary}

\theoremstyle{definition}
\newtheorem{definition}{Definition}[section]
\newtheorem{remark}[proposition]{Remark}

\DeclareMathAlphabet{\mathcalligra}{T1}{calligra}{m}{n}
\DeclareFontShape{T1}{calligra}{m}{n}{<->s*[2.2]callig15}{}

\newcommand{\Domain}{\mathcal{R}}
\newcommand{\Int}{\text{int}}
\newcommand{\tderivative}{\frac{d}{dt}}

\newcommand{\vort}{\omega}

\newcommand{\speedone}{c_1}
\newcommand{\speedtwo}{c_2}

\newcommand{\variables}{\Phi}
\newcommand{\vvariables}{\vec{\variables}}
\newcommand{\cvariables}{\Psi}
\newcommand{\cvvariables}{\vec{\cvariables}}
\newcommand{\vardiv}{\phi}
\newcommand{\vvardiv}{\vec{\vardiv}}
\newcommand{\varcurl}{\psi}
\newcommand{\vvarcurl}{\vec{\varcurl}}

\newcommand{\norm}[1]{\left\lVert#1\right\rVert}
\newcommand{\onenorm}[4]{\norm{#1}_{L_{#2}^{#3}(#4)}}
\newcommand{\twonorms}[6]{\norm{#1}_{L_{#2}^{#3}L_{#4}^{#5}(#6)}}
\newcommand{\sobolevtwonorms}[6]{\norm{#1}_{L_{#2}^{#3}H_{#4}^{#5}(#6)}}
\newcommand{\threenorms}[8]{\norm{#1}_{L_{#2}^{#3}L_{#4}^{#5}L_{#6}^{#7}(#8)}}
\newcommand{\holdernorm}[5]{\norm{#1}_{C_{#2}^{#3,#4}(#5)}}
\newcommand{\semiholdernorm}[5]{\norm{#1}_{\dot{C}_{#2}^{#3,#4}(#5)}}
\newcommand{\holdertwonorms}[7]{\norm{#1}_{L_{#2}^{#3}C_{#4}^{#5,#6}(#7)}}
\newcommand{\holderthreenorms}[9]{\norm{#1}_{L_{#2}^{#3}L_{#4}^{#5}C_{#6}^{#7,#8}(#9)}}
\newcommand{\sobolevnorm}[3]{\norm{#1}_{H^{#2}(#3)}}
\newcommand{\abs}[1]{\left\vert#1\right\vert}
\newcommand{\sgabs}[1]{\left\vert#1\right\vert_{\gsphere}}

\newcommand{\dive}{\mbox{\upshape div}}
\newcommand{\curl}{\mbox{\upshape curl}}

\newcommand{\antisymmetic}{\epsilon}

\newcommand{\asphere}{\upomega^{A}}
\newcommand{\bsphere}{\upomega^{B}}
\newcommand{\csphere}{\upomega^{C}}
\newcommand{\cartiasphere}{\Theta_{(A)}}

\newcommand{\vcartiasphere}{\vec{\Theta}_{(A)}}

\newcommand{\sangle}{\upomega}

\newcommand{\derinormal}{\frac{\p}{\p w}}
\newcommand{\deriasphere}{\frac{\p}{\p\asphere}}
\newcommand{\deribsphere}{\frac{\p}{\p\bsphere}}
\newcommand{\dericsphere}{\frac{\p}{\p\csphere}}

\newcommand{\Tstar}{T_{*}}
\newcommand{\Trescale}{T_{*;(\lambda)}}

\newcommand{\region}{\mathcal{M}}
\newcommand{\intregion}{\region^{(\text{Int})}}
\newcommand{\extregion}{\region^{(\text{Ext})}}
\newcommand{\Stwo}{\mathbb{S}^2}
\newcommand{\stu}{S_{t,u}}
\newcommand{\sTauu}{S_{\tau,u}}
\newcommand{\Stau}{\Sigma_{\tau}}
\newcommand{\St}{\Sigma_t}
\newcommand{\Stint}{\St^{(\text{Int})}}
\newcommand{\Szero}{\Sigma_0}
\newcommand{\initialSzero}{\Sigma_0^{w(\lambda)}}
\newcommand{\coneu}{\mathcal{C}_{u}}

\newcommand{\tstart}{[u]_+}

\newcommand{\gensmoothfunction}{\mathrm{f}}
\newcommand{\quadsmoothfunction}{\mathscr{Q}}
\newcommand{\linsmoothfunction}{\mathscr{L}}

\newcommand{\lgensmoothfunction}{\mathrm{f}_{(\vLunit)}}

\newcommand{\Ricfour}[2]{\mathbf{Ric}_{#1 #2}}
\newcommand{\Riemfour}[4]{\mathbf{Riem}_{#1 #2 #3 #4}}
\newcommand{\Riemonethree}[4]{\mathbf{Riem}^{#1}_{#2 #3 #4}}

\newcommand{\spherenormal}{N}

\newcommand{\CC}{\mathcal{A}}
\newcommand{\ACC}{\tilde{\mathcal{A}}}
\newcommand{\nulllapse}{b}
\newcommand{\lapse}{a}
\newcommand{\volume}{\upsilon}
\newcommand{\conformalvol}{\tilde{\upsilon}}
\newcommand{\spheresecondfund}{\uptheta}
\newcommand{\modtorsion}{\tilde{\upzeta}}
\newcommand{\uzeta}{\underline{\upzeta}}
\newcommand{\hchi}{\hat{\upchi}}
\newcommand{\uchi}{\underline{\upchi}}
\newcommand{\huchi}{\hat{\uchi}}
\newcommand{\Restrace}{\mytr_{\congsphere}}
\newcommand{\reschi}{\tilde{\upchi}}
\newcommand{\chismall}{z}
\newcommand{\resuchi}{\tilde{\uchi}}

\newcommand{\varphiconformal}{\tilde{\varphi}}
\newcommand{\HJen}[1]{^{(#1)} \mkern-3mu \breve{\mathbf{J}}}
\newcommand{\Jenarg}[2]{{^{(#1)} \mkern-3mu \mathbf{J}^{#2}}}

\newcommand{\energycurrent}[1]{\mkern-1mu ^{(#1)} \mkern-3mu \mathbf{J}}
\newcommand{\modifiedcurrent}[1]{\mkern-1mu ^{(#1)} \mkern-3mu \tilde{\mathbf{J}}}
\newcommand{\energycurrentwlot}[1]{\mkern-1mu ^{(#1)} \mkern-3mu \tilde{\mathbf{J}}}
\newcommand{\deform}[1]{{^{(#1)} \mkern-1mu \pmb{\pi}}}
\newcommand{\Hdeform}[1]{{^{(#1)} \mkern-1mu \breve{\pmb{\pi}}}}

\newcommand{\minkowski}{\eta}

\newcommand{\gfour}{\mathbf{g}}
\newcommand{\invgfour}{\left(\mathbf{g}^{-1}\right)}
\newcommand{\hfour}{\mathbf{h}}
\newcommand{\rescaledgfour}{\widetilde{\mathbf{g}}}
\newcommand{\gt}{g}

\newcommand{\gsphere}{g \mkern-8.5mu / }
\newcommand{\esphere}{e \mkern-8.5mu / }

\newcommand{\congsphere}{\widetilde{g} \mkern-8.5mu / }
\newcommand{\tip}{\mathbf{o}}

\newcommand{\conformalfactor}{\upsigma}
\newcommand{\Chfour}{\pmb{\Gamma}}

\newcommand{\gvol}{\varpi_{\gfour}}
\newcommand{\hvol}{\varpi_{\hfour}}

\newcommand{\conspherevol}{\varpi_{\congsphere}}
\newcommand{\tvol}{\varpi_g}
\newcommand{\htvol}{\varpi_h}

\newcommand{\spherevol}{\varpi_{\gsphere}}
\newcommand{\flatspherevol}{\varpi_{\esphere}}

\newcommand{\cphi}{\ln\left(\rgeo^{-2}v\right)}

\newcommand{\sphereproject}{{\Pi \mkern-12mu / } \, }

\newcommand{\remainder}{\mathfrak{R}}

\newcommand{\mytr}{{\mbox{\upshape tr}}}
\newcommand{\gtr}{\mytr_{\gsphere}}
\newcommand{\ggtr}{\mytr_{g}}
\newcommand{\Timelike}{\mathbf{T}}
\newcommand{\materialderivative}{\mathbf{B}}
\newcommand{\Lunit}{L}
\newcommand{\Conenormal}{\mathbf{V}}
\newcommand{\uLunit}{\underline{L}}
\newcommand{\Lgeo}{L_{(Geo)}}

\newcommand{\p}{\partial}
\newcommand{\pfour}{\pmb\partial}
\newcommand{\vLunit}{\vec{\Lunit}}
\newcommand{\uvLunit}{\vec{\uLunit}}

\newcommand{\Dfour}{\mathbf{D}}
\newcommand{\angD}{ {\Dfour \mkern-14mu / \,}}
\newcommand{\angnabla}{ {\nabla \mkern-14mu / \,} }
\newcommand{\angDarg}[1]{{\angD_{\mkern-3mu #1}}}

\newcommand{\angdiv}{\mbox{\upshape{div} $\mkern-17mu /$\,}}
\newcommand{\angcurl}{\mbox{\upshape{curl} $\mkern-17mu /$\,}}
\newcommand{\anglap}{ {\Delta \mkern-12mu / \, } }

\newcommand{\rgeo}{\widetilde{r}}

\newcommand{\modmass}{\check{\upmu}}
\newcommand{\hodgemass}{{\upmu \mkern-10mu / \,}}
\newcommand{\mass}{\upmu}
\newcommand{\umodmass}{\bar{\modmass}}
\newcommand{\massone}{\hodgemass_{(1)}}
\newcommand{\masstwo}{\hodgemass_{(2)}}

\newcommand{\Lie}{\mathcal{L}}
\newcommand{\SigmatLie}{\underline{\mathcal{L}}}
\newcommand{\angLie}{ { \mathcal{L} \mkern-10mu / } }
\newcommand{\littlewood}{P_\upnu}
\newcommand{\anglittlewood}{{P    \mkern-12mu  / }}
\newcommand{\boxg}{\square_{\gfour}}
\newcommand{\boxgrescale}{\square_{\tilde{\gfour}}}
\newcommand{\resboxg}{\hat{\square}_{\gfour}}
\newcommand{\boxh}{\square_{\hfour}}
\newcommand{\diff}{\mathrm d}

\newcommand{\neighborhood}{\mathfrak{O}}

\newcommand{\zero}{\mathcal{O}}
\newcommand{\average}[1]{\overline{#1}}

\newcommand{\Energyconformal}{\mathfrak{C}}
\newcommand{\iEnergyconformal}{\Energyconformal^{(i)}}
\newcommand{\eEnergyconformal}{\Energyconformal^{(e)}}
\newcommand{\initialenergy}{\mathcal{C}_0}
\newcommand{\geonormals}{\mathbf{n}}
\title{{Low-Regularity Local Well-Posedness\\ for the Elastic Wave System}}

\date{\today}       

\begin{document}
	\author{Xinliang An\thanks{National University of Singapore, Singapore. matax@nus.edu.sg},
		Haoyang Chen\thanks{National University of Singapore, Singapore. hychen@nus.edu.sg} and 
		Sifan Yu\thanks{National University of Singapore, Singapore. sifanyu@nus.edu.sg}}
	\maketitle         
	\begin{abstract}
		
		We study the elastic wave system in three spatial dimensions. For admissible harmonic elastic materials, we prove a desired low-regularity local well-posedness result for the corresponding elastic wave equations. For such materials, we can split the dynamics into the ``divergence-part” and the ``curl-part,” and each part satisfies a distinct coupled quasilinear wave system with respect to different acoustical metrics.
		Our main result is that the Sobolev norm $H^{3+}$ of the ``divergence-part” (the ``faster-wave part”) and the $H^{4+}$ of the ``curl-part” (the ``slower-wave part”) can be controlled in terms of initial data for short times. We note that the Sobolev norm assumption $H^{3+}$ is optimal for the ``divergence-part.” This marks the first favorable low-regularity local well-posedness result for a wave system with multiple wave speeds. Compared to the quasilinear wave equation, new difficulties arise from the multiple wave-speed nature of the system. Specifically, the acoustic metric $\gfour$ of the faster-wave depends on both the faster-wave and slower-wave parts. Additionally, the dynamics of the faster-wave {``divergence-part"} require higher regularity of the ``curl-part”. In particular, the Ricci curvature associated with the faster-wave is one derivative rougher than that of the slower-wave dynamics.	
		This phenomenon {also appears in} the compressible Euler equations (featuring multiple characteristic speeds) and {is} a major {obstacle to} obtaining low-regularity local well-posedness results {for general quasilinear wave systems} if the two parts do not exhibit strong decoupling properties or if the ``curl-part” lacks the structure necessary for better regularity results. For the elastic wave system governing the dynamics of the admissible harmonic elastic materials, we report that we can overcome these difficulties. {For this system, we exploit its geometric structures and find that the ``divergence-part” and ``curl-part” exhibit decoupling properties and both parts show regularity gains.}
		Moreover, we {prove} that the ``divergence-part” maintains to represent the faster-wave {throughout} the entire time of the existence of the solution, ensuring that the characteristic hypersurfaces of the faster-wave are spacelike with respect to the slower-wave. {This implies a crucial} coerciveness for the geometric {cone-flux} energy of the ``curl-part” on such characteristic hypersurfaces of the ``divergence-part." We {furthermore} carefully address all {these} challenges through {spacetime energy estimates}, Strichartz estimates, frequency-localized decay estimates, and conformal energy estimates. In all these estimates, we also precisely trace the impact of the ``curl-part” on the faster-wave dynamics and control the associated geometry {via employing} the vector field method and the Littlewood-Paley theory.\\
		
		\\
	\end{abstract}
	\tableofcontents
	\newpage
	\section{Introduction}\label{S:sectionintro}
	{In this paper, we initiate the study} of low-regularity local well-posedness for elastic waves in three spatial dimensions. For homogeneous isotropic hyperelastic materials, the motion of displacement $\vec{U}=(U^1,U^2,U^3)$ is governed by a \emph{quasilinear} wave system with \emph{multiple wave-speeds}:
	\begin{equation} \label{EQ:wave}
		\partial_t^2 \vec{U}-\speedtwo^2\Delta \vec{U}-(\speedone^2-\speedtwo^2)\nabla(\nabla\cdot \vec{U})=\vec{\mathcal{G}}(\nabla \vec{U},\nabla^2\vec{U}).
	\end{equation}
	Here, $c_1,c_2$ are two constants representing respectively the speeds of longitudinal and transverse waves with $c_1>c_2>0$. The nonlinear form $\vec{\mathcal{G}}(\nabla \vec{U},\nabla^2\vec{U})$ is linear in $\nabla^2\vec{U}$ and is fixed by different materials. The precise expression of $\vec{\mathcal{G}}(\nabla \vec{U},\nabla^2\vec{U})$ will be provided in Sections \ref{SS:elasticwaveeq}-\ref{SS:ModifiedElasticWaves}.
	
	{To derive equations \eqref{EQ:wave}, we start from considering} the motion of a 3D elastic body, {which can be} described by time-dependent orientation-preserving diffeomorphisms, denoted by $\vec{P}:\mathbb{R}^3\times\mathbb{R}\to\mathbb{R}^3$. Here, $\vec{P}=\vec{P}(\vec{x},t)$ verifies $\vec{P}(\vec{x},0)=\vec{x}=(x^1,x^2,x^3)$. The deformation gradient is further defined as $F:=\nabla P$  with  $F^i_{j}:=\partial P^i/\partial x^j.$ For a homogeneous isotropic hyperelastic material, {a physical quantity called} the stored energy $W$ depends solely on the principal invariants of $FF^T$. {Here, $FF^T$ is called the Cauchy-Green strain tensor. With $\vec{P}$ and $W$, we can define} the action functional $\mathcal{S}$:
	\begin{align}\label{EQ:Lagr}
		\mathcal{S}:=\int\int_{\mathbb{R}^3}\left\{\frac12|\partial_t\vec{P}|^2-W(FF^T)\right\}\diff x\diff t.
	\end{align}
	{Applying the least action principle, }the resulting Euler-Lagrange equations are given by:
	\begin{equation}\label{EQ:el0}
		\frac{\partial^2P^i}{\partial t^2 }-\frac{\partial}{\partial{x^l}}\frac{\partial W(FF^T)}{\partial F^i_{l}}=0.
	\end{equation}
	Let $\vec{U}:=\vec{P}-\vec{x}$ denote the displacement. Via (\ref{EQ:el0}), we then derive the elastic wave system (\ref{EQ:wave}).

	\subsection{Main Result}
	{In this paper, we focus on the elastic wave equations for the admissible harmonic elastic materials}, which belong to the class of harmonic elastic materials {\emph{described by John} in \cite{John1960,John1966} and are associated with} remarkable decoupling properties. The equations {that govern the motion of the admissible harmonic materials} take the form:
	\begin{align}\label{EQ:elasticwaves0.0}
		\p_t^2\vec{U}-\speedtwo^2\Delta \vec{U}-(\speedone^2-\speedtwo^2)\nabla(\nabla\cdot \vec{U})=\dive \left\{	G(\p\vec{U})I+b\det(F)(F^T)^{-1}\right\}.
	\end{align}
	Here, $I$ is the identity matrix, $b$ is a constant related to the material, $F^T$ is the transpose of $F$, $\det(F)$ is the determinant of $F$, and $G(\p\vec{U})$ is a smooth scalar function of $\p\vec{U}$. See Section \ref{SS:ModifiedElasticWaves} for detailed discussions. Both \eqref{EQ:wave} and \eqref{EQ:elasticwaves0.0} are multi-wave-speeds systems, see Section \ref{SS:admissible} for discussions on the special structure of the admissible harmonic elastic wave equations. For equation \eqref{EQ:elasticwaves0.0}, we now state the main theorem of this article.

	\begin{theorem}[Main theorem]\label{TH:IntroMainTheorem}
		We study the admissible harmonic elastic wave system \eqref{EQ:elasticwaves0.0} with initial data $\vec{U}|_{\Sigma_0}:=\vvardiv|_{\Sigma_0}+\vvarcurl|_{\Sigma_0}$. Here, $\vvardiv$ corresponds to the ``divergence-part" and $\vvarcurl$ corresponds to the ``curl-part" of the displacement\footnote{Precise forms of $\vvardiv$ and $\vvarcurl$ are defined in Definition \ref{DE:DefPhiPsi}}. Let $D$ be a fixed positive constant and $i=1,2,3$. Let $\breve{\Domain}$, the region of hyperbolicity, be defined as in Definition \ref{DE:Hyperbolicity}. For any real number $3<N<7/2$, assume that $\vvardiv$ and $\vvarcurl$ satisfy the following initial conditions:
		\begin{enumerate}
			\item $\sobolevnorm{\vardiv^i}{N}{\Sigma_0}+\sobolevnorm{\varcurl^i}{N+1}{\Sigma_0}\leq D$, with $D$ being a fixed positive constant and $i=1,2,3$.
			\item The initial data functions lie within the interior of $\breve{\Domain}$. 
		\end{enumerate}
		Then, for \eqref{EQ:elasticwaves0.0}, the solution's time of classical existence $T>0$ can be bounded from below in terms of $D$. Furthermore, the Sobolev regularity\footnote{See Section \ref{SS:norms} for the definitions of the mixed-norms.} of the data is preserved during the propagation of the solution throughout the time slab of classical existence. Specifically, this means that\footnote{In this article, $A\lesssim B$ means $A\leq C\cdot B$ for some universal constant $C$.}
		\begin{subequations}\label{ES:introenergy}
			\begin{align}
				\label{ES:introenergy1}\sobolevtwonorms{\vardiv^i}{t}{\infty}{x}{N}{[0,T]\times\mathbb{R}^3}\lesssim&1,\\
				\label{ES:introenergy2}\sobolevtwonorms{\varcurl^i}{t}{\infty}{x}{N+1}{[0,T]\times\mathbb{R}^3}\lesssim&1.
			\end{align}
		\end{subequations} 
		
		In addition, $\vardiv$ and $\varcurl$ satisfy the following Strichartz estimates\footnote{We denote the schematic spatial partial derivatives by $\p$ and the schematic space-time partial derivatives by $\pfour$.}:
		\begin{subequations}\label{ES:introStrichartz}
			\begin{align}
				\twonorms{\pfour^2\vardiv^i}{t}{2}{x}{\infty}{[0,T]\times\mathbb{R}^3}\lesssim&1,\\
				\label{ES:introstrichartz2}\twonorms{\pfour^3\varcurl^i}{t}{2}{x}{\infty}{[0,T]\times\mathbb{R}^3}\lesssim&1.
			\end{align}
		\end{subequations}
	\end{theorem}
\begin{remark}\label{RM:insshockformation}
		{For general 3D elastic wave system, in \cite{an2020low}, An-Chen-Yin proved the $H^3$ ill-posedness, which is driven by instantaneous shock formation. The ``divergence-part" in \cite{an2020low} is in $H^{3}$, and the ``curl-part" in \cite{an2020low} is smooth.} Therefore, the\footnote{Here, $H^{3+}$ means $H^{3+\varepsilon}$ for arbitrarily small $\varepsilon$.} $H^{3+}$ local well-posedness (for the divergence part) obtained in Theorem \ref{TH:IntroMainTheorem} is the desired result for the admissible harmonic materials.
	\end{remark}	
	 Theorem \ref{TH:IntroMainTheorem} establishes optimal low-regularity local well-posedness for the admissible harmonic elastic wave system. This will be fundamental for studying this model both mathematically and numerically. Under minimal low-regularity initial conditions, the above theorem shows that singularity formation, including shock formation from sound wave compression, can be avoided at least for short times. Proving the above theorem requires a thorough understanding of the solution dynamics and the associated acoustic geometries of both the “divergence part” and the “curl part,” which feature different characteristic speeds. Meanwhile, although this is a local result, to work with optimal regularity, we need to rescale the spacetime with respect to a large frequency $\lambda$ and prove a decay estimate over the long rescaled time. This decay estimate is essential for proving a frequency-localized Strichartz estimate, which yields \eqref{ES:introStrichartz} via Duhamel's principle.\\
	\color{black}
    
	To establish Theorem \ref{TH:IntroMainTheorem}, we first explore the notable structures of the admissible harmonic elastic materials, which yield decoupling properties (not necessarily present in elastic materials in general\footnote{See Section \ref{SS:admissible} for further discussion of harmonic elastic materials.}) and regularity gains for both the faster-wave part and the slower-wave part.
	
	Compared to previous low-regularity results\footnote{See Section \ref{SS:Overviewofpreviousresults} for a discussion of these results.}, we face new challenges arising from the multiple wave-speed nature of the system. Specifically, the acoustic metric $\gfour$ of the faster-wave depends on both the faster-wave part and the slower-wave parts, requiring us to carefully analyze the corresponding geometric quantities associated with $\gfour$, such as the connection coefficients. Additionally, controlling the dynamics of the faster-wave part requires a higher regularity for the ``curl-part". In particular, the Ricci curvature associated with the faster-wave is one derivative rougher than that of the slower-wave dynamics. This phenomenon is common for physical quasilinear wave systems (e.g., the compressible Euler equations) featuring multiple characteristic speeds and would be a major obstacle to proving low-regularity local well-posedness results if the two parts with different traveling speeds do not exhibit strong (enough) decoupling properties or if the ``curl-part" lacks the necessary structure for better regularity gains. 
	
	For the admissible harmonic elastic materials, we overcome those difficulties by fully exploiting the associated geometric structure. Specifically, we show that the ``divergence-part" maintains to represent the faster-wave in the entire time of the existence of the solution. This property provides us that characteristic hypersurfaces of the faster-wave are spacelike with respect to the slower-wave. {For the ``curl-part", this results in a crucial coercive flux energy along the geometric cone of the faster-wave divergence part.} Moreover, we carefully go through\footnote{See Section \ref{SSS:Blueprint} for the blueprint for proving Theorem \ref{TH:IntroMainTheorem}.} the spacetime energy estimates, a rescaled version of Strichartz estimate, a frequency-localized decay estimate, and conformal energy estimates. Meanwhile, we precisely trace and address the impact of the ``curl-part” on the faster-wave {``divergence-part"} dynamics by a careful control of the geometry associated with the faster-wave. This includes computing complicated geometric structure equations satisfied by various appropriate geometric quantities and deriving mixed-norm estimates for those quantities.

	\subsection{History of Previous Low-Regularity Results in Quasilinear Wave Equations, Compressible Euler Equations, and Related Systems}\label{SS:Overviewofpreviousresults} 
	In the past two decades, significant advancements have been made for studying the Cauchy problem with low-regularity initial data for the following quasilinear wave equations\footnote{$\quadsmoothfunction[A](B,C)$ denotes any scalar-valued function that is quadratic in the components of $B$ and $C$ with coefficients that are functions of the components of $A$.} in three space-dimensions and for the related compressible fluid systems:
	\begin{align}\label{quasi}
		\square_{\gfour(\psi)}(\psi)=\quadsmoothfunction(\psi)[\pfour\psi,\pfour\psi].
	\end{align} 
	Note that in \eqref{quasi}, the Lorentzian metric $\gfour(\psi)$ depends on the unknown $\psi$, and $\quadsmoothfunction(\psi)[\pfour\psi,\pfour\psi]$ is quadratic of $\pfour\psi$.
	
	Classical local well-posedness with initial data in $H^{\frac{5}{2}+}(\mathbb{R}^3)$ for the quasilinear wave system (\ref{quasi}) can be obtained by employing standard energy estimates and Sobolev embedding, as shown by Kato \cite{kato}. Progress of local well-posedness below $H^{(5/2)}$ began in the 1990s. For quasilinear wave equations of the form (\ref{quasi}), Bahouri-Chemin \cite{Equationsd'ondesquasilineairesetestimationdeStrichartz} and Tataru \cite{Strichartzestimatesforoperatorswithnonsmooth} independently showed that local well-posedness holds with $H^{(2+\frac{1}{4})+}$ Cauchy data. These improvements rely on Strichartz estimates via applying Fourier integral parametrix representations. Bahouri-Chemin subsequently improved their earlier result to $H^{(2+\frac{1}{5})+}$ in \cite{Equationsd'ondesquasilinearireseteffetdispersif}. Tataru and Klainerman independently pushed the regularity down to $H^{(2+\frac{1}{6})+}$ in \cite{Strichartzestimatesforsecondorderhyperbolicoperators,ACommutingvector fieldsApproachtoStrichartz}. Klainerman-Rodnianski further achieved an $H^{(2+\frac{2-\sqrt{3}}{2})+}$ local well-posedness in \cite{ImprovedLocalwellPosedness}. They also proved the $H^{2+}$ local well-posedness for Einstein equations in wave coordinates in \cite{KR3}. Finally, the bounded $L^2$ curvature conjecture ($H^2$ local well-posedness in Yang-Mills frame) for Einstein vacuum equations was proved by Klainerman-Rodnianski-Szeftel in \cite{L2curvatureconjecture} with accompanying contributions by Szeftel in \cite{szeftel2012parametrix1,szeftel2012parametrix2,szeftel2012parametrix3,szeftel2012parametrix4}. For generic quasilinear wave equations, the optimal $H^{2+}$ local well-posedness result was first achieved by Smith-Tataru in \cite{Sharplocalwell-posednessresultsforthenonlinearwaveequation} by using methods of wave-packets. Recently, Wang \cite{AGeoMetricApproach} gave an alternative proof of the Smith-Tataru \cite{Sharplocalwell-posednessresultsforthenonlinearwaveequation} result via applying an enhanced version of Klainerman-Rodnianski's vector field method in \cite{ImprovedLocalwellPosedness} and via employing a novel conformal energy method. 
	
	For the 3D compressible Euler flow with vorticity and entropy\footnote{We note that irrotational and isentropic compressible Euler equation can be reduced to the quasilinear wave equation of form \eqref{quasi}.}, a geometric formulation is derived in Luk-Speck \cite{ThehiddennullstructureofthecompressibleEulerequationsandapreludetoapplications} and Speck \cite{ANewFormulationofthe3DCompressibleEulerEquationswithDynamicEntropy}. The 3D Compressible Euler equations can be rewritten as coupled quasilinear wave equations for the ``wave-part" and transport equations for the ``transport-part". Then, under $H^{2+}$ initial data assumption on the ``wave-part", and a more regular initial data assumption on the ``transport-part" of the data, local well-posedness result was proved by Disconzi-Luo-Mazzone-Speck \cite{3DCompressibleEuler}, and Wang \cite{Roughsolutionsofthe3-DcompressibleEulerequations} with a sharper regularity assumption for the ``transport-part". We also refer to Zhang-Anderson \cite{Alterproof} for a different proof using wave packets. In a similar vein, for the relativistic fluid, Disconzi and Speck \cite{TheRelativistivEulerEquations} derived a geometric formulation of the relativistic Euler equations with dynamic vorticity and entropy; thereafter, Yu \cite{yu2022rough} and Zhang \cite{zhang2024wellposedness} proved the desired low-regularity local well-posedness results. We remark that, for the general quasilinear wave equations of the form (\ref{quasi}), $H^{2+}$ local well-posedness is the optimal result. Specifically, Lindblad provided an example of $H^2$ ill-posedness for a specific quasilinear wave equation in \cite{counterexamplestolocalexistenceforquasilinearwaveequations}. An-Chen-Yin \cite{an2021low,xinliangAn} later extended the corresponding $H^2$ ill-posedness to the compressible Euler equations and to the ideal magnet-hydrodynamics (MHD) equations. 
	\newline

	For quasilinear wave equations of the following form:
	\begin{align}\label{quasi2}
		\square_{\gfour(\pfour\psi)}(\psi)=\quadsmoothfunction(\psi)[\pfour\psi,\pfour\psi],
	\end{align} 
	the critical Sobolev exponent in 3D is $H^{5/2}$. By differentiating (\ref{quasi2}) and changing the variable $\pfour\psi\rightarrow\phi$, (\ref{quasi2}) can be reduced to (\ref{quasi}). Hence, for general nonlinearities of $\quadsmoothfunction(\psi)[\pfour\psi,\pfour\psi]$, the standard local well-posedness result is $H^{\frac{7}{2}+}$ in 3D. Moreover, as shown by \cite{AGeoMetricApproach,Sharplocalwell-posednessresultsforthenonlinearwaveequation}, with general nonlinear forms, the optimal low-regularity local well-posedness result is for Cauchy data belonging to $H^{3+}$. On the other hand, with more favorable nonlinear structures of $\quadsmoothfunction(\psi)[\pfour\psi,\pfour\psi]$, it was conjectured that the associated local well-posedness can be obtained with initial regularity assumption below $H^3$; see \cite{Tataru2002ICM}. For the timelike minimal surface equation in Minkowski space, Ai-Ifrim-Tataru \cite{ai2021timelike} proved such a result and they showed that the Cauchy problem is locally well-posed in $H^{\left(\frac{5}{2}+\frac{1}{4}\right)+}$. 
	
	The elastic wave equations (\ref{EQ:wave}) are closely related to the quasilinear wave equations (\ref{quasi2}). In 3D, the classical local well-posedness result in $H^{\frac{7}{2}+}$ was obtained by Hughes-Kato-Marsden \cite{HughesKatoMarsden1976}. Without additional assumptions on the  nonlinearity $\mathcal{G}(\nabla \vec{U},\nabla^2\vec{U})$ in (\ref{EQ:wave}), An-Chen-Yin \cite{AnJMP,an2020low,an2022elasticwaves} proved ill-posedness results for elastic wave system. They showed that elastic wave equations are ill-posed for $H^3$ Cauchy data in $3D$ and are ill-posed for $H^{\frac{11}{4}}$ Cauchy data in $2D$, due to shock formation. They also analyzed the dynamics of the system by decomposing the elastic waves into multiple characteristic strips and precisely proved that the shock formation occurs for the fastest characteristic wave. Meanwhile, limited research has been done on the low-regularity local well-posedness for the elastic wave system. To our knowledge, this article marks the \textit{first} work of obtaining the low-regularity local well-posedness for wave systems featuring multiple wave speeds.

	\subsection{History of Previous Small Data Results for the Elastic Wave Equations}\label{SS:SmallData}

	The study of elastic wave system \eqref{EQ:wave} was pioneered by John. In 3D,  with (smooth) small initial data, he proved that singularities can arise in the radially symmetric case \cite{F.John} (see also \cite{John82}). Without assumptions on symmetry, for this system, the existence of almost-global solutions was later established by John in \cite{john88} and by Klainerman-Sideris in \cite{klainerman}. In the case that the nonlinearities $\mathcal{G}(\nabla \vec{U},\nabla^2\vec{U})$ satisfy the null conditions (as introduced by Klainerman in \cite{Klainerman86}), with small initial data, global existence of solutions to the Cauchy problem of \eqref{EQ:wave} are achieved by Agemi \cite{agemi} and Sideris \cite{Sideris,Sideris00}.
	
	It is worth noting that, to prove global existence results under small data assumptions, one only needs to consider the nonlinearity up to the quadratic terms\footnote{In such small data regime, the nonlinearity can be determined explicitly. See \cite{Sideris}}. This is because the smallness of $\nabla\vec{U}$ can be propagated via using energy estimates and weighted $L^\infty-L^2$ estimates. However, for the study of local well-posedness with general (possibly large) initial data, such a truncated perspective based on $|\nabla \vec{U}|$ being small is no longer valid. Furthermore, for general initial data, the characteristic hypersurfaces of the faster-speed wave and the slower-speed wave may intersect or even intertwine. Additionally, the system (\ref{EQ:wave}) may fail to remain hyperbolic for some large $\nabla\vec{U}$. Consequently, the general strategies and techniques used in the small data regime are not applicable to the low-regularity well-posedness problem studied in this paper. 
	
	\subsection{Multi-wave-speed Properties of the Elastic Wave Equations and the Admissible Harmonic Elastic Wave System}\label{SS:admissible}
	The elastic wave equations (\ref{EQ:wave}) form a multi-wave-speed quasilinear wave system. After taking divergence and curl of the equation (\ref{EQ:wave}) respectively, we arrive at
	\begin{subequations}\label{EQ:decouplewave1}
		\begin{align}
			\p_t^2\dive \vec{U}-\speedone^2\Delta\dive \vec{U}=&\dive \left(\mathcal{G}(\nabla \vec{U},\nabla^2\vec{U})\right),\\
			\p_t^2\curl \vec{U}-\speedtwo^2\Delta\curl \vec{U}=&\curl\left(\mathcal{G}(\nabla \vec{U},\nabla^2\vec{U})\right).
		\end{align}
	\end{subequations}
	Notice that we have $0<\speedtwo^2<\speedone^2$ and system (\ref{EQ:decouplewave1}) is quasilinear since the terms on the right-hand side (RHS) of (\ref{EQ:decouplewave1}) are of the same order as those on the left-hand side (LHS). In the past studies\footnote{See Section \ref{SS:Overviewofpreviousresults} for discussions.} of the low-regularity local well-posedness results for non-relativistic/relativistic compressible Euler equations, which are multi-speed systems, the following properties of the Euler system are crucially used:
	\begin{itemize}
		\item Regularity gains for both/all parts of the system.
		\item Decoupling of ``transport-part" and ``wave-part". (e.g., in the compressible Euler equations, the ``transport-part" (``slower-speed part") remains trivial if it initially vanishes.)
	\end{itemize}
	For elastic wave system, the above properties were not fully explored. To understand the low-regularity problem for the elastic wave equations, it forces us to dive deep into the different parts of the system's dynamics. We split the dynamics of the elastic wave into the ``divergence-part" and the ``curl-part", and investigate the coupling properties between these two parts. 
	
For general elastic material, in the expression of $\mathcal{G}(\nabla \vec{U},\nabla^2\vec{U})$, the ``divergence-part" and ``curl-part" of the dynamics of the elastic waves are intricately coupled. In particular, the dynamic vorticity $\curl \vec{U}$ of the displacement $\vec{U}$ can emerge even if it vanishes initially, meaning that pseudo-irrotationality (defined in Definition \ref{DE:pseudoirrotational}) is not preserved in general. This raises a natural question: for which materials do pseudo-irrotationality exhibit compatibility? More precisely, which materials allow motions with pseudo-irrotational initial conditions to remain pseudo-irrotational in the absence of body forces? John \cite{John1966} provided a precise classification of such materials, namely, the \textit{harmonic elastic materials}. That is the corresponding strain energy function $W$ (used in (\ref{EQ:Lagr}) for constructing the Lagrangian of the system) satisfies the following form:
	\begin{align}
		W=f(r)+as+bt,
	\end{align}
	where $a, b$ are constants, $f$ is an arbitrary function, and
 \begin{align}
     r:=&e_1+e_2+e_3,&
     s:=&e_1e_2+e_2e_3+e_3e_1,&
     t:=&e_1e_2e_3,
 \end{align} where $e_1,e_2,e_3$ are the eigenvalues of the positive definite symmetric matrix\footnote{This unique matrix $E$ is called the principal square root of $FF^T$, and $e_1,e_2,e_3$ are known as the singular values of the matrix $F$.} $E$ such that $E^2=FF^T$ with the Jacobian $F$ defined as $F_{ij}:=\p_iU^j$.
	
	In the regime of harmonic elastic wave equations, assuming $\curl \vec{U}=\vec{0}$ at $t=0$, the vorticity $\curl\vec{U}$ will stay $\vec{0}$ in later time and system (\ref{EQ:wave}) reduces to a single-wave-speed quasilinear wave equation (\ref{quasi2}) (see 
	\cite[Section 6]{John1966}). However, even for the harmonic elastic materials, provided non-zero vorticity for initial data, the ``divergence-part" and the ``curl-part"  are still deeply coupled at the top-order level. This indicates that the evolution of the harmonic elastic material with non-trivial vorticity is more \textit{complicated}\footnote{For the compressible Euler equations, the vorticity $\vort:=\curl v$ satisfies a transport-div-curl system such that $\p\vort\lesssim\p v$.} than the evolution of the compressible Euler equations. In this paper, we find a subset of the harmonic elastic materials, which we call \textit{admissible harmonic elastic materials}. For these materials, the displacement $\vec{U}$ satisfies:
		\begin{align}\label{EQ:elasticwaves0.1}
		\p_t^2U^i-\speedtwo^2\Delta U^i-(\speedone^2-\speedtwo^2)\p_i(\nabla\cdot \vec{U})=\sum\limits_{j,k=1}^3G^\prime(\p_jU^k)\p_i\p_jU^k.
	\end{align}
	This system allows dynamical vorticity and its local evolution is controlled in this paper. To get the form of equation (\ref{EQ:elasticwaves0.1}), we first rewrite elastic wave equation \eqref{EQ:el0} as follows:
	\begin{align}\label{EQ:introelasticwaves1.1}
		\p_t^2\vec{U}-\speedtwo^2\Delta \vec{U}-(\speedone^2-\speedtwo^2)\nabla(\nabla\cdot \vec{U})=\dive \left\{E(\p \vec{U})\right\},
	\end{align}
	where $E(\p\vec{U})$ is a $2$-vector field depending on $\p\vec{U}$. We note that $\dive E$ is the nonlinear-part of (\ref{EQ:el0}). 
	Secondly, we impose the assumption for the admissible harmonic elastic material. That is, we require its deformation $\vec{U}$ to satisfy equation \eqref{EQ:introelasticwaves1.1} with
	\begin{align}\label{DE:introH}
		E(\p \vec{U})=&G(\p\vec{U})I+b\det(F)(F^T)^{-1},
	\end{align}
	where $G(\p\vec{U})$ is a smooth scalar function of $\p\vec{U}$. Then, using the fact that $\dive\left\{\det(F)(F^T)^{-1}\right\}=0$, equation \eqref{EQ:elasticwaves0.1} follows from \eqref{EQ:introelasticwaves1.1} and \eqref{DE:introH}. We refer readers to Section \ref{SS:ModifiedElasticWaves} for a more detailed derivation.
\color{black}

	{After taking the divergence and the curl separately, the system (\ref{EQ:elasticwaves0.1}) exhibits remarkable decoupling properties below}:
	\begin{subequations}\label{EQ:introEW2}
		\begin{align}
			\label{EQ:introEW2.1}\hat{\square}_{\gfour_E(\p\vvardiv,\p\vvarcurl)}(\dive\vec{\vardiv})=&P_{(\dive\vec{\vardiv})}(\p^2\vvardiv,\p^3\vvarcurl),\\
			\label{EQ:introEW2.2}\square_{\hfour_E}(\curl\vec{\varcurl})^i=&0.
		\end{align}
	\end{subequations}
	Here, $\vardiv$ and $\varcurl$ denote the ``divergence-part" and the ``curl-part" of the solution $\vec{U}=\vec{\vardiv}+\vec{\varcurl}$, respectively; $\gfour_E$ and $\hfour_E$ are the Lorentzian metrics, with their inverses $\gfour_E^{-1}, \hfour_{E}^{-1}$ defined as follows:
	\begin{subequations}\label{EQ:introEW3}
		\begin{align}\label{EQ:introEW3.1}
			(\gfour^{-1})(\p\vvardiv,\p\vvarcurl):=&-\p_t\otimes\p_t+\sum\limits_{j,k=1}^3\left\{\delta^{jk}\speedone^2+\frac{1}{2}G^\prime\left(\p_j\vardiv^k+\p_j\varcurl^k\right)+\frac{1}{2}G^\prime\left(\p_k\vardiv^j+\p_k\varcurl^j\right)\right\}\p_j\otimes\p_k,\\
			\label{EQ:introEW3.2}(\hfour^{-1}_E):=&-\p_t\otimes\p_t+\speedtwo^2\sum\limits_{a=1}^3\p_a\otimes\p_a.
		\end{align}
	\end{subequations} 
Readers are referred to Section \ref{S:GeometricFormulation} for the detailed definitions of $\vardiv$ and $\varcurl$, as well as the precise form of (\ref{EQ:introEW2}).

\begin{remark}
	In this article, compared with equation \eqref{EQ:introEW2.2}, we study a more general case, where $\vvarcurl$ is allowed to satisfy the following quasilinear wave equation:
	\begin{align}\label{EQ:introquasiwave}
		\square_{\hfour(\p\vvarcurl)}(\curl\vec{\varcurl})^i=P^i_{(\curl\vvarcurl)}(\p^2\vvarcurl,\p^2\vvarcurl).
	\end{align}
Here, $\hfour^{-1}$ is of the form:
\begin{align}
	\label{EQ:introEW3.3}(\hfour^{-1}):=&-\p_t\otimes\p_t+\sum\limits_{j,k=1}^3\left\{\delta^{jk}\speedtwo^2+H^{jk}(\p\vvarcurl)\right\}\p_j\otimes\p_k.
\end{align}
\end{remark}

 According to (\ref{EQ:introEW2}) and \eqref{EQ:introquasiwave}, we can see that the evolution of $\varcurl$ is decoupled from the dynamics of the ``divergence-part". In addition, given that $\speedone>\speedtwo$, \textbf{under the hyperbolicity assumption \eqref{EQ:symmetry}, we can show that $\vardiv$ (strictly) represents the faster-speed wave, while $\varcurl$ (strictly) represents the slower-speed wave in the system}. Assuming one extra regularity on $\varcurl$, we are in the shape to utilize the machinery developed in the works mentioned in Section \ref{SS:Overviewofpreviousresults} for  quasilinear wave equations.
	
	\subsection{Proof Strategy of Theorem \ref{TH:IntroMainTheorem} and New Difficulties}\label{SS:SummaryofMainIdeas}
	\subsubsection{Blueprint for obtaining Theorem \ref{TH:IntroMainTheorem}}\label{SSS:Blueprint}
	In this article, {we carry out the blueprint in Figure \ref{steps} below, which is developed to achieve low-regularity local well-posedness for quasilinear wave equations as mentioned in Section \ref{SS:Overviewofpreviousresults}, and is appropriately modified for our system with multi-wave speeds.} Theorem \ref{TH:IntroMainTheorem} provides a priori estimates for smooth solutions, which are essential for establishing the local well-posedness. The remaining aspects of a full proof of local well-posedness could be shown by deriving uniform estimates for sequences of smooth solutions and their differences. We refer readers to \cite[Sections 2-3]{Sharplocalwell-posednessresultsforthenonlinearwaveequation} for the proof of local well-posedness based on a priori estimates in the low-regularity setting.\\
	
	Considering the admissible harmonic elastic wave system \eqref{EQ:elasticwaves0.0}, by employing the standard energy estimate (see, for example, \cite[Section 7]{john88}) and Littlewood-Paley theory, {for $N>3$}, one can derive the following estimate:
	\begin{align}\label{ES:introsdenergy0}
		\sobolevnorm{\vec{U}}{N}{\St}\lesssim\sobolevnorm{\vec{U}}{N}{\Szero}+\int_{0}^{t}\left(\onenorm{\pfour\p\vec{U}}{x}{\infty}{\Stau}+1\right)\sobolevnorm{\vec{U}}{N}{\Stau}\diff\tau.
	\end{align}
	To prove low-regularity result, building upon \eqref{ES:introsdenergy0}, we hope to further derive the following Strichartz estimate:
	\begin{align}\label{ES:introStrichartz0}
		\twonorms{\pfour\p\vec{U}}{t}{2}{x}{\infty}{[0,\Tstar]\times\mathbb{R}^3}^2\leq&1. 
	\end{align}
	
	The major efforts in this article are dedicated to prove the Strichartz-type estimates \eqref{ES:introStrichartz0}.
	In the picture below, we present the logic chain of the proofs in this paper, and it is based on a bootstrap argument. The colored steps involve our main modifications different from the single-wave case and we need to perform analysis tailored for the multi-speed elastic wave system. The techniques in the uncolored steps are mostly established in previous works such as \cite{3DCompressibleEuler,AGeoMetricApproach,ImprovedLocalwellPosedness,Equationsd'ondesquasilineairesetestimationdeStrichartz,Equationsd'ondesquasilinearireseteffetdispersif,Sharplocalwell-posednessresultsforthenonlinearwaveequation,Strichartzestimatesforsecondorderhyperbolicoperators,ACommutingvector fieldsApproachtoStrichartz,yu2022rough} on the low-regularity problems, and we carry out all the details with addressing extra terms arising from the elastic wave system.
	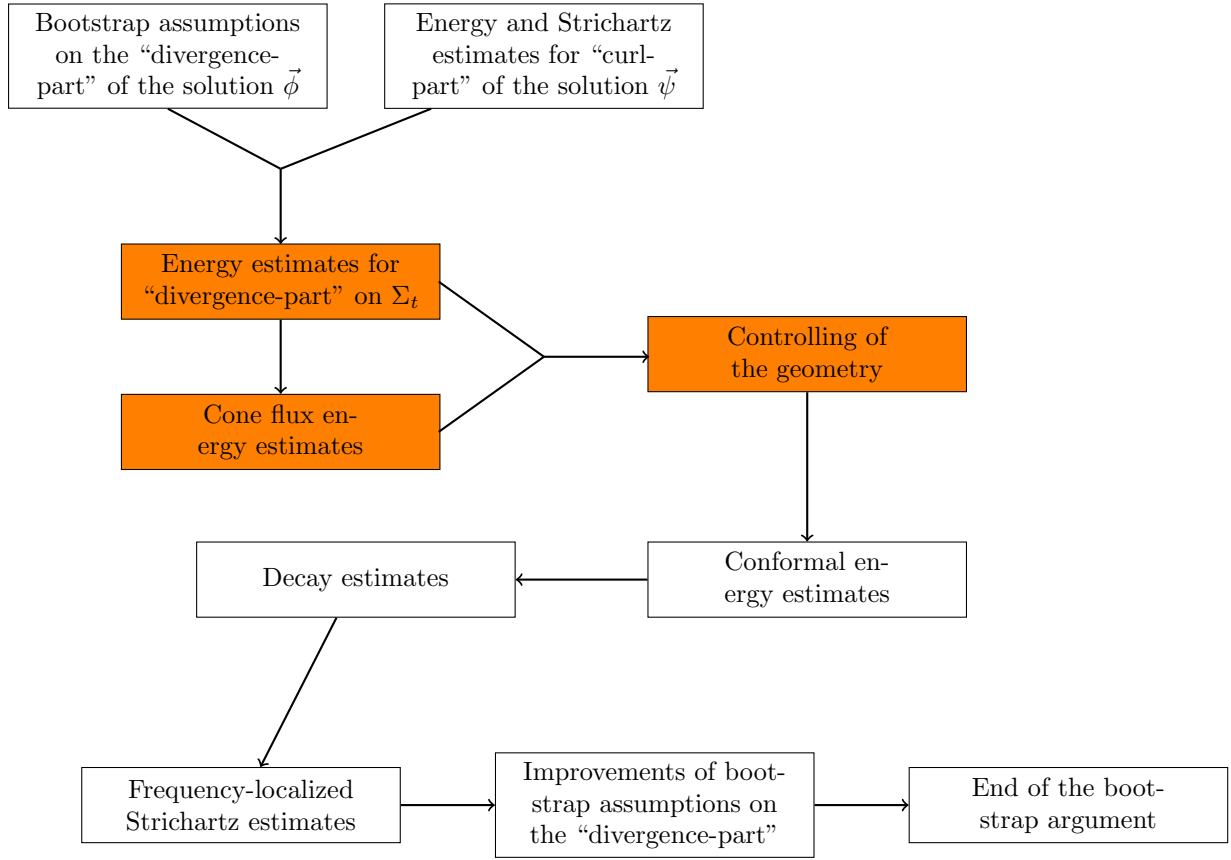
\begin{figure}[htbp]
		\begin{tikzpicture}
			\tikzstyle{process} = [rectangle, minimum width=3cm, minimum height=1cm, text centered, text width=4cm, draw=black]
			\node(ba)[process,xshift=0cm]{\hyperlink{BA}{Bootstrap assumptions} on the ``divergence-part" of the solution $\vvardiv$};
			\node(ee)[process,fill={rgb:red,1,;yellow,1},below of=ba, xshift=1.5cm, yshift=-2cm] {Energy estimates for ``divergence-part" on $\St$};
			\node(cfe)[process,fill={rgb:red,1,;yellow,1},below of=ee, xshift=0cm, yshift=-1cm] {Cone flux energy estimates};
			\node (cg)[process,fill={rgb:red,1,;yellow,1}, above of=cfe,right of=cfe, yshift=1, xshift=6cm]{Controlling of the geometry};
			\node(cee)[process,below of=cg, yshift=-2cm] {Conformal energy estimates};
			\node(de) [process,left of=cee, xshift=-5cm]{Decay estimates};
			\node(se)[process,left of=de, below of=de, xshift=-15,yshift=-2cm] {Frequency-localized Strichartz estimates};
			\node(te) [process,right of=ba, xshift=4cm]{Energy and Strichartz estimates for ``curl-part" of the solution $\vvarcurl$};
			\node(iw) [process,right of=se, xshift=4.5cm]{Improvements of bootstrap assumptions on the ``divergence-part"};
			\node(end)[process,right of=iw, xshift=4.5cm] {End of the bootstrap argument};
			\draw[thick,-](0,-0.7)--(1.5,-1.5);
			\draw[thick,-](3.5,-0.7)--(1.5,-1.5);
			\draw[thick,->](1.5,-1.5)--(1.5,-2.5);
			
			\draw[thick,->](ee)--(cfe);
			
			\draw[thick,-](3.6,-3)--(5,-4);
			\draw[thick,-](3.6,-5)--(5,-4);
			\draw[thick,->](5,-4)--(6.4,-4);
			
			\draw[thick,->](cg)--(cee);
			\draw[thick,->](cee)--(de);
			\draw[thick,->](de)--(se);
			\draw[thick,->](se)--(iw);
			\draw[thick,->](iw)--(end);
		\end{tikzpicture}
		\caption{Blueprint of proving low-regularity local well-posedness for admissible harmonic elastic materials} \label{steps}
	\end{figure}
	
	In below, we state some crucial steps listed in Figure \ref{steps}. More detailed discussions are provided in Section \ref{SS:mainidea}.\\
	
	\noindent	\textbf{\underline{Bootstrap assumptions on the ``divergence-part" of the solution $\vvardiv$:}}
	\begin{subequations}\label{BA:intro}
		\begin{align}
			\label{BA:intro0}\left(\p\vvardiv,\p\vvarcurl\right)\subset&\Domain,\\	
			\label{BA:introDiv0}\twonorms{\pfour\p\vvardiv}{t}{2}{x}{\infty}{[0,\Tstar]\times\mathbb{R}^3}^2+\sum\limits_{\upnu\geq2}\upnu^{2\delta_0}\twonorms{\littlewood\pfour\p\vvardiv}{t}{2}{x}{\infty}{[0,\Tstar]\times\mathbb{R}^3}^2\leq&1,
		\end{align}
	\end{subequations}
	Here, $\Domain$ is the defined hyperbolicity\footnote{Initially, we assume the system to be hyperbolic, and have strictly different wave speeds on $\Sigma_{0}$. And we will prove via a bootstrap argument that, the hyperbolicity, as well as the separation of wave speeds, can be maintained, for at least a short time, in the propagation of solutions.} region (see Section \ref{SS:Data}), $\littlewood$ is the Littlewood-Paley operator (defined in Section \ref{SS:definitionlittlewood}), the small parameter $\delta_0$ is defined in Section \ref{SS:ChoiceofParameters}, and $\Tstar$ is the bootstrap time.\\

	\noindent	\textbf{\underline{Energy estimates for $\vvardiv,\vvarcurl$, and Strichartz estimates for ``curl-part" of the solution $\vvarcurl$:}}\\
	
	In this step, we first establish the following theorem for the ``curl-part" of the solution $\vvarcurl$:
	\begin{theorem}[Energy estimates and Strichartz estimates for the ``curl-part"]\label{TH:introLinearStrichartz}
		Let $\delta_0$ and $\delta$ be defined as in (\ref{DE:ChoiceofParameters}). $\Tstar$ is the  existence time (see Section \ref{SS:bootstrap}). Then, for the the solution's ``curl-part" $\vvarcurl$, the following estimates hold:
		\begin{subequations}\label{ES:introlinearcurl}
			\begin{align}
				\label{ES:introenergyestimates1}&\sum\limits^3_{k=0}\sobolevnorm{\p_t^k\p\vvarcurl}{N-k}{\St}\lesssim D+1,\\
				&\label{ES:introlinearStri}\twonorms{\pfour^2\p\vvarcurl}{t}{2}{x}{\infty}{[0,\Tstar]\times\mathbb{R}^3}+\sum\limits_{\upnu\geq2}\upnu^{2\delta_0}\twonorms{\littlewood\pfour^2\p\vvarcurl}{t}{2}{x}{\infty}{[0,\Tstar]\times\mathbb{R}^3}^2\lesssim\Tstar^{2\delta},\\
				&\label{ES:introLinearSobolev}\onenorm{\p\vvarcurl,\pfour\p\vvarcurl}{}{\infty}{[0,\Tstar]\times\mathbb{R}^3}\lesssim 1.
			\end{align}
		\end{subequations}
		Here $D$ is the size of the initial data (defined in Theorem \ref{TH:IntroMainTheorem}).
	\end{theorem}
	We then further derive
	\begin{proposition}[Energy estimates for ``divergence-part"]\label{PR:introEnergyEstimatesforDiv}
		Under the initial data assumption in Theorem \ref{TH:IntroMainTheorem} and bootstrap assumptions \eqref{BA:intro}, for any $3<N<7/2$ and $t\in[0,\Tstar]$, solutions $\vvardiv$ to the wave equations (\ref{EQ:boxg}) satisfy the following estimates:
		\begin{align}
			\sum\limits^2_{k=0}\sobolevnorm{\pfour^k\p\vvardiv}{N-k-1}{\St}\lesssim&D+1,
		\end{align}
		where $D$ is defined in Theorem \ref{TH:maintheorem}.
	\end{proposition}
	
	 More discussions about energy estimates can be found in Section \ref{SSS:introEnergy}-\ref{SSS:introBA}. The detailed proofs are referred to Section \ref{S:EnergyEstimates}.\\
	
	\noindent	\textbf{\underline{Controlling of the geometry:}}\\

	For a fixed large $\lambda$, we rescaled the spacetime $(t,x)\mapsto(\lambda(t-t_k),\lambda x)$ as in Section \ref{SS:SelfRescaling}. We present the following control of the geometric quantities with respect to the rescaled metric $\gfour_{(\lambda)}$:
	\begin{proposition}[The main estimates for the eikonal function quantities, heuristic version]\label{PR:introMain} Let $\coneu$ denote\footnote{We note that $\coneu$ is the level set of the eikonal function $u$} the null cone of the rescaled metric $\gfour_{(\lambda)}$. Let $\intregion:=\{t\geq0,u\geq0\}$ be the interior region. Under the bootstrap assumptions \eqref{BA:intro}, we derive mixed space-time norm estimates for all the geometric quantities. In summary, we obtain the following estimates for the connection coefficients $\CC:=(\chismall,\hchi,\upzeta)$ (defined in Definition \ref{DE:DEFSOFCONNECTIONCOEFFICIENTS}- \ref{DE:conformalchangemetric})  under mixed-norms:
	\begin{align}
    \twonorms{\CC,\rgeo\angD_{\Lunit}\CC,\rgeo\angnabla\CC}{t}{2}{\sangle}{p}{\coneu},\twonorms{\rgeo^{1/2}\CC}{t}{\infty}{\sangle}{p}{\coneu},\holdertwonorms{\CC}{t}{2}{\sangle}{0}{\delta_0}{\coneu}\lesssim\lambda^{-1/2}.
	\end{align}	
Moreover, letting $\varepsilon_0:=\frac{N-3}{10}$, we obtain the following improved estimates for $\CC$, the conformal factor $\conformalfactor$ (defined in Definition \ref{DE:conformalchangemetric}), the modified mass aspect function $\modmass$ and the modified torsion $\modtorsion$ (defined in Definition \ref{DE:mass}) in the interior region $\intregion$:
\begin{subequations}
	\begin{align}
	\holderthreenorms{\angnabla\conformalfactor}{u}{2}{t}{2}{\sangle}{0}{\delta_0}{\intregion},\threenorms{\rgeo\modmass,\rgeo\angnabla\modtorsion}{u}{2}{t}{2}{\sangle}{p}{\intregion},\threenorms{\rgeo^{\frac{3}{2}}\modmass}{u}{2}{t}{\infty}{\sangle}{p}{\intregion}\lesssim&\lambda^{-4\varepsilon_0},\\ \twonorms{\rgeo^{1/2}\angnabla\conformalfactor}{t}{\infty}{\sangle}{p}{\coneu}, \twonorms{\angnabla\conformalfactor}{t}{2}{\sangle}{p}{\coneu}\lesssim&\lambda^{-1/2},\\
	\holderthreenorms{\chismall,\hchi}{t}{2}{u}{\infty}{\sangle}{0}{\delta_0}{\intregion},\twonorms{\upzeta}{t}{2}{x}{\infty}{\intregion}\lesssim&\lambda^{-1/2-4\varepsilon_0}.
\end{align}
\end{subequations}
	\end{proposition}
	
	We refer readers to Proposition \ref{PR:mainproof} for the full list of results for our control of geometric quantities. See Section \ref{SSS:introcontrol} for an introduction to the proof of Proposition \ref{PR:introMain}. Also, see Section \ref{S:connectioncoefficientsandpdes} for details of cone-flux energy estimates along null-hypersurfaces, and Section \ref{S:causalgeometry}-Section \ref{S:ControloftheGeometry} for details of the control of the various geometric quantities.\\
	\color{black}
    
	\noindent	\textbf{\underline{Conformal energy estimates:}}\\
	
	Building upon the appropriate control of the geometry in Proposition \ref{PR:introMain}, we proceed to derive the following estimate of the conformal energy $\Energyconformal[\varphi](t)$, which is defined in Definition \ref{DE:definitionofconformalenergy}:
	\begin{theorem}[Boundedness {of conformal energy}, heuristic version]\label{TH:introBoundnessTheorem}
		Let $\varphi$ be any solution of 
		\begin{align}\label{EQ:introCEE1}
			\boxg\varphi=0
		\end{align}
		in $[0,\Trescale]\times\St$ with $\varphi(1)$ supported in $B_R\subset\{u\geq 0\}\bigcap\Sigma_{1}$, where $B_R$ is defined in Theorem \ref{TH:Spatiallylocalizeddecay}. Then, under bootstrap assumptions \eqref{BA:intro}, for any $t\in[1,\Trescale]$, the conformal energy $\Energyconformal[\varphi]$ obeys
		\begin{align}
			\Energyconformal[\varphi](t)\lesssim(1+t)^{2\varepsilon}\left(\norm{\pfour\varphi}^2_{L^2(\Sigma_{1})}+\norm{\varphi}^2_{L^2(\Sigma_{1})}\right).
		\end{align}
		Here, $\varepsilon>0$ is an arbitrary small number and the constants in $``\lesssim"$ can depend on $\varepsilon$.
	\end{theorem}
	
	We refer readers to Section \ref{SSS:introconformal} for more introductions and Section \ref{SS:Sectionconformalenergy} for further discussions.\\
	
	\noindent	\textbf{\underline{Decay estimates and Strichartz estimates:}}\\
	
	{With the aid of conformal energy estimates in Theorem \ref{TH:introBoundnessTheorem}, in Section \ref{S:sectionreduction}} we prove the following decay estimate:
 
	\begin{theorem}[Decay estimate, heuristic version]\label{TH:introSpatiallylocalizeddecay}
	Let $\varphi$ be the solution to the equation $\boxg\varphi=0$ on $[0,\Trescale]\times\mathbb{R}^3$ with $\varphi(1,x)$ supported in the Euclidean ball $B_{R}$, where $B_{R}$ is defined in Theorem \ref{TH:Spatiallylocalizeddecay}.
	Then, there exist a large number $\Lambda$, a number $q>2$ sufficiently close to $2$, and a function $d(t)$ satisfying
		\begin{align}
			\norm{d}_{L_t^{\frac{q}{2}}([0,\Trescale])}\lesssim1,
		\end{align}	
		such that the following estimate holds for any $\lambda\geq\Lambda$ and $t\in[0,\Trescale]$:
		\begin{align}
			\norm{P_1\p_t\varphi}_{L_x^\infty(\St)}\lesssim\left(\frac{1}{(1+\abs{t-1})^{\frac{2}{q}}}+d(t)\right)\left(\norm{\pfour\varphi}_{L_x^2(\Sigma_1)}+\norm{\varphi}_{L_x^2(\Sigma_1)}\right).
		\end{align}
	\end{theorem}
	\color{black}
{Note that in Theorem \ref{TH:introSpatiallylocalizeddecay},  with respect to the rescaled metric, we provide a decay estimate within a rescaled time interval $[0,\Trescale]$. See Section \ref{SSS:reductionofstrichartz} for more discussions about the rescaling. }With the help of Theorem \ref{TH:introSpatiallylocalizeddecay}, we close our bootstrap argument by proving the following Strichartz estimates:
	\begin{theorem}[Improvement of the Strichartz-type bootstrap assumptions]\label{TH:introMainEstimates}
		Let $\delta>0$ be sufficiently small as {required} in Section \ref{SS:ChoiceofParameters}. Under the initial data assumptions stated in Theorem \ref{TH:IntroMainTheorem}, and the bootstrap assumptions \eqref{BA:intro},  with {$8\delta_0<\delta_1<N-3$} (see (\ref{DE:ChoiceofParameters})), the following estimates holds for $\vvardiv$:
		\begin{align}
			\label{ES:introestimateswavevariables}\twonorms{\pfour\dive\vvardiv}{t}{2}{x}{\infty}{[0,\Tstar]\times\mathbb{R}^3}^2+\sum\limits_{\upnu\geq2}\upnu^{2\delta_1}\twonorms{\littlewood\pfour\dive\vvardiv}{t}{2}{x}{\infty}{[0,\Tstar]\times\mathbb{R}^3}^2\lesssim\Tstar^{2\delta}.
		\end{align}
	\end{theorem}
	We note that for $\Tstar$ small, (\ref{ES:introestimateswavevariables}) is a strict improvement of (\ref{BA:introDiv}). {This can be seen via applying Schauder's estimates and noting $\curl\vvardiv=0$. See Section \ref{SSS:reductionofstrichartz} for more introductions and Section \ref{S:sectionreduction} for the details.}	
	
	\subsubsection{Difficulties in proving Theorem \ref{TH:IntroMainTheorem}}\label{SS:MainDifficulties}
	{In this section, we list the main new difficulties in proving Theorem \ref{TH:IntroMainTheorem} and state our strategies to overcome these difficulties. We refer readers to Section \ref{SS:mainidea} for detailed discussions of new ingredients and main ideas in the proof of Theorem \ref{TH:IntroMainTheorem}.}\\

	\noindent\textbf{\underline{Difficulties to find miraculous decoupling properties of elastic material:}}\\
	{For 3D elastic wave system \eqref{EQ:wave} with general nonlinearities, under planar symmetry,   the dynamics can be algebraically decomposed into characteristic waves featuring different speeds. In the ill-posedness results for elastic wave equations, An-Chen-Yin \cite{AnJMP,an2020low,an2022elasticwaves} analyzed the dynamics by decomposing the elastic waves into multiple characteristic strips and precisely proved that the plane-symmetric shock formation occurs for the fastest characteristic wave. However, without any symmetry assumptions, the techniques in proving low-regularity local well-posedness based on such decomposition are not developed. In fact, the system \eqref{EQ:wave} may not even be hyperbolic for some large data. Even if there is a regime where we have such decomposition, the interactions of the waves featuring different speeds is mysterious.
	
	In this paper, we focus on the admissible harmonic materials and we discover miraculous structures of these materials. Our discovery enables us to carry out the blueprint and to close the bootstrap arguments listed in Figure \ref{steps}. In particular, we first find that there presents a decoupling that the faster-wave does not get involved in the slower-wave dynamics. This phenomena reveals a stronger decoupling property than what presents for the harmonic elastic materials (There we have the preservation of pseudo-irrotationality). This stronger decoupling property allows us to derive (separate) higher regularity results for the ``curl-part". In addition, for the admissible harmonic elastic materials, we can describe the system's dynamics via employing two separate quasilinear wave systems, which are associated to distinct acoustic metrics. In our analysis, ``divergence-part" strictly represents the ``faster-wave" and ``curl-part" strictly represents the ``slower-wave". We emphasize again that such pleasant decomposition is customarily absent for the general elastic wave system.}\\
	
	\noindent\textbf{\underline{Difficulties connected to and different from the compressible Euler equations:}}\\	
	{In this article, in proving Theorem \ref{TH:IntroMainTheorem}, we require higher regularity for the ``curl-part" of the solution $\vvarcurl$ compared with the ``divergence-part". This phenomena also arises in studies of the compressible Euler equations and the relativistic Euler equations. In those cases, in the wave equations for the ``divergence-part", a higher regularity of the ``curl part" is involved. Thanks to the good structures of the Euler equations that the ``curl-part" of the solution satisfies a transport-div-curl system. Applying transport-div-curl estimates on such system, one can gain one derivative of the ``curl-part" back to solve the loss of regularity problem. This approach was carried out in \cite{3DCompressibleEuler,Roughsolutionsofthe3-DcompressibleEulerequations,yu2022rough,zhang2024wellposedness}. The optimal regularity assumptions on the ``curl-part" for the compressible Euler equations and the relativistic Euler equations still remains open. To our knowledge, however, such remarkable structure is not present in the elastic wave equations in general. We remark that, the elastic wave system is a multi-wave-speeds system, while the compressible Euler/relativistic Euler equations are wave-transport systems. It remains open whether low-regularity local well-posedness results can be obtained for general harmonic elastic materials.
	
	For the admissible harmonic elastic materials, we find and utilize the decoupling property to derive a higher regularity result for the ``curl-part", which enables us to reach our local well-posedness result under the similar regularity assumptions as in \cite{Roughsolutionsofthe3-DcompressibleEulerequations} for Euler equations. Following the blueprint in Figure \ref{steps}, we further derive desired Strichartz estimates with these new difficulties. In particular, we precisely address the impact of the ``curl-part" on the faster-wave dynamics and we employ the higher regularity estimates of ``curl-part" to overcome further potential difficulties in dealing with dangerous Ricci curvature terms.}\\
	
	\noindent\textbf{\underline{Difficulties in the energy estimates:}}\\
	{As mentioned before, for initial data, compared with the ``divergence part", we require one-extra regularity on ``curl-part" $\varcurl$, i.e.,}
	\begin{align}
		\sobolevnorm{\vardiv^i}{N}{\Sigma_0}+\sobolevnorm{\varcurl^i}{N+1}{\Sigma_0}\leq D, \quad \text{with} \, N>3.
	\end{align}
{Due to the decoupling property of our system, we find that the higher-order regularity of the ``curl-part" can be propagated by along the evolution of the solutions. 
	
	In addition, we also need to bound the cone-flux energy. Besides similar aforementioned difficulties, we encounter the obstacles that we have to show the cone-flux energy of the ``curl-part" is coercive along the null hypersurfaces $\coneu$ with respect to $\gfour$. This is equivalent to answering:  ``\textit{Are $\coneu$ always spacelike with respect to $\hfour$?}" We note that this is the important question that has to be answered for studying system with multiple traveling speeds. We give a confirmative answer to this question by:}
	\begin{itemize}
		\item In this article, we only consider the initial data in a region where hyperbolicity condition is satisfied. That is, $\left(\p\vvardiv,\p\vvarcurl\right)(\Sigma_{0})\subset\Int\,\breve{\Domain}\subset\breve{\Domain}\subset\Int\,\Domain$, where $\Domain,\breve{\Domain}$ are defined in Section \ref{SS:Data}. This ensure that, at the level of initial data, \textbf{1)} \eqref{EQ:elasticwaves0.1} is a hyperbolic system, and \textbf{2)} the ``divergence-part" represents the strictly faster-speed wave.
		\item To propagate the hyperbolicity condition we discussed above, we rely on {the Strichartz estimates, which is the main body of this article. We prove Strichartz estimates via a bootstrap argument.} Specifically, Theorem \ref{TH:introMainEstimates} provides an improvement of the bootstrap assumption \eqref{BA:introDiv0}. Finally, the bootstrap assumptions on the hyperbolicity condition \eqref{BA:intro0} can be improved by \eqref{BA:introDiv0} and the fundamental theorem of calculus.
	\end{itemize}
	
	\noindent\textbf{\underline{Difficulties in control of the geometry:}}\\
{For the wave equation of the ``divergence part", its Lorentz metric $\gfour$ verifies $\gfour=\gfour(\p\vvardiv,\p\vvarcurl)$. To control the geometry associated with $\gfour$ requires the estimates of $\p\vvardiv,\p\vvarcurl$ and their derivatives. This suggests that the ``curl-part" profoundly influences the dynamics of the ``divergence-part". Hence, we need to carefully take account of such impact when deriving the estimates for various geometric quantities related to $\gfour$. In particular, in Proposition \ref{PR:introMain} and Proposition \ref{PR:mainproof},} we use null structure equations and curvature decomposition to derive proper estimates for geometric quantities, where the ``curl-part" of the solution $\vvarcurl$ interacts in two scenarios:
	\begin{itemize}
		\item The first derivative of the metric $\gfour$ depends on $\p^2\vvarcurl$. In particular, the second fundamental form $k$ and the Christoffel symbol $\Chfour$ depend on $\pfour\cvvariables$. Those quantities enter as the below-top-order terms in the null structure equations and curvature decomposition, e.g., \eqref{introRaychau}, \eqref{introRicci}. This incorporation introduces low-order terms of the ``curl-part" into the system. In this scenario, one does not require higher regularity for $\vvarcurl$.
		\item {The curvature components depend on $\p^3\vvarcurl$.} A major difficulty arises in the curvature decomposition, e.g., for the first term on the RHS of equation \eqref{introRicci}, one can only gain extra regularity for the ``divergence-part" by using wave equation (\ref{EQ:introEW2.1}) and elliptic estimates, while one should view $\boxg\cvvariables\sim\pfour^2\p\vvarcurl$. Moreover, we gain extra regularity for the ``divergence-part" by paying the price of losing one regularity of the ``curl-part", since $\p^3\vvarcurl$ terms appear on the RHS of (\ref{EQ:introEW2.1}). Therefore, in controlling the geometry of the faster-wave, for the first term on the RHS of \eqref{introRicci}, we have:
		\begin{align}\label{ES:introboxgg}
			\boxg \gfour_{\alpha\beta}(\p\vvardiv,\p\vvarcurl)\sim\gensmoothfunction(\p^2\vvardiv,\p^3\vvarcurl),
		\end{align}
		which indicates that improvements are obtained only on the ``divergence-part". Such loss of regularity cause difficulties in energy estimates (as we discussed in the previous paragraph), and a dangerous rough Ricci curvature of the acoustic geometry of the faster wave. This necessitates deriving better (separate) results for the ``curl-part", and understand precisely how the ``curl-part" influences the ``divergence-part".
		To understand this loss of regularity phenomena, we go through the details in deriving the mix-norm estimates for various geometric quantities, where Ricci and Riemann curvature tensor components plays a crucial role as they appear as main error terms in the null structure equations (see, e.g., \eqref{introRaychau}).
	\end{itemize}
	\subsection{Main Ideas and New Ingredients}\label{SS:mainidea}
	
	{Facing the difficulties stated above, in this section we illustrate main ideas and new ingredients for the proof of Theorem \ref{TH:IntroMainTheorem}}.

	\subsubsection{Energy estimates on constant-time hypersurfaces}\label{SSS:introEnergy}
	To obtain the classical local well-posedness for the 3D elastic wave equations, the classical argument requires that initial data $\vec{U}$ belong to $H^{N}(\Sigma_{0})$ with $N>\frac{7}{2}$. To derive standard energy estimates, we first rewrite the elastic wave equations \eqref{EQ:el0} as:
	\begin{align}
		\mathcal{L}_{ik}U^k=0,
	\end{align}
where the operator matrix $\mathcal{L}_{ik}$ is defined as:

\begin{align}
	\mathcal{L}_{ik}:=&\square_{ik}-C_{ik}^{rs}(\p\vec{U})\p_r\p_s,
    \end{align}
  with
    \begin{align}
	\square_{ik}:=&\delta_{ik}\p_0^2-\left[\speedtwo^2\delta^{rs}\delta_{ik}+(\speedone^2-\speedtwo^2)\delta^r_{i}\delta^s_{k}\right]\p_r\p_s.
\end{align}

Then, in the spacetime region, we integrate the following identity:
	\begin{align}\label{EQ:introDQ}
		D_\alpha Q^\alpha(\p\vec{v})=2(D_0\vec{v})(\mathcal{L}\vec{v})+q,
	\end{align}
where
\begin{subequations}
	\begin{align}
		Q^0:=&\sum\limits_{i,k=1,2,3}\left\{(\p_0v^i)(\p_0v^i)+\speedtwo^2(\p_kv^i)(\p_kv^i)+(\speedone^2-\speedtwo^2)(\dive v)^2+(\p_r v^i)C_{ik}^{rs}(\p_sv^k)\right\},\\
		Q^r:=&\sum\limits_{i,k=1,2,3}\left\{-2\speedtwo^2(\p_0v^i)(\p_rv^i)-2(\speedone^2-\speedtwo^2)(\p_0v^r)(\dive v)-2(\p_0v^i)C^{rs}_{ik}(\p_sv^k)\right\},\\
		q:=&(\p_rv^i)C^{rs}_{ik}(\p_sv^k)-2(\p_0v^i)(\p_rC^{rs}_{ik})(\p_sv^k).
	\end{align}
\end{subequations}
 Substituting $\vec{v}$ by $\p^{k}\vec{U}$ with $k\leq N-1$, one then derive the following standard energy estimates:
	\begin{align}\label{ES:introsdenergy}
		\sobolevnorm{\vec{U}}{N}{\St}\lesssim\sobolevnorm{\vec{U}}{N}{\Szero}+\int_{0}^{t}\left(\onenorm{\pfour\p\vec{U}}{x}{\infty}{\Stau}+1\right)\sobolevnorm{\vec{U}}{N}{\Stau}\diff\tau.
	\end{align}
We refer reader to \cite[Section 7]{john88} for a detailed calculation of \eqref{EQ:introDQ} and the derivation of \eqref{ES:introsdenergy}. 
\color{black}
{We can further bound  $\twonorms{\p^2\vec{U}}{t}{\infty}{x}{\infty}{[0,\Tstar]\times\mathbb{R}^3}$ via applying Sobolev embedding $H^{\frac{3}{2}+}(\mathbb{R}^3)\hookrightarrow L^\infty(\mathbb{R}^3)$ on each constant-time hypersurface $\Sigma_{t}$. Together with a Gr\"onwall's argument, this gives the local well-posedness in $H^{N}$ with $N>\frac{7}{2}$. Note that back to (\ref{ES:introsdenergy}) we only need that $\twonorms{\p^2\vec{U}}{t}{1}{x}{\infty}{[0,\Tstar]\times\mathbb{R}^3}$ is bounded. Below the threshold $H^{7/2}$, the absence of Sobolev embedding necessitates new analytic methods to bound the RHS of (\ref{ES:introsdenergy}). A Strichartz-type estimates based on spacetime energy estimates is the answer and was introduced and used in the aforementioned works of Section \ref{SS:Overviewofpreviousresults}.
	
	Meanwhile, to utilize the decoupling property of the ``admissible harmonic elastic material" and to prove one-derivative higher top-order energy estimates (due to the reason explained in Section \ref{SSS:introBA}), for the ``curl-part", without making extra-regularity bootstrap assumptions, we derive energy estimates for ``divergence-part" and ``curl-part" directly based on their respective quasilinear wave equations (\ref{EQ:introEW2}). These spacetime energy estimates using wave equations \eqref{EQ:introEW2} and associated elliptic estimates are summarized into Proposition \ref{PR:EnergyEstimatesforDiv} and Theorem \ref{TH:LinearStrichartz}.}
	
	\subsubsection{Bootstrap assumptions on the ``divergence-part" and standard estimates for the ``curl-part"}\label{SSS:introBA}
	Thanks to the decoupling property with initial data stated in Theorem \ref{TH:IntroMainTheorem}, we are able to derive spacetime energy estimates (\ref{ES:introenergy2}) and linear Strichartz estimate (\ref{ES:introstrichartz2}). Note that we require an additional one-derivative regularity\footnote{We remark that similar requirement is also imposed for the compressible Euler equations and for the relativistic Euler equations as in \cite{3DCompressibleEuler,Roughsolutionsofthe3-DcompressibleEulerequations,yu2022rough,zhang2024wellposedness}. There the the local well-posedness has been proved for $(v^i,(\curl \vec{v})^i)\in H^{2+}$, where $\vec{v}$ is the velocity vector field with $v^i$ denoting its component.} for the ``curl-part" $\varcurl$ compared with ``divergence-part" $\vardiv$. This is because, utilizing equation (\ref{EQ:introEW2.1}), we have that $\vardiv$ depends on one more derivative of $\varcurl$. That is, $\p^k\vardiv$ is bounded by
$\gensmoothfunction(\p^k\vardiv,\p^{k+1}\varcurl)$, and $\gensmoothfunction$ is a certain smooth function. 
	
	Our main focus in this article is to control the ``divergence-part" with influence from the ``curl-part". In the beginning, we impose the following bootstrap assumptions:
	\begin{align}
		\label{BA:introDiv}\twonorms{\pfour\p\vvardiv}{t}{2}{x}{\infty}{[0,\Tstar]\times\mathbb{R}^3}^2+\sum\limits_{\upnu\geq2}\upnu^{2\delta_0}\twonorms{\littlewood\pfour\p\vvardiv}{t}{2}{x}{\infty}{[0,\Tstar]\times\mathbb{R}^3}^2&\leq1,
	\end{align}
	where $\littlewood$ is the Littlewood-Paley projection (defined in Section \ref{SS:definitionlittlewood}). A large part of this article is to recover the bootstrap assumptions. {The assumption (\ref{BA:introDiv}) is recovered and improved by the following Strichartz estimate:}
	\begin{align}
		\label{ES:introStrichartz1.1}\twonorms{\pfour\dive\vvardiv}{t}{2}{x}{\infty}{[0,\Tstar]\times\mathbb{R}^3}^2+\sum\limits_{\upnu\geq2}\upnu^{2\delta_1}\twonorms{\littlewood\pfour\dive\vvardiv}{t}{2}{x}{\infty}{[0,\Tstar]\times\mathbb{R}^3}^2&\leq\Tstar^{2\delta}.
	\end{align}
{Note that $\Tstar$ is the bootstrap time (or time of the local well-posedness) and small parameters $\delta_0,\delta_1$ satisfy $\delta_0>0$, $8\delta_0<\delta_1<N-3$ (see Section \ref{SS:ChoiceofParameters} for the precise definition of $\delta_0$ and $\delta_1$). For the divergence part, inequality (\ref{ES:introStrichartz1.1}) provides a strict improvement of (\ref{BA:introDiv}) since $\Tstar>0$ is small and $\delta_1>8\delta_0$. }
	\subsubsection{Reductions of the Strichartz type estimate}\label{SSS:reductionofstrichartz}
	
	As explained in Section \ref{SSS:introEnergy}, our argument crucially relies on the control of $\twonorms{\pfour\p\vvardiv}{t}{1}{x}{\infty}{\region}$. Here, we explain how we improve the bootstrap assumption (\ref{BA:introDiv}) and how to establish the Strichartz estimate (\ref{ES:introStrichartz1.1}). {Our argument is achieved through a series of reductions.} This reduction proceeds through the following steps: improvement of bootstrap assumptions$\longleftarrow$Strichartz estimate$\longleftarrow$Decay estimates$\longleftarrow$Conformal energy estimates$\longleftarrow$Controlling of the geometry. Note that here each left arrow signifies that the latter estimate implies the former. For clarity, we refer readers to the logic diagram in Figure \ref{steps}.
	\begin{itemize}
		\item \textit{Reduction to dyadic Strichartz estimates.} The first step in the proof of (\ref{ES:introStrichartz1.1}) is to reduce the Strichartz estimates to a dyadic version. Specifically, for a fixed large dyadic frequency $\lambda$, we partition $[0,\Tstar]$ into disjoint unions of sub-intervals $I_{k}:=[t_{k-1},t_k]$ and {the total number of $I_k$ being $\lesssim\lambda^{8\varepsilon_0}$ with $|I_k|\lesssim\lambda^{-8\varepsilon_0}$.} By Littlewood-Paley decomposition and Duhamel principle, the proof of (\ref{ES:introStrichartz1.1}) is then reduced to deriving the following dyadic Strichartz estimate:
		\begin{align}\label{ES:introDyadicStrichartz}
			\twonorms{P_\lambda\pfour\varphi}{t}{q}{x}{\infty}{[\tau,t_{k+1}]\times\mathbb{R}^3}\lesssim\lambda^{\frac{3}{2}-\frac{1}{q}}\norm{\pfour\varphi}_{L_x^2(\Sigma_{\tau})}.
		\end{align}
		Here, $\varphi$ is the solution of the covariant linear wave equation $\boxg\varphi=0$ on the time interval $I_k$. In (\ref{ES:introDyadicStrichartz}), we note that $q>2$ is any real number, which is sufficiently close to $2$ and we also require $\tau\in[t_{k},t_{k+1}]$. In this part, the reason we focus on large frequencies is that the control with small frequency is much easier due to the Bernstein inequalities.
		
		\item \textit{Self-Rescaling.} We then rescale our solution in terms of frequency $\lambda$. In particular, for $\vvariables:=\p\vvardiv,\cvvariables:=\p\vvarcurl$ (as in Definition \ref{DE:vvariablecvvariables}), we operate the following rescaling:
		\begin{subequations}
			\begin{align}
				\variables^i_{(\lambda)}:=&\variables^i(t_k+\lambda^{-1}t,\lambda^{-1}x),\\
				\cvariables^i_{(\lambda)}:=&\cvariables^i(t_k+\lambda^{-1}t,\lambda^{-1}x).
			\end{align}
		\end{subequations}
		With such rescaling, we can track the size of {the unknowns} and geometric quantities in terms of $\lambda$ (see Proposition \ref{PR:mainproof}). Meanwhile, we also have reduced the large-data small-time problem to a small-data long-time problem and we can comfortably explore the decay of the solutions over the rescaled (long) bootstrap time.
		\item \textit{Reduction to decay estimates.} For a large frequency $\lambda$, by rescaling the coordinates and applying an abstract $\mathcal{T}\mathcal{T}^*$ argument, we can reduce the proof of the dyadic Strichartz estimate (\ref{ES:introDyadicStrichartz}) to the establishment of the following $L^2-L^\infty$ decay estimates at any $t\in[0,\Trescale]$ with $\Trescale$ being the rescaled bootstrap time (see Section \ref{SS:SelfRescaling} for the definition of $\Trescale$):
		\begin{align}\label{introdecayestimates}
			\norm{P_1\Timelike\varphi}_{L_x^\infty(\St)}\lesssim\left(\frac{1}{(1+\abs{t-1})^{\frac{2}{q}}}+d(t)\right)\left(\norm{\pfour\varphi}_{L_x^2(\Sigma_1)}+\norm{\varphi}_{L_x^2(\Sigma_1)}\right).
		\end{align}
{Here, the timelike vector field $\Timelike$ (defined in \eqref{DE:Timelike}) is $\gfour$-unit normal to $\St$, and  $\varphi$ represents a solution to the equation $\boxg\varphi=0$ on the time interval $[0,\Trescale]\times\mathbb{R}^3$, satisfying that $\varphi(1,x)$ is supported in the Euclidean ball $B_{R}$ (see Theorem \ref{TH:Spatiallylocalizeddecay} for detailed definition of $R$).} Moreover, the function $d(t)$ satisfies
		\begin{align}
			\norm{d}_{L_t^{\frac{q}{2}}([0,\Trescale])}\lesssim1,
		\end{align}
		for $q>2$ sufficiently close to 2.
		\item \textit{Reduction to conformal energy estimates.} By employing product estimates and Littlewood-Paley theory, we further reduce the proof of (\ref{introdecayestimates}) to establishing the estimate below for the conformal energy $\Energyconformal[\varphi](t)$:
		\begin{align}\label{introconformalestimates}
			\Energyconformal[\varphi](t)\lesssim(1+t)^{2\varepsilon}\left(\norm{\pfour\varphi}^2_{L_x^2(\Sigma_{1})}+\norm{\varphi}^2_{L_x^2(\Sigma_{1})}\right),
		\end{align}
		{Here, the precise definition of the conformal energy $\Energyconformal[\varphi](t)$ is given in Section \ref{SS:SectionSetupconformal}. The constant $\varepsilon>0$ is an arbitrary small number, and $\varphi$ is an arbitrary solution to the equation $\boxg\varphi=0$, with $\varphi(1,x)$ supported in $B_R\subset\intregion\bigcap\Sigma_{1}$ (see Section \ref{SS:opticalfunction} for the definition of $\intregion$).}
	\end{itemize}
	
	The reduction of the Strichartz estimate is a very useful tool and was established before. We refer readers to \cite{AGeoMetricApproach,ImprovedLocalwellPosedness,Equationsd'ondesquasilineairesetestimationdeStrichartz,Equationsd'ondesquasilinearireseteffetdispersif,Sharplocalwell-posednessresultsforthenonlinearwaveequation,Strichartzestimatesforsecondorderhyperbolicoperators,ACommutingvector fieldsApproachtoStrichartz} for a comprehensive development of the underlying ideas and techniques. In this paper, {we emphasize that both the reduction of (\ref{introdecayestimates}) to (\ref{introconformalestimates}), and the very definition of $\Energyconformal[\varphi]$ rely crucially on our geometric framework. The sharp control on geometry is essential for deriving (\ref{introconformalestimates}).} We describe our approach for obtaining such geometric control in Section \ref{SSS:introgeometry}-Section \ref{SSS:introcontrol}.
	The full details of the aforementioned reduction process are provided in Section \ref{S:connectioncoefficientsandpdes}-Section \ref{S:ControloftheGeometry}.
	
	\subsubsection{Structures for the faster-wave geometry}\label{SSS:introgeometry}
	In order to reduce the decay estimates to the conformal energy estimates and to establish the conformal energy estimate (\ref{introconformalestimates}), one needs sharp information from the geometry. In this section, we discuss the geometric framework that is crucial for our analysis. The original idea of this part dates back to \cite[Section 9]{GlobalStabiliityOfMinkowski} by Christodoulou-Klainerman. The central object of this geometric framework is the eikonal function $u$, which is defined as the solution of the following eikonal equation:
	\begin{align}\label{EQ:Eikonal}
		\left(\gfour^{-1}\right)^{\alpha\beta}\p_\alpha u\p_\beta u=0.
	\end{align} 
	Here, $\gfour^{-1}$ represents the inverse of the Lorentzian metric $\gfour$. We denote the level set of $u$ by $\coneu$ and it corresponds to forward truncated null cones.
	
	Based on the null cones, a null pair $\Lunit, \uLunit$ (such that $\gfour(\Lunit,\uLunit)=-2$) is introduced. We then derive transport and Hodge type equations for the Ricci coefficients. Among them, an important example is the following Raychaudhuri equation (originally introduced in \cite{Raychaudhuri}):
	\begin{equation}\label{introRaychau}
		\Lunit\gtr\upchi+\frac{1}{2}(\gtr\upchi)^2=-\sgabs{\hchi}^2-k_{\spherenormal\spherenormal}\gtr\upchi-\Ricfour{\Lunit}{\Lunit}.
	\end{equation}
	{With $\Chfour_\alpha$ being the contracted Christoffel symbol (see Definition \ref{DE:Christoffelsymbols}) and $\gensmoothfunction$ being a smooth function of its arguments, we have
	\begin{align}\label{introRicci}
		\Ricfour{\alpha}{\beta}=-\frac{1}{2}\boxg \gfour_{\alpha\beta}(\vvariables,\cvvariables)+\frac{1}{2}\left(\Dfour_\alpha\Chfour_\beta+\Dfour_\beta\Chfour_\alpha\right)+\gensmoothfunction(\vvariables,\cvvariables)\cdot\pfour\gfour\cdot\pfour\gfour.
	\end{align}
	With the help of this remarkable decomposition of the Ricci curvature tensor\footnote{See \cite{ImprovedLocalwellPosedness} for original ideas.}, and the Bianchi identities, we substitute $\boxg \gfour_{\alpha\beta}(\vvariables,\cvvariables)$ in (\ref{introRicci}) using the elastic wave equations (\ref{EQ:introEW2}), and obtain estimates for the Ricci coefficients, which will be utilized in the subsequent control on geometry.}
	
	We reiterate the key distinction between this article and the previous works stated in Section \ref{SS:Overviewofpreviousresults} on quasilinear wave equations. In our paper, the metric $\gfour$ depends on both the ``divergence-part" ($\vvariables:=\p\vvardiv$) and the ``curl-part" ($\cvvariables:=\p\vvarcurl$) of the solution. {While, in the previous works for quasilinear wave equations, the principal part depends only on one single unknown, without nonlinear coupling with another wave.} In addition, the term $\boxg \gfour_{\alpha\beta}(\vvariables,\cvvariables)$ only gains regularity for the ``divergence-part" but would cause the loss of derivative for the ``curl-part".\\

    Originating\footnote{Also, see \cite{GlobalStabiliityOfMinkowski} for this framework applied to Einstein equations.} from \cite{ACommutingvector fieldsApproachtoStrichartz,ImprovedLocalwellPosedness}, the introduction of the geometric framework in the proof of low-regularity local well-posedness is motivated by the decay and conformal energy estimates. 
    \color{black}
    The advantage of using the above geometric framework in this low-regularity setting is that it reveals the dispersive properties of solutions to the wave equations. Specifically, for a solution $\varphi$ of wave equation $\boxg\varphi=0$, derivatives tangent to the characteristic null cones $\coneu$ have better decay than the transversal derivatives. {We now point out the importance of controlling the geometric quantities in below}: 
	\begin{itemize}
		\item The geometric framework that we set up must be well-defined. Particularly, we need to rule out the possibility of short-time shock formation arising from the intersection of distinct null cones. One also needs to overcome the derivative loss in defining the geometric frames, i.e., $\Lunit\sim \pfour\gfour$ as a result of \eqref{EQ:Eikonal}.
		\item To bound weighted energy for further deriving decay estimates, we apply the multiplier vector field method. The multipliers we use are related to both $\Lunit$ and $\Sigma_{t}$-tangent sphere normal vector $\spherenormal$. {Since both $\Lunit$ and $\spherenormal$ depend on the wave variables $\p\vvardiv,\p\vvarcurl$, the eikonal function $u$, and their first derivatives, the weighted energy estimate whence requires control on relative derivatives of these above quantities.  }
	\end{itemize}
	\subsubsection{Cone flux energy estimates}\label{SSS:ConeFluxEnergyEstimates}
	In Section \ref{S:connectioncoefficientsandpdes}, utilizing the spacetime energy estimates based on constant-time hypersurfaces for both ``divergence-part" and ``curl-part", we derive the cone flux energy estimates along the null hypersurfaces (with respect to $\gfour$) for rescaled variables as follows:
	\begin{subequations}\label{ES:introconeenergy}
		\begin{align}
			\label{ES:introconeenergy1}\mathcal{F}_{(1)}[\pfour\vvariables;\coneu]+\sum\limits_{\upnu>1}\upnu^{2(N-2)}\mathcal{F}_{(1)}[\littlewood\pfour\vvariables;\coneu]\lesssim\lambda^{-1},\\
			\label{ES:introconeenergy2}\mathcal{F}_{(2)}[\pfour\cvvariables;\coneu]+\sum\limits_{\upnu>1}\upnu^{2(N-3)}\mathcal{F}_{(2)}[\littlewood \pfour\cvvariables;\coneu]\lesssim\lambda^{-1},\\
			\label{ES:introconeenergy3}\mathcal{F}_{(2)}[\pfour^2\cvvariables;\coneu]+\sum\limits_{\upnu>1}\upnu^{2(N-3)}\mathcal{F}_{(2)}[\littlewood \pfour^2\cvvariables;\coneu]\lesssim\lambda^{-3}.
		\end{align}
	\end{subequations}
	Here, we have 
	\begin{subequations}\label{DE:introconeenergy}
		\begin{align}
			\mathcal{F}_{(1)}[\varphi;\coneu]&:=\int_{\coneu}\left((\Lunit\varphi)^2+\sgabs{\angnabla\varphi}^2\right)\diff\spherevol\diff t,\\
			\label{DE:introconeenergy2}\mathcal{F}_{(2)}[\varphi;\coneu]&\simeq\int_{\coneu}\abs{\pfour\varphi}^2\diff\spherevol\diff t
		\end{align}
	\end{subequations}
	with $\diff\spherevol$ being the volume form of the induced metric $\gsphere$ on the $\stu$-sphere from $\gfour$ and $\angnabla$ being the Levi-Civita connection on $\stu$ with respect to $\gsphere$. The precise form of $\mathcal{F}_{(2)}[\varphi;\coneu]$ is referred to Definition \ref{DE:acouticnullfluxes}.
	It is important to note that {the coerciveness of \eqref{DE:introconeenergy2} crucially relies on the fact that the ``divergence-part" (strictly) represents the faster-speed wave.} In particular, in view of the wave equation \eqref{EQ:introEW2.2} of the ``curl-part", in this paper we have that null hypersurfaces $\coneu$ are spacelike with respect to the metric $\hfour$. Combining energy estimates \eqref{ES:introsdenergy}, \eqref{ES:introconeenergy}, the linear Strichartz estimate for the ``curl-part", and bootstrap assumptions \eqref{BA:introDiv}, we are also able to control mixed-norms of $\boxg\gfour(\vvariables,\cvvariables)$ on the $\gfour$-null hypersurfaces as well as in the entire space-time as below:
	\begin{subequations}\label{ES:introfv}
		\begin{align}
			\label{introfv9}\holdertwonorms{\boxg \gfour}{t}{1}{x}{0}{\delta_0}{\region}\lesssim&\lambda^{-1-8\varepsilon_0},\\
			\label{introfv11}\twonorms{\rgeo(\angnabla,\angD_{\Lunit})\boxg\gfour}{t}{1}{\sangle}{p}{\coneu}\lesssim&\lambda^{-1-4\varepsilon_0},\\
			\label{introfv12}\threenorms{\rgeo\pfour\boxg\gfour}{t}{1}{u}{2}{\sangle}{p}{\region}\lesssim&\lambda^{-1/2-8\varepsilon_0}.
		\end{align}
	\end{subequations}
	We refer to Proposition \ref{estimateoffluidvariables} for more details and proofs. We emphasize that
	\eqref{ES:introfv} is crucial in controlling the geometry associated with metric $\gfour$, since $\boxg \gfour$ are the main terms to be bounded, when we study the Ricci curvatures.
	\subsubsection{Control of the conformal energy}\label{SSS:introconformal}
	A crucial part in the reduction steps to establish Strichartz estimates is to derive the decay estimates. As discussed in Section \ref{SSS:reductionofstrichartz}, we apply the conformal energy method introduced by Wang in \cite{AGeoMetricApproach}. We need to consider both equation $\boxg\varphi=0$ and the conformal wave equation with LHS being $\square_{\rescaledgfour}(e^{-\conformalfactor}\varphi)$. Here we have $\rescaledgfour:=e^{2\conformalfactor}\gfour$ with $\conformalfactor$ being the conformal factor defined in Definition \ref{DE:conformalchangemetric}). We are interested in these two equations for the following two reasons:
	\begin{itemize}
		\item We consider $\boxg\varphi=0$ since we need to prove the Strichartz estimate for solution $\varphi$ of the geometric equation $\boxg\varphi=0$.
		\item Following Wang's discovery in \cite{AGeoMetricApproach}, in our paper, we consider $\square_{\rescaledgfour}(e^{-\conformalfactor}\varphi)=\cdots$ because the Raychaudhuri equation with respect to the metric $\rescaledgfour$ reveals better geometric structure. To demonstrate this, we re-normalize the Raychaudhuri equation (\ref{introRaychau}) to overcome the difficulty that $\Dfour\Chfour$ terms on the RHS of \eqref{introRicci} are of top order. 
  Specifically, we consider the ``perturbative-part" $\chismall$ of the null expansion scalar $\Restrace\reschi$ for $\rescaledgfour$:
  \begin{align}\label{DE:introchismall}
		\chismall:=\gtr\upchi+\Chfour_{\Lunit}-\frac{2}{\rgeo}=\gtr\upchi+2\Lunit\conformalfactor-\frac{2}{\rgeo}=\Restrace\reschi-\frac{2}{\rgeo}.
		\end{align}
	 Then, we use a modified Raychaudhuri equation satisfied by $\chismall$:
		\begin{align}\label{EQ:introRaychau2}
			\Lunit\chismall+\frac{2}{\rgeo}\chismall=\lgensmoothfunction\cdot\boxg\gfour+\lgensmoothfunction\cdot(\chismall,\hchi)\cdot(\chismall,\hchi,\pfour\gfour)+\lgensmoothfunction\cdot\rgeo^{-1}\cdot\pfour\gfour.
		\end{align}		
	\end{itemize} 
	\color{black}
	Since the metric $\rescaledgfour$ is smoother only along null hypersurfaces, we first employ the original wave equation $\boxg\varphi=0$ and choose $X=f\spherenormal$ (see Section \ref{SSS:BoundnessTheorem} for the definition of $f$ and Section \ref{SS:geometricquantities} for the precise form of sphere normal vector $\spherenormal$) as the multiplier in a modified current. Applying the divergence theorem for this modified current in an appropriate region, we then obtain a Morawetz-type energy estimate and this leads to a uniform bound of the standard energy of $\varphi$, {along the union of a portion of the null cones and a portion of the constant-time hypersurfaces.}
	
	We then proceed to consider the conformally changed wave equations:
	\begin{subequations}\label{EQ:introboxtildeg}
		\begin{align}
			\label{EQ:introboxgrescaledvarphiconformal1}
			e^{2\conformalfactor}\boxgrescale\varphiconformal
			=&\anglap\varphiconformal-\uLunit\left(\Lunit\varphiconformal+\frac{1}{2}\Restrace\reschi\varphiconformal\right)-\left(\frac{1}{2}\Restrace\resuchi-k_{\spherenormal\spherenormal}\right)\left(\Lunit\varphiconformal+\frac{1}{2}\Restrace\reschi\varphiconformal\right)\\
			\notag&+2\left(\upzeta+\angnabla\conformalfactor\right)\cdot\angnabla\varphiconformal+\frac{1}{2}\left(\uLunit\Restrace\reschi+\frac{1}{2}\Restrace\reschi\Restrace\resuchi-\Restrace\reschi k_{\spherenormal\spherenormal}\right)\varphiconformal,\\
			\label{EQ:introboxgrescaledvarphiconformal2}
			e^{3\conformalfactor}\boxgrescale\varphiconformal=&\boxg\varphi-\left(\boxg\conformalfactor+\Dfour^\alpha\conformalfactor\Dfour_\alpha\conformalfactor\right)\varphi.
		\end{align}
	\end{subequations}
	We apply the multiplier approach using $\rgeo^p\Lunit$ type vector fields in the region $\{\tau_1\leq u\leq\tau_2\}\bigcap\{\rgeo\geq R\}$, where $1\leq\tau_1<\tau_2<\Trescale$, to control the conformal energy in the exterior region and to provide energy decay for each null slice. Finally, we control the conformal energy in the interior region with the help of the argument in \cite{Anewphysical-spaceapproachtodecayforthewaveequationwith}, ensuring energy decay in each spatial-null slice.
	
	The very definition of the conformal energy, as well as its analysis, requires delicate and precise control of the geometry. The boundedness theorem of the conformal energy is stated in Theorem \ref{TH:BoundnessTheorem}.
	
	\subsubsection{Control of the geometry} \label{SSS:introcontrol}{To control the conformal energy, it requires the control of the geometry. Klainerman-Rodnianski \cite{ImprovedLocalwellPosedness} and Wang \cite{AGeoMetricApproach} developed a novel approach to bound the geometric quantities in the low-regularity settings. Compared with \cite{AGeoMetricApproach,ImprovedLocalwellPosedness}, where $\gfour=\gfour(\vvariables)$, in this paper we encounter the metric $\gfour=\gfour(\vvariables,\cvvariables)$. In addition, the Ricci curvatures in this paper are less regular in the ``curl-part" due to the structure of equations \eqref{introRicci} and \eqref{ES:introboxgg}. \\

    Here, we state a way to control our geometric quantities with $\gfour=\gfour(\vvariables,\cvvariables)$. We first derive the PDEs that are satisfied by the geometric quantities. These include the geometric transport equations and the div-curl system for the connection coefficients. These equations depend on our geometric formalism and are independent of the elastic wave equations. We also derive the estimates for certain Ricci and Riemann curvature components via employing the decomposition of the Ricci curvatures and the Bianchi identities. It is at this step that the structure of the elastic wave equations is used. Specifically, we can substitute $\boxg \gfour_{\alpha\beta}(\vvariables,\cvvariables)$ in (\ref{introRicci}) via appealing to the elastic wave system (\ref{EQ:introEW2}). Then, by combining the geometric transport equations, the aforementioned estimates for Ricci and Riemann curvature components, and estimates \eqref{ES:introfv}, we further derive and analyze the equations for various quantities, including the important modified mass aspect function $\modmass$ (see \eqref{DE:modmass} for the definition of $\modmass$) and the conformal factor $\conformalfactor$ (see Definition \ref{DE:conformalchangemetric} for $\conformalfactor$). }Finally, we derive mixed space-time norm estimates, as in Proposition \ref{PR:introMain} for various geometric quantities, which are essential for deriving the conformal energy estimates. The key in proving Proposition \ref{PR:introMain} is to show\footnote{In \eqref{ES:Frame} $\Dfour$ denotes the covariant derivative with respect to $\gfour$.}:
	\begin{subequations}\label{ES:Frame}
		\begin{align}
			(\Dfour_\Lunit,\Dfour_A)\Dfour(\Lunit,\uLunit,e_A)\lesssim&(\Lunit,\angnabla)\pfour\gfour,& \text{ on $\gfour$-null hypersurfaces}&,\\
			\Dfour^2(\Lunit,\uLunit,e_A)\lesssim&\pfour^2\gfour,& \text{ on constant-time hypersurfaces}&,
		\end{align}
	\end{subequations}under appropriate mixed space-time norms, so that we can control geometry using energy estimates on constant-time hypersurfaces \eqref{ES:introsdenergy} and on $\gfour$-null hypersurfaces (\ref{ES:introconeenergy}).
Moreover, in proving Proposition \ref{PR:introMain}, we carefully bound several renormalized geometric quantities $\modmass,\modtorsion,\chismall$ and the rough conformal factor $\conformalfactor$ via establishing transport estimates along null geodesics and via conducting Hodge estimates on $\stu$. We refer readers to Proposition \ref{PR:mainproof} for the detailed version of Proposition \ref{PR:introMain}. These estimates are necessary for proving conformal energy estimates in Section \ref{SS:Sectionconformalenergy}.
\color{black}	
	
	\subsection{Paper Outline}
	The logical structure of this paper is summarized as below:
	\begin{itemize}
		\item In Section \ref{S:StandardElasticWave}, we introduce 3D elastic wave system and our focus: \emph{admissible harmonic elastic materials}.
		\item In Section \ref{S:GeometricFormulation}, we derive the geometric formulation of the admissible harmonic elastic wave system.
		\item In Section \ref{S:PDEs}, we present our decompositions of the elastic wave system, that are used throughout the paper.
		\item In Section \ref{S:sectionmainthm}, we define the Littlewood-Paley projections, which are frequently employed in our analysis. We also state our main theorem and our bootstrap assumptions.
		\item In Section \ref{S:preliminary}, we provide preliminary estimates relying on Littlewood-Paley calculus and standard elliptic estimates. Then, we present the frequency-projected version of ``admissible elastic wave equations" in Lemma \ref{LE:lwpequations}. 
		\item In Section \ref{S:EnergyEstimates}, we deduce energy estimates on constant-time hypersurfaces and linear Strichartz estimates for the ``curl-part" in Theorem \ref{TH:LinearStrichartz}. Then, under the bootstrap assumptions, we derive the energy estimates for the ``divergence-part" on constant-time hypersurfaces. Also, we prove $L^2$-elliptic estimates for $\boxg\vvariables$ in Lemma \ref{LE:ellipticforboxg}
		\item In Section \ref{S:sectionreduction}, following the approach of Tataru \cite{Strichartzestimatesforsecondorderhyperbolicoperators} and Wang \cite{AGeoMetricApproach}, we rescale the variables and reduce the proof of the Strichartz estimates to deriving the spatially localized decay estimates. 
		\item {In Section \ref{S:sectiongeometrysetup}, we establish the geometric framework by constructing the eikonal function $u$ and setting up its geometry, including constructing an appropriate null frame.} Finally, we define the conformal energy and state its boundedness theorem as in Theorem \ref{TH:BoundnessTheorem}, which plays a crucial role in deriving the decay estimates stated in Theorem \ref{TH:Spatiallylocalizeddecay}.
		\item In Section \ref{S:connectioncoefficientsandpdes}, we define the null fluxes in Section \ref{energyalongnullhypersurfaces} and prove the null flux energy estimates along the $\gfour$-null hypersurfaces in Section \ref{SS:EnergyFluxes}.
		\item In Section \ref{S:causalgeometry}, we define additional geometric quantities, including the connection coefficients with respect to the null frame, the conformal factor, mass aspect functions, and curvature tensor components. We state the null structure equations in Section \ref{SS:PDES} and the curvature decomposition in Section \ref{SS:Curvature}. {In Proposition \ref{PR:InitialCT1} and Proposition \ref{PR:InitialCT2}, we also restate the estimates from \cite{AGeoMetricApproach,3DCompressibleEuler}, which yield control over geometry along the initial data hypersurface.} Finally, in Section \ref{SS:mainversionPDES} and Section \ref{SS:proofofmainversion}, we state and derive the main PDEs that will be employed to control the geometry.
		\item In Section \ref{S:ControloftheGeometry}, we state the bootstrap assumptions for the rescaled ``divergence-part", the linearized Strichartz estimates satisfied by the ``curl-part",  as well as the bootstrap assumptions for geometric quantities in Section \ref{SS:BAGeo}. {Then, we state the main estimates for the geometric quantities in Proposition \ref{PR:mainproof}. The detailed proof of this proposition is provided in Section \ref{SS:controlacoustic}. In the proof, we appeal to the preliminary geometric estimates in Section \ref{SSS:TraceEstimates}, and estimates of $\vvariables,\cvvariables,\boxg\gfour$ in mixed-norms in Proposition \ref{estimateoffluidvariables}.}
	\end{itemize}
	
	\subsection{Notations}\label{SS:Notations}
	We use the following notations throughout the paper:
	\begin{itemize}
		\item Greek ``space-time" indices range over $0,1,2,3$, while Latin ``spatial" indices take on the values $1,2,3$. 
		\item We adopt the Einstein summation convention.
		\item $\antisymmetic_{\alpha\beta\gamma\delta}$ and $\antisymmetic_{ijk}$ denote the fully antisymmetric symbol, normalized by $\antisymmetic_{0123}=1$ and $\antisymmetic_{123}=1$.
		\item We raise and lower indices for tensor fields with respect to metric\footnote{See Definition \ref{DE:Defg} for the definition of $\gfour$ and its inverse.} $\gfour$ and its inverse. For example, $V^\alpha=(\gfour^{-1})^{\alpha\beta}V_\beta$ and $V_\beta=\gfour_{\alpha\beta}V^\alpha$.
		\item We denote $\St:=\{(t^\prime,x^1,x^2,x^3)\in\mathbb{R}^{1+3}|t\equiv t^\prime\}$ as the standard constant-time slice.
		\item $\linsmoothfunction[A](B)$ denotes any scalar-valued function, that is linear in the component of $B$, with coefficients depending on the components of $A$.
		\item $\quadsmoothfunction[A](B,C)$ denotes any scalar-valued function quadratic in the components of $B$ and $C$, with coefficients that are a function of the components of $A$.
		\item $\gensmoothfunction(A)$ represents a generic smooth function in terms of $A$.
		\item We denote the schematic spatial partial derivatives by $\p$ and the schematic space-time partial derivatives by $\pfour$.
		\item $A\lesssim B$ means $A\leq C\cdot B$ for some universal constant $C$.
		\item For a vector field $\vec{V}\in\mathbb{R}^3$, we use $\curl V^i:=\antisymmetic_{iab}\p_a V^b$ to denote the curl of $\vec{V}$, and we set $\dive V:=\p_iV^i$ to represent the divergence of $\vec{V}$.
	\end{itemize}
	We also refer readers to Appendix \ref{AP:notations} for the full list of notations that we use throughout this article.

	\subsection{Acknowledgment} The authors would like to thank Jared Speck for pointing out reference \cite{John1960,John1966}. XA is supported by MOE Tier 1 grants A-0004287-00-00, A-0008492-00-00 and MOE Tier 2 grant A-8000977-00-00. HC acknowledges the support of MOE Tier 2 grant A-8000977-00-00 and NSFC (Grant No. 12171097). SY acknowledges the support of MOE Tier 2 grant A-8000977-00-00.

	\section{The Elastic Wave Equations}\label{S:StandardElasticWave}
	We study the low-regularity well-posedness problem for the elastic wave equations in three spatial
	dimensions\footnote{The ideas and techniques also apply to the 2D case.}.	In this section, we derive the general 3D elastic wave equations, and introduce the (admissible) harmonic materials that are the focus of our study.
	\subsection{The Standard Elastic Wave Equations and the Truncated Elastic Wave Equations}\label{SS:elasticwaveeq}
	The motion of a 3D elastic body is described by time-dependent orientation-preserving diffeomorphisms, denoted by $\vec{P}:\mathbb{R}^3\times\mathbb{R}\to\mathbb{R}^3$, where $\vec{P}=\vec{P}(\vec{x},t)$ satisfies $\vec{P}(\vec{x},0)=\vec{x}=(x^1,x^2,x^3)$. The deformation gradient is defined as $F=\nabla P$  with  $F_{ij}:=\partial P^i/\partial x^j.$ For a homogeneous isotropic hyperelastic material, the stored energy $W$ depends solely on the principal invariants of $FF^T$ (the Cauchy-Green strain tensor). The equations of motion can be derived by applying the principle of least action to the action functional $\mathcal{S}$:
	\begin{align}\label{EQ:Lagr1}
		\mathcal{S}:=\int\int_{\mathbb{R}^3}\left\{\frac12|\partial_t\vec{P}|^2-W(FF^T)\right\}\diff x\diff t.
	\end{align}
	The resulting Euler-Lagrange equations are as follows:
	\begin{equation}\label{EQ:el}
		\frac{\partial^2P^i}{\partial t^2 }-\frac{\partial}{\partial{x^l}}\frac{\partial W(FF^T)}{\partial F^i_{l}}=0.
	\end{equation}
	With $\vec{U}=\vec{P}-\vec{x}$ being the displacement, we define $G:=\nabla U=F-I$ to be the displacement gradient, which leads to $FF^T=I+G+G^T+GG^T$. Then, as shown in \cite{Sideris}, the stored energy function $W$ can be expressed as $W=\hat{W}(k_1,k_2,k_3)$. Here, $k_1,k_2,k_3$ denote the principal invariants of the strain matrix $C=G+G^T+GG^T$, given by
	\begin{subequations}
		\begin{align}
			k_1&=\mu_1+\mu_2+\mu_3=\text{tr} C,\\
			k_2&=\mu_1\mu_2+\mu_1\mu_3+\mu_2\mu_3=\frac12\{(\text{tr} C)^2-\text{tr} C^2\},\\
			k_3&=\mu_1\mu_2\mu_3=\frac16\{(\text{tr} C)^3-3(\text{tr} C)(\text{tr} C^2)+2\text{tr} C^3\},
		\end{align}
	\end{subequations}
	with $\mu_1,\ \mu_2,\ \mu_3$ being the eigenvalues of $C$. Applying Taylor expansion to $\hat{W}$ near $k_i=0$, we have
	\begin{equation}\label{EQ:taylor}
		\hat{W}=\sigma_0+\sigma_1k_1+\frac12\sigma_{11}k_1^2+\sigma_2 k_2+\frac16\sigma_{111}k_1^3+\sigma_{12}k_1k_2+\sigma_3k_3+\cdots
	\end{equation}
	where the constant coefficients $\sigma_0,\sigma_1,\sigma_{11}$, etc., {correspond to the partial derivatives of $\hat{W}$ at $k_i=0$ ($i=1,2,3$).} Here, the material under consideration is assumed to satisfy the stress-free\footnote{We refer readers to Tahvildar-Zadeh in \cite{shadi} for a study on the $\sigma_1\neq 0$ case, where the small data global existence under certain null condition assumptions was shown.} reference state condition, i.e., $\sigma_1=0$. Then, we set $4(\sigma_{11}+\sigma_2)$ and $-2\sigma_2$ to be the positive Lam\'e constants. For more details on the derivation, interested readers are referred to \cite{agemi,Sideris}, which studied the small-data problem for the following truncated (at the cubic order) elastic wave system:
	{\begin{proposition}\label{PR:ElasticWaves}\cite[Truncated elastic wave equations]{Sideris}.
		For homogeneous isotropic hyperelastic materials, the motion of the
		displacement $\vec{U} = (U^1, U^2, U^3)$ is governed by the following quasilinear wave system:
		\begin{align}\label{EQ:elasticwaves1}
			\p_t^2U^i-\speedtwo^2\Delta U^i-(\speedone^2-\speedtwo^2)\p_i(\nabla\cdot \vec{U})=\mathcal{G}^i(\p \vec{U},\p^2\vec{U}),
		\end{align}
		where the two constants $\speedone>\speedtwo>0$, representing the speeds of the longitudinal and transverse waves, are defined as
		\begin{align}
			c_1^2=4\sigma_{11},\quad c_2^2=-2\sigma_2.
		\end{align}
		Moreover, the nonlinearity $\mathcal{G}(\p \vec{U},\p^2\vec{U})$ is of a quadratic form, and can be expressed as follows:
		\begin{align}
			\mathcal{G}^i(\p \vec{U},\p^2 \vec{U}):=&\sum\limits_{j,l,m=1}^3C^{lm}_{ij}(\p\vec{U})\p_l\p_m U^j,\\
			C^{lm}_{ij}(\p\vec{U}):=&\sum\limits_{k,n=1}^3C^{lmn}_{ijk}\p_n U^k.
		\end{align} 
		The coefficients also verify $C^{lm}_{ij}(\p\vec{U})=C^{ml}_{ij}(\p \vec{U})$, see \cite[Section 4]{Sideris}.
	\end{proposition}}
	\begin{remark}
		We emphasize that, for the local well-posedness problems, proving a truncated version of the result is not indicative of the full elastic wave system's behaviour; see Section \ref{SS:SmallData} for the discussion.
	\end{remark}
	\subsection{Harmonic Elastic Materials and Admissible Harmonic Elastic Materials}\label{SS:ModifiedElasticWaves}
	Assuming favorable evolution of $\curl U$, we can control $\dive U$ by taking divergence of the elastic wave equations (\ref{EQ:elasticwaves1}). However, it is unknown whether such favorable evolution exists for all elastic wave equation. Hence, we investigate the following special elastic material, which we call admissible harmonic elastic materials.
	
	Before introducing the admissible harmonic elastic materials, we first recall the definition of harmonic elastic material in \cite{John1966,John1960}:
	\begin{definition}\cite[Pseudo-irrotationality and harmonic elastic materials]{John1966}\label{DE:pseudoirrotational}.
		We call the deformation $\vec{U}$ \textbf{pseudo-irrotational} (at a particular moment) if 
		\begin{align}\label{EQ:pseudoirrotational}
			\p_k U^i-\p_i U^k=0.
		\end{align}
		Pseudo-irrotationality of the deformation is equivalent to symmetry of the
		Jacobian matrix F (defined in Section \ref{SS:elasticwaveeq}), or, at least locally, to the existence of a "deformation potential",
		i.e., a scalar $\tilde{\phi}$ such that $U^i=\p_i\tilde{\phi}$.
		
		We call the material \textbf{harmonic} if and only if
			for such material,
			the pseudo-irrotational initial conditions in the absence of body forces always lead to pseudo-irrotational motions. That is, the material is harmonic if and only if (\ref{EQ:pseudoirrotational}) on $\Sigma_0$ implies (\ref{EQ:pseudoirrotational}) on $\Sigma_t$ for all $t$.
        \end{definition}
	For hyperelastic materials satisfying (\ref{EQ:el}), the following equation holds for deformation $\vec{U}$:
	\begin{align}\label{EQ:elasticwaves1.1}
		\p_t^2\vec{U}-\speedtwo^2\Delta \vec{U}-(\speedone^2-\speedtwo^2)\nabla(\nabla\cdot \vec{U})=\dive \left\{E(\p \vec{U})\right\},
	\end{align}
	where $E(\p\vec{U})$ is a $2$-vector field depending on $\p\vec{U}$. We note that $\dive E$ is the nonlinear-part\footnote{In \cite{agemi}, Agemi considered the truncated version of the elastic wave. Nevertheless, \eqref{EQ:elasticwaves1.1} holds true for generic elastic wave equation with $E$ not necessarily being quadratic.} of (\ref{EQ:el}); see \cite[(2.8)]{agemi}. Now we define the admissible harmonic elastic materials as below.
	
	\begin{definition}[Admissible harmonic elastic materials]
		Let $I$ be the identity matrix, $b$ be a constant, $F$ be the deformation gradient defined in Section \ref{SS:elasticwaveeq}, $\det(F)$ be the determinant of $F$. We call the elastic material \textit{admissible harmonic} if the deformation $\vec{U}$ satisfies equation \eqref{EQ:elasticwaves1.1} with
		\begin{align}\label{DE:H}
			E(\p \vec{U})=&G(\p\vec{U})I+b\det(F)(F^T)^{-1}.
		\end{align}
	Here, $G(\p\vec{U})$ is a smooth scalar function of $\p\vec{U}$. In other words,  we say that $\vec{U}$ satisfies the following admissible harmonic elastic wave system if $\vec{U}$ satisfies:
	\begin{align}\label{EQ:elasticwaves1.3}
		\p_t^2\vec{U}-\speedtwo^2\Delta \vec{U}-(\speedone^2-\speedtwo^2)\nabla(\nabla\cdot \vec{U})=\dive \left\{G(\p\vec{U})I+b\det(F)(F^T)^{-1}\right\}.
	\end{align}
	\end{definition}

	\begin{definition}[Hyperbolicity of admissible elastic wave equations]
		We say the admissible elastic wave system \eqref{EQ:elasticwaves1.3} is hyperbolic if $G(\p\vec{U})$ satisfies\footnote{In this article, we only consider the system that are hyperbolic. Hyperbolicity assumption (\ref{EQ:symmetry}) can be imposed on initial data and be propagated (for at least a short time); see Section \ref{SS:Data}-\ref{SS:bootstrap}}:
		\begin{align}\label{EQ:symmetry}
			C^{-1}\abs{\xi}^2\leq\sum\limits_{j,k=1}^3\left\{(\speedone^2-\speedtwo^2)\delta_{jk}+G^\prime(\p_jU^k)\right\}\xi_j\xi_k\leq C\abs{\xi}^2,
		\end{align} 
		for some $C>0$ and all $\xi\in\mathbb{R}^3$.
	\end{definition} 
Meanwhile, we have that
\begin{proposition}[Equations satisfied by the deformation $\vec{U}$ of the admissible harmonic elastic materials]\label{PR:admissible2}
	Let $\vec{U}$ be the solution of the admissible harmonic elastic wave equations \eqref{EQ:elasticwaves1.3}. Then, $\vec{U}$ solves the following equation:
	\begin{align}\label{EQ:elasticwaves1.2}
		\p_t^2U^i-\speedtwo^2\Delta U^i-(\speedone^2-\speedtwo^2)\p_i(\nabla\cdot \vec{U})=\sum\limits_{j,k=1}^3G^\prime(\p_jU^k)\p_i\p_jU^k,
	\end{align}
	where $G^\prime(\p_jU^k):=\frac{\p G}{\p(\p_jU^k)}$ is a smooth scalar function of $\p_jU^k$.
\end{proposition}
\begin{proof}[Proof of Proposition \ref{PR:admissible2}]	
	We consider the RHS of (\ref{EQ:elasticwaves1.3}). First, we compute $\dive\left\{\det(F)(F^T)^{-1}\right\}$. Since $F_{ij}:=\p_j P_i$ as listed in Section \ref{SS:elasticwaveeq}, the following identity holds:
\begin{align}\label{EQ:pF}
	\p_kF_{ab}=\p_bF_{ak}.
\end{align}
Then, by \eqref{EQ:pF} and Jacobi's formula, we derive:
	\begin{align}\label{EQ:DivFt}
		\dive\left\{\det(F)(F^T)^{-1}\right\}=&\p_k\left[\det(F)(F^{-1})^{ki}\right]\\
		\notag=&\det(F)\text{Tr}\left[(F^{-1})^{ba}\cdot(\p_kF_{ac})\right]\cdot(F^{-1})^{ki}-\det(F)\cdot(F^{-1})^{ka}\cdot(\p_kF_{ab})\cdot(F^{-1})^{bi}\\
		\notag=&\det(F)\cdot(F^{-1})^{ba}\cdot(\p_kF_{ab})\cdot(F^{-1})^{ki}-\det(F)\cdot(F^{-1})^{ka}\cdot(\p_bF_{ak})\cdot(F^{-1})^{bi}\\
		\notag=&0.
	\end{align}	

For the first term $\dive\{G(\p\vec{U})I\}$ on the RHS of (\ref{EQ:elasticwaves1.3}), there also holds
	\begin{align}\label{EQ:DivG}
		\left(\dive\{G(\p\vec{U})I\}\right)^i=\p_i(G(\p\vec{U}))=\sum\limits_{j,k=1}^3G^\prime(\p_jU^k)\p_i\p_jU^k.
	\end{align}	
	Gathering \eqref{EQ:elasticwaves1.3}, \eqref{EQ:DivFt} and (\ref{EQ:DivG}), we then obtain the desired equations \eqref{EQ:elasticwaves1.2} for admissible harmonic elastic materials.
\end{proof}
\color{black}
	\begin{proposition}\label{PR:admissible}
		The admissible harmonic elastic material defined in (\ref{DE:H}) is harmonic in the sense of Definition \ref{DE:pseudoirrotational}.
	\end{proposition}
	\begin{proof}[Proof of Proposition \ref{PR:admissible}]
		By taking curl of (\ref{EQ:elasticwaves1.2}), we derive
		\begin{align}
			\p_t^2\curl \vec{U}-\speedtwo^2\Delta\curl\vec{U}=0.
		\end{align}
		Hence, $\curl\vec{U}=0$ at some time $t$ implies $\p_t^2\curl \vec{U}=0$ at the same time. We thus obtained the desired result.
	\end{proof}
	\section{A Geometric Formulation of the Admissible Harmonic Elastic Materials}\label{S:GeometricFormulation}
	
	{In this section, we first perform a div-curl decomposition for our admissible harmonic elastic wave system. Then, we present a geometric formulation for the decomposed waves.}\\

    {Let the scalar function $\tilde{\phi}$ be the solution of the following elliptic equation:
	\begin{align}
		\Delta\tilde{\phi}=\dive \vec{U}.
	\end{align} 
	Then, we can decompose $\vec{U}$ in the following way:
	\begin{align}\label{EQ:Helmholtz}
		\vec{U}^i=\p_i\tilde{\phi}+\vvarcurl^i
	\end{align}
	with $\vvarcurl$ being a divergence free vector field.
	
	\begin{definition}\label{DE:DefPhiPsi}
		Let $\tilde{\phi},\tilde{\psi}$ be as in (\ref{EQ:Helmholtz}). For $i=1,2,3$, we define $\vvardiv=(\vardiv^1,\vardiv^2,\vardiv^3)$ and $\vvarcurl=(\varcurl^1,\varcurl^2,\varcurl^3)$ as follows:
		\begin{align}\label{DE:vardivvarcurl}
			\vardiv^i:=&\p_i\tilde{\phi},&
			\varcurl^i:=&\vvarcurl^i.
		\end{align} 
	\end{definition}}
	Applying the decomposition in Definition \ref{DE:DefPhiPsi} to (\ref{EQ:elasticwaves1.2}), we have:
	\begin{align}\label{EQ:EW1}
		(\p_t^2\vardiv^i-\speedone^2\Delta\vardiv^i)+(\p_t^2\varcurl^i-\speedtwo^2\Delta\varcurl^i)=\sum\limits_{j,k=1}^3G^\prime\left(\p_j\vvardiv^k+\p_j\vvarcurl^k\right)\cdot\left(\p_i\p_j\vvardiv^k+\p_i\p_j\vvarcurl^k\right),
	\end{align}
	where $a_i$ is a constant, and $G^\prime(\p_jU^k):=\frac{\p G}{\p(\p_jU^k)}$ is a smooth scalar function of $\p_jU^k$.
	
		\begin{definition}[Wave operators]\label{DE:waveoperators}
		For the metric $\gfour$ and a smooth function $\varphi$, the standard Laplace-Beltrami operator is defined as 
		\begin{align}
			\boxg\varphi:=\frac{1}{\sqrt{\abs{\det\gfour}}}\p_\alpha\left(\sqrt{\abs{\det\gfour}}\invgfour^{\alpha\beta}\p_\beta\varphi\right).
		\end{align}
		In addition, we define the reduced wave operator with respect to $\gfour$ as
		\begin{align}
			\resboxg\varphi:=\invgfour^{\alpha\beta}\p_\alpha\p_\beta\varphi.
		\end{align}
	\end{definition}
	
	\begin{proposition}[A geometric formulation of the modified elastic wave equations]\label{PR:GeoformulationofModifiedElasticwaves}
		Taking divergence and curl of the equations (\ref{EQ:EW1}), we obtain the following geometric wave equations for $\dive\vardiv$ and $\curl\varcurl$:
		\begin{subequations}\label{EQ:EW2}
			\begin{align}
				\label{EQ:EW2.1}\hat{\square}_{\gfour_{E}(\p\vvardiv,\p\vvarcurl)}(\dive\vardiv)=&P_{(\dive\vvardiv)}(\p^2\vvardiv,\p^3\vvarcurl),\\
				\label{EQ:EW2.2}\square_{\hfour_{E}}(\curl\varcurl)^i=&P^i_{(\curl\vvarcurl)}(\p^2\vvarcurl,\p^2\vvarcurl),
			\end{align}
		\end{subequations}
		{Here, $\gfour_{E}, \hfour_{E}$ are Lorentzian metrics with their inverses $\gfour_{E}^{-1}, \hfour_{E}^{-1}$ defined as follows:}
	\begin{subequations}\label{EQ:EW3}
			\begin{align}\label{EQ:EW3.1}
				(\gfour_{E}^{-1})(\p\vvardiv,\p\vvarcurl):=&-\p_t\otimes\p_t+\sum\limits_{j,k=1}^3\left\{\delta^{jk}\speedone^2+\frac{1}{2}G^\prime\left(\p_j\vardiv^k+\p_j\varcurl^k\right)+\frac{1}{2}G^\prime\left(\p_k\vardiv^j+\p_k\varcurl^j\right)\right\}\p_j\otimes\p_k,\\
				\label{EQ:EW3.2}\hfour_{E}^{-1}:=&-\p_t\otimes\p_t+\speedtwo^2\sum\limits_{a=1}^3\p_a\otimes\p_a.
			\end{align}
		\end{subequations}
		{And also recall the definition of the reduced wave operator:}
		\begin{align*}
\resboxg:=\invgfour^{\alpha\beta}\p_\alpha\p_\beta.
		\end{align*}
		The nonlinear terms on the RHS of (\ref{EQ:EW2.1}) and (\ref{EQ:EW2.2}) are given by:
		\begin{subequations}
			\begin{align}
				\label{EQ:EW3.3}P_{(\dive\vvardiv)}(\p^2\vvardiv,\p^3\vvarcurl)=&\sum\limits_{j,k=1}^3G^\prime\left(\p_j\vardiv^k+\p_j\varcurl^k\right)\cdot\Delta\p_j\varcurl^k+\sum\limits_{i,j,k=1}^3G^{\prime\prime}\left(\p_j\vardiv^k+\p_j\varcurl^k\right)\cdot\left(\p_i\p_j\vardiv^k+\p_i\p_j\varcurl^k\right)^2,\\
				\label{EQ:EW3.4}P^i_{(\curl\vvarcurl)}(\p^2\vvarcurl,\p^2\vvarcurl)=&0.
			\end{align}
		\end{subequations}
	\end{proposition}
	\begin{proof}[Proof of Proposition \ref{PR:GeoformulationofModifiedElasticwaves}]
		Equation (\ref{EQ:EW2.1}) follows directly from taking divergence of (\ref{EQ:EW1}) and using the facts $\p_i\vardiv^j=\p_j\vardiv^i$ and $\sum\limits_{j,k=1}^3G^\prime(\p_j\vardiv^k+\p_j\varcurl^k)\p_j\p_kf=\sum\limits_{j,k=1}^3G^\prime(\p_k\vardiv^j+\p_k\varcurl^j)\p_j\p_kf$ for any function $f$. {Next,}  equation (\ref{EQ:EW2.2}) is a direct result by taking curl of (\ref{EQ:EW1}) and using the fact $\curl(\nabla f)=0$ for any function/tensor field $f$.
	\end{proof}
	\begin{remark}
		We observe that a similar low-regularity local well-posedness result can be derived for system (\ref{EQ:EW2}) with (\ref{EQ:EW3.4}) substituted by $P^i_{(\curl\vvarcurl)}(\p^2\vvarcurl,\p^2\vvarcurl)=\gensmoothfunction(\varcurl)\cdot\p^2\vvarcurl\cdot\p^2\vvarcurl$. However, it is unknown that whether such system can be derived from the elastic wave equations.
	\end{remark}
	\section{PDEs under Study and Relations to the Elastic Wave Equations}\label{S:PDEs}
    {This section is devoted to specifying the geometric wave system to be studied in this paper and demonstrating its relations to the elastic wave equations.\\}

Let $\gensmoothfunction$ denote a generic smooth function that is free to vary from line to line.	Let $\gfour,\hfour$ be symmetric Lorentzian metrics and their inverses, $\gfour^{-1},\hfour^{-1}$, have the following schematic\footnote{See (\ref{EQ:EW3.1}) for the precise definition of $\gfour^{-1}_{E}$, as one candidate for $\gfour^{-1}$, and (\ref{EQ:EW3.3}) for inhomogeneous term.} form:
\begin{subequations}
		\begin{align}\label{DE:Defg}
		\gfour^{-1}:=&-\p_t\otimes\p_t+\speedone^2\sum\limits_{a=1}^3\p_a\otimes\p_a+\sum\limits_{j,k=1}^3\gensmoothfunction(\p\vvardiv,\p\vvarcurl)\p_j\otimes\p_k,\\
\label{DE:Defh}
		\hfour^{-1}:=&-\p_t\otimes\p_t+\speedtwo^2\sum\limits_{a=1}^3\p_a\otimes\p_a+\sum\limits_{j,k=1}^3\gensmoothfunction(\p\vvarcurl)\p_j\otimes\p_k.
	\end{align}
\end{subequations}
	 Let $\speedone>\speedtwo>0$, and $i=1,2,3$. We recall the wave operators from Definition \ref{DE:waveoperators}. In the rest of the article, we study the schematic form of the system:

\begin{align}\label{EQ:EW4}
	\hat{\square}_{\gfour}\vardiv^i+\hat{\square}_{\hfour} \varcurl^i=\gensmoothfunction(\p^2\vvarcurl,\p\vvardiv,\p\vvarcurl),
\end{align}
such that, after taking the divergence and curl of \eqref{EQ:EW4}, we obtain the following multi-wave-speed coupled system of $\dive\vvardiv$ and $\curl\vvarcurl$:\\

\noindent\textbf{Faster-wave system:}
\begin{subequations}
	\begin{align}\label{EQ:boxg}
		\hat{\square}_{\gfour(\p\vvardiv,\p\vvarcurl)}(\dive\vvardiv)=&\gensmoothfunction(\p\vvardiv,\p\vvarcurl)\cdot(\p^2\vvardiv,\p^2\vvarcurl)\cdot(\p^2\vvardiv,\p^2\vvarcurl)+\gensmoothfunction(\p\vvardiv,\p\vvarcurl)\cdot(\p\vvardiv,\p\vvarcurl)\cdot\p^3\vvarcurl=:\tilde{Q}(\p^2\vvardiv,\p^3\vvarcurl),\\
		\curl\vvardiv=&0.
	\end{align}
\end{subequations}

\noindent\textbf{Slower-wave system:}
\begin{subequations}
	\begin{align}\label{EQ:boxh}
		\square_{\hfour(\p\vvarcurl)}(\curl\vvarcurl)^i=&\gensmoothfunction(\p\vvarcurl)\p^2\vvarcurl\cdot\p^2\vvarcurl=:\breve{Q}(\p^2\vvarcurl\cdot\p^2\vvarcurl),\\
		\dive\vvarcurl=&0.
	\end{align}
\end{subequations}
	\color{black}
	In this article, for Greek and Latin indices, we use the following conventions for lowering and raising indices of any vector field or one-form $V$: 
    \begin{itemize}
    \item With respect to the metric $\hfour_{\alpha\beta}$ and its inverse, we lower and raise indices using the notations $(V_\flat)_\beta:=\hfour_{\alpha\beta}V^\alpha$ and $(V^\sharp)^\beta:=(\hfour^{-1})^{\alpha\beta}V_\alpha$. \item Similarly, we lower and raise indices with the metric $\gfour_{\alpha\beta}$ and its inverse, using $V_\beta:=\gfour_{\alpha\beta}V^\alpha$ and $V^\beta:=(\gfour^{-1})^{\alpha\beta}V_\alpha$. 
    \end{itemize}
    These notations apply to all tensor fields as well.
	
	\begin{definition}[Ellipticity of the induced metric $\gt$]\label{DE:Ellipticity}
		For the metrics $\gfour_E,\hfour_{E}$ defined in (\ref{EQ:EW3}), we denote $\gt,h$ to be the induced metrics on $\St$ respectively. In this article, we assume the ellipticity\footnote{Noting that $\speedone>\speedtwo$, therefore, (\ref{AS:ellipticity}) holds when the size of $\p U$ is less than $\speedone-\speedtwo$.} of $\gt-h$ as follows:
		\begin{subequations}\label{AS:ellipticity}
			\begin{align}
				C^{-1}\abs{\xi}^2\leq (\gt^{-1}-h^{-1})^{ij}\xi_i\xi_j\leq C\abs{\xi}^2,\\
				-C\abs{\breve{\xi}}^2\leq (\gt-h)_{ij}\breve{\xi}^i\breve{\xi}^j\leq -C^{-1}\abs{\breve{\xi}}^2.
			\end{align}
		\end{subequations}
	\end{definition}
	We then have the following proposition:
	\begin{proposition}[Equivalence criterion]\label{PR:Equivalence}
		Let $\vvardiv,\vvarcurl$ be the solution of (\ref{EQ:EW4}). Therefore, $\vvardiv$ satisfies (\ref{EQ:boxg}), and $\vvarcurl$ satisfies the solution of (\ref{EQ:boxh}). With $\gfour_{E}$ defined in (\ref{EQ:EW3.1}), and $\hfour_{E}$ defined in (\ref{EQ:EW3.2}), we assume that the following conditions hold:
	\begin{itemize}
			\item $\vardiv^i=\p_i\tilde{\phi}$ for some scalar function $\tilde{\phi}$,
			\item  $\varcurl^i=\curl\tilde{\psi}$ for some $\St$-tangent vector field $\tilde{\psi}$,
			\item $\tilde{Q}(\p^2\vvardiv,\p^3\vvarcurl)=P_{(\dive\vvardiv)}(\p^2\vvardiv,\p^3\vvarcurl), \breve{Q}=P_{(\curl\vvarcurl)}(\p^2\vvardiv,\p^2\vvarcurl)$, and $\gfour=\gfour_{E}$, $\hfour=\hfour_{E}$ coincides with theirs in Proposition \ref{PR:GeoformulationofModifiedElasticwaves}.
		\end{itemize}
		Then, $\vec{U}:=\vvardiv+\vvarcurl$ is the solution to the admissible harmonic elastic wave equations (\ref{EQ:elasticwaves1.2}). In particular, for $\vvarcurl=0$, equation (\ref{EQ:elasticwaves1.2}) belongs to the category of harmonic elastic material; see Proposition \ref{PR:admissible}.
	\end{proposition}
	\begin{definition}[Christoffel symbols]\label{DE:Christoffelsymbols}
		With respect to the metric $\gfour$ (defined in (\ref{EQ:EW3.1})), we define the Christoffel symbols $\Chfour_{\alpha\kappa\lambda}$ and $\Chfour^\beta_{\kappa\lambda}$:
		\begin{subequations}
			\begin{align}
				\Chfour_{\alpha\kappa\lambda}&:=\frac{1}{2}\left(\p_\kappa\gfour_{\alpha\lambda}+\p_\lambda\gfour_{\alpha\kappa}-\p_\alpha\gfour_{\kappa\lambda}\right),\\
				\Chfour^\beta_{\kappa\lambda}&:=\invgfour^{\alpha\beta}\Chfour_{\alpha\kappa\lambda}.
			\end{align}
		\end{subequations}
		Moreover, the contracted Christoffel symbol is defined as below:
		\begin{align}
			\Chfour_{\alpha}:=\Chfour_{\alpha\kappa\lambda}\invgfour^{\kappa\lambda}.
		\end{align}
	\end{definition}
	\begin{definition}\label{DE:vvariablecvvariables}
		We define PDE quantities for analysis on $\vvariables=\{\variables^i\}\in\mathbb{R}^{9}$ and $\cvvariables=\{\cvariables^i\}\in\mathbb{R}^{9}$ with $i=1,2\dots9$, as follows:
		\begin{align}\label{DE:vvariables}
			\vvariables:=&\p\vvardiv,&
			\cvvariables:=&\p\vvarcurl.
		\end{align}
	\end{definition}
	In conclusion, we obtain the following quasilinear wave system which is used to establish the Strichartz estimates:
	\begin{proposition}[Geometric wave equations for Strichartz estimates]\label{PR:analysisPDE}
		Let $\vvardiv$ be the solution of (\ref{EQ:boxg}), and $\vvarcurl$ be the solution of (\ref{EQ:boxh}). Suppose that the conditions in Proposition \ref{PR:Equivalence} are satisfied. Then, $\variables^i$ and $\cvariables$ (defined in (\ref{DE:vvariables})) satisfy the following (schematic\footnote{This schematic form is sufficient for us to derive the desired estimates.}) wave equation:
		\begin{subequations}
			\begin{align}\label{EQ:boxgPHI}
				\hat{\square}_{\gfour(\vvariables,\cvvariables)}(\dive\vvardiv)=&\gensmoothfunction(\vvariables,\cvvariables)\cdot(\p\vvariables,\p\cvvariables)\cdot(\p\vvariables,\p\cvvariables)+\gensmoothfunction(\vvariables,\cvvariables)\cdot\p^2\cvvariables,\\
				\label{EQ:boxhPSI}\hat{\square}_{\hfour}\curl\vvarcurl=&\gensmoothfunction(\cvvariables)\p\cvvariables\cdot\p\cvvariables.
			\end{align}
		\end{subequations}
	\end{proposition}
	\begin{remark}
		For wave equations in Proposition \ref{PR:analysisPDE}, we can apply the machinery on the quasilinear wave equation where the metric components depend on the solution itself. Note that since $\curl\vvardiv=0$, $\dive\vvardiv$ represents the non-trivial part of $\p\vvardiv$, or $\vvariables$. Similarly, $\curl\vvarcurl$ represents the non-trivial part of $\p\vvarcurl$, or $\cvvariables$.
	\end{remark}
	\section{Norms, Littlewood-Paley Projections, Statement of Main Results and Bootstrap Assumptions}\label{S:sectionmainthm}
	In this section, we define various norms and the standard Littlewood-Paley projections used throughout the analysis. Then, we state our main results of the paper along with the bootstrap assumptions.	
	\subsection{Norms}\label{SS:norms}
	In this article, for functions $f,g$ on a normed space $(X,\norm{\cdot}_X)$, we use the notation $\norm{f,g}_X:=\norm{f}_X+\norm{g}_X$. Similarly, for an array of functions $\vec{V}=(V^1,V^2,\dots,V^k)$, we define $\norm{\vec{V}}_X:=\sum\limits_{a=1}^k\norm{V^a}_X$. In particular, we set $\abs{\vec{V}}:=\sqrt{\sum\limits_{i=1}^{k}(V^i)^2}$. For functions $f$ and arrays $\vec{g}$, we also use $\norm{f,g}_X:=\norm{f}_X+\norm{\vec{g}}_X$.
	
	Since the volume form on the constant-time hypersurface $\St$ induced by Minkowski metric $\minkowski$ is $\diff x^1\diff x^2\diff x^3$, by identifying $(t,x^1,x^2,x^3)\in\St$ with $(x^1,x^2,x^3)\in\mathbb{R}$, we define the standard Sobolev norm on $\St$ for $s\in\mathbb{R}$: $\sobolevnorm{F}{s}{\St}:=\norm{\langle\xi\rangle^s\hat{F}(\xi)}_{L_x^2(\St)}$, where $\langle\xi\rangle:=(1+\abs{\xi}^2)^{1/2}$ and $\hat{F}(\xi):=\int_{\mathbb{R}^3}e^{-2\pi ix\cdot\xi}F(x)\diff x$ is the Fourier transform of $F$.

	{With respect to the flat metric on constant-time hypersurface $\St$, for $0<\beta<1$, we denote the standard H\"older semi-norm $\dot{C}_x^{0,\beta}$ and H\"older norm $C_x^{0,\beta}$  of a function $F$ by }
	\begin{subequations}
		\begin{align}
			\semiholdernorm{F}{x}{0}{\beta}{\St}:=&\sup\limits_{x\neq y\in\St}\frac{\abs{F(x)-F(y)}}{\abs{x-y}^\beta},\\
			\holdernorm{F}{x}{0}{\beta}{\St}:=&\sup\limits_{x\in\St}\abs{F(x)}+\sup\limits_{x\neq y\in\St}\frac{\abs{F(x)-F(y)}}{\abs{x-y}^\beta}.
		\end{align}
	\end{subequations}
	
	We also use the following mixed norms for a function $F: \mathbb{R}^3\rightarrow\mathbb{R}$, where $1\leq q_1<\infty$, $1\leq q_2\leq\infty$, and $I$ is an interval of time:
	\begin{subequations}
		\begin{align}
			\twonorms{F}{t}{q_1}{x}{q_2}{I\times\St}&:=\left\{\int_I\onenorm{F}{x}{q_2}{\Sigma_{\tau}}^{q_1}\diff\tau\right\}^{1/q_1},&\twonorms{F}{t}{\infty}{x}{q_2}{I\times\St}&:=\text{ess}\sup\limits_{\tau\in I}\onenorm{F}{x}{q_2}{\Sigma_{\tau}},\\
			\holdertwonorms{F}{t}{q_1}{x}{0}{\beta}{I\times\St}&:=\left\{\int_I\holdernorm{F}{x}{0}{q_2}{\Sigma_{\tau}}^{q_1}\diff\tau\right\}^{1/q_1},&\holdertwonorms{F}{t}{\infty}{x}{0}{q_2}{I\times\St}&:=\text{ess}\sup\limits_{\tau\in I}\holdernorm{F}{x}{0}{q_2}{\Sigma_{\tau}}.
		\end{align}
	\end{subequations}
	
	If $\{F_\upnu\}_{\upnu\in 2^\mathbb{N}}$ is a dyadic-indexed sequence of functions on $\St$, we define
	\begin{align}
		\norm{F_{\upnu}}_{l^2_\upnu L_x^2(\St)}:=\left(\sum\limits_{\upnu\geq 1}\onenorm{F_{\upnu}}{x}{2}{\St}^2\right)^{1/2}.
	\end{align}
	
	\subsection{Littlewood-Paley Projections}\label{SS:definitionlittlewood} We fix a smooth function $\psi=\psi(\abs{\xi}): \mathbb{R}^3\rightarrow [0,1]$ supported on the frequency space annulus $\{\xi\in\mathbb{R}^3|1/2\leq\abs{\xi}\leq2\}$ such that for $\xi\neq0$, we have $\sum\limits_{k\in\mathbb{Z}}\psi(2^k\xi)=1$. For dyadic frequencies $\upnu=2^k$ with $k\in\mathbb{Z}$, we define the standard Littlewood-Paley projection $\littlewood$, which acts on scalar functions $F:\mathbb{R}\rightarrow\mathbb{C}$, as follows:
	\begin{align}
		\littlewood F(x):=\int_{\mathbb{R}^3}e^{2\pi ix\cdot\xi}\psi(\upnu^{-1}\xi)\hat{F}(\xi)\diff\xi,
	\end{align}
	where $\hat{F}(\xi):=\int_{\mathbb{R}^3}e^{-2\pi ix\cdot\xi}F(x)\diff x$ is the Fourier transform of $F$. If $I\subset 2^\mathbb{Z}$ is an interval of dyadic frequencies, then $P_IF:=\sum\limits_{\upnu\in I}\littlewood F$, and $P_{\leq\upnu}F:=P_{[-\infty,\upnu]}F$.
	For functions $f,g$, we use the schematic notation $\littlewood(f,g)$ to denote a linear combination of $\littlewood f$ and $\littlewood g$, namely, $\littlewood f+\littlewood g$.
	
	\subsection{Assumptions on the Initial Data}\label{SS:Data}
	In this section, we precisely state our assumptions on the data.
	\begin{definition}[Regime of hyperbolicity]\label{DE:Hyperbolicity}
		We define $\Domain$ to be
		\begin{align}\label{DE:Hyperbolicity1}
			\Domain:=\left\{(\p\vvardiv,\p\vvarcurl)\in\mathbb{R}^{18}|\ \gt-h \text{ is elliptic}\right\}.
		\end{align}
		See Definition \ref{DE:Ellipticity} for ellipticity of $\gt-h$.
	\end{definition}
	With $N$ as given in Section \ref{SS:ChoiceofParameters} and $D>0$, we assume along $\Sigma_0$ that 
	\begin{subequations}
		\begin{align}
			\textbf{``Divergence-part"}&&\label{AS:divpart}\sum\limits_{k=0}^{2}\sobolevnorm{\pfour^{k}\p\vardiv^i}{N-k-1}{\Sigma_0}&\leq D,\\
			\textbf{``Curl-part"}&&\label{AS:curlpart}\sum\limits_{l=0}^{3}\sobolevnorm{\pfour^{l}\p\varcurl^{i}}{N-l}{\Sigma_0}&\leq D.
		\end{align}
	\end{subequations}
	Assumptions (\ref{AS:divpart}) and (\ref{AS:curlpart}) correspond to regularity assumptions on the ``divergence-part" and ``curl-part" of the data, respectively. 
	
	Let $\Int\,{\mathcal{U}}$ denote the interior of the set $\mathcal{U}$. Assume that there is a compact subset $\breve{\Domain}$ such that 
	\begin{align}
		\left(\p\vvardiv,\p\vvarcurl\right)(\Sigma_{0})\subset\Int\,\breve{\Domain}\subset\breve{\Domain}\subset\Int\,\Domain,
	\end{align}
	where $\Domain$ is defined in (\ref{DE:Hyperbolicity1}).

	\subsection{Statement of Main Theorem}
	Using the notations introduced in the previous sections, we now state the main theorem.
	\begin{theorem}[Main theorem]\label{TH:maintheorem}
		Consider the solution $\vvardiv$ to the wave equations (\ref{EQ:boxg}), and the solution  $\vvarcurl$ to (\ref{EQ:boxh}). {For any real numbers $3<N<7/2$, $D>0$ and $i=1,2,3$, suppose that} the initial data of $\vvardiv,\vvarcurl$ satisfy the following assumptions\footnote{Recall that in Section \ref{SS:Notations}, the notation $\p$ denotes spatial partial derivatives and $\pfour$ denotes the space-time partial derivatives.}:
		\begin{enumerate}
			\item $\sum\limits_{k=0}^{2}\sobolevnorm{\pfour^{k}\p\vardiv^i}{N-k-1}{\Sigma_0}+\sum\limits_{l=0}^{3}\sobolevnorm{\pfour^{l}\p\varcurl^{i}}{N-l}{\Sigma_0}\leq D$.
			\item The image of data functions lies within the interior of $\breve{\Domain}$ (See Definition \ref{DE:Hyperbolicity} for definition of $\breve{\Domain}$). 
		\end{enumerate}
		Then, the time of classical existence $T>0$ for the solution can be bounded from below in terms of $D$. Moreover, the Sobolev regularity of the initial data is propagated by the solution throughout the slab of classical existence. In addition, $\vardiv$ satisfies the following Strichartz estimate:
		\begin{align}
			\twonorms{\pfour\p\vardiv^i}{t}{2}{x}{\infty}{[0,T]\times\mathbb{R}^3}\lesssim1.
		\end{align}
		Furthermore, if the equivalence conditions in Proposition \ref{PR:Equivalence} are verified, then $\vec{U}:=\vvardiv+\vvarcurl$ is the solution to the admissible harmonic elastic wave equations (\ref{EQ:elasticwaves1.2}).
	\end{theorem}
	
    \subsection{Choice of Parameters}\label{SS:ChoiceofParameters}In this section, we introduce several parameters, each of which either measures the regularity or plays a role in our analysis. Let $3<N<7/2$. For the purpose of analysis, we select positive numbers $q,\varepsilon_0,\delta_0,\delta$ and $\delta_1$ that satisfy the following conditions:
	\begin{subequations}\label{DE:ChoiceofParameters}
		\begin{align}
			&2<q<\infty,\\
			&0<\varepsilon_0:=\frac{N-3}{10}<\frac{1}{10},\\
			&\delta_0:=\varepsilon_0^2,\\
			&0<\delta:=\frac{1}{2}-\frac{1}{q}<\varepsilon_0,\\
			&\delta_1:=N-3-4\varepsilon_0-\delta(1-8\varepsilon_0)>8\delta_0>0.
		\end{align}
	\end{subequations}
	More precisely, we consider $N,\varepsilon_0$ and $\delta_0$ to be fixed throughout the paper, while $q,\delta$ and $\delta_1$ will be treated as parameters.

	\subsection{Bootstrap Assumptions}\label{SS:bootstrap}
     In below, we provide our bootstrap assumptions to be used in the proof of Theorem \ref{TH:maintheorem}.\\
     
     Throughout the article, $0<\Tstar\ll1$ denotes a bootstrap time that depends only on the initial data. Let $\Domain$ be defined as in (\ref{DE:Hyperbolicity1}). We assume that $\vvardiv$ and $\vvarcurl$ are smooth\footnote{By smooth, we mean sufficiently smooth for the analysis arguments to hold. Meanwhile, all of our quantitative estimates depend only on the Sobolev and H\"older norms.} 
	solutions to the equation (\ref{EQ:boxg}) and (\ref{EQ:boxh}), and the following estimates hold for $\delta_0$ (defined in Section \ref{SS:ChoiceofParameters}):
	\begin{subequations}\label{BA:ba}
		\begin{align}
    \label{BA:databa}&\left(\p\vvardiv,\p\vvarcurl\right)\subset\Domain,\\
			\label{BA:Div}&\twonorms{\pfour\p\vvardiv}{t}{2}{x}{\infty}{[0,\Tstar]\times\mathbb{R}^3}^2+\sum\limits_{\upnu\geq2}\upnu^{2\delta_0}\twonorms{\littlewood\pfour\p\vvardiv}{t}{2}{x}{\infty}{[0,\Tstar]\times\mathbb{R}^3}^2\leq1.
		\end{align}
	\end{subequations}

    Assuming $\Tstar$ is sufficiently small, in Theorem \ref{TH:MainEstimates}, we derive an improvement of (\ref{BA:Div}). By the fundamental theorem of calculus, equation (\ref{BA:databa}) follows directly from (\ref{BA:Div}).
	
	\section{Preliminary Littlewood-Paley Estimates, Standard Elliptic Estimates, and Commuted Wave Equations}\label{S:preliminary}
	In this section, we provide preliminary estimates and commuted wave equations that we will use throughout the article.
	\begin{proposition}\cite[A.1. Properties of the Littlewood-Paley projection]{pseudodifferentialoperatorsandnonlinearpde}.
		For a function $F$, with $\littlewood$ denoting the Littlewood-Paley projection defined in Section \ref{SS:definitionlittlewood}, the standard results in Littlewood-Paley theory give the following:
		\begin{align}
			\label{EQ:Littlewoodsobolev}\sobolevnorm{F}{s}{\St}&\approx\onenorm{F}{x}{2}{\St}+\left(\sum\limits_{\upnu>1}\upnu^{2s}\onenorm{\littlewood F}{x}{2}{\St}^2\right)^{1/2},\\
			\label{EQ:Littlewoodholder}\holdernorm{F}{x}{0}{s}{\St}&\approx\norm{ F}_{L_x^\infty(\St)}+\sup\limits_{\upnu\geq2}\upnu^s\norm{\littlewood F}_{L_x^\infty(\St)},
		\end{align}
		where $H^s$ is the standard Sobolev norm and $C^{0,s}$ is the standard H\"older norm. One can refer to \cite[Section 1]{AGeometricApproachtotheLittlewoodPaleyTheory} and \cite[A.1]{pseudodifferentialoperatorsandnonlinearpde} for the above results.
		
		In addition, for scalar functions $\phi$ and $\psi$ on $\St$ and a smooth function $\gensmoothfunction$, we have
		\begin{align}\label{holdercal1}
			\holdernorm{\phi\cdot\psi}{x}{0}{\delta_0}{\St}\lesssim&\norm{\phi}_{L_x^\infty(\St)}\semiholdernorm{\psi}{x}{0}{\delta_0}{\St}+\norm{\psi}_{L_x^\infty(\St)}\semiholdernorm{\phi}{x}{0}{\delta_0}{\St}\\
			\notag\lesssim&\holdernorm{\psi}{x}{0}{\delta_0}{\St}\holdernorm{\phi}{x}{0}{\delta_0}{\St}.
		\end{align}
		\begin{align}\label{holdercal2}
			\holdernorm{\gensmoothfunction\circ\phi}{x}{0}{\delta_0}{\St}\lesssim 1+\holdernorm{\phi}{x}{0}{\delta_0}{\St}.
		\end{align}
	\end{proposition}

	\begin{lemma}\cite[Lemma 5.3. Preliminary product estimates]{3DCompressibleEuler}.
		For $3<N<7/2$, let $F,G$ and $\varphi$ be functions on $\St$, $\gensmoothfunction$ be a smooth function of its arguments, and $f^\prime$ be the derivative of $f$ with respect to its arguments. Then, the following frequency-summed product estimates hold for dyadic frequencies $\upnu\geq 1$:
		\begin{subequations}
			\begin{align}
				\label{ES:LPproduct1}\norm{\upnu^{N-3}\littlewood(F\cdot\p G)}_{l^2_\upnu L_x^2(\St)}\lesssim&\onenorm{F}{x}{\infty}{\St}\sobolevnorm{\p G}{N-3}{\St}+\onenorm{G}{x}{\infty}{\St}\sobolevnorm{\p F}{N-3}{\St},\\
				\label{ES:LPproduct2}\norm{\upnu^{N-3}\littlewood(\gensmoothfunction(F)\cdot G)}_{l^2_\upnu L_x^2(\St)}\lesssim&\onenorm{\gensmoothfunction(F)}{x}{\infty}{\St}\sobolevnorm{G}{N-3}{\St}+\sobolevnorm{G}{N-3}{\St}\sobolevnorm{\p F}{1}{\St}.
			\end{align}
		\end{subequations}
	\end{lemma}
	\begin{lemma}\cite[Lemma 2.4. Preliminary commutator estimates]{AGeoMetricApproach}.
		For $3<N<7/2$, with the functions $F,G,\varphi$ defined on $\St$, and $f$ being a smooth function of its arguments, we have the following frequency-summed commutator estimates for dyadic frequencies $\upnu\geq 1$:
		\begin{subequations}
			\begin{align}
				\label{ES:LPcommutator1}\norm{\upnu^{N-3}\left[f\circ\varphi-P_{\leq\upnu}(f\circ\varphi)\right]\cdot\littlewood \p G}_{l^2_\upnu L^2_x(\St)}\lesssim&\onenorm{\p\varphi}{x}{\infty}{\St}\sobolevnorm{G}{N-3}{\St},\\
				\label{ES:LPcommutator2}\norm{\upnu^{N-3}\left\{\littlewood\left[f\circ\varphi\cdot\p^2 G\right]-P_{\leq\upnu}(f\circ\varphi)\cdot\littlewood\p^2 G\right\}}_{l^2_\upnu L^2_x(\St)}\lesssim&\onenorm{\p\varphi}{x}{\infty}{\St}\sobolevnorm{G}{N-2}{\St}\\
				\notag&+\onenorm{ G}{x}{\infty}{\St}\sobolevnorm{\p\varphi}{N-2}{\St}.
			\end{align}
		\end{subequations}	
	\end{lemma}
	
	\begin{lemma}\cite[Lemma 4.5. Standard $L^2$-elliptic estimates]{3DCompressibleEuler}.\label{LE:elliptic}
		Let $\vec{V}$ be a vector field in $\mathbb{R}^3$. Then, by standard elliptic estimates, for\footnote{(\ref{ES:Elliptic0}) holds for all $k\geq0$. Here, we are only concerned with the range $0\leq k\leq5/2$.} $0\leq k\leq5/2$, the following estimate holds:
		\begin{align}\label{ES:Elliptic0}
			\onenorm{\p^{k+1} \vec{V}}{x}{2}{\mathbb{R}^3}\lesssim\onenorm{\p^{k}\dive \vec{V}}{x}{2}{\mathbb{R}^3}+\onenorm{\p^{k}\curl \vec{V}}{x}{2}{\mathbb{R}^3}.
		\end{align}
		In particular, for the curl-free part $\vec{V}=\vvardiv$ (as in Proposition \ref{PR:Equivalence}), there holds
		\begin{align}\label{ES:Elliptic}
			\onenorm{\p^{k+1}\vvardiv}{x}{2}{\mathbb{R}^3}\lesssim\onenorm{\p^{k}\dive\vvardiv}{x}{2}{\mathbb{R}^3}.
		\end{align}
	\end{lemma}
	\begin{proof}[Proof of Lemma \ref{LE:elliptic}]
		Desired estimates follow directly from integrating the following identity on $\mathbb{R}^3$:
		\begin{align}\label{EQ:elliptic}
			\sum\limits_{a,b=1}^3(\p_aV^b)^2=(\dive V)^2+\abs{\curl V}^2+\p_a\left(V^b\p_b V^a\right)-\p_a\left(V^a\dive V\right).
		\end{align}
	\end{proof}
	
	\begin{lemma}\cite[Lemma 8.2. Schauder estimates]{3DCompressibleEuler}.\label{LE:schauder1}
		Let $V$ be a vector field in $\mathbb{R}^3$, and $\delta_1>0$ be the parameter in \eqref{DE:ChoiceofParameters}. Then, the following estimate holds:
		\begin{align}
			\holdernorm{\p V}{x}{0}{\delta_1}{\mathbb{R}^3}\lesssim\holdernorm{\dive V}{x}{0}{\delta_1}{\mathbb{R}^3}+\holdernorm{\curl V}{x}{0}{\delta_1}{\mathbb{R}^3}+\sobolevnorm{V}{2}{\mathbb{R}^3}.
		\end{align}
		In particular, for $V=\vvardiv$ (as in Proposition \ref{PR:Equivalence}), we have:
		\begin{align}\label{ES:Schauder}
			\holdernorm{\p \vvardiv}{x}{0}{\delta_1}{\mathbb{R}^3}\lesssim\holdernorm{\dive \vvardiv}{x}{0}{\delta_1}{\mathbb{R}^3}+\sobolevnorm{\vvardiv}{2}{\mathbb{R}^3}.
		\end{align}
	\end{lemma}
	
	\subsection{Commuted Wave Equations}\label{SS:commutedwaves}
	The following two lemmas provide the commuted version of the equations. Lemma \ref{LE:lemmacommute} addresses the commutation of $\boxg$ with $\pfour$, which are needed for below-top order estimates. In Lemma \ref{LE:lwpequations}, we  commute $\boxg$ with $\littlewood\pfour$, which is required for the top order estimates.
	
	\begin{lemma}[Commuted equations satisfied by one derivative of the solution variables]\label{LE:lemmacommute}
		
		We consider the solutions to the equations in Section \ref{S:PDEs}, specifically, $\vardiv^i$ solves equation (\ref{EQ:boxg}) and $\varcurl^i$ solves equation (\ref{EQ:boxh}) for $i=1,2,3$. Assume that the conditions in Proposition \ref{PR:Equivalence} hold. Then, we obtain the following schematic\footnote{We use the schematic equation as it is sufficient in our analysis.} equations:
		\begin{align}
			\label{EQ:boxg3}\resboxg\pfour\dive\vvardiv=&\gensmoothfunction(\p\vvardiv,\p\vvarcurl)\cdot(\pfour\p\vvardiv,\pfour\p\vvarcurl)\cdot\left(\pfour\p^2\vvardiv,\pfour\p^2\vvarcurl\right)+\gensmoothfunction(\p\vvardiv,\p\vvarcurl)\cdot\pfour\p^3\vvarcurl\\
			\notag&+\gensmoothfunction(\p\vvardiv,\p\vvarcurl)\cdot(\pfour\p\vvardiv,\pfour\p\vvarcurl)\cdot(\pfour\p\vvardiv,\pfour\p\vvarcurl)\cdot(\pfour\p\vvardiv,\pfour\p\vvarcurl)=:Q_{(\pfour^2\vardiv)},
		\end{align}
		where $\gensmoothfunction$ denotes a generic smooth function that may vary from line to line.
	\end{lemma}
	\begin{proof}[Sketch of the proof of Lemma \ref{LE:lemmacommute}]
		By commuting (\ref{EQ:boxg}) with $\pfour$, (\ref{EQ:boxg3}) can be derived through straightforward computations.
	\end{proof}
	The following lemma provides the commuted equations with the Littlewood-Paley projections.
	\begin{lemma}\label{LE:lwpequations}
		For solutions $\vvardiv$ to (\ref{EQ:boxg}) and $\vvarcurl$ to (\ref{EQ:boxh}), provided that the conditions in Proposition \ref{PR:Equivalence} are satisfied, the following equations hold:
		\begin{subequations}
			\begin{align}
				\label{EQ:lwpbox}\boxg\littlewood\dive\vvardiv&=\remainder_{(Q_{(\pfour\vardiv)});\upnu},\\
				\label{EQ:lwpboxg}\boxg\littlewood\pfour\dive\vvardiv&=\remainder_{(Q_{(\pfour^2\vardiv)});\upnu},
			\end{align}
		\end{subequations}
		where $\Chfour_{\alpha}$ is defined in Definition \ref{DE:Christoffelsymbols}, and 
		\begin{subequations}
			\begin{align}
				\label{EQ:remainder0}\remainder_{(Q_{(\pfour\vardiv)});\upnu}=&\littlewood \tilde{Q}-\Chfour^{\alpha}\littlewood\p_\alpha\dive\vvardiv\\\
				\notag&+
				\left[\invgfour^{ab}-P_{\leq\upnu}\invgfour^{ab}\right]\littlewood\p_a\p_{b}\dive\vvardiv\\
				\notag&+\left\{P_{\leq\upnu}\invgfour^{ab}\littlewood\p_a\p_{b}\dive\vvardiv-\littlewood[\invgfour^{ab}\p_a\p_b\dive\vvardiv]\right\},\\
				\label{EQ:remainder1}\remainder_{(Q_{(\pfour^2\vardiv)});\upnu}=&\littlewood Q_{(\pfour^2\vardiv)}-\Chfour^{\alpha}\littlewood\p_\alpha\pfour\dive\vvardiv\\\
				\notag&+
				\left[\invgfour^{ab}-P_{\leq\upnu}\invgfour^{ab}\right]\littlewood\p_a\p_{b}\pfour\dive\vvardiv\\
				\notag&+\left\{P_{\leq\upnu}\invgfour^{ab}\littlewood\p_a\p_{b}\pfour\dive\vvardiv-\littlewood[\invgfour^{ab}\p_a\p_b\pfour\dive\vvardiv]\right\}.
			\end{align}
		\end{subequations}
		
		Moreover, the following estimates are valid for the remainders:
		\begin{subequations}
			\begin{align}
				\label{ES:remainderestimates0}\norm{\upnu^{N-2}\remainder^i_{(Q_{(\pfour\vardiv)});\upnu}}_{l^2_\upnu L_x^2(\St)}\lesssim&\onenorm{\p\vvardiv,\pfour\p\vvardiv,\pfour\p\vvarcurl}{x}{\infty}{\St}\cdot\left(\sum\limits_{k=1}^{2}\sobolevnorm{\pfour^k\p\vvardiv}{N-k-1}{\St}+\sum\limits_{j=2}^{4}\sobolevnorm{\pfour^{j}\vvarcurl}{N+1-j}{\St}\right),\\
				\label{ES:remainderestimates}\norm{\upnu^{N-3}\remainder^i_{(Q_{(\pfour^2\vardiv)});\upnu}}_{l^2_\upnu L_x^2(\St)}\lesssim&\onenorm{\p\vvardiv,\pfour\p\vvardiv,\pfour\p\vvarcurl}{x}{\infty}{\St}\cdot\left(\sum\limits_{k=1}^{2}\sobolevnorm{\pfour^k\p\vvardiv}{N-k-1}{\St}+\sum\limits_{j=2}^{4}\sobolevnorm{\pfour^{j}\vvarcurl}{N+1-j}{\St}\right),
			\end{align}
		\end{subequations}
        where the $l^2_\upnu$-seminorm is taken over dyadic frequencies with $\upnu>1$.
	\end{lemma}
	\begin{proof}[Discussion of the proof of Lemma \ref{LE:lwpequations}]
		{Equations} (\ref{EQ:lwpbox})-(\ref{EQ:lwpboxg}) are derived via straightforward computations by commuting $\littlewood$ with \eqref{EQ:boxg} and (\ref{EQ:boxg3}), using the following connection between the geometric wave equation and reduced wave equation:
		\begin{align}\label{EQ:diffwave}
			\boxg\varphi=\resboxg\varphi-\Chfour^{\alpha}\p_\alpha\varphi,
		\end{align}
		where $\Chfour$ is defined in Definition \ref{DE:Christoffelsymbols}.
		
		Now we prove (\ref{ES:remainderestimates}). We consider (\ref{EQ:remainder1}), first, by (\ref{ES:LPproduct1}), (\ref{EQ:boxg3}), and the fact that $\Chfour=\pfour\p\vvardiv,\pfour\p\vvarcurl$, we have that the following quantity is bounded by the RHS of (\ref{ES:remainderestimates}):
		\begin{align}
			\norm{\upnu^{N-3}\left(\littlewood Q^i_{(\pfour^2\vardiv)}-\Chfour^{\alpha}\p_\alpha\pfour\dive\vvardiv\right)}_{l^2_\upnu L_x^2(\St)}.
		\end{align} 
		Then, applying (\ref{ES:LPcommutator1}) for $\varphi=(\p\vvardiv,\p\vvarcurl)$, $f=(\gfour^{-1})^{ab}$, and $G=\pfour^2 \dive\vvardiv$, we obtain that
		\begin{align}
			\norm{\upnu^{N-3}\left[\invgfour^{ab}-P_{\leq\upnu}\invgfour^{ab}\right]\littlewood\p_a\p_{b}\pfour\dive\vvardiv}_{l^2_\upnu L_x^2(\St)}
		\end{align} is bounded by the RHS of (\ref{ES:remainderestimates}).
		Finally, applying (\ref{ES:LPcommutator2}) for $\varphi=(\p\vvardiv,\p\vvarcurl)$, $f=(\gfour^{-1})^{ab}$ and $G=\pfour\dive \vvardiv$, we can bound the following norm by the RHS of (\ref{ES:remainderestimates}):
		\begin{align}
			\norm{\upnu^{N-3}\left(P_{\leq\upnu}\invgfour^{ab}\littlewood\p_a\p_{b}\pfour\dive\vvardiv-\littlewood[\invgfour^{ab}\p_a\p_b\pfour\dive\vvardiv]\right)}_{l^2_\upnu L_x^2(\St)}.
		\end{align} 
		Combining the above results, we obtain the desired estimate (\ref{ES:remainderestimates}). The proof of (\ref{ES:remainderestimates0}) follows in a similar fashion. We refer readers to \cite[Section 2]{AGeoMetricApproach} and \cite[Section 5]{3DCompressibleEuler} for the detailed proofs.
	\end{proof}
	
	\section{Energy Estimates on Constant-time Hypersurfaces}\label{S:EnergyEstimates}
	In this section, we first present energy estimates and Strichartz estimates for the ``curl-part" of the solution $\vvarcurl$. Next, we derive the energy estimates for the ``divergence-part" of the solution, $\vvardiv$, along constant-time hypersurfaces. 
	
	\subsection{Energy and Strichartz Estimates for the ``Curl-Part" on Constant-time Hypersurfaces}
	In this section, we provide standard result of energy and Strichartz estimates for ``curl-part" of the solution $\vvarcurl$.
	\begin{theorem}[Energy estimates and Strichartz estimates for the ``curl-part"]\label{TH:LinearStrichartz}
		Let $\vvarcurl$ be the solution of the wave equation (\ref{EQ:boxh}) with initial data as specified in Section \ref{SS:Data}. The parameters $\delta_0$ and $\delta$ are defined in (\ref{DE:ChoiceofParameters}), and $\Tstar$ denotes the time of existence/bootstrap time (see Section \ref{SS:bootstrap}). Then, the following estimates hold:
		\begin{subequations}\label{ES:linearcurl}
			\begin{align}
				\label{ES:energyestimates1}&\sum\limits^3_{k=0}\sobolevnorm{\p_t^k\p\vvarcurl}{N-k}{\St}\lesssim D+1,\\
				&\label{ES:linearStri}\twonorms{\pfour^2\p\vvarcurl}{t}{2}{x}{\infty}{[0,\Tstar]\times\mathbb{R}^3}+\sum\limits_{\upnu\geq2}\upnu^{2\delta_0}\twonorms{\littlewood\pfour^2\p\vvarcurl}{t}{2}{x}{\infty}{[0,\Tstar]\times\mathbb{R}^3}^2\lesssim\Tstar^{2\delta},\\
				&\label{ES:LinearSobolev}\onenorm{\p\vvarcurl,\pfour\p\vvarcurl}{}{\infty}{[0,\Tstar]\times\mathbb{R}^3}\lesssim 1.
			\end{align}
		\end{subequations}
		where $D$ is defined in Theorem \ref{TH:maintheorem}.
	\end{theorem}
	\begin{proof}[Discussion of the proof for Theorem \ref{TH:LinearStrichartz}]
	Theorem \ref{TH:LinearStrichartz} follows from \cite[Theorem 1.1]{AGeoMetricApproach}, elliptic estimates Lemma \ref{LE:elliptic}-\ref{LE:schauder1}, and the fact that $\dive\vvarcurl=0$. The details of the proof are similar to, and simpler\footnote{In fact, since the regularity is above the classical level ($H^{\frac{7}{2}+}(\mathbb{R}^3)$), \eqref{ES:energyestimates1} can be obtained by standard energy estimates and Sobolev embedding in $\mathbb{R}^3$. Meanwhile, there is still a tiny derivative loss in the Strichartz estimates due to the rough metric $\hfour$ (below $C^2$); see \cite{Smith98,SmithSogge94}.} than, those of Proposition \ref{PR:EnergyEstimatesforDiv} and Theorem \ref{TH:MainEstimates}, since ``curl-part" is decoupled from the ``divergence-part". Therefore, we omit the details of the proof. We emphasize that the proof of Theorem \ref{TH:LinearStrichartz} is independent of the analysis of the ``divergence-part" due to this decoupling property.
	\end{proof}
\begin{remark}
	For the admissible harmonic elastic wave system \eqref{EQ:elasticwaves1.2}, where the ``curl-part" satisfies the linear wave equation \eqref{EQ:EW2.2}, Theorem \ref{TH:LinearStrichartz} follows from the linear Strichartz estimates and energy estimates. We refer to \cite[Section 1]{Strichartzestimatesforsecondorderhyperbolicoperators} for a more detailed discussion of the linear Strichartz estimates.
\end{remark}
\color{black}
	\subsection{Energy Estimates for ``Divergence-Part" on Constant-time Hypersurfaces}\label{SS:EnergyforDiv}
	The following proposition is the main result of the energy estimates.
	\begin{proposition}[Energy estimates for ``divergence-part"]\label{PR:EnergyEstimatesforDiv}
		Given the initial data and bootstrap assumptions in Section \ref{S:sectionmainthm}, the smooth solutions $\vvardiv$ to the wave equations (\ref{EQ:boxg}) satisfy the following estimates for $3<N<7/2$ and $t\in[0,\Tstar]$:
		\begin{align}\label{ES:energyestimates2}
			\sum\limits^2_{k=0}\sobolevnorm{\pfour^k\p\vvardiv}{N-k-1}{\St}\lesssim&D+1,
		\end{align}
		where $D$ is defined in Theorem \ref{TH:maintheorem}.
	\end{proposition}
	The proof of Proposition \ref{PR:EnergyEstimatesforDiv} will be provided in Section \ref{SSS:proofofEnergy}. {Before the proof, we first derive the standard energy estimates for wave equations.}
	\subsubsection{The basic energy inequality for wave equations}\label{SSS:energymethodforwave}
	In this section, we provide the basic energy inequality for the wave equations.
	\begin{definition}[Energy-momentum tensor, energy current, and deformation tensor]\label{DE:energymomentum} We define the energy-momentum tensor $Q_{\mu\nu}[\varphi]$ associated to a scalar function $\varphi$ in the following way:
		\begin{align}\label{DE:Qmunu}
			Q_{\mu\nu}[\varphi]:=\p_\mu\varphi\p_\nu\varphi-\frac{1}{2}\gfour_{\mu\nu}(\gfour^{-1})^{\alpha\beta}\p_\alpha\varphi\p_\beta\varphi.
		\end{align} 
		Then, given $\varphi$ and any multiplier vector field $\mathbf{X}$, the corresponding energy current $\energycurrentwlot{\mathbf{X}}^{\alpha}[\varphi]$ vector field is defined to be:
		\begin{align}\label{DE:energycurrent}
			\energycurrentwlot{\mathbf{X}}^{\alpha}[\varphi]:=Q^{\alpha\beta}[\varphi]\mathbf{X}_{\beta}-\varphi^2\mathbf{X}^\alpha.
		\end{align}
		Moreover, we define the deformation tensor of the vector field $\mathbf{X}$ as follows:
		\begin{align}\label{DE:deformtensor}
			\deform{\mathbf{X}}_{\alpha\beta}:=\Dfour_{\alpha}\mathbf{X}_\beta+\Dfour_{\beta}\mathbf{X}_\alpha,
		\end{align}
		where $\Dfour$ is the Levi-Civita connection with respect to $\gfour$.
        \end{definition}
        
		For the energy current defined above, there holds the following well-known divergence identity:
		\begin{align}\label{DE:DalphaJalpha}
			\Dfour_\alpha\energycurrentwlot{\mathbf{X}}^{\alpha}[\varphi]=\boxg\varphi(\mathbf{X}\varphi)+\frac{1}{2}Q^{\mu\nu}[\varphi]\deform{\mathbf{X}}_{\mu\nu}-2\varphi\mathbf{X}\varphi-\frac{1}{2}\varphi^2(\gfour^{-1})^{\mu\nu}\deform{\mathbf{X}}_{\mu\nu}.
		\end{align}
		
		Now let $\Timelike$ be the future-directed unit $\gfour$-normal (to $\St$) vector field:
		\begin{align}\label{DE:Timelike}
			\Timelike^\alpha:=\delta^{\alpha}_0,
		\end{align}
		where $\delta$ is the Kronecker delta. Then, we define the energy $\mathbb{E}[\varphi](t)$:
		\begin{align}\label{DE:Energy1}
			\mathbb{E}[\varphi](t):=\int_{\St}\energycurrentwlot{\Timelike}^{\alpha}\Timelike_\alpha\diff\tvol=\int_{\St}\left(Q^{00}[\varphi]+\varphi^2\right)\diff\tvol,
		\end{align}
		where $\diff\tvol$ is the volume form on $\St$ with respect to the induced metric $g$. In the following, we provide the coerciveness lemma for the energy $\mathbb{E}$.
	\begin{lemma}[Coerciveness of $\mathbb{E}$]\label{LE:coeciveness} Under the bootstrap assumptions in Section \ref{SS:bootstrap}, the following estimate holds for $t\in[0,\Tstar]$:
		\begin{align}\label{ES:Coecive}
			\mathbb{E}[\varphi](t)\approx\norm{\varphi}_{H^1(\St)}^2+\norm{\p_t\varphi}_{ L_x^2(\St)}^2.
		\end{align}
	\end{lemma}
	\begin{proof}[{Discussion of the proof of Lemma \ref{LE:coeciveness}}]
		\eqref{ES:Coecive} follows from \eqref{DE:Energy1}, \eqref{DE:Qmunu}, \eqref{DE:Defg} and bootstrap assumption \eqref{BA:databa}.
	\end{proof}
	\begin{lemma}[Basic energy inequality for the wave equations]\label{LE:basicenergy}
		Let $\varphi$ be smooth on $[0,\Tstar]\times\mathbb{R}^3$. Under the bootstrap assumptions in Section \ref{SS:bootstrap}, the following inequality holds for $t\in[0,\Tstar]$: 
		\begin{align}\label{ES:Energy1}
			\mathbb{E}[\varphi](t)\lesssim &\mathbb{E}[\varphi](0)
			+\int_{0}^{t}\norm{\pfour\p\vvardiv,\pfour\p\vvarcurl}_{L_x^\infty(\Sigma_{\tau})}\mathbb{E}[\varphi](\tau)\diff\tau\\
			\notag&+\int_{0}^{t}\norm{\boxg\varphi}_{L_x^2(\Sigma_{\tau})}\norm{\pfour\varphi}_{L_x^2(\Sigma_{\tau})}\diff\tau.
		\end{align}
		Moreover, we have the following estimate:
		\begin{align}\label{ES:Energy2}
			\onenorm{\pfour\varphi}{x}{2}{\St}\lesssim\onenorm{\pfour\varphi}{x}{2}{\Sigma_{0}}+\int_{0}^{t}\onenorm{\boxg\varphi}{x}{2}{\Stau}\diff\tau.
		\end{align}
	\end{lemma}
	
	\begin{proof}[Proof of Lemma \ref{LE:basicenergy}]
		We apply the divergence theorem on the space-time region $[0,t]\times\mathbb{R}^3$, relative to the volume form $\diff\gvol=\sqrt{\det g}\diff x^1\diff x^2\diff x^3\diff t=\diff\tvol\diff t$. Note that $\Timelike$ (defined in (\ref{DE:Timelike})) is the future-directed $\gfour$-unit normal (to $\St$) vector. By (\ref{DE:deformtensor}), (\ref{DE:DalphaJalpha} and (\ref{DE:Energy1}), with $\mathbf{X}:=\Timelike$, there holds
		\begin{align}
			\mathbb{E}[\varphi](t)=&\mathbb{E}[\varphi](0)\\
			\notag&-\int_{0}^{t}\int_{\Sigma_{\tau}}\left(\boxg\varphi(\Timelike\varphi)+\frac{1}{2}Q^{\mu\nu}[\varphi]\deform{\Timelike}_{\mu\nu}-2\varphi\Timelike\varphi-\frac{1}{2}\varphi^2\invgfour^{\mu\nu}\deform{\Timelike}_{\mu\nu}\right)\diff\tvol\diff\tau.
		\end{align}
		Using (\ref{DE:deformtensor}), we have 
		\begin{align}\label{ES:Deformaton}
			\abs{\deform{\Timelike}_{\mu\nu}}\lesssim\abs{\pfour\p\vvardiv,\pfour\p\vvarcurl}.
		\end{align} 
		Combining  (\ref{ES:Deformaton}) with $\abs{\Timelike\varphi}\lesssim\abs{\pfour\varphi}$ and $\abs{Q^{\mu\nu}[\varphi]}\lesssim\abs{\pfour\varphi}^2$, by Cauchy-Schwarz inequality along $\Sigma_{\tau}$ and (\ref{ES:Coecive}), we get the desired estimate (\ref{ES:Energy1}). Using bootstrap assumptions (\ref{BA:Div}), the proof of (\ref{ES:Energy2}) follows in a similar fashion.
	\end{proof}

	\subsubsection{Proof of Proposition \ref{PR:EnergyEstimatesforDiv}}\label{SSS:proofofEnergy}
	\begin{proof}[Proof of Proposition \ref{PR:EnergyEstimatesforDiv}]
		For $N$ defined as in Section \ref{SS:ChoiceofParameters}, we let
		\begin{align}\label{DE:totalenergy}
			P_N(t):=\sum\limits_{k=0}^2\sobolevnorm{\pfour^k\p\vvardiv}{N-k-1}{\St}^2+\sum\limits_{k=0}^3\sobolevnorm{\pfour^k\p\vvarcurl}{N-k}{\St}^2.
		\end{align}
		
		In this proof, we derive integral inequalities for $\sum\limits_{k=0}^2\sobolevnorm{\pfour^k\p\vvardiv}{N-k-1}{\St}^2$ in $P_N(t)$, namely (\ref{PF:proof6.1.2}) and (\ref{PF:proof6.1.7}). We then apply Gr\"onwall's inequality to all the terms in $P_N(t)$ collectively.
		
		The proof of Proposition \ref{PR:EnergyEstimatesforDiv} for $\vvardiv$ combines the vector field multiplier method and Littlewood-Paley theory. More precisely, to derive the energy estimates at the top order, we integrate (\ref{DE:DalphaJalpha}), and applies the divergence theorem, using the energy current $\energycurrentwlot{\Timelike}^{\alpha}[\pfour\dive\vvardiv]:=Q^{\alpha\beta}[\pfour\dive\vvardiv]\Timelike_\beta-\Timelike^\alpha(\pfour\dive\vvardiv)^2$ and  $\energycurrentwlot{\Timelike}^{\alpha}[\littlewood\pfour\dive\vvardiv]:=Q^{\alpha\beta}[\littlewood\pfour\dive\vvardiv]\Timelike_\beta-\Timelike^\alpha(\littlewood\pfour\dive\vvardiv)^2$, on the space-time region bounded by $\Sigma_{0}$ and $\Sigma_{t}$, where $\Timelike$ is defined in (\ref{DE:Timelike}). Then, using Lemma \ref{LE:basicenergy} with $\pfour\dive\vvardiv$ and $\littlewood\pfour\dive\vvardiv$ in the role of $\varphi$, we derive, respectively, the following estimates:
		\begin{subequations}
			\begin{align}
				\mathbb{E}[\pfour\dive\vvardiv](t)\lesssim &\mathbb{E}[\pfour\dive\vvardiv](0)
				+\int_{0}^{t}\norm{\pfour\p\vvardiv,\pfour\p\vvarcurl}_{L_x^\infty(\Sigma_{\tau})}\mathbb{E}[\pfour\dive\vvardiv](\tau)\diff\tau\\
				\notag&+\int_{0}^{t}\norm{\boxg\pfour\dive\vvardiv}_{L^2_x(\Sigma_{\tau})}\left[\mathbb{E}[\pfour\dive\vvardiv](\tau)\right]^{1/2}\diff\tau,\\
				\label{PF:proof6.1.1}\mathbb{E}[\littlewood\pfour\dive\vvardiv](t)\lesssim &\mathbb{E}[\littlewood\pfour\dive\vvardiv](0)
				+\int_{0}^{t}\norm{\pfour\p\vvardiv,\pfour\p\vvarcurl}_{L_x^\infty(\Sigma_{\tau})}\mathbb{E}[\littlewood\pfour\dive\vvardiv](\tau)\diff\tau\\
				\notag&+\int_{0}^{t}\norm{\boxg\littlewood\pfour\dive\vvardiv}_{L^2_x(\Sigma_{\tau})}\left[\mathbb{E}[\littlewood\pfour\dive\vvardiv](\tau)\right]^{1/2}\diff\tau.
			\end{align}
		\end{subequations}
		Then, we substitute $\boxg\pfour\dive\vvardiv$ using equation (\ref{EQ:boxg3}) , and we use equation (\ref{EQ:lwpboxg}) to replace $\boxg\littlewood\pfour\dive\vvardiv$ on the right-hand side of (\ref{PF:proof6.1.1}). Multiplying (\ref{PF:proof6.1.1}) by $\upnu^{2(N-3)}$, summing over $\upnu$, and using (\ref{EQ:Littlewoodsobolev}), estimates (\ref{ES:remainderestimates}), (\ref{ES:Elliptic}), (\ref{ES:linearcurl}) and H\"older's inequality, we obtain:
		\begin{align}\label{PF:proof6.1.2}
			\sobolevnorm{\pfour^2\p\vvardiv}{N-3}{\St}^2\lesssim\sobolevnorm{\pfour^2\dive\vvardiv}{N-3}{\St}^2\lesssim&\sobolevnorm{\pfour^2\p\vvardiv}{N-3}{\Sigma_{0}}^2+\int_{0}^{t}\norm{\p\vvardiv,\pfour\p\vvardiv,\pfour\p\vvarcurl}_{L_x^\infty(\Sigma_{\tau})}\sobolevnorm{\pfour^2\p\vvardiv}{N-3}{\Sigma_{\tau}}^2\diff\tau\\
			\notag\lesssim& P_N(0)+\int_{0}^{t}\norm{\p\vvardiv,\pfour\p\vvardiv,\pfour\p\vvarcurl}_{L_x^\infty(\Sigma_{\tau})}\diff\tau\\
			\notag&+\int_{0}^{t}\left(\norm{\p\vvardiv,\pfour\p\vvardiv,\pfour\p\vvarcurl}_{L_x^\infty(\Sigma_{\tau})}+1\right)P_N(\tau)\diff\tau.
		\end{align} 
		For $\sobolevnorm{\pfour\p\vvardiv}{N-2}{\St}$, first, there holds
		\begin{align}\label{PF:proof6.1.2.4}
			\sobolevnorm{\pfour\p\vvardiv}{N-2}{\St}^2\lesssim\norm{\pfour\p\vvardiv}_{L^2_x(\St)}^2+\sobolevnorm{\pfour^2\p\vvardiv}{N-3}{\St}^2.
		\end{align}
		Then, by the fundamental theorem of calculus in time, the Minkowski integral inequality, and the smallness of $\Tstar$, we deduce:
		\begin{align}\label{PF:proof6.1.2.5}
			\norm{\pfour\p\vvardiv}_{L^2_x(\St)}^2=\int_{\St}\left\{\pfour\p\vvardiv(0,x)+\int_0^t\p_t\pfour\p\vvardiv(\tau,x)\diff\tau\right\}^2\diff x\lesssim P_N(0)+{\int_{0}^{t} P_N(\tau)\diff\tau.}
		\end{align}
		Similarly, we have:
		\begin{align}\label{PF:proof6.1.2.6}
			\sobolevnorm{\p\vvardiv}{N-1}{\St}\lesssim P_N(0)+{\int_{0}^{t} P_N(\tau)\diff\tau.}
		\end{align}
		Therefore, by adding (\ref{PF:proof6.1.2})-(\ref{PF:proof6.1.2.6}), we obtain:
		\begin{align}\label{PF:proof6.1.7}
			\sum\limits_{k=0}^2\sobolevnorm{\pfour^k\p\vvardiv}{N-k-1}{\St}^2\lesssim& P_N(0)+\int_{0}^{t}\norm{\p\vvardiv,\pfour\p\vvardiv,\pfour\p\vvarcurl}_{L^\infty_x(\Sigma_{\tau})}\diff\tau\\
			\notag&+\int_{0}^{t}\left(\norm{\p\vvardiv,\pfour\p\vvardiv,\pfour\p\vvarcurl}_{L^\infty_x(\Sigma_{\tau})}+1\right)P_N(\tau)\diff\tau.
		\end{align}
		
		Then, by (\ref{PF:proof6.1.7}), there holds
		\begin{align}\label{PF:proof6.1.8}
			P_N(t)\lesssim& P_N(0)+\int_{0}^{t}\norm{\p\vvardiv,\pfour\p\vvardiv,\pfour\p\vvarcurl}_{L^\infty_x(\Sigma_{\tau})}\diff\tau\\
			\notag&+\int_{0}^{t}\left(\norm{\p\vvardiv,\pfour\p\vvardiv,\pfour\p\vvarcurl}_{L^\infty_x(\Sigma_{\tau})}+1\right)P_N(\tau)\diff\tau.
		\end{align}
		Via using the bootstrap assumptions (\ref{BA:Div}), energy and Strichartz estimates (\ref{ES:linearcurl}), H\"older inequality in time, Fundamental Theorem of Calculus in time, and the Gr\"onwall's inequality, we conclude the desired result (\ref{ES:energyestimates2}).
	\end{proof}
	\begin{lemma}[$L^2$-elliptic estimates for $\boxg\variables$]\label{LE:ellipticforboxg}
		For $\vvariables$ defined in (\ref{DE:vvariables}), by bootstrap assumptions (\ref{BA:ba}) and energy estimates (\ref{ES:linearcurl}), (\ref{ES:energyestimates2}), $\boxg\variables$ satisfies the following estimate:
		\begin{subequations}
			\begin{align}\label{ES:elliptic1}
				\twonorms{\boxg\variables}{t}{2}{x}{2}{[0,\Tstar]\times\St}\lesssim&\twonorms{\p\vvariables,\p\cvvariables,\p^2\cvvariables}{t}{\infty}{x}{2}{[0,\Tstar]\times\St}\cdot\twonorms{\vvariables,\p\vvariables,\p\cvvariables}{t}{2}{x}{\infty}{[0,\Tstar]\times\St}\lesssim1,\\
				\label{ES:elliptic2}	\twonorms{\boxg\pfour\variables}{t}{2}{x}{2}{[0,\Tstar]\times\St}\lesssim&\twonorms{\pfour\p\vvariables,\pfour\p\cvvariables,\pfour\p^2\cvvariables}{t}{\infty}{x}{2}{[0,\Tstar]\times\St}\cdot\twonorms{\vvariables,\p\vvariables,\p\cvvariables}{t}{2}{x}{\infty}{[0,\Tstar]\times\St}\\
				\notag&+\twonorms{\pfour\vvariables,\pfour\cvvariables}{t}{\infty}{x}{2}{[0,\Tstar]\times\St}^2\cdot\twonorms{\pfour\vvariables,\pfour\cvvariables}{t}{2}{x}{\infty}{[0,\Tstar]\times\St}
				\lesssim1.
			\end{align}
		\end{subequations}
	\end{lemma}
	\begin{proof}[Proof of Lemma \ref{LE:ellipticforboxg}]
		We only prove (\ref{ES:elliptic2}) here, (\ref{ES:elliptic1}) can be proved in a similar fashion. Applying the elliptic estimates (\ref{ES:Elliptic0}) with $\vec{V}:=\boxg\pfour\vvardiv$, we have:
		\begin{align}\label{ES:elliptic3}
			\onenorm{\p\boxg\pfour\vvardiv}{x}{2}{\St}\lesssim&\onenorm{\dive(\boxg\pfour\vvardiv)}{x}{2}{\St}+\onenorm{\curl(\boxg\pfour\vvardiv)}{x}{2}{\St}.
		\end{align}
		Then, by using (\ref{ES:elliptic3}), $\p\gfour=\p(\vvariables,\cvvariables)$, wave equation (\ref{EQ:boxg3}), and the facts that $\vvariables=\p\vvardiv$ and $\curl\vvardiv=0$, we derive:
		\begin{align}\label{ES:elliptic4}
			\onenorm{\boxg\pfour\variables}{x}{2}{\St}\lesssim&\onenorm{\p\boxg\pfour\vvardiv}{x}{2}{\St}+\onenorm{(\p\vvariables,\p\cvvariables)\cdot(\pfour\p\vvariables,\pfour\p\cvvariables)}{x}{2}{\St}\\
			&\notag+\onenorm{(\pfour\vvariables,\pfour\cvvariables)\cdot(\pfour\vvariables,\pfour\cvvariables)\cdot(\pfour\vvariables,\pfour\cvvariables)}{x}{2}{\St}\\
			\notag\lesssim&\onenorm{(\p\vvariables,\p\cvvariables)\cdot(\pfour\p\vvariables,\pfour\p\cvvariables)}{x}{2}{\St}+\onenorm{\vvariables\cdot\pfour\p^2\cvvariables}{x}{2}{\St}\\
			&\notag+\onenorm{(\pfour\vvariables,\pfour\cvvariables)\cdot(\pfour\vvariables,\pfour\cvvariables)\cdot(\pfour\vvariables,\pfour\cvvariables)}{x}{2}{\St}\\
			\notag\lesssim&\onenorm{\vvariables,\p\vvariables,\p\cvvariables}{x}{\infty}{\St}\cdot\onenorm{\pfour\p\vvariables,\pfour\p\cvvariables,\pfour\p^2\cvvariables}{x}{2}{\St}\\
			\notag&+\onenorm{\pfour\vvariables,\pfour\cvvariables}{x}{\infty}{\St}\cdot\onenorm{\pfour\vvariables,\pfour\cvvariables}{x}{2}{\St}^2.
		\end{align}
		Desired estimate (\ref{ES:elliptic2}) follows then by taking $L^2_t$ norm of (\ref{ES:elliptic4}).
	\end{proof}
	\section{Reduction of Strichartz Estimates and the Rescaled Solution}\label{S:sectionreduction}
	In this section, we state our main estimates as Theorem \ref{TH:MainEstimates}, which provide an improvement of the bootstrap assumption (\ref{BA:Div}). Following that, we then present a series of analytic reductions, from the Strichartz estimates of Theorem \ref{TH:MainEstimates} to the decay estimates of Theorem \ref{TH:Spatiallylocalizeddecay}. This section is presented concisely, since the full proofs of these reductions are lengthy and technically complicated, though standard. We refer readers to \cite{3DCompressibleEuler,AGeoMetricApproach,ImprovedLocalwellPosedness} for the detailed proofs.
	\begin{theorem}[Improvement of the Strichartz-type bootstrap assumptions for the wave variables]\label{TH:MainEstimates}
		Let $\vvariables$ be defined as in (\ref{DE:vvariables}), and let $\vvardiv$ be the solution of (\ref{EQ:boxgPHI}). Given a sufficiently small $\delta>0$, as defined in Section \ref{SS:ChoiceofParameters}, then under the initial data and bootstrap assumptions of Section \ref{S:sectionmainthm}, the following estimates hold for $8\delta_0<\delta_1<N-3$ (see (\ref{DE:ChoiceofParameters})):
		\begin{align}
			\label{ES:estimateswavevariables}\twonorms{\pfour\dive\vvardiv}{t}{2}{x}{\infty}{[0,\Tstar]\times\mathbb{R}^3}^2+\sum\limits_{\upnu\geq2}\upnu^{2\delta_1}\twonorms{\littlewood\pfour\dive\vvardiv}{t}{2}{x}{\infty}{[0,\Tstar]\times\mathbb{R}^3}^2\lesssim\Tstar^{2\delta}.
		\end{align}
	\end{theorem}
	
	Recall the bootstrap assumption (\ref{BA:Div}):
	\begin{equation*}
		\twonorms{\pfour\p\vvardiv}{t}{2}{x}{\infty}{[0,\Tstar]\times\mathbb{R}^3}^2+\sum\limits_{\upnu\geq2}\upnu^{2\delta_0}\twonorms{\littlewood\pfour\p\vvardiv}{t}{2}{x}{\infty}{[0,\Tstar]\times\mathbb{R}^3}^2\leq1.
	\end{equation*}
	Note that for $\Tstar$ small\footnote{$\Tstar$ is in fact small as we are focusing on the local well-posedness of the solution.} enough, (\ref{ES:estimateswavevariables}) is a strict improvement of (\ref{BA:Div}). To see this, first note that by (\ref{EQ:Littlewoodholder}) and the choice of parameters that $\delta_1>8\delta_0$ (see Section \ref{SS:ChoiceofParameters}), we have:
	\begin{align}\label{ES:improvement1}
		\sum\limits_{\upnu\geq2}\upnu^{\delta_0}\onenorm{\littlewood\pfour\p\vvardiv}{x}{\infty}{\mathbb{R}^3}\lesssim\sum\limits_{\upnu\geq2}\upnu^{-\delta_1+\delta_0}\holdernorm{\pfour\p\vvardiv}{x}{0}{\delta_1}{\mathbb{R}^3}\lesssim\holdernorm{\pfour\p\vvardiv}{x}{0}{\delta_1}{\mathbb{R}^3}.
	\end{align}
	Then, by (\ref{ES:Schauder}) and (\ref{ES:improvement1}), we obtain
	\begin{align}\label{ES:improvement2}
		\holdernorm{\pfour\p\vvardiv}{x}{0}{\delta_1}{\mathbb{R}^3}\lesssim\holdernorm{\pfour\dive\vvardiv}{x}{0}{\delta_1}{\mathbb{R}^3}+\sobolevnorm{\pfour\vvardiv}{2}{\mathbb{R}^3}.
	\end{align}
	Finally, by (\ref{ES:improvement1}), (\ref{ES:improvement2}), energy estimates (\ref{ES:energyestimates2}), and (\ref{EQ:Littlewoodholder}), for a small $\Tstar$, (\ref{ES:estimateswavevariables}) is indeed a strict improvement of (\ref{BA:Div}).
	
	\subsection{Partitioning of the Bootstrap Time Interval}\label{SS:partitioning}
{In the following, we first reduce the proof of Theorem \ref{TH:MainEstimates} to the proof of Strichartz estimates on small time intervals.}\\

	Let $\vvariables$ be defined as in (\ref{DE:vvariables}), and let $\lambda$ be a fixed large number, $0<\varepsilon_0<\frac{N-3}{5}$ be a fixed number as mentioned in Section \ref{SS:ChoiceofParameters}. By the bootstrap assumptions (\ref{BA:Div}), we can\footnote{The existence of such a partition follows easily from the bootstrap assumptions, see \cite[Remark 1.3]{ImprovedLocalwellPosedness}.} partition $[0,\Tstar]$ into a disjoint union of sub-intervals $I_k:=[t_{k-1},t_k]$ with total number $\lesssim\lambda^{8\varepsilon_0}$. Moreover, for each $I_k$, we have the properties that $\abs{I_k}\leq\lambda^{-8\varepsilon_0}\Tstar$ and 
	\begin{align}
		\twonorms{\pfour\p\vvardiv}{t}{2}{x}{\infty}{I_k\times\mathbb{R}^3}^2+\sum\limits_{\upnu\geq2}\upnu^{2\delta_1}\twonorms{\littlewood\pfour\p\vvardiv}{t}{2}{x}{\infty}{I_k\times\mathbb{R}^3}^2&\lesssim\lambda^{-8\varepsilon_0}.
	\end{align}
	Here, $\delta_1$ satisfies that $8\delta_0<\delta_1<N-3$ (see Section \ref{SS:ChoiceofParameters} for the precise definition of $\delta_1$).
	Now we reduce Theorem \ref{TH:MainEstimates} to a frequency localized estimate.
	\begin{theorem}[Frequency localized Strichartz estimate]\label{TH:FrequencyStrichartz}
		Let $\varphi$ be a solution of
		\begin{align}
			\boxg\varphi=0
		\end{align}
		on the time interval $I_k$. Then, for any $q>2$ sufficiently close to $2$ and any $\tau\in[t_{k},t_{k+1}]$, under the bootstrap assumptions, there holds the following estimate:
		\begin{align}\label{ES:Stri2}
			\twonorms{P_\lambda\pfour\varphi}{t}{q}{x}{\infty}{[\tau,t_{k+1}]\times\St}\lesssim\lambda^{\frac{3}{2}-\frac{1}{q}}\norm{\pfour\varphi}_{L_x^2(\Stau)}.
		\end{align}
	\end{theorem}
	
	\subsection{ Reducing the proof of Theorem \ref{TH:MainEstimates} to Theorem \ref{TH:FrequencyStrichartz}}
	The proof that Theorem \ref{TH:MainEstimates} follows from Theorem \ref{TH:FrequencyStrichartz} is similar to the proof of \cite[Section 3]{AGeoMetricApproach}. The reason is that the reduction relies only on: \textbf{1)} Duhamel's principle; \textbf{2)} our top-order energy estimates (\ref{ES:energyestimates2}); \textbf{3)} Littlewood-Paley estimates for the inhomogeneous terms in a frequency-projected version of the wave equations.
	
	Let $W(t,s)(f_0,f_1)$ denote the solution $\varphi(t)$ of $\boxg\varphi=0$ with initial conditions $\varphi(s)=f_0,\p_t\varphi(s)=f_1$. Define $\tilde{P}_\lambda:=\sum\limits_{\frac{1}{2}\leq\frac{\upnu}{\lambda}\leq 2}\littlewood$, so we have $P_\lambda=\tilde{P}_\lambda P_{\lambda}$. Then, by using (\ref{EQ:lwpbox}) and Duhamel's principle, we have:
	\begin{align}\label{EQ:Duhamel1}
		P_{\lambda}\dive\vvardiv(t)=W(t,t_k)[P_{\lambda}\dive\vvardiv(t_k),P_{\lambda}\p_t\dive\vvardiv(t_k)]+\int_{t_k}^{t}W(t,s)[0,\remainder_{(Q_{(\pfour\vardiv)});\lambda}(s)]\diff s.
	\end{align}
	Differentiating (\ref{EQ:Duhamel1}) and applying the Littlewood-Paley projection $\tilde{P}_\lambda$, we get:
	\begin{align}\label{EQ:Duhamel2}
		P_{\lambda}\pfour\dive\vvardiv(t)=&\tilde{P}_\lambda\pfour W(t,t_k)[P_{\lambda}\dive\vvardiv(t_k),P_{\lambda}\p_t\dive\vvardiv(t_k)]+\int_{t_k}^{t}\tilde{P}_\lambda\pfour W(t,s)[0,\remainder_{(Q_{(\pfour\vardiv)});\lambda}(s)]\diff s.
	\end{align}
	Consider the first term on the RHS of (\ref{EQ:Duhamel2}), employing Theorem \ref{TH:FrequencyStrichartz}, we obtain:
	\begin{align}\label{ES:Duhamel1}
		&\twonorms{\tilde{P}_\lambda\pfour W(t,t_k)[P_{\lambda}\dive\vvardiv(t_k),P_{\lambda}\p_t\dive\vvardiv(t_k)]}{t}{2}{x}{\infty}{[t_k,t]\times\Stau}\\
		\notag\lesssim&\abs{I_k}^{\frac{1}{2}-\frac{1}{q}}\twonorms{\tilde{P}_\lambda\pfour W(t,t_k)[P_{\lambda}\dive\vvardiv(t_k),P_{\lambda}\p_t\dive\vvardiv(t_k)]}{t}{q}{x}{\infty}{[t_k,t]\times\Stau}\\
		\notag\lesssim&\lambda^{\frac{3}{2}-\frac{1}{q}}\abs{I_k}^{\frac{1}{2}-\frac{1}{q}}\onenorm{P_\lambda\pfour \dive\vvardiv}{x}{2}{\Sigma_{t_k}}.
	\end{align}
	Let $\delta=\frac{1}{2}-\frac{1}{q}$. By using the fact that $\abs{I_k}\lesssim \lambda^{-8\varepsilon_0}\Tstar$ (see Section \ref{SS:partitioning}) and the energy estimate (\ref{ES:Energy2}), we derive:
	\begin{align}\label{ES:Duhamel2}
		&\twonorms{\tilde{P}_\lambda\pfour W(t,t_k)[P_{\lambda}\dive\vvardiv(t_k),P_{\lambda}\p_t\dive\vvardiv(t_k)]}{t}{2}{x}{\infty}{[t_k,t]\times\Stau}\\
		\notag\lesssim&\lambda^{1+(1-8\varepsilon_0)\delta}\Tstar^\delta\left(\onenorm{\pfour P_{\lambda}\dive\vvardiv}{x}{2}{\Sigma_{0}}+\twonorms{\remainder_{(Q_{(\pfour\vardiv)});\lambda}}{t}{1}{x}{2}{[0,\Tstar]\times\Stau}\right).
	\end{align}
	For the second term on the RHS of (\ref{EQ:Duhamel2}), via using the Minkowski inequality, Theorem \ref{TH:FrequencyStrichartz}, the energy estimate (\ref{ES:Energy2}), and the fact that $W(s,s)(0,\remainder_{(Q_{(\pfour\vardiv)});\upnu}(s))=0$, we deduce:
	\begin{align}\label{ES:Duhamel3}
		&\twonorms{\int_{t_k}^{t}\tilde{P}_\lambda\pfour W(t,s)[0,\remainder_{(Q_{(\pfour\vardiv)});\lambda}(s)]\diff s}{t}{2}{x}{\infty}{[t_k,t_{k+1}]\times\Stau}\\
		\notag\lesssim&\int_{t_k}^{t_{k+1}}\twonorms{\tilde{P}_\lambda\pfour W(t,s)[0,\remainder_{(Q_{(\pfour\vardiv)});\lambda}(s)]}{t}{2}{x}{\infty}{[s,t_{k+1}]\times\Stau} \diff s\\
		\notag\lesssim&\lambda^{1+(1-8\varepsilon_0)\delta}\Tstar^\delta\twonorms{\remainder_{(Q_{(\pfour\vardiv)});\lambda}}{t}{1}{x}{2}{[t_k,t_{k+1}]\times\Stau}.
	\end{align} 
	Together with (\ref{EQ:Duhamel2}), (\ref{ES:Duhamel2}), this implies
	\begin{align}\label{ES:Duhamel4}
		\twonorms{P_{\lambda}\pfour\dive\vvardiv}{t}{2}{x}{\infty}{[t_k,t_{k+1}]\times\Stau}\lesssim&\lambda^{1+(1-8\varepsilon_0)\delta}\Tstar^\delta\left(\onenorm{\pfour P_{\lambda}\dive\vvardiv}{x}{2}{\Sigma_{0}}+\twonorms{\remainder_{(Q_{(\pfour\vardiv)});\lambda}}{t}{1}{x}{2}{[0,\Tstar]\times\Stau}\right).
	\end{align}
	Then, we sum over $I_k$ (with total number $\lesssim\lambda^{8\varepsilon_0}$). Recalling (\ref{DE:ChoiceofParameters}) for $N$, we thus obtain:
	\begin{align}\label{ES:Duhamel5}
		\lambda^{2\delta_1}\twonorms{P_{\lambda}\pfour\dive\vvardiv}{t}{2}{x}{\infty}{[0,\Tstar]\times\St}^2\lesssim&\lambda^{2\delta_1}\lambda^{8\varepsilon_0}\lambda^{2+2(1-8\varepsilon_0)\delta}\Tstar^{2\delta}\left(\onenorm{\pfour P_{\lambda}\dive\vvardiv}{x}{2}{\Sigma_{0}}^2+\twonorms{\remainder_{(Q_{(\pfour\vardiv)});\lambda}}{t}{1}{x}{2}{[0,\Tstar]\times\St}^2\right)\\
		\notag\lesssim&\Tstar^{2\delta}\left(\onenorm{\lambda^{N-2}\pfour P_{\lambda}\dive\vvardiv}{x}{2}{\Sigma_{0}}^2+\twonorms{\lambda^{N-2}\remainder_{(Q_{(\pfour\vardiv)});\lambda}}{t}{1}{x}{2}{[0,\Tstar]\times\St}^2\right).
	\end{align} 
	Summing over large frequencies, and employing (\ref{EQ:Littlewoodsobolev}), H\"older's inequality, remainder estimates (\ref{ES:remainderestimates0}), and energy estimates (\ref{ES:energyestimates2}), we have:
	\begin{align}\label{ES:Duhamel6}
		&\sum\limits_{\upnu\geq\lambda}\upnu^{2\delta_1}\twonorms{\littlewood\pfour\dive\vvardiv}{t}{2}{x}{\infty}{[0,\Tstar]\times\mathbb{R}^3}^2\\
		\notag\lesssim&\sum\limits_{\upnu\geq\lambda}\left\{\Tstar^{2\delta}\left(\onenorm{\upnu^{N-2}\pfour \littlewood\dive\vvardiv}{x}{2}{\Sigma_{0}}^2+\twonorms{\upnu^{N-2}\remainder_{(Q_{(\pfour\vardiv)});\upnu}}{t}{1}{x}{2}{[0,\Tstar]\times\St}^2\right)\right\}\\
		\notag\lesssim&\Tstar^{2\delta}.
	\end{align}
	The desired estimate (\ref{ES:estimateswavevariables}) follows from (\ref{ES:Duhamel6}) and the following low-frequency estimate:
	\begin{align}
		\twonorms{\littlewood\pfour\dive\vvardiv}{t}{2}{x}{\infty}{[0,\Tstar]\times\mathbb{R}^3}\lesssim&\Tstar^{\frac{1}{2}}\onenorm{\pfour\dive\vvardiv}{x}{2}{\St}\lesssim\Tstar^{\frac{1}{2}},& 2\leq\upnu\leq&\lambda,
	\end{align}
{	which can be derived by using the Bernstein inequality and the energy estimates.}
	\subsection{Rescaled Quantities and Rescaled Elastic Wave Equations}\label{SS:SelfRescaling}
	In this section, to facilitate further reductions, we perform the following coordinate transformation $(t,x)\mapsto(\lambda(t-t_k),\lambda x)$. Define
	\begin{align}\label{DE:Trescale}
		\Trescale:=\lambda(t_{k+1}-t_k).
	\end{align}
	Noting that, by construction, it holds that
	\begin{align}\label{ES:Trescale}
		0\leq\Trescale\leq\lambda|I_k|\leq\lambda^{1-8\varepsilon_0}\Tstar.
	\end{align}
	\begin{definition}[Rescaled quantities]\label{DE:rescaledquantities}
		First, we define the following rescaled variables:
		\begin{subequations}
			\begin{align}
				\variables^i_{(\lambda)}:=&\variables^i(t_k+\lambda^{-1}t,\lambda^{-1}x),\\
				\cvariables^i_{(\lambda)}:=&\cvariables^i(t_k+\lambda^{-1}t,\lambda^{-1}x),\\
				\vardiv^i_{(\lambda)}:=&\vardiv^i(t_k+\lambda^{-1}t,\lambda^{-1}x),\\
				\varcurl^i_{(\lambda)}:=&\varcurl^i(t_k+\lambda^{-1}t,\lambda^{-1}x),\\
				(\dive\vvardiv)_{(\lambda)}:=&(\dive\vvardiv)(t_k+\lambda^{-1}t,\lambda^{-1}x).
			\end{align}
		\end{subequations}
		Then, we accordingly define the rescaled tensor fields:
		\begin{subequations}
			\begin{align}
				\label{frequencyrescaledg}(\gfour_{(\lambda)})_{\alpha\beta}(t,x):=&\gfour_{\alpha\beta}(\vvariables_{(\lambda)}(t,x),\cvvariables_{(\lambda)}(t,x)),\\
				\label{frequencyrescaledgt}(g_{(\lambda)})_{\alpha\beta}(t,x):=&g_{\alpha\beta}(\vvariables_{(\lambda)}(t,x),\cvvariables_{(\lambda)}(t,x)).
			\end{align}
		\end{subequations}
	\end{definition}
	\begin{remark}
		Note that $(\dive\vvardiv)_{(\lambda)}=\lambda\dive\left(\vvardiv_{(\lambda)}\right)$.
	\end{remark}
	\begin{remark}\label{newsigma0}
		After rescaling, the new initial constant-time hypersurface $\Sigma_{0}$ corresponds to the constant-time hypersurface previously denoted as $\Sigma_{t_k}$ for some $k$ in Section \ref{S:sectionintro}-\ref{S:EnergyEstimates}.
	\end{remark}
	The following proposition provides the equations satisfied by the rescaled quantities. The proof is straightforward and therefore omitted.
	\begin{proposition}[The rescaled geometric wave-transport formulation of the elastic wave equations] Consider a solution to (\ref{EQ:EW4}) and (\ref{EQ:boxg}). Let $\{\vvardiv_{(\lambda)},\vvariables_{(\lambda)},\cvvariables_{(\lambda)}\}$ be the rescaled quantities from Definition \ref{DE:rescaledquantities}. Then, the rescaled quantities satisfy the following equations:
		\begin{subequations}
			\begin{align}
				\label{EQ:EW5}\square_{\gfour_{(\lambda)}}\vvariables_{(\lambda)}=&\gensmoothfunction(\vvariables_{(\lambda)},\cvvariables_{(\lambda)})\cdot\p\vvariables_{(\lambda)}\cdot(\p\vvariables_{(\lambda)},\p\cvvariables_{(\lambda)})+(\gensmoothfunction(\vvariables_{(\lambda)},\cvvariables_{(\lambda)})+1)\cdot\pfour^2\cvvariables_{(\lambda)},\\
				\label{EQ:boxg4}\square_{\gfour_{(\lambda)}}(\dive\vvardiv)_{(\lambda)}=&\gensmoothfunction(\vvariables_{(\lambda)},\cvvariables_{(\lambda)})\cdot\p\vvariables_{(\lambda)}\cdot(\p\vvariables_{(\lambda)},\p\cvvariables_{(\lambda)})+\gensmoothfunction(\vvariables_{(\lambda)},\cvvariables_{(\lambda)})\cdot\p^2\cvvariables_{(\lambda)}.
			\end{align}
		\end{subequations}
	\end{proposition}
	\begin{remark}
		For notational convenience,  in the remainder of the article, we drop the sub and super scripts $(\lambda)$, except for the rescaled time $\Trescale$.
	\end{remark}
	
	\subsection{Consequences of the Bootstrap Assumptions and Strichartz Estimates}
	After rescaling in Section \ref{SS:SelfRescaling}, assuming bootstrap assumptions (\ref{BA:Div}), standard computations  based on Littlewood-Paley calculus yield the following consequences of the bootstrap assumptions for $\tau\in[t_k,t_{k+1}]$: 
	\\
	
	\noindent\textbf{Estimates by using bootstrap assumptions of ``divergence-part"}
	\begin{subequations}\label{BA:bt1}
		\begin{align}
			\twonorms{\pfour\vvariables}{t}{2}{x}{\infty}{[t_k,t_{k+1}]\times\Stau}+\lambda^{\delta_0}\sqrt{\sum\limits_{\upnu>2}\upnu^{2\delta_0}\twonorms{\littlewood\pfour\vvariables}{t}{2}{x}{\infty}{[t_k,t_{k+1}]\times\Stau}^2}\lesssim&\lambda^{-1/2-4\varepsilon_0},\\
			\holdertwonorms{\pfour\vvariables}{t}{2}{x}{0}{\delta_0}{[t_k,t_{k+1}]\times\Stau}\lesssim&\lambda^{-1/2-4\varepsilon_0}
			.\end{align}
	\end{subequations}
	
	Using the Strichartz estimates (\ref{ES:linearcurl}), we have the following rescaled version of estimates for the ``curl-part":\\
	
	\noindent\textbf{Estimates by using Strichartz estimates for ``curl-part"}
	\begin{subequations}
		\begin{align}
			\label{BA:bt2}\twonorms{\pfour\cvvariables}{t}{2}{x}{\infty}{[t_k,t_{k+1}]\times\Stau}+\lambda^{\delta_0}\sqrt{\sum\limits_{\upnu>2}\upnu^{2\delta_0}\twonorms{\littlewood\pfour\cvvariables}{t}{2}{x}{\infty}{[t_k,t_{k+1}]\times\Stau}^2}\lesssim&\lambda^{-1/2-4\varepsilon_0},\\
			\label{BA:bt2.1}\holdertwonorms{\pfour\cvvariables}{t}{2}{x}{0}{\delta_0}{[t_k,t_{k+1}]\times\Stau}\lesssim&\lambda^{-1/2-4\varepsilon_0},\\
			\label{BA:bt3}\twonorms{\pfour^2\cvvariables}{t}{2}{x}{\infty}{[t_k,t_{k+1}]\times\Stau}+\lambda^{\delta_0}\sqrt{\sum\limits_{\upnu>2}\upnu^{2\delta_0}\twonorms{\littlewood\pfour^2\cvvariables}{t}{2}{x}{\infty}{[t_k,t_{k+1}]\times\Stau}^2}\lesssim&\lambda^{-3/2-4\varepsilon_0},\\
			\label{BA:bt3.1}\holdertwonorms{\pfour^2\cvvariables}{t}{2}{x}{0}{\delta_0}{[t_k,t_{k+1}]\times\Stau}\lesssim&\lambda^{-3/2-4\varepsilon_0}
			.\end{align}
	\end{subequations}
	\subsection{Further Reduction of the Strichartz Estimates}
	By the rescaling in Section \ref{SS:SelfRescaling} and direct computation, to prove Theorem \ref{TH:FrequencyStrichartz}, it is equivalent to establishing the following Strichartz estimate on $[0,\Trescale]$, with respect to the LP projection onto the frequency domain $\{1/2\leq\abs{\xi}\leq 2\}$.
	\begin{theorem}\label{TH:RescaledStrichartz}
		Under the bootstrap assumptions, for any solution $\varphi$ of $\boxg\varphi=0$ on the slab $[0,\Trescale]\times\mathbb{R}^3$, the following estimate holds:
		\begin{align}
			\twonorms{P_1\pfour\varphi}{t}{q}{x}{\infty}{[0,\Trescale]\times\mathbb{R}^3}\lesssim\norm{\pfour\varphi}_{L_x^2(\Sigma_{0})},
		\end{align}
		where $q>2$ is sufficiently close to 2, and $\gfour$ here represents the rescaled metric $\gfour_{(\lambda)}$ defined in (\ref{frequencyrescaledg}).
	\end{theorem}
	The proof of Theorem \ref{TH:RescaledStrichartz} crucially relies on the following decay estimate. 
	\begin{theorem}[Decay estimate]\label{TH:DecayEstimates}
		There exists a sufficiently large number $\Lambda$, such that for any $\lambda\geq\Lambda$ and any solution $\varphi$ of the equation $\boxg\varphi=0$ on $[0,\Trescale]\times\mathbb{R}^3$, there is a function $d(t)$ satisfying
		\begin{align}
			\norm{d}_{L_t^{\frac{q}{2}}([0,\Trescale])}\lesssim1,
		\end{align}
		such that for $q>2$ sufficiently close to 2, and for any $t\in[0,\Trescale]$, there holds the following decay estimate:
		\begin{align}
			\norm{P_1\p_t\varphi}_{L_x^\infty(\St)}\lesssim\left(\frac{1}{(1+t)^{\frac{2}{q}}}+d(t)\right)\left(\sum\limits_{m=0}^3\norm{\p^m\varphi}_{L_x^1(\Sigma_{0})}+\sum\limits_{m=0}^2\norm{\p^m\p_t\varphi}_{L_x^1(\Sigma_{0})}\right).
		\end{align}
		
	\end{theorem}

	The proof of Theorem \ref{TH:RescaledStrichartz} using Theorem \ref{TH:DecayEstimates} is based on a $\mathcal{T}\mathcal{T}^*$ argument\footnote{This argument comes from functional analysis, and it does not rely on the specific structure of the elastic wave equations.} (see \cite[Section 8.6]{ImprovedLocalwellPosedness}, \cite[Section 8.30]{KR3}, and \cite[Appendix B]{AGeoMetricApproach}). The argument in the proof goes
 through for the elastic wave equations, {since under the bootstrap assumptions (\ref{BA:Div}),} we share the same estimate as in \cite{KR3,ImprovedLocalwellPosedness,AGeoMetricApproach}:
	\begin{align}
		\twonorms{\pfour\gfour}{t}{2}{x}{\infty}{[0,\Trescale]\times\mathbb{R}^3}\lesssim\lambda^{-1/2-4\varepsilon_0}.
	\end{align}
	
	Theorem \ref{TH:DecayEstimates} can be further reduced to the following spatially localized version.
	\begin{theorem}[Spatially localized version of decay estimate]\label{TH:Spatiallylocalizeddecay}
		Let $R$ be a fixed radius\footnote{The radius $R$ will be used in the following sections as well with the same definition. The existence of such an $R$ is guaranteed by properties of $\gfour$ (\ref{DE:Defg}), Strichartz estimates (\ref{ES:linearStri}), bootstrap assumptions (\ref{BA:Div}), which ensure $g$ is comparable to the Euclidean metric on $\St$ under the bootstrap assumptions.} such that
		\begin{align}
			B_R(p)&\subset B_{\frac{1}{2}}(p,g_{(\lambda)}),&
			&\text{for any}\ p\in\Sigma_t, \, 0\leq t\leq\Trescale,
		\end{align} 
		where $B_{\rho}(p,g_{(\lambda)})$ is the geodesic ball\footnote{The notation $B_{\rho}(p,g_{(\lambda)})$ will be used in the remainder of the article. This is consistent with the notation that is used in \cite{3DCompressibleEuler} and \cite{AGeoMetricApproach}} centered at $p$ with radius $\rho$, and $g_{(\lambda)}$ is the rescaled induced metric of $\gfour$ on $\St$ (defined in (\ref{frequencyrescaledgt})). {Then, there exist
		a function $d(t)$ satisfying}
		\begin{align}
			\norm{d}_{L_t^{\frac{q}{2}}([0,\Trescale])}\lesssim1,
		\end{align}	
		{and a large number $\Lambda$ such that for any $\lambda\geq\Lambda$,}  
		there is a $q>2$ sufficiently close to 2, ensuring that the following estimate holds for any $t\in[0,\Trescale]$:
		\begin{align}\label{Decay}
			\norm{P_1\p_t\varphi}_{L_x^\infty(\St)}\lesssim\left(\frac{1}{(1+\abs{t-1})^{\frac{2}{q}}}+d(t)\right)\left(\norm{\pfour\varphi}_{L_x^2(\Sigma_1)}+\norm{\varphi}_{L_x^2(\Sigma_1)}\right).
		\end{align}
    {Here, $\varphi$ is a solution of the equation $\boxg\varphi=0$ on $[0,\Trescale]\times\mathbb{R}^3$, and  $\varphi(1,x)$ is supported in the Euclidean ball $B_{R}$.}
	\end{theorem}
	
	The proof of Theorem \ref{TH:DecayEstimates} using Theorem \ref{TH:Spatiallylocalizeddecay} can be done via an approach involving the Bernstein inequalities of LP projection, the partition of unity\footnote{We take a sequence of Euclidean balls $\{B_I\}$ of radius $R$ such that their union covers $\mathbb{R}^3$. Each ball $B_R$ is centered at $\upgamma_\tip(1)$ (defined in  Section \ref{SSS:solutionofu}) for some $\tip$.} of $\varphi$, and Sobolev embedding $W^{2,1}\hookrightarrow L^2$. We 
	refer readers to \cite[Section 8.5]{ImprovedLocalwellPosedness} for the detailed proof.\\
	
	To summarize, in this section we have reduced Theorem \ref{TH:MainEstimates} to Theorem \ref{TH:Spatiallylocalizeddecay}. Before proceeding, we need to introduce the geometric setup. We will discuss the proof of Theorem \ref{TH:Spatiallylocalizeddecay} in Section \ref{SSS:reductiondecay}.
	
	\section{Geometric Foliations and Conformal Energy}\label{S:sectiongeometrysetup}
	In this section, we first construct the geometry, which is intricately coupled to the elastic wave equations via the metric $\gfour$ and plays a crucial role in our analysis. In Section \ref{SS:Sectionconformalenergy}, building on the constructed geometry, we introduce the conformal energy, and derive estimates for it. This is a fundamental ingredient in the vector field method that we use to derive Strichartz estimates.
	
	\subsection{Construction of the Optical Function}\label{SS:opticalfunction}
	
	The goal of this section is to construct the geometry based on a solution $u$ to the following eikonal equation\footnote{For more details of the geometric construction of $u$, we refer readers to \cite[Chapter 9 and Chapter 14]{GlobalStabiliityOfMinkowski}.}:
	\begin{align}\label{EQ:eikonalequation}
		(\gfour^{-1})^{\alpha\beta}\p_\alpha u\p_\beta u=0.
	\end{align}
	
	Many previous results address the construction of the geometry based on (\ref{EQ:eikonalequation}) for the quasilinear wave equations, the Einstein equations, the compressible Euler equations, and the relativistic Euler equations; see \cite{3DCompressibleEuler,AGeoMetricApproach,ImprovedLocalwellPosedness,yu2022rough,Roughsolutionsofthe3-DcompressibleEulerequations}. We adopt the same approach for the elastic wave equations. Note that, by direct computation, $\Timelike:=\p_t$ is geodesic.
	\subsubsection{Point $\tip$ and integral curve $\upgamma_\tip(t)$}\label{SSS:pointtip}
	Recall the center $p$ of the Euclidean ball $B_R$ on $\Sigma_1$ in Theorem \ref{TH:Spatiallylocalizeddecay}. 
	Let $\tip\in\Sigma_{0}$ be the point\footnote{Note that in the original spacetime $[0,\Tstar]\times\mathbb{R}^3$, $\tip$ lies in $\Sigma_{t_k}$ for some $k$.} in the rescaled space-time $[0,\Trescale]\times\mathbb{R}^3$, where $\Trescale$ is defined in (\ref{DE:Trescale}), such that the integral curve $\upgamma_\tip(t)$ of the future-directed vector field $\Timelike$ emanating from the point $\tip$ satisfies 
	\begin{align}
		\upgamma_\tip(1)=p.
	\end{align}
	We refer to $\upgamma_\tip(t)$ as the cone-tip axis. Notice that the estimates, constants and parameters in Section \ref{S:sectiongeometrysetup}-Section \ref{S:connectioncoefficientsandpdes} are independent of $\tip$.
	\subsubsection{Parameterize a family of null vectors $\Lunit_{\sangle}$ on $\tip$}\label{SSS:familyofnullvectors}
	Let $w_*=\frac{4}{5}\Trescale$. Using the arguments in \cite{szeftel2012parametrix1}, we can guarantee\footnote{The existence of $w$-foliation for $w\in[0,\varepsilon]$ with a small $\varepsilon>0$ can be proven by the Nash-Moser implicit function theorem, and such foliation can be extended to $w_*$ by an argument of continuity (see \cite{szeftel2012parametrix1}).} that there is a neighborhood $\neighborhood\in\Sigma_{0}$ contained in the geodesic ball $B_{\Trescale}(\tip,g_{(\lambda)})$, where $g_{(\lambda)}$ is the frequency rescaled induced metric of $\gfour$ on $\Sigma_{0}$ defined in (\ref{frequencyrescaledgt}). This neighborhood $\neighborhood$ can be foliated by the level sets $S_{0,-w}$ of a positive function $w$, taking all values in $[0,w_*]$ with $w(\tip)=0$. Each $S_{0,-w}$ for positive $w$ is diffeomorphic to $\mathbb{S}^2$. Fix a diffeomorphism $\sangle\rightarrow\Phi_\sangle(w_*)$ from $\Stwo$ to $S_{0,-w_*}$. Then, for each point $\textbf{p}=\Phi_\sangle(w_*)$, we denote the lapse
	\begin{align}
		\lapse:=&\left((g^{-1})^{cd}\p_cw\p_dw\right)^{-\frac{1}{2}},&\lapse(\tip)=&1.
	\end{align} 
	We note that $\lapse\approx1$; see Proposition \ref{PR:InitialCT1} (we note that the proof of Proposition \ref{PR:InitialCT1} is independent of the construction of $u$). Then, there exists a unique integral curve\footnote{The existence and uniqueness of such integral curve are guaranteed by the estimates on $\Sigma_{0}$ in Proposition \ref{PR:InitialCT1}. We refer readers to \cite[Section 9.4.2]{3DCompressibleEuler} for details.} 
	$w\rightarrow\Phi_\sangle(w)$ of the vector field $\lapse^2(g^{-1})^{ic}\p_cw$ in $\Sigma_{0}$ with $\Phi_\sangle(w_*)=\textbf{p}$, and such that this integral curve can be extended to $\tip$, i.e. $\Phi_{\sangle}(0)=\tip$ (the extendibility of $\Phi_\sangle$ to $\tip$ follows from estimates (\ref{ES:pNatSigma0}) and the fundamental theorem of Calculus). Denoting $\dot{\Phi}_\sangle:=\frac{\diff}{\diff w}\Phi_{\sangle}$, we then define
	\begin{align}
		\spherenormal_\sangle|_\tip:=&\dot{\Phi}_\sangle(0),&
		\left.\Lunit_\sangle\right|_\tip:=&\spherenormal_\sangle|_\tip+\Timelike|_\tip.
	\end{align}
	Note that the diffeomorphism $\sangle\rightarrow\left.\Lunit_\sangle\right|_\tip$ is uniquely determined by the vector field $\lapse^2(g^{-1})^{ic}\p_cw$.
	\subsubsection{The interior and exterior solution $u$}\label{SSS:solutionofu}
	\paragraph{The interior solution $u$.}
	We let $\{\left.\Lunit_\sangle\right|_\tip\}_{\sangle\in\Stwo}$ be the family of future-directed null vectors (parameterized by $\Stwo$) in the tangent space $T_{\tip}\region$ such that $\Lunit_\sangle^0=1$, as defined in Section \ref{SSS:familyofnullvectors}. To propagate $\Lunit_\sangle$ along the cone-tip axis $\upgamma_\tip(t)$, for any $\textbf{p}\in\upgamma_\tip(t)$ (as defined in Section \ref{SSS:pointtip}) and $\sangle\in\Stwo$, we define $\left.\Lunit_\sangle\right|_{\textbf{p}}$ by solving the following parallel transport equation:
	\begin{align}
		\Dfour_\Timelike\Lunit_{\sangle}=0,
	\end{align} with initial conditions $\left.\Lunit_\sangle\right|_\tip$. Note that for any $\textbf{p}\in\upgamma_\tip(t)$, we have $\langle\left.\Lunit_\sangle\right|_{\textbf{p}},\Timelike\rangle=-1$ since $\Timelike$ is geodesic. Then, for each $u\in[0,\Trescale]$ and $\sangle\in\Stwo$, there exists a unique null geodesic $\Upsilon_{u,\sangle}(t)$, where $t\in[u,\Trescale]$, emanating from $\textbf{p}=\upgamma_\tip(u)$ with $\left.\frac{\diff}{\diff t}\Upsilon_{u,\sangle}\right|_{t=u}=\Lunit_\sangle$ and $\Upsilon_{u,\sangle}^0(t)=t$. Specifically, $\Upsilon_{u,\sangle}(t)$ is constructed by solving the following ``geodesic" ODE system\footnote{See \cite[Section 9.4.1]{3DCompressibleEuler} for the detailed explanation of this ODE system, where the correct analog of $\Timelike$ is denoted by $\materialderivative$ there.}:
	\begin{subequations}
		\begin{align}
			\label{EQ:geodesiceq}\ddot{\Upsilon}_{u,\sangle}^\alpha(t)=&-\left.\Chfour^\alpha_{\kappa\lambda}\right|_{\Upsilon_{u,\sangle}(t)}\dot{\Upsilon}_{u,\sangle}^\kappa(t)\dot{\Upsilon}_{u,\sangle}^\lambda(t)\\
			\notag&+\frac{1}{2}\left.[\Lie_\Timelike\gfour]_{\kappa\lambda}\right|_{\Upsilon_{u,\sangle}(t)}(\dot{\Upsilon}_{u,\sangle}^\kappa(t)-\left.\Timelike^\kappa\right|_{\Upsilon_{u,\sangle}(t)})(\dot{\Upsilon}_{u,\sangle}^\lambda(t)-\left.\Timelike^\lambda\right|_{\Upsilon_{u,\sangle}(t)})\dot{\Upsilon}_{u,\sangle}^\alpha(t),\\
			\Upsilon_{u,\sangle}^\alpha(u)=&\upgamma_\tip^\alpha(u), \ \ \ 
			\dot{\Upsilon}_{u,\sangle}^\alpha(u)=\Lunit^\alpha_{\sangle},
		\end{align} 
	\end{subequations}
	where $\ddot{\Upsilon}_{u,\sangle}^\alpha:=\frac{\diff^2}{\diff t^2}\Upsilon_{u,\sangle}^\alpha$ and $\dot{\Upsilon}_{u,\sangle}^\alpha:=\frac{\diff}{\diff t}\Upsilon_{u,\sangle}^\alpha$, $\Chfour$ denote the Christoffel symbol of $\gfour$, and $\Lie_\Timelike\gfour$ is the Lie derivative of $\gfour$ with respect to $\Timelike$. The curve $t\rightarrow\Upsilon_{u,\sangle}(t)$ is a non-affinely parameterized null geodesic, with the vector field $\Lunit_{u,\sangle}^\alpha:=\frac{\diff}{\diff t}\Upsilon_{u,\sangle}^\alpha$ being null and normalized by $\Lunit_{u,\sangle}^0=1$. In fact, this vector field coincides (in the interior region) with the vector field $\Lunit$ defined in (\ref{DE:DefinitionofL}) below. The truncated null cone $\coneu^{(\text{int})}$ is defined as
	\begin{align}\label{DE:coneuint}
		\coneu^{(\text{int})}:=\bigcup\limits_{\sangle\in\Stwo,t\in[u,\Trescale]}\Upsilon_{u,\sangle}(t).
	\end{align} 
	Now we define the optical function $u$ by asserting that its level sets $\{u=u^\prime\}$ are $\mathcal{C}^{(\text{int})}_{u^\prime}$. For $u\in[0,\Trescale]$ and $t\in[u,\Trescale]$, we define $\stu$ as follows: 
	\begin{align}
		\stu:=\coneu^{(\text{int})}\bigcap\St.
	\end{align}
	For $u\neq t$, $\stu$ is a smooth surface diffeomorphic\footnote{$\stu$ being diffeomorphic to $\Stwo$ is a direct result of Proposition \ref{PR:mainproof}} to $\Stwo$. Then, we set the interior region $\intregion$ to be
	\begin{align}
		\intregion:=\bigcup\limits_{t\in[0,\Trescale], 0\leq u\leq t}\stu.
	\end{align}
	For each $\sangle\in\Stwo$ and $u\in[0,\Trescale]$, we define the angular coordinate functions $\{\sangle^A\}_{A=1,2}$ to be constant along each fixed null geodesics $\Upsilon_{u,\sangle}(t)$. These coordinates coincide with the standard angular coordinates on $\Stwo$ at the tip $\textbf{p}$, which corresponds to $t=u$.
	
	\paragraph{The exterior solution $u$.}
	
	Now we extend the foliation of space-time by null hypersurfaces to a neighborhood of $\bigcup\limits_{t\in[0,\Trescale], 0\leq u\leq t}\stu$ in $[0,\Trescale]\times\mathbb{R}^3$.
	Recall $w_*$ and $\Phi_\sangle(w)$ in Section \ref{SSS:familyofnullvectors}. For each $w\in(0,w_*]$ and any $\textbf{p}\in S_{0,-w}$, there exists a $\sangle\in\Stwo$ such that $\textbf{p}=\Phi_\sangle(w)$. The outward unit normal (in $\Sigma_{0}$) to $S_{0,-w}$ at $\textbf{p}$ is $\spherenormal_{w,\sangle}:=\dot{\Phi}_\sangle(w)$. We set
	\begin{align}
		\Lunit_{w,\sangle}:=\spherenormal_{w,\sangle}+\Timelike|_{\Phi_\sangle(w)},
	\end{align}
	which is a null vector in $T_\textbf{p}\region$. Then, with $u=-w$, there is a unique (non-affinely parameterized) null geodesic $\Upsilon_{u,\sangle}$ emanating from $\textbf{p}$ and solving (\ref{EQ:geodesiceq}) with the initial condition $\left.\frac{\diff}{\diff t}\Upsilon_{u,\sangle}\right|_{t=0}=\Lunit_{-u,\sangle}$ and $\Upsilon_{u,\sangle}(0)=\textbf{p}$. For $-w_*\leq u<0$, we define the null cone $\coneu^{(\text{ext})}$ as follows:
	\begin{align}\label{DE:coneuext}
		\coneu^{(\text{ext})}:=\bigcup\limits_{\sangle\in\Stwo,t\in[0,\Trescale]}\Upsilon_{u,\sangle}(t).
	\end{align}
	The optical function $u$ is then given by asserting that its level sets $\{u=u^\prime\}$ are $\mathcal{C}^{(\text{ext})}_{u^\prime}$. For $t\in[0,\Trescale]$ and $u\in[-w_*,0)$, we define $\stu$ to be
	\begin{align}
		\stu:=\coneu^{(\text{ext})}\bigcap\St.
	\end{align} 
	Note that $\coneu^{(\text{ext})}$'s are the outgoing truncated null cones, specifically, $\coneu^{(\text{ext})}=\bigcup\limits_{t\in[0,\Trescale]}\stu$. The exterior region $\extregion$ is defined as below:
	\begin{align}
		\extregion:=\bigcup\limits_{t\in[0,\Trescale], u\in[-w_*,0)}\stu.
	\end{align}
	{Note that on $\Sigma_{0}\setminus\{\tip\}$, by the construction in Section \ref{SSS:familyofnullvectors}, the angular coordinate functions $\{\sangle^A\}_{A=1,2}$ are constant along the integral curves of the vector field $\lapse^2(g^{-1})^{ic}\p_cw\p_i$. For each $\sangle\in\Stwo$ and $u\in[-w_*,0)$, we similarly define the angular coordinate functions $\{\sangle^A\}_{A=1,2}$ to be constant along null geodesics $\Upsilon_{u,\sangle}(t)$. The constant values coincide with the constructed angular coordinates $\{\sangle^A\}_{A=1,2}$ on $\Sigma_{0}$.} 
	
	Now let the space-time region $\region$ be 
	\begin{align}
		\region=\intregion\bigcup\extregion.
	\end{align}
	
	For convenience, we denote 
	\begin{align}\label{DE:coneu}
		\coneu:=
		\begin{cases}
			\coneu^{(\text{ext})},&  u\in[-w_*,0),\\
			\coneu^{(\text{int})},&  u\in[0,\Trescale].
		\end{cases}
	\end{align}
	
	By the constructions above, we have established the geometric coordinates $(t,u,\sangle^A)$ in $\region$.
	
	\begin{definition}[Norms of $\stu$-tangent tensor fields]\label{DE:norms}
		Recall we define various norms in Section \ref{SS:norms}. Now we introduce norms on $2$-sphere. Let $\esphere=\esphere(\sangle)$ be the canonical metric on $\Stwo$, and let $\{\asphere\}_{A=1,2}$ denote local angular coordinates on $\Stwo$. For $p\in[1,\infty)$, we define the following Lebesgue norms $L^p_{\sangle}$ and $L^p_{\gsphere}$ for $\stu$-tangent tensor fields $\xi$:
		\begin{align}
			\onenorm{\xi}{\sangle}{p}{\stu}:=&\left(\int_{\sangle\in\Stwo}\sgabs{\xi}^p\diff\flatspherevol\right)^{1/p},&
			\onenorm{\xi}{\gsphere}{p}{\stu}:=&\left(\int_{\sangle\in\Stwo}\sgabs{\xi}^p\diff\spherevol\right)^{1/p},\\
			\onenorm{\xi}{\sangle}{\infty}{\stu}:=&\mathrm{ess}\sup\limits_{\sangle\in\Stwo}\sgabs{\xi}.
		\end{align}
		We use the same notation as in Section \ref{SS:norms} for mixed-norms, for example:
		\begin{align}
			\twonorms{\xi}{t}{q}{\sangle}{p}{\coneu}=\left\{\int_{\max\{u,0\}}^{\Trescale}\onenorm{\xi}{\sangle}{p}{\stu}^q\diff t\right\}^{1/q}.
		\end{align}
	\end{definition}
	
	\subsection{Geometric Quantities}\label{SS:geometricquantities}
	
	\begin{definition}[The radial variable]\label{DE:rescaledquantity}
		Recall that
		\begin{align}
			\label{ES:lengtht}0\leq\Trescale\leq\lambda^{1-8\varepsilon_0}\Tstar.
		\end{align}
		We define the geometric radial variable $\rgeo$ as follows:
		\begin{align}
			\label{DE:rgeo}\rgeo:=t-u.
		\end{align}
		
		{Recall $w_*=\frac{4}{5}\Trescale$. Since $t\in[0,\Trescale]$ and $u\in[-w_*,t]$ in $\region$} we have
		\begin{align}\label{ES:estimatesradial}
			0&\leq\rgeo<2\Trescale=2\lambda^{1-8\varepsilon_0}\Tstar,&
			-\frac{4}{5}\lambda^{1-8\varepsilon_0}T_*&\leq u\leq\lambda^{1-8\varepsilon_0}\Tstar.
		\end{align}
		We will silently use estimates (\ref{ES:estimatesradial}) throughout the paper.
	\end{definition}
	\begin{definition}[vector fields and scalar functions corresponding to $\gfour$]\label{DE:Nullframe}
		
		We define the null vector field
		\begin{align}\label{DE:Lgeo}
			\Lgeo&:=-(\gfour^{-1})^{\alpha\beta}\partial_\beta u\partial_\alpha.
		\end{align}	
		Note that by (\ref{EQ:eikonalequation}) and straightforward computations, we have $\Dfour_{\Lgeo}\Lgeo=0$.	
		
		Recall the induced metric $\gt$ on $\St$ defined in Definition \ref{DE:Ellipticity}. We define the null lapse $\nulllapse$ as follows:
		\begin{align}
			\nulllapse &:=(\sqrt{(\gt^{-1})^{ij}\partial_i u\partial_j u})^{-1}.
		\end{align}
		
		Then, the vector field $\spherenormal$ is defined in the following manner:
		\begin{align}\label{DE:spherenormal}
			\spherenormal&:=-\nulllapse(\gt^{-1})^{ij}\partial_iu\partial_j.
		\end{align}
		Note that $\left.\spherenormal\right|_{\Sigma_{0}}=\spherenormal_\sangle$, where $\spherenormal_\sangle$ is defined in Section \ref{SS:opticalfunction}.
		
		Next, we define the principal null vector fields $\Lunit$ and $\uLunit$ to be
		\begin{align}\label{DE:DefinitionofL}
			\Lunit&:=\Timelike+\spherenormal,
			&
			\uLunit&:=\Timelike-\spherenormal.
		\end{align}	
		Notice that by (\ref{EQ:eikonalequation}) and (\ref{DE:spherenormal}), we obtain the identity $\Timelike u=\abs{\nabla u}_g=\frac{1}{\nulllapse}$. Then, we have
		\begin{align}\label{DE:Lunit}
			\Lunit=\nulllapse\Lgeo.
		\end{align}
		
		Moreover, the following basic properties hold:
		\begin{align}
			\gfour(\Timelike,\Timelike)&=-1,
			&
			\gfour(\spherenormal,\spherenormal)&=1,\\
			\gfour(\Timelike,\spherenormal)&=0,
			&
			\gfour(\Lunit,\Lunit)&=0,\\
			\gfour(\Lunit,\uLunit)&=-2,
			&
			\gfour(\uLunit,\uLunit)&=0.
		\end{align}
		
		We define $\gsphere$ to be the inverse of $\gsphere^{-1}$, with $\gsphere^{-1}$ taking the following form:
		\begin{align}\label{DE:decompositionofmetric}
			(\gsphere^{-1})^{\alpha\beta}:=(\gfour^{-1})^{\alpha\beta}+\frac{1}{2}\Lunit^\alpha\uLunit^\beta+\frac{1}{2}\uLunit^\alpha\Lunit^\beta.
		\end{align}
		It is straightforward to verify that $\gsphere$ is the induced metric of $\gfour$ on $\stu$. We fix a pair of locally defined, unit orthogonal spherical vector fields\footnotemark
		\footnotetext{In the rest of the paper, we automatically sum over $A$ if there are two $A$ appearing in the expression.} on $\stu$, denoted by $\{e_A\}_{A=1,2}$, such that $(\gsphere^{-1})^{\alpha\beta}=\sum\limits_{A=1,2}e_A^\alpha e_A^\beta$. We refer to, $\Lunit,\uLunit,e_1,e_2$, as a null frame for the geometry. \\
		
		We denote the Levi-Civita connection on $\stu$ with respect to $\gsphere$ by $\angnabla$.\\
		
		As discussed in Section \ref{SS:opticalfunction}, the angular coordinate functions $\{\sangle^A\}_{A=1,2}$ satisfy the equation $\Lunit(\sangle^A)=0$ along null cones. Consequently, with respect to the geometric coordinates $(t,u,\sangle^A)$, we have
		\begin{align}\label{DE:framework}
			\Lunit&=\frac{\p}{\p t},&  \spherenormal&=-\nulllapse^{-1}\frac{\p}{\p u}+Y^A\deriasphere.
		\end{align}
		By the construction in Section \ref{SSS:solutionofu}, it holds that $Y^A=0$ on $\Sigma_{0}$. 
	\end{definition}
	
	\begin{lemma}[$\hfour$-spacelike property of $\coneu$]\label{LE:Lspacelike}
		{The vector field $\Lunit$ defined in (\ref{DE:Lunit}) is spacelike with respect to $\hfour$,} that is:
		\begin{align} \label{L spacelike}
			\hfour(\Lunit,\Lunit):=H>0.
		\end{align}
		Furthermore, the $\gfour$-null cone $\coneu$ is $\hfour$-spacelike.
	\end{lemma}
	\begin{proof}[Proof of Lemma \ref{LE:Lspacelike}]
		By the definition of $\gfour$ in (\ref{DE:Defg}) and $\hfour$ in (\ref{DE:Defh}), we have:
		\begin{align}\label{EQ:hLL}
			\hfour(\Lunit,\Lunit)=\gfour(\Lunit,\Lunit)-\left(\gt_{ij}-h_{ij}\right)\Lunit^i\Lunit^j=-\left(\gt_{ij}-h_{ij}\right)\Lunit^i\Lunit^j.
		\end{align}
		The desired result \eqref{L spacelike} follows from (\ref{EQ:hLL}) and the ellipticity assumption (\ref{AS:ellipticity}).

		Next, in order to prove that $\coneu$ is spacelike with respect to $\hfour$, it suffices to prove that the normal vector $^{(\hfour)}\nabla u$ is $\hfour$-timelike, where $^{(\hfour)}\nabla u$ denotes the $\hfour$-gradient of $u$. By using definitions of $\gfour$ in (\ref{DE:Defg}) and $\hfour$ in (\ref{DE:Defh}), along with the eikonal equation \eqref{EQ:eikonalequation}, we derive:
		\begin{align}\label{EQ:h spacelike}
			\hfour(^{(\hfour)}\nabla u,^{(\hfour)}\nabla u)&=-(\partial_t u)^2+(h^{-1})^{ij}\partial_i u\partial_j u\\
			\nonumber=&\left(h^{-1}-\gt^{-1}\right)^{ij}\partial_i u\partial_j u.
		\end{align}
		The desired result follows from (\ref{EQ:h spacelike}) and the ellipticity assumption (\ref{AS:ellipticity}).
	\end{proof}
	\begin{definition}[$\hfour$-unit normal vector field of $\coneu$ and a decomposition of $\hfour$ and $\Timelike$]\label{DE:Conenormal}
		Notice that $\coneu$ is $\hfour$-spacelike. We define $\Conenormal$ to be the (timelike) future-directed $\hfour$-unit normal vector field to $\coneu$:
		\begin{align}
			\hfour(\Conenormal,\Lunit)=&0,&
			\hfour(\Conenormal,e_A)=&0,&
			\hfour(\Conenormal,\Conenormal)=&-1.
		\end{align}
		By such construction, we have:
		\begin{align}
			\hfour(\Conenormal,\Timelike)=&-\Conenormal^0,\\\label{unit V}
			(\Conenormal^0)^2=&\sum\limits_{i=1}^3\speedtwo^{-2}(\Conenormal^i)^2+1.
		\end{align}
			\end{definition}
			\begin{remark} \label{V formula}
				Noting that $\coneu$ denotes the level set of $u$, we indeed work with $\Conenormal=-\frac{^{(\hfour)}\nabla u}{\sqrt{-\langle^{(\hfour)}\nabla u,^{(\hfour)}\nabla u\rangle_\hfour}}$.
			\end{remark}
			\begin{definition}[$\stu$-tangent tensor fields]\label{D:stutensor fields}
				
				We define the $\gfour$-orthogonal projection $\sphereproject$ onto $\stu$ as follows:
				\begin{align}
					\sphereproject^\alpha_\beta:=\delta^\alpha_\beta+\frac{1}{2}\Lunit^\alpha\uLunit_\beta+\frac{1}{2}\uLunit^\alpha\Lunit_\beta,
				\end{align}
				where $\delta^\alpha_\beta$ is the Kronecker delta.

				We use the notation $\sgabs{\upxi}$ to denote the norm of the $\stu$-tangent tensor field $\upxi=\sphereproject\xi$, with respect to $\gsphere$, that is,
				\begin{align}
					\sgabs{\upxi}^2:=\gsphere(\upxi,\upxi).
				\end{align}
				
				Then, we denote $\gtr\upxi$ as the trace of a $\binom{0}{2}$ $\stu$-tangent tensor field $\upxi$, with respect to $\gsphere$:
				\begin{align}
					\gtr\upxi:=(\gsphere^{-1})^{AB}\upxi_{AB}.
				\end{align}
				
				Furthermore, we define $\hat{\upxi}$ to be the trace-free part of the $\binom{0}{2}$ $\stu$-tangent tensor field $\upxi$:
				\begin{align}
					\hat{\upxi}_{AB}:=\upxi_{AB}-\frac{1}{2}(\gtr\upxi)\gsphere_{AB}.
				\end{align}
			\end{definition}
			\subsection{Conformal Energy}\label{SS:Sectionconformalenergy}
			In this section, we first reduce the proof of localized decay estimates in Theorem \ref{TH:Spatiallylocalizeddecay} to the conformal energy estimates in Theorem \ref{TH:BoundnessTheorem}. Then, appealing to the control of the geometry established in Proposition \ref{PR:mainproof}, we demonstrate the boundedness of the conformal energy in Theorem \ref{TH:BoundnessTheorem}. These steps are crucial for concluding the Strichartz estimates (\ref{ES:estimateswavevariables}). We omit the detailed proof of these steps, since they follow identically from \cite[Section 4,7]{AGeoMetricApproach}.
			
			\subsubsection{Setup of the Conformal Energy}\label{SS:SectionSetupconformal}
			In order to give the definition of our conformal energy, we fix two smooth non-negative cut-off functions $\underline{\upvarpi}$ and $\upvarpi$, depending only on $t$ and $u$, such that,
			\begin{align}\label{DE:region1}
				\underline{\upvarpi}&=
				\begin{cases}
					1& \text{on } 0\leq u\leq t\\
					0& \text{on } u\leq-\frac{t}{4}
				\end{cases} , &
				\upvarpi&=
				\begin{cases}
					1& \text{on } 0\leq u\leq \frac{1}{2}t\\
					0& \text{if } u\geq\frac{3}{4}t\ \text{or}\ u\leq-\frac{t}{4}
				\end{cases}.
			\end{align}
			Also, see Figure \ref{FI:figureCEE1} for an illustration of the regions.
			
			\begin{definition}[Conformal energy]\label{DE:definitionofconformalenergy}
				For any scalar $\varphi$ vanishing outside $\intregion$ (defined in Section \ref{SS:opticalfunction}), we define the conformal energy $\mathfrak{C}[\varphi]$ as follows:
				\begin{align}
					\Energyconformal[\varphi](t)=\iEnergyconformal[\varphi](t)+\eEnergyconformal[\varphi](t),
				\end{align}
				where
				\begin{subequations}
					\begin{align}
						\iEnergyconformal[\varphi](t)&=\int_{\St}(\underline{\upvarpi}-\upvarpi)t^2\left(\abs{\Dfour\varphi}^2+\abs{\rgeo^{-1}\varphi}^2\right)\diff\tvol,\\
						\eEnergyconformal[\varphi](t)&=\int_{\St}\upvarpi\left(\rgeo^2\abs{\Dfour_\Lunit\varphi}^2+\rgeo^{2}\abs{\angnabla\varphi}^2+\abs{\varphi}^2\right)\diff\tvol.
					\end{align}
				\end{subequations}
			\end{definition}
			\label{FI:figureCEE1}\begin{figure}
				\centering
				\includegraphics[width=12cm,height=10cm]{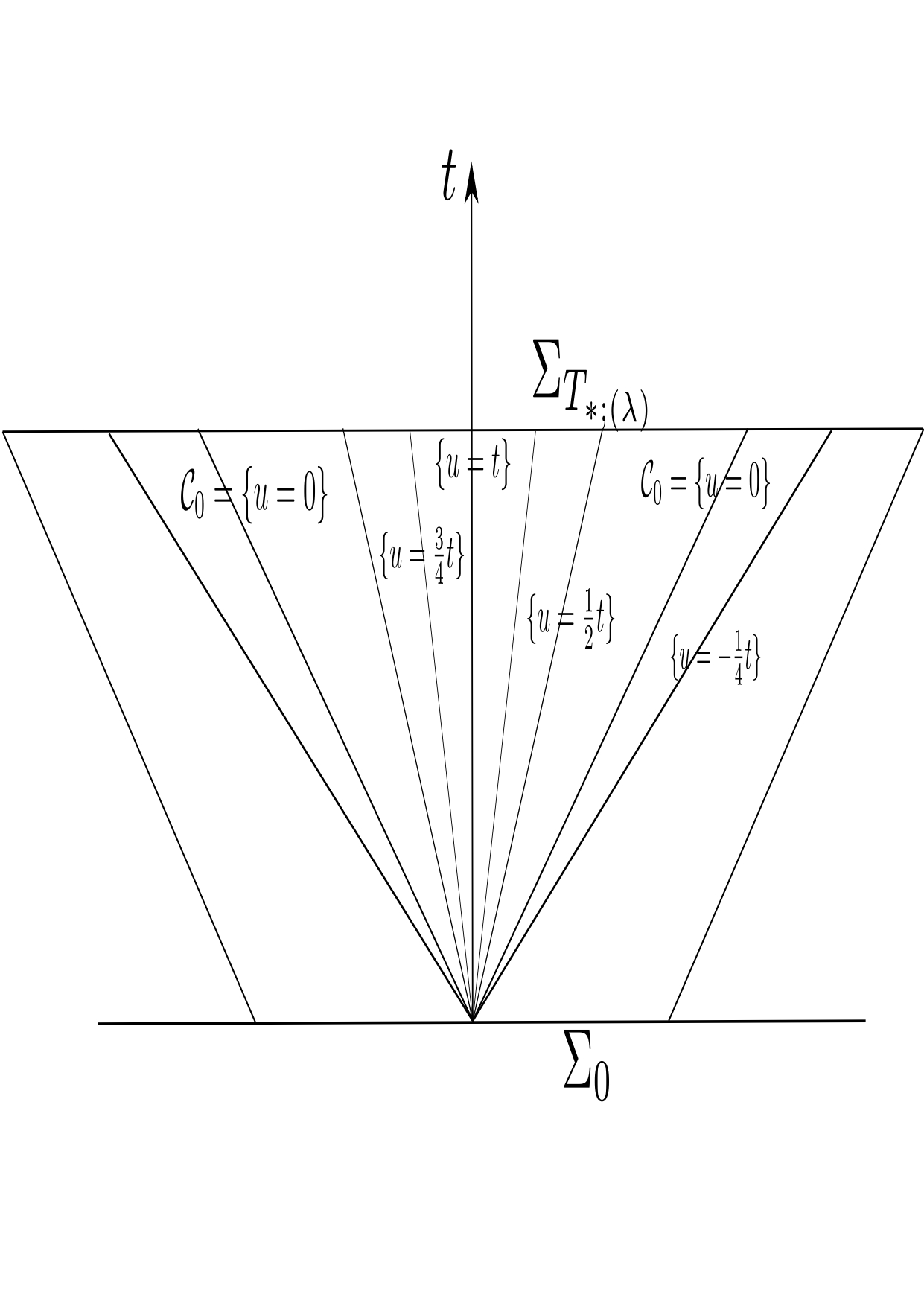}
				\caption{The illustration of the regions in (\ref{DE:region1})}
			\end{figure}
			In this section, we prove the following theorem.
			\begin{theorem}[Boundedness theorem]\label{TH:BoundnessTheorem}
				Let $\varphi$ be any solution of 
				\begin{align}\label{EQ:CEE1}
					\boxg\varphi=0
				\end{align}
				on $[0,\Trescale]\times\St$, with $\varphi(1)$ supported in $B_R\subset\intregion\bigcap\Sigma_{1}$, where $B_R$ is defined in Theorem \ref{TH:Spatiallylocalizeddecay}. Then, under the bootstrap assumptions, for $t\in[1,\Trescale]$, the conformal energy $\Energyconformal[\varphi]$ satisfies the estimate:
				\begin{align}
					\Energyconformal[\varphi](t)\lesssim(1+t)^{2\varepsilon}\left(\norm{\pfour\varphi}^2_{L^2(\Sigma_{1})}+\norm{\varphi}^2_{L^2(\Sigma_{1})}\right),
				\end{align}
				where $\varepsilon>0$ is an arbitrary small number, and the constants in $``\lesssim"$ can depend on $\varepsilon$.
			\end{theorem}
			Before proving Theorem \ref{TH:BoundnessTheorem}, we first show that the conformal energy estimates imply the decay estimates in Theorem \ref{TH:Spatiallylocalizeddecay}.

			\subsubsection{Theorem \ref{TH:BoundnessTheorem} implies Theorem \ref{TH:Spatiallylocalizeddecay} }\label{SSS:reductiondecay}
			The proof of Theorem \ref{TH:Spatiallylocalizeddecay} by using Theorem \ref{TH:BoundnessTheorem} is achieved via product estimates and the Bernstein inequality from Littlewood-Paley theory. We refer readers to \cite[Section 8]{ImprovedLocalwellPosedness} and \cite[Section 4]{AGeoMetricApproach} for the detailed proof.
			
			\subsubsection{Discussion of the Proof of Theorem \ref{TH:BoundnessTheorem} using Proposition \ref{PR:mainproof}}\label{SSS:BoundnessTheorem}
			Let $\initialenergy:=\norm{\pfour\varphi}_{L^2(\Sigma_{1})}+\norm{\varphi}_{L^2(\Sigma_{1})}$. For $\varepsilon_0$ defined as in Section \ref{SS:ChoiceofParameters}, we make the following assumptions:
			\begin{subequations}\label{BA:conformalenergy}
				\begin{align}
					\label{BA:conformalenergy1}\Energyconformal[\varphi](t)\leq&\lambda^{2\varepsilon_0}\initialenergy^2,\\
					\label{BA:conformalenergy2}\twonorms{\varphi}{u}{2}{\sangle}{2}{\St}\leq&(t+1)^{-1}\lambda^{\varepsilon_0}\initialenergy.
				\end{align}
			\end{subequations}
			
			Then, we restate Theorem \ref{TH:BoundnessTheorem} and improve bootstrap assumptions \eqref{BA:conformalenergy} in the following theorem:
			\begin{theorem}[Boundedness of the conformal energy]\label{TH:Boundnesstheorem}
				The following three estimates hold for solutions to equation (\ref{EQ:CEE1}):
				\begin{subequations} \label{ES:conformalenergy0}
					\begin{align}
						\label{ES:conformalenergy1}\int_{\St}\abs{(t+1)\angnabla\varphi}^2+\abs{(t+1)\Lunit\varphi}^2\diff\tvol\lesssim&(t+1)^{2\varepsilon}\initialenergy^2,\\
						\label{ES:conformalenergy2}\int_{\St\bigcap\{u\geq\frac{t}{2}\}}\energycurrent{\Timelike}_\alpha[\varphi]\geonormals^{\alpha}\diff\tvol\lesssim&(1+t)^{-2}\initialenergy^2,\\
						\label{ES:conformalenergy3}\twonorms{\varphi}{u}{2}{\sangle}{2}{\St}\leq&(t+1)^{-1+\varepsilon}\initialenergy,
					\end{align}
				\end{subequations}
				where $\geonormals$ denotes the normal vector field of the space-like hypersurface(s) with respect to $\gfour$.

				The estimates in \eqref{ES:conformalenergy0} can then imply the following improvement of the bootstrap assumptions:
				\begin{align}\label{ES:conformalenergy}
					\Energyconformal[\varphi](t)\lesssim(t+1)^{2\varepsilon}\initialenergy^2,
				\end{align}
				where $\varepsilon>0$ can be arbitrarily small, and the constants in $``\lesssim"$ can depend on $\varepsilon$. Since $ t\leq\Trescale\leq\lambda^{1-8\varepsilon_0}\Tstar$ as in (\ref{ES:lengtht}) and $\varepsilon>0$ can be arbitrarily close to $0$ (in particular, $\varepsilon<\varepsilon_0$), it follows that (\ref{ES:conformalenergy}) is an improvement of (\ref{BA:conformalenergy1}).
			\end{theorem}
			
\begin{figure} \centering    
				\subfigure[$\Sigma_{\tau_1,R}^\tau$ and $\region_{\tau_1,R}^\tau$.] {
					\label{FI:FigureCEE2}     
					\includegraphics[width=7cm,height=7cm]{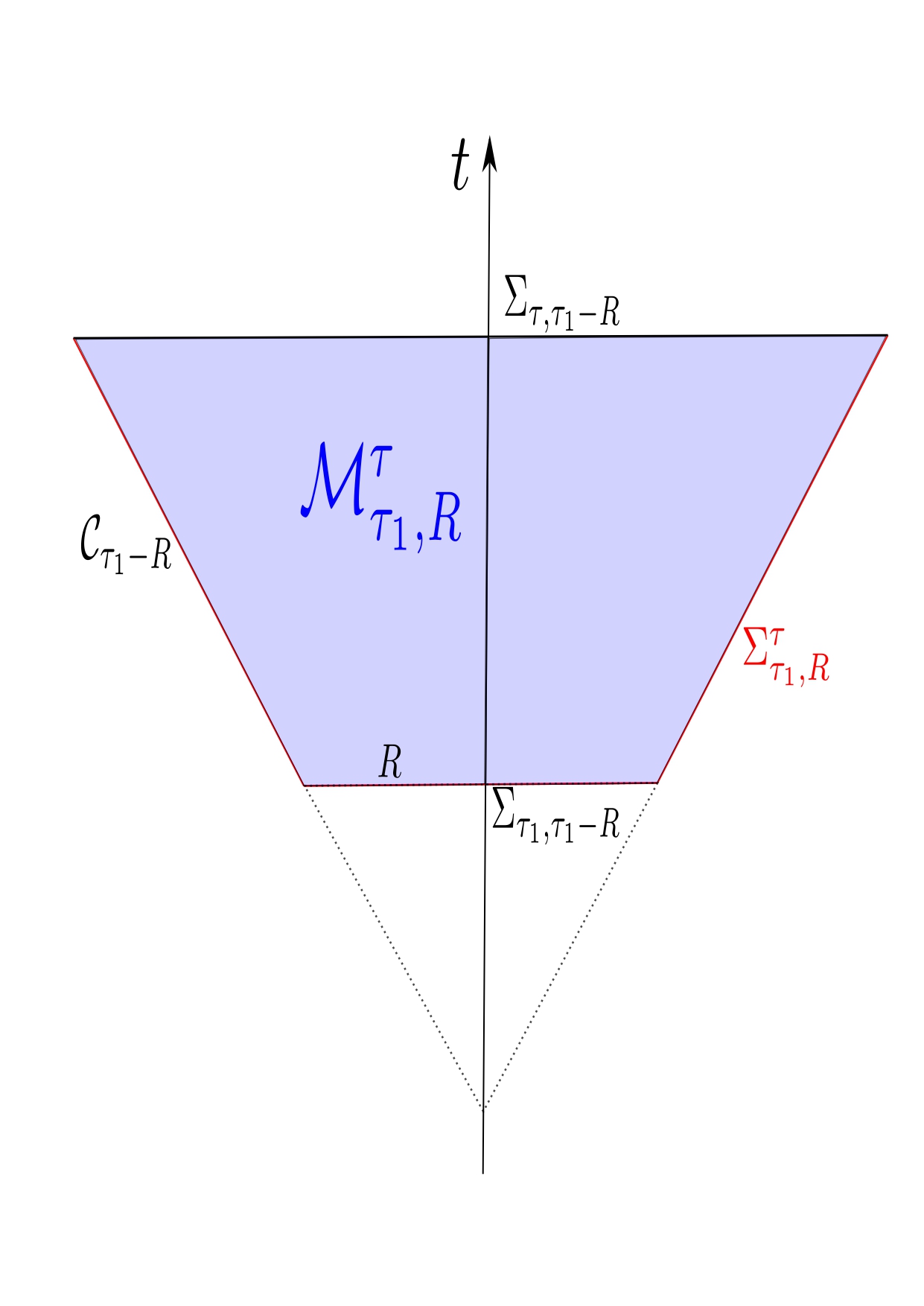}  
				}     
				\subfigure[$D_{\tau_1,R}^{\tau_2,\tau}$.] {
					\label{FI:FigureCEE3}     
					\includegraphics[width=7cm,height=7cm]{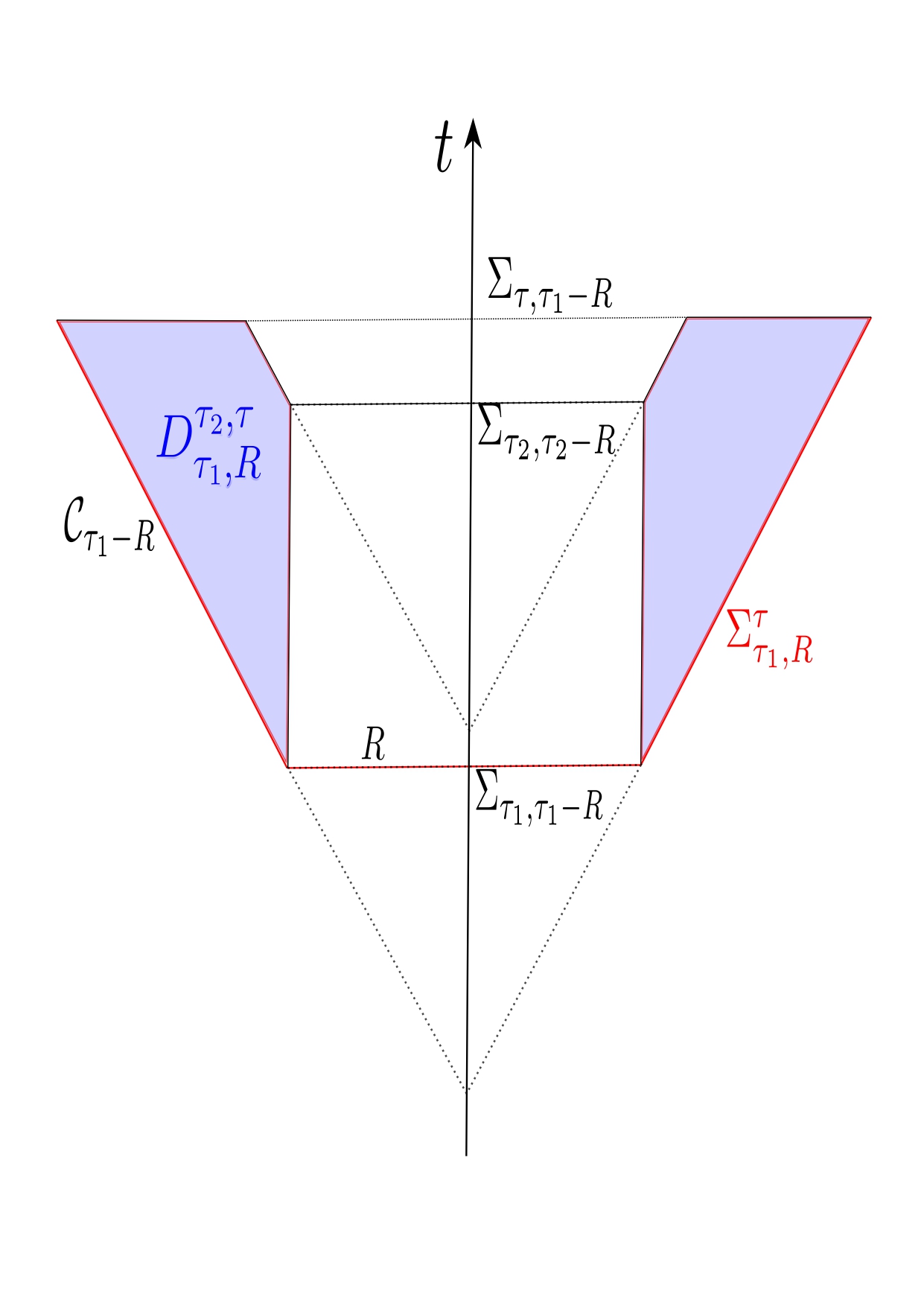}  
				}  
				\caption{The regions where the divergence theorem is applied in conformal energy estimates.}     
			\end{figure}

			Now we discuss the proof of the conformal energy estimates. First, we consider equation $\boxg\varphi=0$. Choose $X=f\spherenormal$, with $f$ defined as follows:
			\begin{align}
				f:=&\beta-\frac{\beta}{(1+\rgeo)^\alpha},
			\end{align}
			where $\rgeo$ is defined in (\ref{DE:rgeo}), $\alpha=2\varepsilon_0$ and $\beta=\frac{2}{\alpha}$. Now we let
			\begin{align}
				\label{D:Theta}\Theta:=&\rgeo^{-1}\left(\beta-\frac{\beta}{(1+\rgeo)^\alpha}\right).
			\end{align} We apply divergence theorem on region $\region_{\tau_1,R}^\tau$ (see Figure \ref{FI:FigureCEE2} for $\region_{\tau_1,R}^\tau$) for the following modified current:
			\begin{align}
				\label{DE:modifiedcurrent}\modifiedcurrent{X}_\mu[\varphi]=Q_{\mu\nu}[\varphi]X^\mu+\frac{1}{2}\Theta\p_\mu(\varphi^2)-\frac{1}{2}\varphi^2\p_\mu\Theta.
			\end{align}
			Then, we get the following Morawetz-type energy estimate (which allows us to obtain the uniform bound for a standard energy of $\varphi$) on $\region_{\tau_1,R}^{\tau_2}$ (as in Figure \ref{FI:FigureCEE2}):
			\begin{align}
				&\label{ES:coneenergy1.3}\int_{\region_{\tau_1,R}^{\tau_2}}\left(\frac{(\Timelike\varphi)^2+(\spherenormal\varphi)^2+\abs{\angnabla\varphi}^2}{(1+\rgeo)^{\alpha+1}}+\frac{\varphi^2}{\rgeo(1+\rgeo)^{\alpha+2}}\right)\diff\gvol\\
				\notag&\lesssim\int_{\Sigma_{\tau_1,R}^{\tau_2}}\energycurrent{\Timelike}_\alpha[\varphi]\geonormals^\alpha\diff\gvol+\onenorm{\rgeo^{\frac{1}{2}}\varphi}{\sangle}{2}{S_{\tau_2,\tau_1-R}}^2+\onenorm{\rgeo^{\frac{1}{2}}\varphi}{\sangle}{2}{S_{\tau_1,\tau_1-R}}^2,
			\end{align}
			where $1\leq\tau_1<\tau_2\leq\Trescale$, $R\leq R^\prime<2R$, and $\energycurrent{\Timelike}_\alpha[\varphi]$ is the standard energy current defined in (\ref{DE:energycurrent}).
			
			Next, we consider the conformally transformed wave equation (see Definition \ref{DE:conformalchangemetric} for the conformally changed metric and the unknown):
			\begin{subequations}\label{EQ:boxtildeg}
				\begin{align}
					\label{EQ:boxgrescaledvarphiconformal1}
					e^{2\conformalfactor}\boxgrescale\varphiconformal
					=&\anglap\varphiconformal-\uLunit\left(\Lunit\varphiconformal+\frac{1}{2}\Restrace\reschi\varphiconformal\right)-\left(\frac{1}{2}\Restrace\resuchi-k_{\spherenormal\spherenormal}\right)\left(\Lunit\varphiconformal+\frac{1}{2}\Restrace\reschi\varphiconformal\right)\\
					\notag&+2\left(\upzeta+\angnabla\conformalfactor\right)\cdot\angnabla\varphiconformal+\frac{1}{2}\left(\uLunit\Restrace\reschi+\frac{1}{2}\Restrace\reschi\Restrace\resuchi-\Restrace\reschi k_{\spherenormal\spherenormal}\right)\varphiconformal,\\
					\label{EQ:boxgrescaledvarphiconformal2}
					e^{3\conformalfactor}\boxgrescale\varphiconformal=&\boxg\varphi-\left(\boxg\conformalfactor+\Dfour^\alpha\conformalfactor\Dfour_\alpha\conformalfactor\right)\varphi.
				\end{align}
			\end{subequations}
			We apply the multiplier approach for (\ref{EQ:boxtildeg}) using $\rgeo^p\Lunit$ type vector fields in the region $D_{\tau_1,R}^{\tau_2,\tau}$ (see Figure \ref{FI:FigureCEE3}), where $1\leq\tau_1<\tau_2<\Trescale$, to control the conformal energy in the exterior region and gain energy decay along null slices. Finally, we control the conformal energy in the interior region with the help of the argument in \cite{Anewphysical-spaceapproachtodecayforthewaveequationwith} regarding energy decay in each spatial-null slice.\\
			
			\indent
			The proof of Theorem \ref{TH:BoundnessTheorem} will complete the reduction of Strichartz estimate. One can follow the steps listed in \cite[Section 11]{3DCompressibleEuler} to derive conformal energy estimates, using the control of the Ricci coefficients given. One could go through the details of the argument in \cite[Section 7]{AGeoMetricApproach} or \cite[Section 10]{Yuthesis}. Interested readers could also look into \cite[Section 3]{ACommutingvector fieldsApproachtoStrichartz} by Klainerman and \cite{AGeoMetricApproach} by Wang for initial ideas.
			We emphasize that, with the control of the geometry in Proposition \ref{PR:mainproof} at our disposal, the proof of Theorem \ref{TH:Boundnesstheorem} follows the same argument as in \cite[Section 7]{AGeoMetricApproach}. Note that the proof of Proposition \ref{PR:mainproof} is independent of the conformal energy estimates here.
			\section{Energies along Null Hypersurfaces of Faster Wave}\label{S:connectioncoefficientsandpdes}
			In this section, we derive energy estimates along null hypersurfaces $\coneu$, which are essential for obtaining the mixed-norm estimates in Proposition \ref{PR:mainproof}.
			\subsection{Null Fluxes along Null Hypersurfaces of Faster Wave}\label{energyalongnullhypersurfaces}
			Here, we define the null fluxes on $\gfour$-null hypersurfaces for both the ``divergence-part" and the ``curl-part",
 and show that the null fluxes for ``curl-part" are coercive.  
			\begin{definition}[Energy-momentum tensor with respect to $\hfour$]
				We define the following energy-momentum tensor:
				\begin{align}\label{DE:breveQ}
					\breve{Q}_{\mu\nu}[\varphi]:=\p_\mu\varphi\p_\nu\varphi-\frac{1}{2}\hfour_{\mu\nu}(\hfour^{-1})^{\alpha\beta}\p_\alpha\varphi\p_\beta\varphi,
				\end{align} 
				where $\hfour$ is defined in (\ref{DE:Defh}). We set the energy current $\HJen{\Timelike}$ to be
				\begin{align}\label{DE:breveJ}
					\HJen{\Timelike}^\alpha[\varphi]:=\breve{Q}^{\alpha\beta}[\varphi](\Timelike_\flat)_\beta,
				\end{align} 
				where $Q_{\alpha\beta}$ is the energy momentum tensor defined in (\ref{DE:energymomentum}), and $(V_\flat)_\beta:=\hfour_{\alpha\beta}V^\alpha$. Notice that
				\begin{align}
					\breve{\Dfour}_\alpha\HJen{\Timelike}^{\alpha}[\varphi]=\boxh\varphi(\Timelike\varphi)+\frac{1}{2}\breve{Q}^{\mu\nu}[\varphi]\Hdeform{\Timelike}_{\mu\nu},
				\end{align} 
				where $\breve{\Dfour}$ denotes the Levi-Civita connection with respect to $\hfour$, and $\Hdeform{\Timelike}$ is defined as
				\begin{align}
					\Hdeform{\Timelike}:=\breve{\Dfour}_{\alpha}(\Timelike_\flat)_\beta+\breve{\Dfour}_{\beta}(\Timelike_\flat)_\alpha.
				\end{align}
			\end{definition}
			\begin{definition}[Null fluxes]\label{DE:acouticnullfluxes}
				For a function $\varphi$ on $\coneu$, we define the null fluxes $\mathcal{F}_{(1)}[\varphi;\coneu]$ and $\mathcal{F}_{(2)}[\varphi;\coneu]$ as follows:
				\begin{subequations}\label{DE:coneenergy}
					\begin{align}
						\mathcal{F}_{(1)}[\varphi;\coneu]&:=\int_{\coneu}\left((\Lunit\varphi)^2+\sgabs{\angnabla\varphi}^2\right)\diff\spherevol\diff t,\\
						\mathcal{F}_{(2)}[\varphi;\coneu]&:=\int_{\coneu}\breve{Q}_{\alpha\beta}[\varphi]\Timelike^\alpha\Conenormal^\beta\diff\spherevol\diff t,
					\end{align}
				\end{subequations}
				where $\diff\spherevol$ is the volume form of the induced (from $\gfour$) metric $\gsphere$ on the $\stu$-sphere, and $\angnabla$ is the Levi-Civita connection on $\stu$ with respect to $\gsphere$.
			\end{definition}
			For the null flux $\mathcal{F}_{(2)}$, we prove the following important coerciveness lemma.
			\begin{lemma}[Coerciveness of $\mathcal{F}_{(2)}$]\label{LE:hCoercive}
				For $\mathcal{F}_{(2)}$ defined in (\ref{DE:coneenergy}), the following coerciveness property holds:
				\begin{align}\label{ES:hCoercive}
					\mathcal{F}_{(2)}[\varphi;\coneu]\approx\int_{\coneu}\abs{\pfour\varphi}^2\diff\spherevol\diff t.
				\end{align}
			\end{lemma}
			\begin{proof}[Proof of Lemma \ref{LE:hCoercive}]
				In Cartesian coordinates, using \eqref{EQ:EW3.2} and \eqref{unit V}, we have
				\begin{align}\label{EQ:Coercive1}
					\HJen{\Timelike}^\alpha(\Conenormal_\flat)_\alpha=&\Timelike\varphi\Conenormal\varphi+\frac{\Conenormal^0}{2}\left((\hfour^{-1})^{\alpha\beta}\p_\alpha\varphi\p_\beta\varphi\right)\\
					\notag=&\Conenormal^0(\partial_t \varphi)^2+\Conenormal^i\partial_t\varphi\partial_i\varphi+\frac{\Conenormal^0}{2}\left(-(\partial_t\varphi)^2+(\hfour^{-1})^{ij}\p_i\varphi\p_j\varphi\right)\\
					\notag=&\frac{\Conenormal^0}{2}\left(\p_t\varphi+\frac{\Conenormal^i}{\Conenormal^0}\p_i\varphi\right)^2-\frac{\Conenormal^0}{2}\left(\frac{\Conenormal^i}{\Conenormal^0}\p_i\varphi\right)^2+\frac{\Conenormal^0}{2}\sum\limits_i\speedtwo^2(\p_i\varphi)^2\\
					\notag=&\frac{1}{2\Conenormal^0}\left(\Conenormal\varphi\right)^2+\frac{1}{2\Conenormal^0}\sum\limits_i\left(\sum\limits_{i\neq j}(\Conenormal^j)^2+\speedtwo^2\right)(\p_i\varphi)^2.
				\end{align}
				Recalling the expression of $\Conenormal$ in Remark \ref{V formula}, we can bound $\Conenormal^0$ from above by using the estimate on the eikonal function in \eqref{bt6}. Then, by (\ref{EQ:Coercive1}), we obtain
				\begin{align}
					(\Conenormal\varphi)^2+(\p\varphi)^2\lesssim\HJen{\Timelike}^\alpha(\Conenormal_\flat)_\alpha.
				\end{align}
				We hence derive the desired estimate (\ref{ES:hCoercive}).
				
			\end{proof}

			\subsection{Energy Estimates along Null Hypersurfaces of Faster Wave}\label{SS:EnergyFluxes}
			In this section, we derive the energy estimates for null fluxes along $\gfour$-Null Hypersurfaces. 
			\begin{proposition}[Energy estimates along $\gfour$-null hypersurfaces]\label{PR:EnergyNullFluxes}
				For any $3<N<7/2$, under the initial data assumption in Section \ref{SS:Data}, bootstrap assumptions (\ref{BA:ba}), applying the standard energy estimates Proposition \ref{PR:EnergyEstimatesforDiv}, Theorem \ref{TH:LinearStrichartz}, we derive the following estimates along null hypersurfaces $\coneu$ (defined in (\ref{DE:coneu})) for $u\in[-w_*,\Trescale]$:
				\begin{subequations}\label{ES:coneenergy}
					\begin{align}
						\label{ES:coneenergy1}\mathcal{F}_{(1)}[\pfour\vvariables;\coneu]+\sum\limits_{\upnu>1}\upnu^{2(N-3)}\mathcal{F}_{(1)}[\littlewood\pfour\vvariables;\coneu]\lesssim\lambda^{-1},\\
						\label{ES:coneenergy2}\mathcal{F}_{(2)}[\pfour\cvvariables;\coneu]+\sum\limits_{\upnu>1}\upnu^{2(N-3)}\mathcal{F}_{(2)}[\littlewood \pfour\cvvariables;\coneu]\lesssim\lambda^{-1},\\
						\label{ES:coneenergy3}\mathcal{F}_{(2)}[\pfour^2\cvvariables;\coneu]+\sum\limits_{\upnu>1}\upnu^{2(N-3)}\mathcal{F}_{(2)}[\littlewood \pfour^2\cvvariables;\coneu]\lesssim\lambda^{-3}.
					\end{align}
				\end{subequations}
			\end{proposition}
			\begin{remark}
				Note that in the proof of Proposition \ref{PR:EnergyNullFluxes}, we used the rescaled version of Theorem \ref{TH:LinearStrichartz} and Proposition \ref{PR:EnergyEstimatesforDiv}, where the rescaled quantities are defined in Definition \ref{DE:rescaledquantities}.
			\end{remark}
			
			\begin{proof}[Proof of Proposition \ref{PR:EnergyNullFluxes}]
				
				We first prove (\ref{ES:coneenergy1}). {We define the energy current 
				\begin{align}
					\Jenarg{\Timelike}{\alpha}[\varphi]:=Q^{\alpha\beta}[\varphi]\Timelike_\beta,
				\end{align} 
				with $Q_{\alpha\beta}$ being the energy momentum tensor defined in (\ref{DE:energymomentum}).} Notice that
				\begin{align}
					\Dfour_\alpha\Jenarg{\Timelike}{\alpha}[\varphi]=\boxg\varphi(\Timelike\varphi)+\frac{1}{2}Q^{\mu\nu}[\varphi]\deform{\Timelike}_{\mu\nu},
				\end{align} 
				where $\deform{\Timelike}_{\mu\nu}$ is define in (\ref{DE:deformtensor}). It is straightforward to check that:
				\begin{align}\label{EQ:JLunit}
					\Jenarg{\Timelike}{\alpha}[\varphi]\Lunit_\alpha=\frac12\left((\Lunit\varphi)^2+\sgabs{\angnabla\varphi}^2\right).
				\end{align}
				{Then, applying the divergence theorem over the region bounded by $\coneu$, $\Sigma_{t_0}$ and $\St$ with $t_0:=\max\{0,u\}$,} we have:
				\begin{align}\label{EQ:coneenergy1.1}
					\mathcal{F}_{(1)}[\varphi;\coneu]=&\int_{\St}\abs{\pfour\varphi}^2\diff\tvol-\int_{\Sigma_{t_0}}\abs{\pfour\varphi}^2\diff\tvol\\
					\notag&+\int_{\bigcup\limits_{t\geq u^\prime\geq u}\mathcal{C}_{u^\prime}}\left(\boxg\varphi(\Timelike\varphi)+\frac{1}{2}Q^{\mu\nu}[\varphi]\deform{\Timelike}_{\mu\nu}\right)\diff\gvol.
				\end{align}
				First, for the last term on the RHS of \eqref{EQ:coneenergy1.1}, we have $\abs{\Timelike\varphi}\lesssim\abs{\pfour\varphi}$, $\abs{Q^{\mu\nu}[\varphi]}\lesssim\abs{\pfour\varphi}^2$ and $\abs{\deform{\Timelike}_{\mu\nu}}\lesssim\abs{\pfour\vvariables,\pfour\cvvariables}$, where one should note that $\vvariables,\cvvariables$ are the rescaled quantities in Definition \ref{DE:rescaledquantities}. . 
				
				Next, we substitute $\pfour\vvariables$ for $\varphi$, 
				and consider the term $(\boxg\pfour\vvariables)(\Timelike\pfour\vvariables)$. By Cauchy-Schwarz inequality, and estimate $\int_{\bigcup\limits_{u^\prime\geq u}\mathcal{C}_{u^\prime}}\boxg\varphi(\Timelike\varphi)\diff\gvol\leq\int_0^{\Trescale}\int_{\Sigma_{\tau}}\abs{\boxg\varphi}\abs{\Timelike\varphi}\diff \tvol\diff\tau$, we derive:
				\begin{align}\label{ES:coneenergy1.1}
					\int_{\bigcup\limits_{t\geq u^\prime\geq u}\mathcal{C}_{u^\prime}}(\boxg\pfour\vvariables)(\Timelike\pfour\vvariables)\diff\gvol\lesssim&\int_0^{\Trescale}\int_{\Sigma_{\tau}}\abs{\boxg\pfour\vvariables}\abs{\Timelike\pfour\vvariables}\diff \tvol\diff\tau\\
					\notag\lesssim&\twonorms{\pfour\p\vvariables,\pfour\p\cvvariables,\pfour\p^2\cvvariables}{t}{\infty}{x}{2}{[0,\Trescale]\times\St}^2\cdot\twonorms{\vvariables,\p\vvariables,\p\cvvariables}{t}{2}{x}{\infty}{[0,\Trescale]\times\St}^2.
				\end{align}
				Applying (\ref{EQ:coneenergy1.1}), (\ref{ES:coneenergy1.1}), the bootstrap assumption (\ref{BA:Div}), energy estimates (\ref{ES:linearcurl}), (\ref{ES:energyestimates2}), and rescaling the coordinates back $(\lambda (t-t_k),\lambda x)\rightarrow(t,x)$, we obtained the desired result (\ref{ES:coneenergy1}) for $\pfour\vvariables$.
				
				The proof for $\littlewood \pfour\vvariables$ is of the same fashion. \\
				
				Now we prove (\ref{ES:coneenergy2}). 
				Applying the divergence theorem to the energy current $\HJen{\Timelike}$ defined in (\ref{DE:breveJ}), over the region bounded by $\coneu$, $\Sigma_{t_0}$ and $\St$ with $t_0:=\max\{0,u\}$, we obtain:
				\begin{align}\label{EQ:coneenergy2}
					\mathcal{F}_{(2)}[\varphi;\coneu]=&\int_{\St}\abs{\pfour\varphi}^2\diff\htvol-\int_{\Sigma_{t_0}}\abs{\pfour\varphi}^2\diff\htvol\\
					\notag&+\int_{\bigcup\limits_{t\geq u^\prime\geq u}\mathcal{C}_{u^\prime}}\left(\boxh\varphi(\Timelike\varphi)+\frac{1}{2}\breve{Q}^{\mu\nu}[\varphi]\Hdeform{\Timelike}_{\mu\nu}\right)\diff\hvol.
				\end{align}
				
				We substitute $\pfour\cvvariables$ for $\varphi$. Considering the last integrand on the RHS of (\ref{EQ:coneenergy2}), we note that by (\ref{BA:databa}), on $\coneu$ there holds $\diff\hvol\approx\diff\tvol\diff t$. Therefore, applying the Cauchy-Schwarz inequality, (\ref{EQ:coneenergy2}), energy estimates and Strichartz estimates (\ref{ES:linearcurl}), and rescaling the coordinates back $(\lambda (t-t_k),\lambda x)\rightarrow(t,x)$, we then obtain the desired result (\ref{ES:coneenergy2}) for $\pfour\cvvariables$.
				\color{black}
				The proofs for $\littlewood \pfour\cvvariables$ and (\ref{ES:coneenergy3}) follow in the same fashion. 
			\end{proof}

				\section{Geometric Spacetime Quantities and Their Evolution Equations}\label{S:causalgeometry}
				In this section, we introduce the notation for many geometric quantities, along with their corresponding bootstrap assumptions. The improved estimates are provided in Section \ref{SS:sectionRestatement}. The proof of these estimates for geometric quantities is obtained by transport equation and div-curl estimates for the geometric quantities, decomposition of Ricci curvature components, trace and Sobolev inequalities. We note that these estimates are similar to those for compressible/relativistic Euler and quasilinear wave equations; see \cite[Chapter 9]{Yuthesis}, \cite[Section 10]{3DCompressibleEuler}, \cite[Section 5]{AGeoMetricApproach}, and \cite{ImprovedLocalwellPosedness}. However, due to the special structure of elastic wave equations, we must account for the influence of ``slower-wave" on the geometry of ``faster-wave".
				\subsection{Connection Coefficients}
                {We start with introducing the connection coefficients with respect to the ``faster-wave" geometry.}
				\subsubsection{Levi-Civita connections, angular operators and curvatures}\label{SSS:levi-civita2}
				If $\upxi$ is a space-time tensor, then $\sphereproject\xi$ denotes its $\gfour$-orthogonal projection onto $\stu$. If both $V$ and $\upxi$ are $\stu$-tangent, we denote the covariant derivative on $\stu$ with respect to the induced metric $\gsphere$ as $\angnabla_V\upxi$. Then we define $\angD_V\upxi:=\sphereproject\Dfour_V\upxi$ where $V$ is a vector and $\Dfour_V\upxi$ is the covariant derivative of $\upxi$ in the $V$ direction. Note that $\angD_V\upxi:=\angnabla_V\upxi$ when both $V$ and $\upxi$ are $\stu$-tangent.\\
				\indent
				We let $\Riemfour{\alpha}{\beta}{\gamma}{\delta}$ denote the Riemann curvature tensor of $\gfour$, and define $\Ricfour{\alpha}{\beta}:=(\gfour^{-1})^{\gamma\delta}\Riemfour{\gamma}{\alpha}{\delta}{\beta}$. We use the convention that
				\begin{align}\label{DE:curvature}
					\gfour(\Dfour_X\Dfour_YW-\Dfour_Y\Dfour_XW,Z)=\Riemfour{Z}{W}{X}{Y}+\gfour(\Dfour_{[X,Y]}W,Z),
				\end{align}
				where $X,Y,W,Z$ are vector fields, $[\cdot,\cdot]$ is the Lie bracket.
				\begin{definition}[Connection coefficients]
					\label{DE:DEFSOFCONNECTIONCOEFFICIENTS}
					We define the second fundamental form $k$ of $\Sigma_t$
					to be the $\binom{0}{2}$ $\Sigma_t$-tangent tensor,
					such that the following relation holds for all $\Sigma_t$-tangent 
					vector fields $X$ and $Y$:
					\begin{align} 
						k(X,Y)
						& :=
						-
						\gfour(\Dfour_X \Timelike, Y).
						\label{DE:SECONDFUNDOFSIGMAT} 
					\end{align}
					
					We denote a pair of unit orthogonal spherical vector fields on $\stu$ by $\{e_A\}_{A=1,2}$.
					We define the second fundamental form $\spheresecondfund$ of $\stu$,
					the null second fundamental form $\upchi$ of $\stu$, 
					and $\underline{\upchi}$ 
					to be the following type $\binom{0}{2}$ $\stu$-tangent tensors:
					\begin{subequations}
						\begin{align}  \label{DE:SECONDFUNDOFSPHERESINSIGMAT}
							\spheresecondfund_{AB}
							& := \gfour(\Dfour_A \spherenormal,e_B),
							\\
							\upchi_{AB}
							& := \gfour(\Dfour_A \Lunit, e_B),
							&
							\underline{\upchi}_{AB}
							& := \gfour(\Dfour_A \uLunit, e_B).
						\end{align}
					\end{subequations}
					
					We define the torsion $\upzeta$ and $\underline{\upzeta}$ to be the following
					$\stu$-tangent one-forms:
					\begin{align} \label{DE:TORSION}
						\upzeta_A
						& := \frac{1}{2}
						\gfour(\Dfour_{\uLunit} \Lunit, e_A),
						&
						\underline{\upzeta}_A
						& := \frac{1}{2}
						\gfour(\Dfour_{\Lunit} \uLunit, e_A).
					\end{align}
					
				\end{definition}
				
				Note that $k$, $\spheresecondfund$ $\upchi$, and $\underline{\upchi}$ are symmetric,
				and the following relations hold:
				\begin{subequations}\label{DE:SECONDFORMLIEDIFFERENTIATIONDEF}
					\begin{align}
						k & = \frac{1}{2} \SigmatLie_{\Timelike} g
						= \frac{1}{2} \SigmatLie_{\Timelike} \gfour,
						\\
						\upchi & = \frac{1}{2} \angLie_{\Lunit} \gsphere
						= \frac{1}{2} \angLie_{\Lunit} \gfour,
						& 
						\underline{\upchi}
						& = \frac{1}{2} \angLie_{\uLunit} \gsphere
						= \frac{1}{2} \angLie_{\uLunit} \gfour,
						\label{DE:NULLSECONDFORMSLIEDIFFERENTIATIONDEF}
					\end{align}
					
				\end{subequations}
				where $\Lie$ is the standard Lie derivative, and $\angLie=\sphereproject\Lie$. Also, we have 
				\begin{align}
					\label{EQ:DfourNN}\Dfour_\spherenormal\spherenormal&=-\angnabla \ln b-\frac{1}{2}k_{\spherenormal\spherenormal}\uLunit-\frac{1}{2}k_{\spherenormal\spherenormal}\Lunit,
					&
					\Dfour_A \spherenormal
					& = \spheresecondfund_{AB}e_B-\frac{1}{2}k_{A\spherenormal}\Lunit-\frac{1}{2}k_{A\spherenormal}\uLunit,
				\end{align}
               and the following identities hold
		\begin{subequations}\label{EQ:CC}
					\begin{align}
						\Dfour_A \Lunit
						& = \upchi_{AB} e_B-k_{A \spherenormal}\Lunit,
						&
						\Dfour_A \uLunit
						& = \underline{\upchi}_{AB} e_B+k_{A \spherenormal}\uLunit, 
						\\
						\Dfour_{\Lunit} \Lunit
						& = - k_{\spherenormal \spherenormal} \Lunit,
						&
						\Dfour_{\Lunit} \uLunit
						& = 2 \underline{\upzeta}_A e_A +k_{\spherenormal \spherenormal}\uLunit, 
						\\
						\Dfour_{\uLunit} \Lunit
						& = 2 \upzeta_A e_A +k_{\spherenormal \spherenormal} \Lunit,
						&
						\Dfour_{\Lunit} e_A
						& = \angDarg{\Lunit} e_A+\underline{\upzeta}_A \Lunit,
						\\
						\Dfour_B e_A
						& = \angnabla_B e_A+\frac{1}{2} \upchi_{AB} \uLunit+\frac{1}{2} \underline{\upchi}_{AB} \Lunit,
						&
						\Dfour_{\uLunit} \uLunit 
						& = - 2 (\angnabla_A \ln \nulllapse) e_A
						-
						k_{\spherenormal \spherenormal}
						\uLunit, \\
						\Dfour_{\uLunit}e_A&=\angD_{\uLunit}e_A+\upzeta_A\uLunit-\angnabla_A\ln\nulllapse\Lunit.
					\end{align}
				\end{subequations}
				Moreover, we have the following relations among the connection coefficients:
				\begin{align} \label{E:CONNECTIONCOEFFICIENT}
					\upchi_{AB} 
					& =  \spheresecondfund_{AB}
					-
					k_{AB},
					&
					\underline{\upchi}_{AB} 
					& =  - 
					\spheresecondfund_{AB}
					-
					k_{AB},
					&
					\underline{\upzeta}_A
					& = - k_{A \spherenormal},
					&
					\upzeta_A
					& = \angnabla_A \ln \nulllapse
					+
					k_{A \spherenormal}.
				\end{align}
				
				The computation of the identities (\ref{EQ:DfourNN})-(\ref{EQ:CC}) relies on the following fact:
				Let $X$ be a vector field, then 
				\begin{equation}
					\label{EQ:decomp}X=\gfour (X, e_A)e_A-\frac{1}{2}\gfour (X, \Lunit)\uLunit-\frac{1}{2}\gfour (X, \uLunit)\Lunit.
				\end{equation}
				\subsection{Conformal Metric, Initial Conditions on $\Szero$ and on the Cone-tip Axis for the Eikonal Function $u$}
				\begin{definition}[Conformal factor and the modified null second fundamental form and geometric quantities]\label{DE:conformalchangemetric}
					In the interior region $\intregion$ (defined in Section \ref{SSS:solutionofu}), we define $\conformalfactor$ to be the solution of the following initial value problem for a transport equation:
					\begin{subequations}
						\begin{align}
							\label{DE:Lunitsigma}\Lunit\conformalfactor(t,u,\sangle)&=\frac{1}{2}[\Chfour_{\Lunit}](t,u,\sangle)&
							u\in[0,\Trescale],t\in[u,\Trescale],\sangle\in\Stwo,\\
							\label{DE:initialsigma}\conformalfactor(u,u,\sangle)&=0,&
							u\in[0,\Trescale], \sangle\in\Stwo,
						\end{align}
					\end{subequations}
					where $\Chfour_{\Lunit}:=\Chfour_\alpha {\Lunit}^\alpha$ with the Christoffel symbol $\Chfour_\alpha$ defined in Definition \ref{DE:Christoffelsymbols}.
					
					We set the conformal space-time metric and the Riemannian metric induced on $S_{t,u}$ to be
					\begin{align}\label{DE:rescaledgfour}
						\rescaledgfour&:=e^{2\conformalfactor}\gfour,&\congsphere:=e^{2\conformalfactor}\gsphere
						.\end{align}
					
					We define the null second fundamental forms to be the following symmetric $S_{t,u}$-tangent tensors:
					\begin{align}\label{DE:rescaledchi}
						\reschi&:=\frac{1}{2}\angLie_{\Lunit}\congsphere,& \resuchi:=\frac{1}{2}\angLie_{\uLunit}\congsphere 
						.\end{align}
					
					Using (\ref{DE:rescaledgfour})-(\ref{DE:rescaledchi}), it follows that $\upchi,\uchi$ and $\reschi,\resuchi$ are related by:
					\begin{subequations}\label{EQ:chiuchi}
						\begin{align}
							\reschi&=e^{2\conformalfactor}\left(\upchi+(\Lunit\conformalfactor)\gsphere\right),&\resuchi&=e^{2\conformalfactor}\left(\uchi+(\uLunit\conformalfactor)\gsphere\right),\\
							\Restrace\reschi&=\gtr\upchi+2\Lunit\conformalfactor=\gtr\upchi+\Chfour_{\Lunit},&\Restrace\resuchi&=\gtr\uchi+2\uLunit\conformalfactor,\\
							\upchi&=\frac{1}{2}\left(\Restrace\reschi-\Chfour_{\Lunit}\right)\gsphere+\hchi,&\uchi&=\frac{1}{2}\left(\Restrace\resuchi-2\uLunit\conformalfactor\right)\gsphere+\huchi
							.\end{align}
					\end{subequations}
					
					Then, we let
	\begin{align}\label{DE:tracechismall}
						\chismall:=\gtr\upchi+\Chfour_{\Lunit}-\frac{2}{\rgeo}=\Restrace\reschi-\frac{2}{\rgeo}
						.\end{align}
					We note that, the first equality in (\ref{DE:tracechismall}) holds in the whole region $\region$, while the second equality is valid only within the interior region $\intregion$.
					
					We denote $\CC:=(\chismall,\hchi,\upzeta)$ as the collection of ``small" connection coefficients, and $\ACC:=(\CC,\pfour\vvariables,\pfour\cvvariables)$ as the collection of ``small" order $1$ quantities\footnote{In view of (\ref{DE:tracechismall}), $\rgeo^{-1}$ is also an order $1$ quantity of $\gfour$. However, $\rgeo^{-1}$ is not necessarily small (especially near the tip axis).} of $\gfour$.
				\end{definition}
				
				\begin{definition}[Average values on $S_{t,u}$ and spherevolume ratio] \label{DE:average}For a scalar functions $f$, we define the average value of $f$, denoted by $\average{f}$, as below:
					\begin{align}
						\average{f}&:=\frac{1}{\sgabs{S_{t,u}}}\int_{S_{t,u}}fd\spherevol,&
						\sgabs{S_{t,u}}&:=\int_{S_{t,u}}1d\spherevol
						.\end{align}
					We define $\diff\flatspherevol$ to be the standard volume form on the Euclidean sphere, and we further define the volume form $\diff \spherevol$, $\diff\conspherevol$ and $\diff\gvol$ as follows:
					\begin{align}
						\diff \spherevol:=&\sqrt{\det\gsphere}\diff\flatspherevol,\\
						\diff \conspherevol:=&\sqrt{\det\congsphere}\diff\flatspherevol,\\
						\diff\gvol:=&\nulllapse\diff t\diff u\diff \spherevol.
					\end{align}
					Then, the spherevolume ratios are defined in the following way:
					\begin{align}
						\label{DE:spherevolume}\volume:=&\frac{\sqrt{\det\gsphere}}{\sqrt{\det\esphere}},\\
						\label{DE:conspherevolume}\conformalvol:=&\frac{\sqrt{\det\congsphere}}{\sqrt{\det\esphere}}.
					\end{align}
				\end{definition}
				\begin{proposition}\cite[Proposition 8.8. Evolution equations for the lapse and volume]{Yuthesis}. For geometric quantities $\volume,\nulllapse$, we have the following transport equations:
					\begin{subequations}
						\begin{align}
							\label{EQ:Lunitv}\Lunit \volume&=\volume\gtr\upchi,\\
							\label{EQ:uLunitv}\uLunit \volume&=\volume\gtr\uchi,\\
							\label{EQ:evol eq nulllapse}\Lunit\nulllapse&=-\nulllapse k_{\spherenormal\spherenormal}.
						\end{align}
					\end{subequations}
				\end{proposition}
				\begin{lemma}\cite[Lemma 8.9. Evolution equation for the average value on $S_{t,u}$]{Yuthesis}. For the average of a scalar functions $f$ (see Definition \ref{DE:average}), there holds
					\begin{align}
						\label{EQ:Evolutionaverage}\Lunit\average{f}+\gtr\upchi\average{f}=\left(\gtr\upchi-\average{\gtr\upchi}\right)\average{f}+\average{\Lunit f+\gtr\upchi f}
						.\end{align}
				\end{lemma}

				\begin{definition}[Modified mass aspect function] \label{DE:mass}
					We define the mass aspect function $\mass$ to be the following scalar function:
					\begin{align}\label{DE:massfunction}
						\mass:=\uLunit\gtr\upchi+\frac{1}{2}\gtr\upchi\gtr\uchi
						.\end{align}
					
					We now introduce the modified mass aspect function $\modmass$ as below:
			\begin{align}\label{DE:modmass}
						\modmass:=2\anglap\conformalfactor+\uLunit\gtr\upchi+\frac{1}{2}\gtr\upchi\gtr\uchi-\gtr\upchi k_{\spherenormal\spherenormal}+\frac{1}{2}\gtr\upchi\Chfour_{\uLunit}
						.\end{align}

                    Moreover, in $\intregion$, we define $\hodgemass$ to be the $S_{t,u}$-tangent one-form that satisfies the following Hodge system on $S_{t,u}$:
					\begin{align}
						\label{DE:hodgemass}\angdiv\hodgemass&=\frac{1}{2}(\modmass-\umodmass)&
						\angcurl\hodgemass&=0
						.\end{align}
					
					In $\intregion$, we also define the modified torsion $\modtorsion$ as the following $S_{t,u}$-tangent one-form:
					\begin{align}
						\label{DE:defmodtorsion}\modtorsion:=\upzeta+\angnabla\conformalfactor
						.\end{align}
				\end{definition}
				
				In the following two propositions, we list the estimates of the initial foliation. These estimates are crucial for establishing the well-defined geometric setup in Section \ref{SS:opticalfunction} and for controlling the geometry.
				\begin{proposition}\cite[Proposition 9.8]{3DCompressibleEuler}, \cite[Proposition 4.3. Existence and properties of the initial foliation]{AGeoMetricApproach}.\label{PR:InitialCT1} On $\Szero$, there exists a function $w=w(x)$ on the domain $0\leq w\leq w_*=\frac{4}{5}\Trescale$, such that $w(\tip)=0$ and each level set $S_{0,w}$ is diffeomorphic to $\Stwo$, and there hold
					\begin{align}\label{EQ:initialtracecondition}
						\gtr\spheresecondfund+k_{\spherenormal\spherenormal}&=\frac{2}{\lapse w}+\ggtr k-\Chfour_{\Lunit},& \lapse(\tip)&=1
						.\end{align}
					By (\ref{DE:tracechismall}), the facts that $\rgeo(0,-u)=w$ and $\upchi_{AB}=\spheresecondfund_{AB}-k_{AB}$, (\ref{EQ:initialtracecondition}) is equivalent to 
					\begin{align}
						\label{EQ:zinitial}\chismall|_{\Szero}&=\frac{2(1-\lapse)}{\lapse w},& \text{for}\ 0&\leq w\leq w_*
						,\end{align}
					where we define the lapse $a$ as follows:
					\begin{align}\label{DE:Defofa}
						a=\left(\sqrt{(g^{-1})^{cd}\p_cw\p_dw}\right)^{-1}.
					\end{align}
					Note that $\derinormal=\lapse\spherenormal|_{\Sigma_{0}}$.
					Then, on $\initialSzero:=\bigcup\limits_{0\leq w\leq w_*}S_{0,w}$, for $0<1-\frac{2}{q_{*}}<N-3$, the following estimates hold
					\begin{subequations}
						\begin{align}
							\label{ES:lapseminusone}\abs{\lapse-1}&\lesssim\lambda^{-4\varepsilon_0}\leq\frac{1}{4},& \twonorms{w^{-1/2}(\lapse-1)}{w}{\infty}{\sangle}{\infty}{\initialSzero}&\lesssim\lambda^{-1/2},&
							\volume&\approx w^2
							,\end{align}
						\begin{align}
							\label{ES:initialfoliationtwo}\twonorms{w^{\frac{1}{2}-\frac{2}{q_{*}}}(\hat{\spheresecondfund},\angnabla\ln\lapse)}{w}{\infty}{\gsphere}{q_{*}}{\initialSzero}, \twonorms{\angnabla\ln\lapse}{w}{2}{\sangle}{\infty}{\initialSzero}, \twonorms{\hchi}{w}{2}{\sangle}{\infty}{\initialSzero}&\lesssim\lambda^{-1/2} 
							,\end{align}
						\begin{align}
							\label{ES:initialfoliation3}\max\limits_{A,B=1,2}\norm{w^{-2}\gsphere\left(\deriasphere,\deribsphere\right)-\esphere\left(\deriasphere,\deribsphere\right)}_{L^\infty_x(\initialSzero)}&\lesssim\lambda^{-4\varepsilon_0},\\
							\label{ES:initialfoliation4}\max\limits_{A,B,C=1,2}\onenorm{\dericsphere\left(w^{-2}\gsphere\left(\deriasphere,\deribsphere\right)-\esphere\left(\deriasphere,\deribsphere\right)\right)}{\sangle}{q_{*}}{S_{0,w}}&\lesssim\lambda^{-4\varepsilon_0},\\
							\label{ES:in5}\onenorm{w^{\frac{1}{2}-\frac{2}{q_{*}}}\angnabla\cphi}{\gsphere}{q_{*}}{S_{0,w}}&\lesssim\lambda^{-1/2}
							.\end{align}
						In addition, we have
						\begin{align}
							\label{ES:in6}\norm{w\chismall}_{L^\infty_x\left(\initialSzero\right)}&\lesssim\lambda^{-4\varepsilon_0},\\
							\label{ES:in7}\twonorms{w^{3/2}\angnabla\chismall}{w}{\infty}{\sangle}{p}{\initialSzero}+\holdertwonorms{w^{1/2}\chismall}{w}{\infty}{\sangle}{0}{\delta_0}{\initialSzero}&\lesssim\lambda^{-1/2}
							.\end{align}
						Finally, we obtain
						\begin{align}
							\label{ES:pNatSigma0}\sum\limits_{i,j=1,2,3}\abs{w\sphereproject_j^a\p_a\spherenormal^i-\sphereproject_j^i}=\zero(w)\ \text{as} \ w\downarrow 0
							.\end{align}
					\end{subequations}
				\end{proposition}
				\begin{proposition}\cite[Lemma 9.9. Initial conditions on the cone-tip axis tied to the eikonal function]{3DCompressibleEuler}.\label{PR:InitialCT2} On any null cone $\coneu$, which initiate from a point on the time axis $0\leq t=u\leq \Trescale$, the following estimates hold
					\begin{subequations}
						\begin{align}
							\label{ES:initialcontipone}&\gtr\upchi-\frac{2}{\rgeo}, \rgeo\chismall, \sgabs{\hchi}, \abs{\rgeo\sphereproject_j^a\p_a\Lunit^i-\sphereproject_j^i},\nulllapse-1,\sgabs{\upzeta},\conformalfactor,\\
							\notag&\rgeo\sgabs{\angnabla\gtr\upchi},\rgeo^2\sgabs{\angnabla\chismall},\rgeo\sgabs{\angnabla\hchi},\rgeo\sgabs{\angnabla\nulllapse},\rgeo\sgabs{\angnabla\upzeta}, \rgeo\sgabs{\angnabla\conformalfactor},\\
							\notag&\rgeo^2\anglap\nulllapse,\rgeo^2\anglap\conformalfactor,\rgeo^2\mass,\rgeo^2\modmass\\
							\notag&=\zero(\rgeo) \ \text{as}\ t\downarrow u,\\
							&\lim\limits_{t\downarrow u}\norm{\uzeta,k}_{L^{\infty}(\stu)}<\infty
							.\end{align}
						
						Moreover, we have
						\begin{align}
							\label{ES:initialcontip3}\lim\limits_{t\downarrow u} \rgeo^{-2}\gsphere\left(\deriasphere,\deribsphere\right)&=\esphere\left(\deriasphere,\deribsphere\right),\\
							\label{ES:initialcontip4}\lim\limits_{t\downarrow u} \rgeo^{-2}\dericsphere\gsphere\left(\deriasphere,\deribsphere\right)&=\dericsphere\esphere\left(\deriasphere,\deribsphere\right)
							.\end{align}
						
					\end{subequations}
				\end{proposition}
				
				\begin{proof}[Discussion of the proof of Proposition \ref{PR:InitialCT1} and Proposition \ref{PR:InitialCT2}]
					
					The existence of such initial foliation can be proved by Nash-Moser implicit function theorem (see \cite{szeftel2012parametrix1}). The proof of the estimates in Proposition \ref{PR:InitialCT1} relies on the energy estimates (\ref{ES:energyestimates2}) and (\ref{ES:linearcurl}) due to the rescaling (see Section \ref{SS:SelfRescaling}). We refer the reader to \cite[Appendix C]{Wang2014} for the proof of the estimates in Proposition \ref{PR:InitialCT1} and in Proposition \ref{PR:InitialCT2}. We emphasize that the proofs are exactly the same since we have the same estimates\footnote{In \cite{Wang2014}, where $2<N<5/2$ and the norm in \eqref{ES:pfourgSigma0} is $\sobolevnorm{\pfour^2\gfour}{N-2}{\Sigma_{0}}$ instead.} for the metric, that is, after rescaling\footnote{Here, $\Sigma_{0}$ after rescaling corresponds to $\Sigma_{t_k}$ for some $k$ in the original setting.}, we obtain
					\begin{align}\label{ES:pfourgSigma0}
						\sobolevnorm{\pfour^2\gfour}{N-3}{\Sigma_{0}}\lesssim\lambda^{-\frac{1}{2}}.
					\end{align}
				\end{proof}
				
				\subsection{Null Structure Equations}\label{SS:PDES}
				In order to control the geometry, we need to derive the estimates for various geometric quantities. We provide the transport-Hodge type of PDEs in Lemma \ref{LE:PDE1}. Then we use the elastic wave equations (\ref{EQ:boxg4}) to express the curvature terms on the right hand side of these PDEs in the following section.
				\begin{lemma}\cite[Section 5.4, PDEs verified by connection coefficients]{ImprovedLocalwellPosedness}.\label{LE:PDE1}
					For the Ricci coefficients defined in Definition \ref{DE:DEFSOFCONNECTIONCOEFFICIENTS}, the following evolution equations hold:\\
					
					\noindent\textbf{Equations for $\upchi$:}
					\begin{subequations}
						\begin{align}
							\label{EQ:Ltraceupchi}\Lunit\gtr\upchi+\frac{1}{2}(\gtr\upchi)^2=&-\sgabs{\hchi}^2-k_{\spherenormal\spherenormal}\gtr\upchi-\Ricfour{\Lunit}{\Lunit},\\
							\label{EQ:uLunitgtrupchi}
							\uLunit\gtr\upchi+\frac{1}{2}\gtr\uchi\gtr\upchi=&2\angdiv\upzeta+k_{\spherenormal\spherenormal}\gtr\upchi-\hchi_{AB}\huchi_{AB}+2\sgabs{\upzeta}^2+\Riemfour{A}{\Lunit}{\uLunit}{A},\\
							\label{EQ:angDLunithchi}\angD_\Lunit\hchi_{AB}+(\gtr\upchi)\hchi_{AB}=&-k_{\spherenormal\spherenormal}\hchi_{AB}-\left(\Riemfour{\Lunit}{A}{\Lunit}{B}-\frac{1}{2}\Ricfour{\Lunit}{\Lunit}\delta_{AB}\right),\\
							\label{EQ:angdivhchi}\angdiv\hchi_A+\hchi_{AB}k_{B\spherenormal}=&\frac{1}{2}(\angnabla_A\gtr\upchi+k_{A\spherenormal}\gtr\upchi)+\Riemfour{B}{\Lunit}{B}{A},\\
							\label{EQ:angDuLunithchi}\angD_{\uLunit}\hchi_{AB}+\frac{1}{2}(\gtr\uchi)\hchi_{AB}=&-\frac{1}{2}(\gtr\upchi)\huchi_{AB}+2\angnabla_A\upzeta_B\\
							\notag&-\angdiv\upzeta\delta_{AB}+k_{\spherenormal\spherenormal}\hchi_{AB}+\left(2\upzeta_A\upzeta_B-\sgabs{\upzeta}^2\delta_{AB}\right)\\
							\notag&-\left(\huchi_{AC}\hchi_{BC}-\frac{1}{2}\huchi_{CD}\hchi_{CD}\delta_{AB}\right)+\Riemfour{B}{\Lunit}{\uLunit}{A}-\frac{1}{2}\Riemfour{C}{\Lunit}{\uLunit}{C}\delta_{AB}.
						\end{align}
					\end{subequations}
					
					\noindent\textbf{Equations for $\uchi$:}
					\begin{subequations}
						\begin{align}
							\label{Ltraceuchi}\Lunit\gtr\uchi+\frac{1}{2}(\gtr\upchi)\gtr\uchi=&2\angdiv\uzeta+k_{\spherenormal\spherenormal}\gtr\uchi-\hchi_{AB}\huchi_{AB}+2\sgabs{\uzeta}^2+\Riemfour{A}{\uLunit}{\Lunit}{A},\\
							\label{angDLunithuchi}\angD_{\Lunit}\huchi_{AB}+\frac{1}{2}(\gtr\upchi)\huchi_{AB}=&-\frac{1}{2}(\gtr\uchi)\hchi_{AB}+2\angnabla_A\uzeta_B\\
							\notag&-\angdiv\uzeta\delta_{AB}+k_{\spherenormal\spherenormal}\huchi_{AB}+\left(2\uzeta_A\uzeta_B-\sgabs{\uzeta}^2\delta_{AB}\right)\\
							&\notag-\left(\hchi_{AC}\huchi_{BC}-\frac{1}{2}\huchi_{CD}\hchi_{CD}\delta_{AB}\right)+\Riemfour{B}{\uLunit}{\Lunit}{A}-\frac{1}{2}\Riemfour{C}{\uLunit}{\Lunit}{C}\delta_{AB}.
						\end{align}
					\end{subequations}
					\noindent\textbf{Equations for $\upzeta$:}
					\begin{subequations}
						\begin{align}
							\label{EQ:angDLupzetaA}\angD_\Lunit\upzeta_A+\frac{1}{2}(\gtr\upchi)\upzeta_A=&-(k_{B\spherenormal}+\upzeta_B)\hchi_{AB}-\frac{1}{2}\gtr\upchi k_{A\spherenormal}+\frac{1}{2}\Riemfour{A}{\Lunit}{\Lunit}{\uLunit},\\
							\label{EQ:angdivupzet}\angdiv\upzeta=&\frac{1}{2}(\upmu-k_{\spherenormal\spherenormal}\gtr\upchi-2\sgabs{\upzeta}^2-\sgabs{\hchi}^2-2k_{AB}\hchi_{AB})-\frac{1}{2}\Riemfour{A}{\uLunit}{\Lunit}{A},\\
							\label{EQ:angcurlupzet}\angcurl\upzeta=&\frac{1}{2}\antisymmetic^{AB}\huchi_{AC}\hchi_{BC}-\frac{1}{2}\antisymmetic^{AB}\Riemfour{B}{\Lunit}{\uLunit}{A}
							,\end{align}
					\end{subequations}
					where $\mass$ is defined in Definition \ref{DE:mass}.
					
					\begin{remark}
						The equations in Lemma \ref{LE:PDE1} are a consequence of geometry, where elastic wave equations play no role. The motivation of deriving those equations is to prove:
						\begin{align}
							(\Lunit,\angnabla)\pfour(\Lunit,\uLunit,e_A)\lesssim(\Lunit,\angnabla)\pfour\gfour,
						\end{align}
						so that we can control geometry by using the energy estimates (\ref{ES:linearcurl}), (\ref{ES:energyestimates2}) and (\ref{ES:coneenergy}).
					\end{remark}
				\end{lemma}
				
				\subsection{Curvature Decompositions}\label{SS:Curvature}
				In this section, we decompose the components of the Riemann and Ricci curvature tensor with the help of the geometric structure equations and the elastic wave equations. Then, in Proposition \ref{PR:mainversionPDE}, we rewrite the PDEs from Lemma \ref{LE:PDE1} into the version that will be used in the analysis.
				\begin{lemma}\cite[Lemma 9.5. Identities for the derivatives of some scalar functions]{3DCompressibleEuler}.\label{LE:Scalarfunctions} With $\CC,\ACC$ define in Definition \ref{DE:conformalchangemetric}, and $\diff f$ denoting the space gradient of the scalar function $f$, we have the following identities
					\begin{subequations}
						\begin{align}
							\sphereproject\cdot \diff(\vLunit,\uvLunit,\spherenormal)&=\lgensmoothfunction\cdot(\ACC,\rgeo^{-1}),&
							\sphereproject\cdot \diff\lgensmoothfunction&=\lgensmoothfunction\cdot(\ACC,\rgeo^{-1}),\\
							\diff(\vLunit,\uvLunit,\spherenormal)&=\lgensmoothfunction\cdot(\ACC,\rgeo^{-1}),&
							\diff\lgensmoothfunction&=\lgensmoothfunction\cdot(\ACC,\rgeo^{-1})
							.\end{align}
					\end{subequations}
					Here, $\sphereproject$ denotes the orthogonal projection onto $\stu$ (see Definition \ref{D:stutensor fields}), $\lgensmoothfunction$ are polynomials of Cartesian components of $\Lunit$ with smooth functions of $\vvariables,\cvvariables$ as coefficients. We will use the same notation $\lgensmoothfunction$ for the rest of the article.
				\end{lemma}
				In the proof of Lemma \ref{LE:Scalarfunctions}, we (implicitly) use (\ref{DE:DEFSOFCONNECTIONCOEFFICIENTS}) whenever necessary.
				\begin{proof}[Proof of Lemma \ref{LE:Scalarfunctions}]
					Notice that $\Timelike=\p_t$ and $(k_{ij},\Chfour)=\gensmoothfunction\cdot(\pfour\vvariables,\pfour\cvvariables)$. Also, by definition (see Definition \ref{DE:Nullframe}), we have $(\spherenormal,\uvLunit)=\lgensmoothfunction$ and $\Chfour_{\Lunit}=\lgensmoothfunction\cdot(\pfour\vvariables,\pfour\cvvariables)$. Hence, we derive:
					\begin{align}
						\Dfour_\beta\Lunit^\alpha&=\p_\beta\Lunit^\alpha+\Gamma_{\beta\gamma}^\alpha\Lunit^\gamma,\\
						\Dfour\Lunit^\alpha&=\langle\Dfour\Lunit^\alpha,e_A\rangle e_A-\frac{1}{2}\langle\Dfour\Lunit^\alpha,\Lunit\rangle\uLunit-\frac{1}{2}\langle\Dfour\Lunit^\alpha,\uLunit\rangle\Lunit\\
						\notag&=\left(\Dfour_A\Lunit^\alpha\right) e_A-\left(\frac{1}{2}\Dfour_\Lunit\Lunit^\alpha\right)\uLunit-\left(\frac{1}{2}\Dfour_{\uLunit}\Lunit^\alpha\right)\Lunit\\
						\notag&=\left(\upchi_{AB}e_B^\alpha-k_{A\spherenormal}\Lunit^\alpha\right)e_A+\left(\frac{1}{2}k_{\spherenormal\spherenormal}\Lunit^\alpha\right)\uLunit-\left(\upzeta_Ae_A^\alpha+\frac{1}{2}k_{\spherenormal\spherenormal}\Lunit^\alpha\right)\Lunit\\
						\notag&=\lgensmoothfunction\cdot\left(\pfour\vvariables,\pfour\cvvariables,\chismall,\hchi,\upzeta,\rgeo^{-1}\right)
						,\end{align}
					where we use the fact that
					\begin{align}
						\gfour(\Dfour\Lunit^\alpha,X)=(\gfour^{-1})^{\beta\gamma}\Dfour_\beta\Lunit^\alpha X_\gamma=\Dfour_\beta\Lunit^\alpha X^\beta=\Dfour_X\Lunit^\alpha
						.\end{align}
					
					It follows obviously from above that:
					\begin{align}
						\sphereproject\cdot d\vLunit=\lgensmoothfunction\cdot\left(\pfour\vvariables,\pfour\cvvariables,\chismall,\hchi,\rgeo^{-1}\right)
						.\end{align}
					
					We then obtain the identities that involve $\spherenormal$ and $\uvLunit$ by using chain rule and the fact that $(\spherenormal,\uvLunit)=\lgensmoothfunction$.
				\end{proof}

				\begin{lemma}\cite[Lemma 2.1. Ricci curvature component decompositions]{ImprovedLocalwellPosedness} We decompose the Ricci curvature tensor components as follows, with $\gfour_{\alpha\beta}$ viewed as a scalar function\label{LE:Ricci}:
					\begin{align}\label{EQ:Riccicurvaturedecomposition}
						\Ricfour{\alpha}{\beta}=-\frac{1}{2}\boxg \gfour_{\alpha\beta}(\vvariables,\cvvariables)+\frac{1}{2}\left(\Dfour_\alpha\Chfour_\beta+\Dfour_\beta\Chfour_\alpha\right)+\gensmoothfunction(\vvariables,\cvvariables)\cdot\pfour\gfour\cdot\pfour\gfour,
					\end{align}
					where $\Chfour_{\alpha}:=\Chfour_{\alpha\kappa\lambda}\invgfour^{\kappa\lambda}$.
				\end{lemma}
				\begin{corollary}\cite{ImprovedLocalwellPosedness,AGeoMetricApproach,3DCompressibleEuler,Yuthesis}.\label{CO:Specialricci}
					Recall (\ref{DE:conformalchangemetric}) for $\ACC$. Let $\xi$ denote schematically the $\stu$-tangent one-forms and symmetric $\binom{0}{2}$ $\stu$-tangent tensors with the property\footnote{We use the symbol $\xi$ in the same way in the rest of the article.} $\xi=\lgensmoothfunction\cdot(\pfour\vvariables,\pfour\cvvariables)$. We note that $\vvariables,\cvvariables$ are the rescaled solution variables defined in Definition \ref{DE:rescaledquantities}. Then, with the help of (\ref{EQ:Riccicurvaturedecomposition}), we can decompose the Ricci curvature tensor components as below:
					
					\begin{subequations}
						\begin{align}
							\label{EQ:riccilunitlunit}\Ricfour{\Lunit}{\Lunit}&=\lgensmoothfunction\cdot\boxg \gfour(\vvariables,\cvvariables)+\Lunit\Chfour_\Lunit+k_{\spherenormal\spherenormal}\Chfour_{\Lunit}+\lgensmoothfunction\cdot\pfour\gfour\cdot\pfour\gfour,\\
							\Ricfour{\Lunit}{\uLunit}&=\lgensmoothfunction\cdot\boxg \gfour(\vvariables,\cvvariables)+\frac{1}{2}(\Lunit\Chfour_{\uLunit}+\uLunit\Chfour_{\Lunit})+\lgensmoothfunction\cdot\ACC\cdot\pfour\gfour,\\
							\label{EQ:ricciLA}\Ricfour{\Lunit}{A}&=\lgensmoothfunction\cdot\boxg \gfour(\vvariables,\cvvariables)+\angnabla_A\xi^{(1)}+\angD_{\Lunit}\xi_A^{(2)}+\lgensmoothfunction\cdot(\ACC,\rgeo^{-1})\cdot\pfour\gfour,\\
							\label{EQ：ricciAB}\Ricfour{A}{B}&=\lgensmoothfunction\cdot\boxg \gfour(\vvariables,\cvvariables)+\angnabla_A\xi_B+\lgensmoothfunction\cdot\pfour\gfour\cdot\pfour\gfour
							.\end{align}
					\end{subequations}
				\end{corollary}
				
				\begin{proof}[Discussion of the proof of Lemma \ref{LE:Ricci} and Corollary \ref{CO:Specialricci}]
					The metric $\gfour$ shares similar properties as in \cite{ImprovedLocalwellPosedness}, except that $\gfour=\gfour(\vvariables,\cvvariables)$ in our work, compared to $\gfour=\gfour(\vvariables)$ in \cite{ImprovedLocalwellPosedness}.
					
					In the proof of the identities in Corollary \ref{CO:Specialricci}, we contract (\ref{EQ:Riccicurvaturedecomposition}) with the null frame.
				\end{proof}
				
				\begin{corollary}\cite[Lemma 5.12. Riemann curvature component decomposition]{AGeoMetricApproach}.\label{CO:Specialriemann} Let $\ACC$ be defined as in Definition \ref{DE:conformalchangemetric}. With the help of Corollary \ref{CO:Specialricci}, we can decompose the Riemann curvature tensor components as follows:
					\begin{subequations}
						\begin{align}
							\label{EQ:rLALB}\Riemfour{\Lunit}{A}{\Lunit}{B}&=(\angnabla,\angD_{\Lunit})\xi+\lgensmoothfunction\cdot(\ACC,\rgeo^{-1})\cdot\pfour\gfour,\\
							\Riemfour{C}{A}{\Lunit}{B}&=\angnabla\xi+\lgensmoothfunction\cdot(\ACC,\rgeo^{-1})\cdot\pfour\gfour,\\
							\label{EQ:riemannABAB}\Riemfour{A}{B}{A}{B}&=\angdiv\xi+\lgensmoothfunction\cdot(\ACC,\rgeo^{-1})\cdot\pfour\gfour,\\
							\label{EQ:riemannABBL}\Riemfour{A}{B}{B}{\Lunit}&=\angnabla_A\xi+\angdiv\xi_A+\lgensmoothfunction\cdot(\ACC,\rgeo^{-1})\cdot\pfour\gfour,\\
							\label{EQ:rALLuLunit}\Riemfour{A}{\Lunit}{\Lunit}{\uLunit}&=(\angnabla,\angD_{\Lunit})\xi+\lgensmoothfunction\cdot\boxg \gfour(\vvariables,\cvvariables)+\lgensmoothfunction\cdot(\ACC,\rgeo^{-1})\cdot\pfour\gfour,\\
							\label{EQ:ALuLA}\Riemfour{A}{\Lunit}{\uLunit}{A}&=\angdiv\xi+\lgensmoothfunction\cdot\boxg \gfour(\vvariables,\cvvariables)+\lgensmoothfunction\cdot(\ACC,\rgeo^{-1})\cdot\pfour\gfour,\\
							\label{EQ:antisym riem curv}
							\antisymmetic^{AB}\Riemfour{A}{\Lunit}{\uLunit}{B}&=\angcurl\xi+\lgensmoothfunction\cdot(\ACC,\rgeo^{-1})\cdot\pfour\gfour
							.\end{align}
					\end{subequations}
				\end{corollary}
				
				\subsection{PDEs Verified by the Geometric Quantities}\label{SS:mainversionPDES}
				{Now we state the main equations for the geometric quantities.}
            \begin{proposition}\cite[Section 9.9.3. PDEs verified by the modified geometric quantities]{3DCompressibleEuler}. Recall that the vector fields $\Lunit,\uLunit,\spherenormal,e_A$ and null lapse $\nulllapse$ in Section \ref{SS:geometricquantities}. $\sangle^A$ is the local coordinates in $\Stwo$ (defined in Definition \ref{DE:norms}). Spherevolume ratio $\volume$ is defined in (\ref{DE:spherevolume}). Connection coefficients $\upxi,\uchi,\upzeta$ are defined in Definition \ref{DE:DEFSOFCONNECTIONCOEFFICIENTS}. Conformal factor $\conformalfactor$ and conformally modified connection coefficients $\reschi,\resuchi,\chismall,\CC,\ACC$ are\footnote{Note that in the rest of the article, for convenience, we prefer notation $\ACC,\CC$ to $(\pfour\vvariables,\pfour\cvvariables,\reschi,\resuchi,\chismall)$. Therefore, we may still use $\ACC,\CC$ even if not every term appears; for example, we may denote $(\pfour\vvariables,\pfour\cvvariables,\reschi)$ by $\ACC$, {which does not affect our analysis}.} defined in Definition \ref{DE:conformalchangemetric}. Mass aspect function $\mass$ and its modified version $\modmass,\hodgemass$ and modified torsion $\modtorsion$ are defined in Definition \ref{DE:mass}. The following identities hold where the terms on the right hand sides are displayed schematically:\label{PR:mainversionPDE}\\
					\textbf{Transport equations involving the Cartesian components $\Lunit^i$ and $\spherenormal^i$:}
					\begin{subequations}
						\begin{align}
							\label{EQ:LunitLunit}\Lunit\Lunit^i&=\lgensmoothfunction\cdot\pfour\gfour,&
							\Lunit\spherenormal^i&=\lgensmoothfunction\cdot\pfour\gfour
							.\end{align}
						Moreover, along $\Szero$(where $w=\rgeo=-u$ and $a=b$, see Proposition \ref{PR:InitialCT1}), we have:
						\begin{align}
							\label{EQ:derinormalLunit}\derinormal\Lunit^i&=\lapse\cdot\lgensmoothfunction\cdot\pfour\gfour-\angnabla\lapse,&
							\derinormal\spherenormal^i&=\lapse\cdot\lgensmoothfunction\cdot\pfour\gfour-\angnabla \lapse
							.\end{align}
					\end{subequations}
					\textbf{Transport equations involving the Cartesian components $\cartiasphere^i$:}
					\begin{align}
						\cartiasphere^i:=\frac{1}{\rgeo}\left(\deriasphere \right)^i
						.\end{align}
					Then, the following evolution equations hold in $\region$:
					\begin{align}
						\label{EQ:Lunitcartiasphere}\Lunit\cartiasphere^i&=\lgensmoothfunction\cdot\ACC\cdot\vcartiasphere,\\
						\label{EQ:NomegaA}\spherenormal\left(\deriasphere\right)^i&=\lgensmoothfunction\cdot\pfour\gfour\cdot\deriasphere+\angnabla_A\ln\nulllapse\cdot\spherenormal
						.\end{align}
					Moreover, along $\Szero$, we have
					\begin{align}
						\label{EQ:derinormalcartiasphere}\derinormal\cartiasphere^i=\frac{(1-a)}{w}\cartiasphere^i+\lapse\cdot\lgensmoothfunction\cdot\ACC\cdot\vcartiasphere+\lgensmoothfunction\cdot\angnabla \lapse\cdot\vcartiasphere
						.\end{align}
					\textbf{Transport equations connected to $\chismall$:}
					\begin{subequations}
						\begin{align}
							\label{EQ:Lz}\Lunit\chismall+\frac{2}{\rgeo}\chismall=&\lgensmoothfunction\cdot\boxg\gfour+\lgensmoothfunction\cdot(\chismall,\hchi)\cdot(\chismall,\hchi,\pfour\gfour)+\lgensmoothfunction\cdot\rgeo^{-1}\cdot\pfour\gfour,\\
							\label{EQ:Lunitangnablaz}\angD_{\Lunit}\angnabla\chismall+\frac{3}{\rgeo}\angnabla\chismall=&\lgensmoothfunction\cdot\angnabla\boxg\gfour+\lgensmoothfunction\cdot\boxg\gfour\cdot(\ACC,\rgeo^{-1})\\
							\notag&+\lgensmoothfunction\cdot\angnabla\pfour\gfour\cdot(\ACC,\rgeo^{-1})+\lgensmoothfunction\cdot(\angnabla\chismall,\angnabla\hchi)\cdot\ACC+\lgensmoothfunction\cdot(\ACC,\rgeo^{-1})\cdot(\ACC,\rgeo^{-1})\cdot\pfour\gfour,\\
							\label{EQ:Lv}\Lunit\left(\frac{1}{2}\gtr\reschi \volume\right)&-\frac{1}{4}(\gtr\upchi)^2v+\frac{1}{2}\left(\Lunit\ln\nulllapse\right)\gtr\reschi \volume-\sgabs{\angnabla\conformalfactor}^2\volume\\
							\notag=&\lgensmoothfunction\cdot\boxg\gfour\cdot\volume+\lgensmoothfunction\cdot\ACC\cdot\ACC\cdot \volume-\sgabs{\angnabla\conformalfactor}^2\volume
							.\end{align}
					\end{subequations}
					\textbf{PDEs involving $\hchi$:}
					\begin{subequations}
						\begin{align}
							\label{EQ:angdvhchi}\angdiv\hchi&=\angnabla\chismall+\angdiv\xi^{(1)}+\angnabla\xi^{(2)}+\lgensmoothfunction\cdot(\ACC,\rgeo^{-1})\cdot\pfour\gfour,\\
							\label{EQ:angDhchi}\angD_{\Lunit}\hchi+(\gtr\upchi)\hchi&=(\angnabla,\angD_{\Lunit})\xi+\lgensmoothfunction\cdot\boxg\gfour+\lgensmoothfunction\cdot(\ACC,\rgeo^{-1})\cdot\pfour\gfour
							.\end{align}
					\end{subequations}
					\textbf{The transport equation for $\upzeta$:}
					\begin{align}
						\label{EQ:angDLunitupzeta}\angD_{\Lunit}\upzeta_A+\frac{1}{2}(\gtr\upchi)\upzeta=&(\angnabla,\angD_{\Lunit})\xi+\lgensmoothfunction\cdot\boxg\gfour+\lgensmoothfunction\cdot\left(\ACC,\rgeo^{-1}\right)\cdot\pfour\gfour+\lgensmoothfunction\cdot\upzeta\cdot\hchi
						.\end{align}
					\textbf{The transport equation for $\nulllapse$:}
					\begin{align}
						\label{EQ:Lunitnulllapse}\Lunit\nulllapse=\nulllapse\cdot\lgensmoothfunction\cdot\pfour\gfour
						.\end{align}
					\textbf{Transport equation for $\gsphere$:} Along the integral curves of $\Lunit$, which are parameterized  by $t$, we have that, with $\esphere$ denoting the standard round metric on the Euclidean unit sphere $\mathbb{S}^2$, the following identity:
					\begin{align}
						\label{EQ:dtrminus2gsphere}\tderivative&\left\{\rgeo^{-2}\gsphere\left(\deriasphere,\deribsphere\right)-\esphere\left(\deriasphere,\deribsphere\right)\right\}\\
						\notag=&\left(\chismall-\Chfour_{\Lunit}\right)\left\{\rgeo^{-2}\gsphere\left(\deriasphere,\deribsphere\right)-\esphere\left(\deriasphere,\deribsphere\right)\right\}+\left(\chismall-\Chfour_{\Lunit}\right)\esphere\left(\deriasphere,\deribsphere\right)+\frac{2}{\rgeo^2}\hchi\left(\deriasphere,\deribsphere\right)
						.\end{align}
					\textbf{Transport equations for $\volume$ and $\angnabla \volume$:}
					\begin{subequations}
						\begin{align}
							\label{EQ:Lunitcphi}\Lunit\cphi&=\gtr\upchi-\frac{2}{\rgeo}=\chismall-\Chfour_{\Lunit},\\
							\label{EQ:Lunitangnablacphi}\Lunit\angnabla\cphi+\frac{1}{2}(\gtr\upchi)\angnabla\cphi&=\lgensmoothfunction\cdot\hchi\cdot\angnabla\cphi+\angnabla\chismall-\angnabla(\Chfour_{\Lunit})
							.\end{align}
					\end{subequations}
					\textbf{An algebraic identity for $\mass$:} The mass aspect function $\mass$ verifies the following identity:
					\begin{align}
						\label{EQ:mass}\mass=&\lgensmoothfunction\cdot\boxg\gfour+\angdiv\xi+\lgensmoothfunction\cdot\hchi\cdot\hchi+\lgensmoothfunction\cdot\angnabla\cphi\cdot\left(\pfour\gfour,\upzeta\right)+\lgensmoothfunction\cdot\pfour\gfour\cdot\left(\ACC,\rgeo^{-1}\right)
						.\end{align}
					\textbf{The transport equation for $\angnabla\conformalfactor$:} In the interior region $\intregion$, $\angnabla\conformalfactor$ satisfies the following transport equation:
					\begin{align}
						\label{EQ:Lunitangnablasigma}\Lunit\angnabla\conformalfactor+\frac{1}{2}\gtr\upchi\angnabla\conformalfactor=\frac{1}{2}\angnabla(\Chfour_{\Lunit})-\hchi\cdot\angnabla\conformalfactor
						.\end{align}
					\textbf{The transport equation for $\modmass$:} The modified mass aspect function $\modmass$ verifies the following transport equation:
					\begin{align}
						\label{EQ:Lunitmodmass}\Lunit\modmass+(\gtr\upchi)\modmass=\mathfrak{J}_{(1)}+\mathfrak{J}_{(2)},
					\end{align}
					where
					\begin{subequations}\label{DE:J1J2}
						\begin{align}
							\mathfrak{J}_{(1)}=&\rgeo^{-1}\angdiv\xi+\rgeo^{-2}\xi,\\
							\label{DE:J2}\mathfrak{J}_{(2)}=&\cdot\lgensmoothfunction\cdot\pfour\boxg\gfour+\lgensmoothfunction\cdot\boxg\gfour\cdot(\ACC,\rgeo^{-1})\\
							\notag&+\lgensmoothfunction\cdot\ACC\cdot(\angnabla\modtorsion,\angnabla\chismall,\pfour^2\gfour)+\lgensmoothfunction\cdot\left(\angnabla\pfour\gfour,\angnabla\chismall\right)\cdot\angnabla\conformalfactor\\
							\notag&+\lgensmoothfunction\cdot\left(\ACC,\rgeo^{-1}\right)\cdot\pfour\gfour\cdot\angnabla\conformalfactor+\lgensmoothfunction\cdot\ACC\cdot\ACC\cdot\left(\ACC,\rgeo^{-1}\right)
							.\end{align}
					\end{subequations}
					\textbf{The Hodge system for $\upzeta$:} The torsion $\upzeta$ satisfies the following Hodge system on $\stu$:
					\begin{subequations}
						\begin{align}
							\label{EQ:angdivupzeta}\angdiv\upzeta=&\angdiv\xi+\lgensmoothfunction\cdot\boxg\gfour+\lgensmoothfunction\cdot\left(\upzeta,\hchi\right)\cdot\left(\upzeta,\hchi\right)+\lgensmoothfunction\cdot\angnabla\cphi\cdot\left(\pfour\gfour,\upzeta\right)+\lgensmoothfunction\cdot\pfour\gfour\cdot\left(\ACC,\rgeo^{-1}\right),\\
							\label{EQ:angcurlupzeta}\angcurl\upzeta=&\angcurl\xi+\lgensmoothfunction\cdot\hchi\cdot\hchi+\lgensmoothfunction\cdot\pfour\gfour\cdot\left(\ACC,\rgeo^{-1}\right)
							.\end{align}
					\end{subequations}
					\textbf{The Hodge system for $\modtorsion$:} The modified torsion $\modtorsion$ verifies the following Hodge system on $\stu$:
					\begin{subequations}
						\begin{align}
							\label{EQ:angdivmodtorsion}\angdiv\modtorsion-\frac{1}{2}\modmass=&\angdiv\xi+\lgensmoothfunction\cdot\boxg\gfour+\lgensmoothfunction\cdot\left(\upzeta,\hchi\right)\cdot\left(\upzeta,\hchi\right)+\lgensmoothfunction\cdot\pfour\gfour\cdot\left(\ACC,\rgeo^{-1}\right),\\
							\label{EQ:angcurlmodtorsion}\angcurl\modtorsion=&\angcurl\xi+\lgensmoothfunction\cdot\hchi\cdot\hchi+\lgensmoothfunction\cdot\pfour\gfour\cdot\left(\ACC,\rgeo^{-1}\right)
							.\end{align}
					\end{subequations}
					\textbf{The Hodge system for $\modtorsion-\hodgemass$:} The difference $\modtorsion-\hodgemass$ satisfies the following Hodge system on $S_{t,u}$
					\footnote{The barred terms appear in the equation (\ref{EQ:angdivmodtorsion}) because $\angdiv(\modtorsion-\hodgemass)$ must have vanishing average.}:
					\begin{subequations}
						\begin{align}
							\label{EQ:angdivhodgemass}\angdiv\left(\modtorsion-\hodgemass\right)=&\angdiv\xi+\left(\lgensmoothfunction\cdot\boxg\gfour-\average{\lgensmoothfunction\cdot\boxg\gfour}\right)\\
							\notag&+\left(\lgensmoothfunction\cdot(\upzeta,\hchi)\cdot(\upzeta,\hchi)-\average{\lgensmoothfunction\cdot(\upzeta,\hchi)\cdot(\upzeta,\hchi)}\right)\\
							\notag&+\left(\lgensmoothfunction\cdot\pfour\gfour\cdot(\ACC,\rgeo^{-1})-\average{\lgensmoothfunction\cdot\pfour\gfour\cdot(\ACC,\rgeo^{-1})}\right),\\
							\label{EQ:angcurlhodgemass}\angcurl\left(\modtorsion-\hodgemass\right)=&\angcurl\xi+\lgensmoothfunction\cdot\hchi\cdot\hchi+\lgensmoothfunction\cdot\pfour\gfour\cdot\left(\ACC,\rgeo^{-1}\right)
							.\end{align}
					\end{subequations}
					\textbf{A decomposition of $\hodgemass$ and a Hodge-transport system for the constituent parts:} In $\intregion$, we can decompose $\hodgemass$ as follows:
					\begin{align}
						\label{EQ:hodgemudecomposition}\hodgemass=\massone+\masstwo
						\end{align}
					where $\massone$ and $\masstwo$ verify the following Hodge-transport PDE systems:\
					\begin{subequations}\label{EQ:div1}
						\begin{align}
							\label{EQ:angdiv1}\angdiv\left(\angD_{\Lunit}\massone+\frac{1}{2}\gtr\upchi\massone\right)&=\mathfrak{J}_{(1)}-\average{\mathfrak{J}_{(1)}},\\
							\label{EQ:angcurl1}\angcurl\left(\angD_{\Lunit}\massone+\frac{1}{2}\gtr\upchi\massone\right)&=0
							,\end{align}
					\end{subequations}
					\begin{subequations}
						\begin{align}
							\label{EQ:angdiv2}\angdiv\left(\angD_{\Lunit}\masstwo+\frac{1}{2}\gtr\upchi\masstwo\right)=&\mathfrak{J}_{(2)}-\average{\mathfrak{J}_{(2)}}+\hchi\cdot\angnabla\hodgemass+\left(\angnabla\pfour\gfour,\angnabla\chismall\right)\cdot\hodgemass\\
							\notag&+\left(\ACC,\rgeo^{-1}\right)\cdot\pfour\gfour\cdot\hodgemass+\left(\gtr\upchi-\average{\gtr\upchi}\right)\average{\modmass},\\
							\label{EQ:angcurl2}\angcurl\left(\angD_{\Lunit}\masstwo+\frac{1}{2}\gtr\upchi\masstwo\right)=&\hchi\cdot\angnabla\hodgemass+\left(\angnabla\pfour\gfour,\angnabla\chismall\right)\cdot\hodgemass\\
							\notag&+\left(\ACC,\rgeo^{-1}\right)\cdot\pfour\gfour\cdot\hodgemass,
						\end{align}
					\end{subequations}
					which are subject to the following initial conditions along the cone-tip axis for $u\in[0,\Trescale]$:
					\begin{align}
						\label{EQ:initialhodgemass}\lim\limits_{t\downarrow u}\rgeo\sgabs{\massone-\hodgemass}&=\zero(\rgeo),&
						\lim\limits_{t\downarrow u}\rgeo\sgabs{\masstwo}&=\zero(\rgeo)
						.\end{align}
				\end{proposition}
				\subsection{Proof of Proposition \ref{PR:mainversionPDE}}\label{SS:proofofmainversion}
                {First, we have the following lemma for commutators.}
				\begin{lemma}\cite[Lemma 2.2, Lemma 2.3. Commutator formulas]{ImprovedLocalwellPosedness}.\label{LE:Commutatorformulas} Let $\xi_A$ be a $S_{t,u}$-tangent covector and $f$ be a scalar, then there hold
					\begin{subequations}
						\begin{align}
							\label{EQ:commutatorformula}[\Lunit,\angnabla_A]f&=-\upchi_{AB}\angnabla_Bf,\\
							\label{EQ:commutatorformulaLuL}[\Lunit,\uLunit]&=2\left(\uzeta_A-\upzeta_A\right)e_A+k_{\spherenormal\spherenormal}\uLunit-k_{\spherenormal\spherenormal}\Lunit,\\
							\label{EQ:commutatorformula1}[\spherenormal,\angnabla_A]f&=-\spheresecondfund_{AB}e_Bf+\angnabla_A\ln\nulllapse\cdot\spherenormal f,\\
							\label{EQ:commutatorformulafortensor}\angD_{\Lunit}\angnabla_B\xi_A-\angnabla_B\angD_{\Lunit}\xi_A&=-\upchi_{AB}k_{C\spherenormal}\xi_C+\upchi_{BC}k_{A\spherenormal}\xi_C-\upchi_{BC}\angnabla_C\xi_A+\Riemfour{A}{C}{\Lunit}{B}\xi_C
							.\end{align}
					\end{subequations}
				\end{lemma}
				
				Now, we prove Proposition \ref{PR:mainversionPDE}. Note that we silently use Lemma \ref{LE:Scalarfunctions}, the decomposition of $\chi,\uchi$ (\ref{EQ:chiuchi}) and $\gtr\upchi$ (\ref{DE:tracechismall}) throughout the proofs.
				\begin{proof}[Proof of $\Lunit\Lunit^i=\lgensmoothfunction\cdot\pfour\vvariables$ in (\ref{EQ:LunitLunit})] By (\ref{EQ:CC}) and the definition of $k$ (\ref{DE:SECONDFUNDOFSIGMAT}), we derive:
					\begin{align}
						\Lunit\Lunit^i=\Dfour_\Lunit\Lunit^i-\Chfour_{\Lunit\Lunit}^i=-k_{\spherenormal\spherenormal}\Lunit^i-\Chfour_{\Lunit\Lunit}^i=\lgensmoothfunction\cdot\pfour\gfour
					\end{align}
					with the notation $\Chfour_{\Lunit\Lunit}^i:=\Chfour_{\alpha\beta}^i \Lunit^\alpha \Lunit^\beta$. And in the following, we will similarly define $\Chfour_{XY}^i$ for other vectors $X$ and $Y$.
				\end{proof}
				\begin{proof}[Proof of $\Lunit\spherenormal^i=\lgensmoothfunction\cdot\pfour\gfour$ in (\ref{EQ:LunitLunit})] By (\ref{EQ:CC}) and the definition of $k$ (\ref{DE:SECONDFUNDOFSIGMAT}), we have:
					\begin{subequations}
						\begin{align}
							\Lunit\uLunit^i&=\Dfour_\Lunit\uLunit^i-\Chfour_{\Lunit\uLunit}^i=-2k_{A\spherenormal}e_A^i+k_{\spherenormal\spherenormal}\Lunit^i-\Chfour_{\Lunit\uLunit}^i=\lgensmoothfunction\cdot\pfour\gfour,\\
							\Lunit\spherenormal^i&=\frac{1}{2}\Lunit\Lunit^i-\frac{1}{2}\Lunit\uLunit^i=\lgensmoothfunction\cdot\pfour\gfour
							.\end{align}
					\end{subequations}
				\end{proof}
				\begin{proof}[Proof of $\derinormal\spherenormal^i=a\cdot\lgensmoothfunction\cdot\pfour\gfour-\angnabla a$ in (\ref{EQ:derinormalLunit})]
					Using (\ref{EQ:DfourNN}), we obtain
					
					\begin{align}\label{EQ:EaN}
						\derinormal\spherenormal^i=a\spherenormal\spherenormal^i=&-a\angnabla\ln a-\frac{a}{2}k_{\spherenormal\spherenormal}\Lunit^i-\frac{a}{2}k_{\spherenormal\spherenormal}\uLunit^i-a\Chfour_{\spherenormal\spherenormal}^i,\\
						\notag=&-\angnabla a+a\cdot\lgensmoothfunction\cdot\pfour\gfour
						.\end{align}
					
				\end{proof}
				\begin{proof}[Proof of $\derinormal\Lunit^i=a\cdot\lgensmoothfunction\cdot\pfour\gfour-\angnabla a$ in (\ref{EQ:derinormalLunit})] Note that
					\begin{subequations}\label{EQ:EaL}
						\begin{align}
							\derinormal\Lunit^i&=\derinormal\Timelike^i+\derinormal\spherenormal^i,\\
							\derinormal\Timelike^i&=a\spherenormal\cdot d\Timelike^i=0
							.\end{align}
					\end{subequations}
					Combining (\ref{EQ:EaL}) and (\ref{EQ:EaN}), we conclude the desired equation.
				\end{proof}
				\begin{proof}[Proof of $\Lunit\cartiasphere^i=\lgensmoothfunction\cdot\ACC\cdot\vcartiasphere$ in (\ref{EQ:Lunitcartiasphere})]
					A direct computation gives
					\begin{align}
						\Lunit\cartiasphere^i&=\Lunit\left(\frac{1}{\rgeo}\deriasphere\right)^i=-\frac{1}{\rgeo^2}\Lunit(\rgeo)\left(\deriasphere\right)^i+\frac{1}{\rgeo}\Lunit\left(\deriasphere\right)^i.
					\end{align}
					Note that $\frac{\p}{\p t}=\Lunit$ relative to geometric coordinates. Hence, $\Lunit$ commutes with $\deriasphere$. By (\ref{EQ:CC}), it holds that
					\begin{align}
						\Lunit\left(\deriasphere\right)^i&=\deriasphere\Lunit^i=\upchi_{\deriasphere B}e_B^i-k_{\deriasphere\spherenormal}\Lunit^i+\lgensmoothfunction\cdot\pfour\gfour\cdot\deriasphere\\
						\notag&=\left(\hchi_{\deriasphere B}+\frac{1}{2}\left(\chismall+\frac{2}{\rgeo}-2\Chfour_{\Lunit}\right)\gsphere_{\deriasphere B}\right)e_B^i+\lgensmoothfunction\cdot\pfour\gfour\cdot\deriasphere\\
						\notag&=\frac{1}{\rgeo}\left(\deriasphere\right)^i+\lgensmoothfunction\cdot\ACC\cdot\deriasphere
						.\end{align}
					Combining the above equations and using the fact $\Lunit(\rgeo)=1$, we obtain the desired equation.
				\end{proof}
				\begin{proof}[Proof of $\derinormal\cartiasphere^i=\lapse\cdot\lgensmoothfunction\cdot\ACC\cdot\vcartiasphere+\lgensmoothfunction\cdot\angnabla \lapse\cdot\vcartiasphere$ in (\ref{EQ:derinormalcartiasphere})]
					We are computing on $\Szero$, where $w=\rgeo=-u, a=b$. Using (\ref{EQ:DfourNN}), we deduce
					\begin{subequations}
						\begin{align}
							\label{EQ:deriaspherecartiasphere}\derinormal\cartiasphere^i&=\derinormal\left(\frac{1}{\rgeo}\deriasphere\right)=-\frac{1}{\rgeo^2}\derinormal(\rgeo)\deriasphere+\frac{1}{\rgeo}\derinormal\left(\deriasphere\right)^i,\\
							\label{EQ:derinormalderiasphere}\derinormal\left(\deriasphere\right)^i&=\deriasphere\left(\derinormal\right)^i=\deriasphere\left(a\spherenormal^i\right)=\left(\deriasphere a\right)\spherenormal^i+a\deriasphere\spherenormal^i,\\
							\label{EQ:deriaspheeN}\deriasphere\spherenormal^i&=\spheresecondfund_{\deriasphere B}e_B^i-k_{\deriasphere B}\spherenormal^i+\lgensmoothfunction\cdot\pfour\gfour\cdot\deriasphere\\
							\notag&=\upchi_{\deriasphere B}e_B^i+\lgensmoothfunction\cdot\pfour\gfour\cdot\deriasphere\\
							\notag&=\left(\hchi_{\deriasphere B}+\frac{1}{2}\left(\chismall+\frac{2}{\rgeo}-\Chfour_{\Lunit}\right)\gsphere_{\deriasphere B}\right)e_B^i+\lgensmoothfunction\cdot\pfour\gfour\cdot\deriasphere
							.\end{align}
					\end{subequations}
					We then utilize the properties on the initial foliation in Proposition \ref{PR:InitialCT1}:
					\begin{align}
						\gtr\spheresecondfund+k_{\spherenormal\spherenormal}&=\frac{2}{\lapse w}+\ggtr k-\Chfour_{\Lunit},\\
						\label{EQ:Restracechismall0}\chismall|_{\Szero}&=\frac{2(1-\lapse)}{\lapse w},&
						\text{for} \ 0\leq w\leq w_*
						.\end{align}
					Using (\ref{EQ:Restracechismall0}) and (\ref{EQ:deriaspheeN}), we have:
					\begin{align}
						\label{EQ:deriasphereNi}\deriasphere\spherenormal^i=\frac{1}{\lapse w}\deriasphere+\lgensmoothfunction\cdot\ACC\cdot\deriasphere
						.\end{align}
					Combining (\ref{EQ:deriasphereNi}), (\ref{EQ:derinormalderiasphere}),  (\ref{EQ:deriaspherecartiasphere}), and noticing that $\derinormal \rgeo=1$, we conclude the desired equation.
				\end{proof}
				\begin{proof}[Proof of $\spherenormal\deriasphere^i=\lgensmoothfunction\cdot\pfour\gfour\cdot\deriasphere+\angnabla_A\ln\nulllapse\cdot\spherenormal$ in (\ref{EQ:NomegaA})]
					Employing (\ref{EQ:deriaspheeN}) and (\ref{EQ:commutatorformula1}), we get
					\begin{align}
						\spherenormal\deriasphere^i=k_{\deriasphere B}e_B^i+\angnabla_A\ln\nulllapse\cdot\spherenormal^i
					\end{align}
				\end{proof}
				\begin{proof}[Proof of $\Lunit\chismall+\frac{2}{\rgeo}\chismall$ in (\ref{EQ:Lz})]
					Using (\ref{DE:tracechismall}) and (\ref{EQ:Ltraceupchi}), we derive
					\begin{align}
						&\Lunit\left(\chismall+\frac{2}{\rgeo}-\Chfour_\Lunit\right)+\frac{1}{2}\left(\chismall+\frac{2}{\rgeo}-\Chfour_\Lunit\right)^2=-\sgabs{\hchi}^2-k_{\spherenormal\spherenormal}\left(\chismall+\frac{2}{\rgeo}-\Chfour_\Lunit\right)-\Ricfour{\Lunit}{\Lunit}
						.\end{align}
					Substituting $\Ricfour{\Lunit}{\Lunit}$ above by using (\ref{EQ:riccilunitlunit}), we arrive at
					\begin{align}
						\Lunit\chismall+\frac{2}{\rgeo}\chismall=&-\frac{1}{2}\chismall^2-\sgabs{\hchi}^2-\frac{2}{\rgeo^2}+\frac{2}{\rgeo^2}\Lunit(\rgeo)+\Lunit\Chfour_{\Lunit}-\frac{1}{2}\Chfour_{\Lunit}^2+\frac{2}{\rgeo}\Chfour_\Lunit\\
						\notag&+\chismall\Chfour_\Lunit-k_{\spherenormal\spherenormal}\chismall-k_{\spherenormal\spherenormal}\frac{2}{\rgeo}+k_{\spherenormal\spherenormal}\Chfour_{\Lunit}\\
						\notag&-\left(\lgensmoothfunction\cdot\boxg \gfour_{\alpha\beta}(\vvariables,\cvvariables)+\Lunit\Chfour_\Lunit+k_{\spherenormal\spherenormal}\Chfour_{\Lunit}+\lgensmoothfunction\cdot\pfour\gfour\cdot\pfour\gfour\right)\\
						\notag=&-\frac{1}{2}\chismall^2-\sgabs{\hchi}^2+\lgensmoothfunction\cdot\boxg\gfour+\lgensmoothfunction\cdot\left(\ACC,\rgeo^{-1}\right)\cdot\pfour\gfour
						.\end{align}
					Thus, we obtain the desired equation (\ref{EQ:Lz}).
				\end{proof}
				\begin{proof}[Proof of $\angD_{\Lunit}\angnabla\chismall+\frac{3}{\rgeo}\angnabla\chismall$ in (\ref{EQ:Lunitangnablaz})]
					By (\ref{EQ:chiuchi}), (\ref{DE:tracechismall}) and (\ref{EQ:CC}), we have:
					\begin{align}\label{EQ:DLAtrace1}
						\angD_{\Lunit}\angnabla_A\chismall-\angnabla_A\angD_{\Lunit}\chismall&=\Lunit e_A(\chismall)-\angnabla_{\angD_{\Lunit}e_A}\chismall-e_A\Lunit\chismall\\
						\notag&=\left([\Lunit,\angnabla_A]-\angD_{\Lunit}e_A\right)\chismall\\
						\notag&=\left(\Dfour_\Lunit e_A-\angD_{\Lunit}e_A-\Dfour_{A}\Lunit\right)\chismall\\
						\notag&=\left(-k_{A\spherenormal}\Lunit-\upchi_{AB}e_B+k_{A\spherenormal}\Lunit\right)\chismall\\
						\notag&=\left(-\hchi_{AB}-\frac{1}{2}\left(\chismall+\frac{2}{\rgeo}-\Chfour_{\Lunit}\right)\gsphere_{AB}\right)e_B\chismall\\
						\notag&=\lgensmoothfunction\cdot\ACC\cdot\angnabla_A\chismall-\frac{1}{\rgeo}\angnabla_A\chismall
						.\end{align}
					Now we apply (\ref{EQ:Ltraceupchi}) to $\angnabla_A\angD_{\Lunit}\chismall$ , replace $\Ricfour{\Lunit}{\Lunit}$ above with (\ref{EQ:riccilunitlunit}), and use the fact $\angnabla_A\rgeo=0$ to obtain:
					\begin{align}\label{EQ:DLAtrace2}
						\angnabla_A\left(\Lunit\chismall\right)=&\angnabla_A\left(-\frac{2}{\rgeo}\chismall\right)+\angnabla_A\left(\lgensmoothfunction\cdot\boxg\gfour\right)\\
						\notag&+\angnabla_A\left(\lgensmoothfunction\cdot\left(\pfour\gfour,\chismall,\rgeo^{-1}\right)\cdot\pfour\gfour\right)\\
						\notag&-\angnabla_A\left(\sgabs{\hchi}^2\right)+\angnabla_A\left(\frac{1}{2}\chismall^2\right)\\
						\notag=&-\frac{2}{\rgeo}\angnabla_A\chismall+\lgensmoothfunction\cdot\left(\ACC,\rgeo^{-1}\right)\cdot(\boxg\gfour)\\
						\notag&+\lgensmoothfunction\cdot\angnabla_A(\boxg\gfour)\\
						\notag&+\lgensmoothfunction\cdot\left(\pfour\gfour,\chismall,\hchi,\rgeo^{-1}\right)\cdot\left(\pfour\gfour,\chismall,\rgeo^{-1}\right)\cdot\pfour\gfour\\
						\notag&+\lgensmoothfunction\cdot\left(\angnabla\pfour\gfour,\angnabla\chismall\right)\cdot\pfour\gfour+\lgensmoothfunction\cdot\left(\pfour\gfour,\chismall,\rgeo^{-1}\right)\cdot\angnabla\pfour\gfour\\
						\notag&-2\hchi\cdot\angnabla_A\hchi+\chismall\cdot\angnabla_A\chismall
						.\end{align}
					Combining (\ref{EQ:DLAtrace1})-(\ref{EQ:DLAtrace2}), we conclude the desired result (\ref{EQ:Lunitangnablaz}).
				\end{proof}
				\begin{proof}[Proof of $\Lunit\left(\frac{1}{2}\gtr\reschi \volume\right)$ in (\ref{EQ:Lv})]
					We use (\ref{EQ:Lunitv}), (\ref{EQ:Ltraceupchi}) and (\ref{EQ:riccilunitlunit}) for the computations below:
					\begin{subequations}\label{EQ:Lvolume2}
						\begin{align}
							\Lunit\left(\frac{1}{2}\gtr\reschi \volume\right)=&\frac{1}{2}\Lunit\left\{(\gtr\upchi+\Chfour_{\Lunit})\volume\right\}=\frac{1}{2}\Lunit(\gtr\upchi \volume)+\frac{1}{2}\Lunit\left(\Chfour_{\Lunit}\volume\right),\\
							\Lunit\left(\Chfour_{\Lunit}\volume\right)=&\Lunit\Chfour_{\Lunit} \volume+\Chfour_{\Lunit}\volume\gtr\upchi,\\
							\Lunit(\gtr\upchi \volume)=&\volume\left(\gtr\upchi\right)^2+\volume(\Lunit\gtr\upchi)\\
							\notag=&\volume\left(\gtr\upchi\right)^2+\volume\left(-\frac{1}{2}(\gtr\upchi)^2-\sgabs{\hchi}^2-k_{\spherenormal\spherenormal}\gtr\upchi\right)\\
							\notag&-\volume\left(\lgensmoothfunction\cdot\boxg \gfour_{\alpha\beta}(\vvariables,\cvvariables)+\Lunit\Chfour_\Lunit+k_{\spherenormal\spherenormal}\Chfour_{\Lunit}+\lgensmoothfunction\cdot\pfour\gfour\cdot\pfour\gfour\right)
							.\end{align}
					\end{subequations}
					Combining (\ref{EQ:Lvolume2}) and the fact that 
					\begin{align}
						\volume k_{\spherenormal\spherenormal}\gtr\upchi+\volume k_{\spherenormal\spherenormal}\Chfour_{\Lunit}=\volume k_{\spherenormal\spherenormal}\gtr\reschi=\volume\gtr\reschi\Lunit(\ln\nulllapse)
						,\end{align}
					we get the desired equation (\ref{EQ:Lv}).
				\end{proof}
				\begin{proof}[Proof of $\angdiv\hchi=\angnabla\chismall+\angdiv\xi^{(1)}+\angnabla\xi^{(2)}+\lgensmoothfunction\cdot(\ACC,\rgeo^{-1})\cdot\pfour\gfour$ in (\ref{EQ:angdvhchi})]
					The desired equation follows by combining equations (\ref{EQ:angdivhchi}), (\ref{EQ:riemannABBL}), $\angnabla_A\rgeo=0$, and the following facts:
					\begin{subequations}\label{EQ:angdivhchi1}
						\begin{align}
							\angnabla_A\gtr\upchi&=\angnabla_A\left(\chismall-\Chfour_{\Lunit}+\frac{2}{\rgeo}\right)=\angnabla_A\chismall-\angnabla_A\Chfour_{\Lunit},\\
							\angnabla_A\Chfour_{\Lunit}&=e_A\left(\Chfour_{\Lunit}\right)-\left(\upchi_{AB}e_B^\alpha-k_{A\spherenormal}\Lunit^\alpha-\Chfour_{A\Lunit}^\alpha\right)\Chfour_\alpha\\
							\notag&=\angnabla\xi+\lgensmoothfunction\cdot\left(\pfour\gfour,\chismall,\hchi,\rgeo^{-1}\right)\cdot\pfour\gfour
							.\end{align}
					\end{subequations}
				\end{proof}
				\begin{proof}[Proof of $\angD_{\Lunit}\hchi+(\gtr\upchi)\hchi=(\angnabla,\angD_{\Lunit})\xi+\lgensmoothfunction\cdot\boxg\gfour+\lgensmoothfunction\cdot(\ACC,\rgeo^{-1})\cdot\pfour\gfour$ in (\ref{EQ:angDhchi})]
					The desired result follows from combining equations (\ref{EQ:angDLunithchi}), (\ref{EQ:riccilunitlunit}) and (\ref{EQ:rLALB}), and the following fact:
					\begin{align}
						\Lunit\Chfour_{\Lunit}&=\Lunit(\Chfour_{\Lunit})-k_{\spherenormal\spherenormal}\Chfour_\Lunit=\angD_{\Lunit}\xi-\lgensmoothfunction\cdot\pfour\gfour\cdot\pfour\gfour
						.\end{align}
				\end{proof}
				\begin{proof}[Proof of (\ref{EQ:angDLunitupzeta})]
					We derive the result by combining equations (\ref{EQ:angDLupzetaA}) and (\ref{EQ:rALLuLunit}).
				\end{proof}
				\begin{proof}[Proof of \eqref{EQ:Lunitnulllapse}]
					The equation \eqref{EQ:Lunitnulllapse} follows directly from \eqref{EQ:evol eq nulllapse}.
				\end{proof}
				\begin{proof}[Proof of (\ref{EQ:dtrminus2gsphere})] First, we have:
					\begin{align}
						\upchi\left(\deriasphere,\deribsphere\right)=\hchi\left(\deriasphere,\deribsphere\right)+\frac{1}{2}\left(\chismall+\frac{2}{\rgeo}-\Chfour_{\Lunit}\right)\gsphere\left(\deriasphere,\deribsphere\right)
						.\end{align}
					Using the facts that $\Lunit(\rgeo)=1$ and $\Dfour_\Lunit\deriasphere=\Dfour_\deriasphere\Lunit$, we obtain:
					\begin{align}
						\Lunit\left\{\rgeo^{-2}\gsphere\left(\deriasphere,\deribsphere\right)\right\}&=-2\rgeo^{-3}\gsphere\left(\deriasphere,\deribsphere\right)+2\rgeo^{-2}\upchi\left(\deriasphere,\deribsphere\right)\\
						\notag&=2\rgeo^{-2}\hchi\left(\deriasphere,\deribsphere\right)+\rgeo^{-2}\left(\chismall-\Chfour_{\Lunit}\right)\gsphere\left(\deriasphere,\deribsphere\right)
						.\end{align}
					This finishes the proof.
				\end{proof}
				\begin{proof}[Proof of $\Lunit\cphi=\gtr\upchi-\frac{2}{\rgeo}=\chismall-\Chfour_{\Lunit}$ in (\ref{EQ:Lunitcphi})]
					We use equation (\ref{EQ:Lunitv}) and the fact that $\Lunit(\rgeo)=1$ to deduce:
					\begin{align}
						\Lunit\cphi&=\rgeo^2\volume^{-1}\left(-2\rgeo^{-3}\Lunit(\rgeo)\volume+\rgeo^{-2}\Lunit \volume\right)\\
						\notag&=\rgeo^2\volume^{-1}\left\{-2\rgeo^{-3}\Lunit(\rgeo)\volume+\rgeo^{-2}\volume\left(\chismall-\Chfour_{\Lunit}+\frac{2}{\rgeo}\right)\right\}\\
						\notag&=\chismall-\Chfour_{\Lunit}
						.\end{align}
				\end{proof}
				\begin{proof}[Proof of $\Lunit\angnabla\cphi+\frac{1}{2}(\gtr\upchi)\angnabla\cphi=\lgensmoothfunction\cdot\hchi\cdot\angnabla\cphi+\angnabla\chismall-\angnabla(\Chfour_{\Lunit})$ in (\ref{EQ:Lunitangnablacphi})] Using (\ref{EQ:commutatorformula}), we have:
					\begin{subequations}\label{EQ:Langnablacphi}
						\begin{align}
							\left(\Lunit\angnabla_A-\angnabla_A\Lunit\right)\cphi=&-\upchi_{AB}e_B\cphi,\\
							\angnabla_A\Lunit\cphi=&\angnabla_A\left(\chismall-\Chfour_{\Lunit}\right),\\
							\upchi_{AB}e_B\cphi=&\hchi_{AB}e_B\cphi+\frac{1}{2}\gtr\upchi e_A\cphi
							.\end{align}
					\end{subequations}
					Combining the equations in (\ref{EQ:Langnablacphi}), we conclude the desired equation (\ref{EQ:Lunitangnablacphi}).
				\end{proof}
				\begin{proof}[Proof of the algebraic identity of $\mass$ in (\ref{EQ:mass})]
					We use (\ref{EQ:uLunitv}) and the fact that $\uLunit(\rgeo)=\Lunit(\rgeo)-2\spherenormal(\rgeo)=1-\frac{2}{\nulllapse}$ to get the following equation:
					\begin{align}\label{EQ:Lbarcphi}
						\uLunit\cphi&=\rgeo^2\volume^{-1}\left(-2\rgeo^{-3}\uLunit(\rgeo)\volume+\rgeo^{-2}\uLunit \volume\right)\\
						\notag&=\rgeo^2\volume^{-1}\left(-2\rgeo^{-3}\left(1-\frac{2}{\nulllapse}\right)\volume+\rgeo^{-2}\gtr\uchi \volume\right)\\
						\notag&=\gtr\uchi-\frac{2}{\rgeo}+\frac{4}{\nulllapse\rgeo}
						.\end{align}
					Combining (\ref{EQ:Lbarcphi}) and (\ref{EQ:Lunitcphi}), we have the following equation:
					\begin{align}
						\spherenormal\cphi=\gtr\spheresecondfund-\frac{2}{\nulllapse\rgeo}
						.\end{align}
					Taking $\Lunit$ derivative of (\ref{EQ:Lbarcphi}) and $\uLunit$ derivative of (\ref{EQ:Lunitcphi}), we obtain
					\begin{align}\label{EQ:LuLcphi1}
						\left(\Lunit\uLunit-\uLunit\Lunit\right)\cphi&=\Lunit\left(\gtr\uchi-\frac{2}{\rgeo}+\frac{4}{\nulllapse\rgeo}\right)-\uLunit\left(\gtr\upchi-\frac{2}{\rgeo}\right)\\
						\notag&=\Lunit\gtr\uchi-\uLunit\gtr\upchi+\frac{4}{\rgeo}\Lunit\left(\nulllapse^{-1}\right)\\
						\notag&=\Lunit\gtr\uchi-\uLunit\gtr\upchi+\frac{4}{\nulllapse\rgeo}k_{\spherenormal\spherenormal}
						.\end{align}
					It then follows from (\ref{EQ:LuLcphi1}) and (\ref{EQ:commutatorformulaLuL}) that
					\begin{align}\label{EQ:LuLcphi2}
						\Lunit\gtr\uchi-\uLunit\gtr\upchi+\frac{4}{\nulllapse\rgeo}k_{\spherenormal\spherenormal}=&[\Lunit,\uLunit]\cphi\\
						\notag=&2\left(\uzeta_A-\upzeta_A\right)e_A\cphi+k_{\spherenormal\spherenormal}\uLunit\cphi-k_{\spherenormal\spherenormal}\Lunit\cphi\\
						\notag=&2\left(\uzeta_A-\upzeta_A\right)e_A\cphi-2k_{\spherenormal\spherenormal}\left(\gtr\spheresecondfund-\frac{2}{\nulllapse\rgeo}\right)
						.\end{align}
					Notice that (\ref{EQ:LuLcphi2}) can be rewritten as follows:
					\begin{align}\label{Ltraceuchi1}
						\Lunit\gtr\uchi-\uLunit\gtr\upchi&=2\left(\uzeta_A-\upzeta_A\right)e_A\cphi-2k_{\spherenormal\spherenormal}\gtr\spheresecondfund
						.\end{align}
					Combining (\ref{Ltraceuchi1}) with (\ref{Ltraceuchi}), by the definition of $\mass$ in Definition (\ref{DE:mass}), we can get
					\begin{align}
						\mass=&2\angdiv\uzeta+k_{\spherenormal\spherenormal}\gtr\uchi-\hchi_{AB}\huchi_{AB}+2\sgabs{\uzeta}^2+\Riemfour{A}{\uLunit}{\Lunit}{A}\\
						\notag&-2\left(\uzeta_A-\upzeta_A\right)e_A\cphi+2k_{\spherenormal\spherenormal}\gtr\spheresecondfund\\
						\notag=&2\angdiv\uzeta+k_{\spherenormal\spherenormal}\gtr\upchi-\hchi_{AB}\left(-2\hat{k}_{AB}-\hchi_{AB}\right)+2\sgabs{\uzeta}^2\\
						\notag&+\Riemfour{A}{\uLunit}{\Lunit}{A}-2\left(\uzeta_A-\upzeta_A\right)e_A\cphi
						.\end{align}
					Using (\ref{EQ:ALuLA}), and noticing that $\uzeta_A=-k_{A\spherenormal}$, we conclude the desired equation (\ref{EQ:mass}).
				\end{proof}
				\begin{proof}[Proof of $\Lunit\angnabla\conformalfactor+\frac{1}{2}\gtr\upchi\angnabla\conformalfactor=\frac{1}{2}\angnabla(\Chfour_{\Lunit})-\hchi\cdot\angnabla\conformalfactor$ in (\ref{EQ:Lunitangnablasigma})] By (\ref{EQ:commutatorformula}), there holds
					\begin{align}
						\Lunit\angnabla_A\conformalfactor&=\angnabla_A\Lunit\conformalfactor-\upchi_{AB}\angnabla_B\conformalfactor=\frac{1}{2}\angnabla_A(\Chfour_{\Lunit})-\hchi_{AB}\angnabla_B\conformalfactor-\frac{1}{2}\gtr\upchi\angnabla_A\conformalfactor
						.\end{align}
				\end{proof}
				
				We first introduce several preliminary ingredients before deriving the transport equation for $\modmass$ (\ref{EQ:Lunitmodmass}).
				\begin{proposition}\cite[(6.6). Transport equation for $\anglap\conformalfactor$]{AGeoMetricApproach}. The angular Laplacian of the conformal factor verifies the following transport equation:
					\begin{align}
						\label{EQ:anglapconformalfactor}\Lunit\anglap\conformalfactor+\gtr\upchi\anglap\conformalfactor=&\frac{1}{2}\anglap(\Chfour_{\Lunit})-2\hchi\cdot\angnabla^2\conformalfactor-\angnabla_B\gtr\upchi\angnabla_B\conformalfactor-k_{B\spherenormal}\gtr\upchi\angnabla_B\conformalfactor\\
						\notag&+2\hchi_{AB}k_{A\spherenormal}\angnabla_B\conformalfactor+2\Riemfour{A}{B}{\Lunit}{A}\angnabla_B\conformalfactor\\
						\notag=&\frac{1}{2}\anglap(\Chfour_{\Lunit})-2\hchi\cdot\angnabla^2\conformalfactor+\lgensmoothfunction\cdot\left(\angnabla\pfour\gfour,\angnabla\chismall\right)\cdot\angnabla\conformalfactor\\
						\notag&+\lgensmoothfunction\cdot\left(\ACC,\rgeo^{-1}\right)\cdot\pfour\gfour\cdot\angnabla\conformalfactor
						.\end{align}
				\end{proposition}
				\begin{proof}
					Here, we use the commutator formulas in Lemma \ref{LE:Commutatorformulas} to compute $[\Lunit,\anglap]\conformalfactor$:
					\begin{align}\label{EQ:Langlapsigma1}
						\Lunit\anglap\conformalfactor-\anglap\Lunit\conformalfactor=&\delta^{AB}\left(\angD_{\Lunit}\angnabla_B\angnabla_A\conformalfactor-\angnabla_B\angD_{\Lunit}\angnabla_A\conformalfactor-\angnabla_B\upchi_{AC}\angnabla_C\conformalfactor\right)\\
						\notag=&\left(-\gtr\upchi k_{B\spherenormal}+\upchi_{AB}k_{A\spherenormal}+\Riemfour{A}{B}{\Lunit}{A}\right)\angnabla_B\conformalfactor-\upchi_{AB}\angnabla_B\angnabla_A\conformalfactor-\angnabla_A\upchi_{AB}\angnabla_B\conformalfactor\\
						\notag=&\left(-\gtr\upchi k_{B\spherenormal}+\upchi_{AB}k_{A\spherenormal}+\Riemfour{A}{B}{\Lunit}{A}\right)\angnabla_B\conformalfactor\\
						\notag&-\left(\hchi_{AB}+\frac{1}{2}\gtr\upchi\gsphere_{AB}\right)\angnabla_B\angnabla_A\conformalfactor-\angnabla_A\upchi_{AB}\angnabla_B\conformalfactor\\
						\notag=&\left(-\gtr\upchi k_{B\spherenormal}+\upchi_{AB}k_{A\spherenormal}+\Riemfour{A}{B}{\Lunit}{A}\right)\angnabla_B\conformalfactor\\
						\notag&-\hchi\cdot\angnabla^2\conformalfactor-\frac{1}{2}\gtr\upchi\anglap\conformalfactor-\angnabla_A\upchi_{AB}\angnabla_B\conformalfactor
						.\end{align}
					Next, we calculate $\angnabla_A\upchi_{AB}\angnabla_B\conformalfactor$ employing (\ref{EQ:angdivhchi})
					\begin{align}\label{EQ:Langlapsigma2}
						\angnabla_A\upchi_{AB}\angnabla_B\conformalfactor=&\angnabla_A\left(\hchi_{AB}+\frac{1}{2}\gtr\upchi\gsphere_{AB}\right)\angnabla_B\conformalfactor+\hchi\cdot\angnabla^2\conformalfactor+\frac{1}{2}\gtr\upchi\anglap\conformalfactor\\
						\notag=&\angdiv\hchi_B\angnabla_B\conformalfactor+\frac{1}{2}\angnabla_A\gtr\upchi\angnabla_A\conformalfactor+\hchi\cdot\angnabla^2\conformalfactor+\frac{1}{2}\gtr\upchi\anglap\conformalfactor\\
						\notag=&(\angnabla_B\gtr\upchi)\angnabla_B\conformalfactor+\frac{1}{2}k_{B\spherenormal}\gtr\upchi\angnabla_B\conformalfactor+\Riemfour{A}{\Lunit}{A}{B}\angnabla_B\conformalfactor\\
						\notag&-\hchi_{AB}k_{A\spherenormal}\angnabla_B\conformalfactor+\hchi\cdot\angnabla^2\conformalfactor+\frac{1}{2}\gtr\upchi\anglap\conformalfactor
						.\end{align}
					Combining (\ref{EQ:Langlapsigma1})-(\ref{EQ:Langlapsigma2}), and writing $\upchi_{AB}k_{A\spherenormal}\angnabla_B\conformalfactor=\hchi_{AB}k_{A\spherenormal}\angnabla_B\conformalfactor+\frac{1}{2}(\chismall+\frac{2}{\rgeo}-\Chfour_{\Lunit}) k_{B\spherenormal}\angnabla_B\conformalfactor$, we conclude the desired equation (\ref{EQ:anglapconformalfactor}).
				\end{proof}
				\begin{proposition}\cite[(6.9). Identity for $\uLunit(k_{\spherenormal\spherenormal})$]{AGeoMetricApproach}. The following identity holds:
					\begin{align}
						\uLunit(k_{\spherenormal\spherenormal})+\Lunit(k_{\spherenormal\spherenormal})=\Ricfour{\Lunit}{\uLunit}+4k_{A\spherenormal}\upzeta_A+\Riemfour{A}{\Lunit}{\uLunit}{A}+2k_{A\spherenormal}k_{A\spherenormal}+2k_{\spherenormal\spherenormal}^2
						.\end{align}
				\end{proposition}
				\begin{proof} Via using the definitions of $k$ in (\ref{DE:SECONDFUNDOFSIGMAT}) and the curvature tensor, the identities \eqref{EQ:CC}, and the fact that $\Timelike$ is geodesic, we have:
					\begin{align}
						\uLunit(k_{\spherenormal\spherenormal})+\Lunit(k_{\spherenormal\spherenormal})=&-\gfour(\Dfour_{\uLunit}\Dfour_\spherenormal\Timelike,\spherenormal)-\gfour(\Dfour_{\Lunit}\Dfour_\spherenormal\Timelike,\spherenormal)-2\gfour(\Dfour_\spherenormal\Timelike,\Dfour_\Timelike\spherenormal)\\
						\notag=&-\frac{1}{8}\left(\Riemfour{\uLunit}{\Lunit}{\Lunit}{\uLunit}-\Riemfour{\uLunit}{\Lunit}{\uLunit}{\Lunit}-\Riemfour{\Lunit}{\uLunit}{\Lunit}{\uLunit}+\Riemfour{\Lunit}{\uLunit}{\uLunit}{\Lunit}\right)\\
						\notag&-2\gfour(\Dfour_\spherenormal\Timelike,\Dfour_\Timelike\spherenormal)-2\gfour(\Dfour_{[\Timelike,\spherenormal]}\Timelike,\spherenormal)\\
				\notag=&\Ricfour{\Lunit}{\uLunit}+\Riemfour{A}{\Lunit}{\uLunit}{A}+4\upzeta_Ak_{A\spherenormal}+2k_{A\spherenormal}k_{A\spherenormal}+2k_{\spherenormal\spherenormal}^2
						.\end{align}
				\end{proof}
				\begin{proposition}[Proposition 8.24. Decomposition of the wave operator under null frame]\cite{Yuthesis}. We have the following decompositions for $\boxg f$, $\boxg\Lunit^\alpha$ and $\boxg(\Chfour_{\Lunit})$, with $f$ being a scalar function:
					\begin{subequations}
						\begin{align}
							\label{EQ:boxgf}\boxg f=&-\Lunit\uLunit f+\anglap f-\frac{1}{2}\gtr\upchi\uLunit f-\frac{1}{2}\gtr\uchi\Lunit f-2k_{A\spherenormal}\angnabla_Af+k_{\spherenormal\spherenormal}\uLunit f,\\
							\label{EQ:boxgL}\boxg\Lunit^\alpha=&\left(-\hchi_{AB}k_{B\spherenormal}+\frac{1}{2}(\angnabla_A\gtr\upchi+k_{A\spherenormal}\gtr\upchi)+\Riemfour{B}{\Lunit}{B}{A}\right)e_A^\alpha\\
							\notag&+\upchi_{AB}\left(\frac{1}{2}\upchi_{AB}\uLunit^\alpha+\frac{1}{2}\uchi_{AB}\Lunit^\alpha\right)-\left(\Dfour_A(k_{A\spherenormal})\Lunit^\alpha+k_{AB}\spheresecondfund_{AB}\Lunit^\alpha\right)\\
							\notag&-k_{A\spherenormal}\left(\upchi_{AB}e_B^\alpha-k_{A\spherenormal}\Lunit^\alpha\right)\\
							\notag&-\frac{1}{2}\gtr\upchi(2\upzeta_Ae_A^\alpha+k_{\spherenormal\spherenormal}\Lunit^\alpha)+\frac{1}{2}\gtr\uchi k_{\spherenormal\spherenormal}\Lunit^\alpha\\
							\notag&+\uLunit(k_{\spherenormal\spherenormal})\Lunit^\alpha+2\upzeta_A\left(k_{\spherenormal\spherenormal}e_A^\alpha+\upchi_{AB}e_B^\alpha-k_{A\spherenormal}\Lunit^\alpha\right)-\frac{1}{2}\Riemonethree{\alpha}{\Lunit}{\Lunit}{\uLunit},\\
							\label{EQ:boxgChfourL}\boxg(\Chfour_{\Lunit})=&\uLunit(k_{\spherenormal\spherenormal})\Chfour_{\Lunit}+\rgeo^{-2}\xi\\
							\notag&+\lgensmoothfunction\cdot\ACC\cdot\pfour^2\gfour+\lgensmoothfunction\cdot\pfour\gfour\cdot\angnabla\chismall\\
							\notag&+\lgensmoothfunction\cdot\boxg\gfour\cdot\pfour\gfour+\lgensmoothfunction\cdot\left(\ACC,\rgeo^{-1}\right)\cdot\ACC\cdot\pfour\gfour
							.\end{align}
					\end{subequations}
				\end{proposition}
				We only present the proof of (\ref{EQ:boxgChfourL}), the proofs of (\ref{EQ:boxgf}) and (\ref{EQ:boxgL}) are identical to those in \cite{Yuthesis}. Also, see \cite[Section 6]{AGeoMetricApproach}.
				\begin{proof}[Proof of $\boxg(\Chfour_{\Lunit})$ in (\ref{EQ:boxgChfourL})] Using the chain rule and (\ref{EQ:CC}), we derive:
					\begin{subequations}
						\begin{align}
							\label{EQ:boxgChfourL0}\boxg(\Chfour_{\Lunit})&=\boxg\Chfour\cdot\Lunit+\boxg\Lunit\cdot\Chfour+2(\gfour^{-1})^{\alpha\beta}\Dfour_\alpha\Chfour\cdot\Dfour_\beta\Lunit,\\
							\label{EQ:boxgChfourL1}2(\gfour^{-1})^{\alpha\beta}\Dfour_\alpha\Chfour\cdot\Dfour_\beta\Lunit&=\left(2e_A^\alpha e_A^\beta-\Lunit^\alpha\uLunit^\beta-\uLunit^\alpha\Lunit^\beta\right)\Dfour_\alpha\Chfour\cdot\Dfour_\beta\Lunit\\
							\notag&=2\Dfour_A\Chfour\cdot\Dfour_A\Lunit-\Dfour_\Lunit\Chfour\cdot\Dfour_{\uLunit}\Lunit-\Dfour_{\uLunit}\Chfour\cdot\Dfour_{\Lunit}\Lunit\\
							\notag&=2\Dfour_A\Chfour\left(\upchi_{AB}e_B-k_{A\spherenormal}\Lunit\right)-\Dfour_\Lunit\Chfour\left(2\upzeta_Ae_A+k_{\spherenormal\spherenormal}\Lunit\right)+\Dfour_{\uLunit}\Chfour\left(k_{\spherenormal\spherenormal}\Lunit\right)\\
							\notag&=\Dfour\Chfour\cdot\lgensmoothfunction\cdot\ACC\\
							\notag&=\lgensmoothfunction\cdot\ACC\cdot\pfour^2\vvariables,\\
							\label{EQ:boxgChfourL2}\boxg\Chfour\cdot\Lunit&=\lgensmoothfunction\cdot\boxg\pfour\gfour+\lgensmoothfunction\cdot\ACC\cdot\pfour^2\gfour+\lgensmoothfunction\cdot\boxg\gfour\cdot\pfour\gfour.\end{align}
					\end{subequations}
					For the second term on the RHS of (\ref{EQ:boxgChfourL0}), by (\ref{EQ:boxgL}), we deduce that
					\begin{align}
						\label{EQ:boxgChfourL3}\boxg\Lunit\cdot\Chfour=&\left(\frac12\angnabla_A\gtr\upchi+\frac12k_{A\spherenormal}\gtr\upchi-k_{B\spherenormal}\upchi_{AB}+\Riemfour{B}{\Lunit}{B}{A}\right)\Chfour_A\\
						\notag&+\upchi_{AB}\left(\frac{1}{2}\upchi_{AB}\Chfour_{\uLunit}+\frac{1}{2}\uchi_{AB}\Chfour_{\Lunit}\right)-\left(\Dfour_A(k_{A\spherenormal})+k_{AB}\spheresecondfund_{AB}\right)\Chfour_{\Lunit}\\
						\notag&-k_{A\spherenormal}\left(\upchi_{AB}\Chfour_B-k_{A\spherenormal}\Chfour_{\Lunit}\right)\\
						\notag&-\frac{1}{2}\gtr\upchi(2\upzeta_A\Chfour_A+k_{\spherenormal\spherenormal}\Chfour_{\Lunit})+\frac{1}{2}\gtr\uchi k_{\spherenormal\spherenormal}\Chfour_{\Lunit}\\
						\notag&+\uLunit(k_{\spherenormal\spherenormal})\Chfour_{\Lunit}+2\upzeta_A\left(k_{\spherenormal\spherenormal}\Chfour_A+\upchi_{AB}\Chfour_B-k_{A\spherenormal}\Chfour_{\Lunit}\right)-\frac{1}{2}\Riemonethree{\alpha}{\Lunit}{\Lunit}{\uLunit}\Chfour_\alpha\\
						\notag=&\frac{1}{2}\sgabs{\upchi}^2\Chfour_{\uLunit}+\frac{1}{2}\upchi_{AB}\uchi_{AB}\Chfour_{\Lunit}+\uLunit(k_{\spherenormal\spherenormal})\Chfour_{\Lunit}\\
						\notag&-\frac12\boxg\gfour\cdot\Chfour_{\Lunit}-\frac14\Lunit\Chfour_{\uLunit}\Chfour_{\Lunit}-\frac14\uLunit\Chfour_{\Lunit}\Chfour_{\Lunit}+\frac12\Riemfour{A}{\uLunit}{A}{\Lunit}\Chfour_{\Lunit}-\frac{1}{2}\Riemfour{A}{\Lunit}{\Lunit}{\uLunit}\Chfour_A\\
						\notag&+\lgensmoothfunction\cdot\pfour\gfour\cdot\ACC\cdot\left(\ACC,\rgeo^{-1}\right)+\lgensmoothfunction\cdot\pfour\gfour\cdot\pfour^2\gfour\\
						\notag&+\lgensmoothfunction\cdot\pfour\gfour\cdot\left(\angnabla\chismall,\angnabla\pfour\gfour\right)\\
						\notag=&\uLunit(k_{\spherenormal\spherenormal})\Chfour_{\Lunit}\\
						\notag&+\rgeo^{-2}\xi+\lgensmoothfunction\cdot\pfour\gfour\cdot\left(\ACC,\rgeo^{-1}\right)\cdot\ACC+\lgensmoothfunction\cdot\pfour\gfour\cdot\pfour^2\gfour\\
						\notag&+\lgensmoothfunction\cdot\pfour\gfour\cdot\angnabla\chismall\\
						\notag&+\lgensmoothfunction\cdot\boxg\gfour\cdot\pfour\gfour+\lgensmoothfunction\cdot\left(\ACC,\rgeo^{-1}\right)\cdot\ACC\cdot\pfour\gfour
						,\end{align}
					where the term $\rgeo^{-2}\xi$ in (\ref{EQ:boxgChfourL3}) comes from the terms $\frac{1}{2}\sgabs{\upchi}^2\Chfour_{\uLunit}$ and $\frac{1}{2}\upchi_{AB}\uchi_{AB}\Chfour_{\Lunit}$, and {in the above derivation}, we use the following identity:
					\begin{align}
						-\frac{1}{2}\Riemonethree{\alpha}{\Lunit}{\Lunit}{\uLunit}\Chfour_\alpha=&\frac14\Riemfour{\uLunit}{\Lunit}{\Lunit}{\uLunit}\Chfour_{\Lunit}-\frac{1}{2}\Riemfour{A}{\Lunit}{\Lunit}{\uLunit}\Chfour_A\\
						\notag=&-\frac12\left(\Ricfour{\Lunit}{\uLunit}-\Riemfour{A}{\uLunit}{A}{\Lunit}\right)\Chfour_{\Lunit}-\frac{1}{2}\Riemfour{A}{\Lunit}{\Lunit}{\uLunit}\Chfour_A+\Riemfour{B}{\Lunit}{B}{A}\Chfour_A\\
						\notag=&-\frac12\boxg\gfour\cdot\Chfour_{\Lunit}-\frac14\Lunit\Chfour_{\uLunit}\Chfour_{\Lunit}-\frac14\uLunit\Chfour_{\Lunit}\Chfour_{\Lunit}+\frac12\Riemfour{A}{\uLunit}{A}{\Lunit}\Chfour_{\Lunit}-\frac{1}{2}\Riemfour{A}{\Lunit}{\Lunit}{\uLunit}\Chfour_A
						\\
						\notag&+\lgensmoothfunction\cdot\ACC\cdot\pfour\gfour\cdot  \Chfour_{\Lunit}.\end{align}
					Combining (\ref{EQ:boxgChfourL0}), (\ref{EQ:boxgChfourL1}), (\ref{EQ:boxgChfourL2}) and (\ref{EQ:boxgChfourL3}), we obtain (\ref{EQ:boxgChfourL}.)
				\end{proof}
				
				\subsubsection{Transport equation for $\mass$}
				Now we compute $\Lunit\mass+\gtr\upchi\mass$. Notice that $\huchi_{AC}\hchi_{BC}\hchi_{AB}=0$. 
				By using the definition of $\mass$ in (\ref{DE:massfunction}), and the commutator formula (\ref{EQ:commutatorformulaLuL}), there holds
				\begin{align}\label{EQ:Lmu1}
					\Lunit\mass+\gtr\upchi\mass=&\Lunit\uLunit\gtr\upchi+\frac{1}{2}\Lunit(\gtr\upchi\gtr\uchi)+\gtr\upchi\left(\uLunit\gtr\upchi+\frac{1}{2}\gtr\upchi\gtr\uchi\right)\\
					\notag=&[\Lunit,\uLunit]\gtr\upchi+\uLunit\Lunit\gtr\upchi+\frac{1}{2}\Lunit(\gtr\upchi\gtr\uchi)+\gtr\upchi\uLunit\gtr\upchi+\frac{1}{2}(\gtr\upchi)^2\gtr\uchi\\
					\notag=&2(\uzeta_A-\upzeta_A)\angnabla_A\gtr\upchi+k_{\spherenormal\spherenormal}(\uLunit-\Lunit)\gtr\upchi\\
					\notag&+\uLunit\Lunit\gtr\upchi+\frac{1}{2}\Lunit(\gtr\upchi\gtr\uchi)+\gtr\upchi\uLunit\gtr\upchi+\frac{1}{2}(\gtr\upchi)^2\gtr\uchi.
				\end{align}
				Together with (\ref{EQ:Ltraceupchi}), this then implies
				\begin{align}\label{EQ:Lmu1.1}
					\Lunit\mass+\gtr\upchi\mass=&2(\uzeta_A-\upzeta_A)\angnabla_A\gtr\upchi+k_{\spherenormal\spherenormal}\uLunit\gtr\upchi-k_{\spherenormal\spherenormal}\Lunit\gtr\upchi\\
					\notag&+\uLunit\left(-\frac{1}{2}(\gtr\upchi)^2-\sgabs{\hchi}^2-k_{\spherenormal\spherenormal}\gtr\upchi-\Ricfour{\Lunit}{\Lunit}\right)+\frac{1}{2}\Lunit(\gtr\upchi\gtr\uchi)+\gtr\upchi\uLunit\gtr\upchi\\
					\notag&+\frac{1}{2}(\gtr\upchi)^2\gtr\uchi\\
					\notag=&2(\uzeta_A-\upzeta_A)\angnabla_A\gtr\upchi+k_{\spherenormal\spherenormal}\uLunit\gtr\upchi-k_{\spherenormal\spherenormal}\Lunit\gtr\upchi\\
					\notag&+\left(-\gtr\upchi\uLunit\gtr\upchi-2\hchi\cdot\angD_{\uLunit}\hchi-k_{\spherenormal\spherenormal}\uLunit\gtr\upchi-\gtr\upchi\uLunit k_{\spherenormal\spherenormal}-\uLunit\Ricfour{\Lunit}{\Lunit}\right)\\
					\notag&+\frac{1}{2}\Lunit(\gtr\upchi\gtr\uchi)+\gtr\upchi\uLunit\gtr\upchi+\frac{1}{2}(\gtr\upchi)^2\gtr\uchi.
				\end{align}
				On the RHS of equation (\ref{EQ:Lmu1.1}), substitute $\Lunit(\gtr\upchi\gtr\uchi)$ by (\ref{EQ:angDLupzetaA}) and (\ref{EQ:Ltraceupchi}), $\uLunit\gtr\upchi$ by (\ref{EQ:uLunitgtrupchi}),  and $\angD_{\uLunit}\hchi$ by (\ref{EQ:angDuLunithchi}).  Then, using the definition of $\modtorsion$ in (\ref{DE:defmodtorsion}), we deduce that
				\begin{align}
					\label{EQ:Transportmu}\Lunit\mass+\gtr\upchi\mass=&2(\uzeta_A-\upzeta_A)\angnabla_A\gtr\upchi-k_{\spherenormal\spherenormal}\Lunit\gtr\upchi+\gtr\uchi\sgabs{\hchi}^2+\gtr\upchi\huchi_{AB}\hchi_{AB}-4\angnabla_A\upzeta_B\hchi_{AB}\\
					\notag&-2k_{\spherenormal\spherenormal}\sgabs{\hchi}^2-4\upzeta_A\upzeta_B\hchi_{AB}+2\huchi_{AC}\hchi_{BC}\hchi_{AB}-2\hchi_{AB}\Riemfour{B}{\Lunit}{\uLunit}{A}\\
					\notag&-\gtr\upchi\uLunit k_{\spherenormal\spherenormal}-\uLunit\Ricfour{\Lunit}{\Lunit}+\frac{1}{2}(\gtr\upchi)^2\gtr\uchi\\
					\notag&+\frac{1}{2}\left(-\sgabs{\hchi}^2\gtr\uchi-\gtr\uchi\Ricfour{\Lunit}{\Lunit}\right.\\
					\notag&\left.+\gtr\upchi\left(2\angdiv\uzeta-\hchi_{AB}\huchi_{AB}+2\sgabs{\uzeta}^2+\Riemfour{A}{\uLunit}{\Lunit}{A}\right)-(\gtr\upchi)^2\gtr\uchi\right)\\
					\notag=&2(\uzeta_A-\upzeta_A)\angnabla_A\gtr\upchi-k_{\spherenormal\spherenormal}\Lunit\gtr\upchi+\frac{1}{2}\gtr\uchi\sgabs{\hchi}^2\\
					\notag&-4\angnabla_A\upzeta_B\hchi_{AB}-4\upzeta_A\upzeta_B\hchi_{AB}-2\hchi_{AB}\Riemfour{B}{\Lunit}{\uLunit}{A}\\
					\notag&-\gtr\upchi\left(\Ricfour{\Lunit}{\uLunit}-\Lunit(k_{\spherenormal\spherenormal})+\Riemfour{A}{\Lunit}{\uLunit}{A}+4k_{A\spherenormal}\upzeta_A+2k_{A\spherenormal}k_{A\spherenormal}+2k_{\spherenormal\spherenormal}^2\right)\\
					\notag&-\uLunit\Ricfour{\Lunit}{\Lunit}-\frac{1}{2}\gtr\uchi\Ricfour{\Lunit}{\Lunit}+\frac{1}{2}\gtr\upchi\left(2\angdiv\uzeta+\hchi_{AB}\huchi_{AB}+2\sgabs{\uzeta}^2+\Riemfour{A}{\uLunit}{\Lunit}{A}\right)\\
					\notag=&-\gtr\upchi\Ricfour{\Lunit}{\uLunit}-\uLunit\Ricfour{\Lunit}{\Lunit}-\frac{1}{2}\gtr\uchi\Ricfour{\Lunit}{\Lunit}\\
					\notag&-k_{\spherenormal\spherenormal}\Lunit\gtr\upchi+\gtr\upchi\Lunit k_{\spherenormal\spherenormal}+4\angnabla^2\conformalfactor\cdot\hchi+\rgeo^{-1}\angdiv\xi+\rgeo^{-2}\xi\\
					\notag&+\lgensmoothfunction\cdot\ACC\cdot(\angnabla\chismall,\angnabla\modtorsion,\pfour^2\gfour)\\
					\notag&+\lgensmoothfunction\cdot\pfour\boxg\gfour+\lgensmoothfunction\cdot\boxg\gfour\cdot\left(\ACC,\rgeo^{-1}\right)\\
					\notag&+\lgensmoothfunction\cdot\ACC\cdot\ACC\cdot\left(\ACC,\rgeo^{-1}\right)
					,\end{align}
				where the term $\rgeo^{-2}\xi$ in (\ref{EQ:Transportmu}) comes from term $\gtr\upchi\Riemfour{A}{\uLunit}{\Lunit}{A}$, and $\rgeo^{-1}\angdiv\xi$ comes from the term $\gtr\upchi\angdiv\upzeta$.\\
				
				Let
				\begin{align}
					R(\mu)=-\gtr\upchi\Ricfour{\Lunit}{\uLunit}-\uLunit\Ricfour{\Lunit}{\Lunit}-\frac{1}{2}\gtr\uchi\Ricfour{\Lunit}{\Lunit}
					.\end{align}
				Recalling (\ref{EQ:boxgChfourL}), by the curvature decomposition in Corollary \ref{CO:Specialricci}, the following identity holds
				\begin{align}
					\label{EQ:Rmu}R(\mu)=&-\Lunit\uLunit\Chfour_{\Lunit}-\frac{1}{2}\gtr\upchi\uLunit\Chfour_{\Lunit}-\frac{1}{2}\gtr\uchi\Lunit\Chfour_{\Lunit}+k_{\spherenormal\spherenormal}\uLunit\Chfour_{\Lunit}\\
					\notag&+\lgensmoothfunction\cdot\pfour\boxg\gfour+\lgensmoothfunction\cdot(\ACC,\rgeo^{-1})\cdot\boxg\gfour\\
					\notag&+\lgensmoothfunction\cdot\pfour\gfour\cdot\pfour^2\gfour+\lgensmoothfunction\cdot\pfour\gfour\cdot\ACC\cdot(\ACC,\rgeo^{-1})\\
					\notag&-\frac{1}{2}\gtr\upchi\Lunit\Chfour_{\uLunit}+2\left(\uzeta_A-\upzeta_A\right)e_A\Chfour_{\Lunit}-k_{\spherenormal\spherenormal}\Lunit\Chfour_{\Lunit}-\uLunit k_{\spherenormal\spherenormal}\Chfour_{\Lunit}-\frac{1}{2}\gtr\uchi k_{\spherenormal\spherenormal}\Chfour_{\Lunit}\\
					\notag=&\boxg\Chfour_{\Lunit}-\anglap\Chfour_{\Lunit}+2k_{A\spherenormal}\angnabla_A\Chfour_{\Lunit}\\
					\notag&+\lgensmoothfunction\cdot\pfour\boxg\gfour+\lgensmoothfunction\cdot(\ACC,\rgeo^{-1})\cdot\boxg\gfour\\
					\notag&+\lgensmoothfunction\cdot\pfour\gfour\cdot\pfour^2\gfour+\lgensmoothfunction\cdot\pfour\gfour\cdot\ACC\cdot(\ACC,\rgeo^{-1})\\
					\notag&-\frac{1}{2}\gtr\upchi\Lunit\Chfour_{\uLunit}+2\left(\uzeta_A-\upzeta_A\right)e_A\Chfour_{\Lunit}-k_{\spherenormal\spherenormal}\Lunit\Chfour_{\Lunit}-\uLunit k_{\spherenormal\spherenormal}\Chfour_{\Lunit}-\frac{1}{2}\gtr\uchi k_{\spherenormal\spherenormal}\Chfour_{\Lunit}
					.\end{align}
				Inserting (\ref{EQ:Rmu}) into (\ref{EQ:Transportmu}) and rearranging the terms, we obtain
				\begin{align}
					\label{EQ:Lmu2}\Lunit\mass+\gtr\upchi\mass=&\boxg\Chfour_{\Lunit}-\anglap\Chfour_{\Lunit}\\
					\notag&-\frac{1}{2}\gtr\upchi\Lunit\Chfour_{\uLunit}-\uLunit k_{\spherenormal\spherenormal}\Chfour_{\Lunit}-k_{\spherenormal\spherenormal}\Lunit\gtr\upchi+\gtr\upchi\Lunit k_{\spherenormal\spherenormal}\\
					\notag&+4\angnabla^2\conformalfactor\cdot\hchi+\rgeo^{-1}\angdiv\xi+\rgeo^{-2}\xi\\
					\notag&+\lgensmoothfunction\cdot\ACC\cdot(\angnabla\chismall,\angnabla\modtorsion,\pfour^2\gfour)\\
					\notag&+\lgensmoothfunction\cdot\pfour\boxg\gfour+\lgensmoothfunction\cdot\boxg\gfour\cdot\left(\ACC,\rgeo^{-1}\right)\\
					\notag&+\lgensmoothfunction\cdot\ACC\cdot\ACC\cdot\left(\ACC,\rgeo^{-1}\right)
					.\end{align}
				
				With the help of (\ref{EQ:Lmu2}), (\ref{EQ:boxgChfourL}) and (\ref{EQ:anglapconformalfactor}), we can now derive the transport equation (\ref{EQ:Lunitmodmass}) for $\modmass$. 
				
				\begin{proof}[Proof of $\Lunit\modmass+\gtr\upchi\modmass$ in (\ref{EQ:Lunitmodmass})] By using the definition of $\modmass$ in (\ref{DE:modmass}), the equations (\ref{EQ:anglapconformalfactor}), (\ref{EQ:boxgChfourL}) and (\ref{EQ:Lmu2}), there holds
					\begin{align}
						\label{EQ:Transportmodmass1}\Lunit\modmass+\gtr\upchi\modmass=&\Lunit\mass+\gtr\upchi\mass+2\left(\Lunit\anglap\conformalfactor+\gtr\upchi\anglap\conformalfactor\right)\\
						\notag&-\left\{\Lunit(k_{\spherenormal\spherenormal}\gtr\upchi)+k_{\spherenormal\spherenormal}(\gtr\upchi)^2\right\}+\frac{1}{2}\left\{\Lunit(\gtr\upchi\Chfour_{\Lunit})+\Chfour_{\Lunit}(\gtr\upchi)^2\right\}\\
						\notag=&-\frac{1}{2}\gtr\upchi\Lunit\Chfour_{\uLunit}-\uLunit (k_{\spherenormal\spherenormal})\Chfour_{\Lunit}-k_{\spherenormal\spherenormal}\Lunit\gtr\upchi+\gtr\upchi\Lunit k_{\spherenormal\spherenormal}\\
						\notag&-\Lunit k_{\spherenormal\spherenormal}\gtr\upchi-k_{\spherenormal\spherenormal}\Lunit\gtr\upchi-k_{\spherenormal\spherenormal}(\gtr\upchi)^2\\
						\notag&+\frac{1}{2}\Lunit\gtr\upchi\Chfour_{\uLunit}+\frac{1}{2}(\gtr\upchi)^2\Chfour_{\uLunit}+\frac{1}{2}\gtr\upchi\Lunit\Chfour_{\uLunit}+\uLunit(k_{\spherenormal\spherenormal})\Chfour_{\Lunit}\\
						\notag&+\rgeo^{-2}\xi+\rgeo^{-1}\angdiv\xi\\
						\notag&+\lgensmoothfunction\cdot\left(\angnabla\pfour\gfour,\angnabla\chismall\right)\cdot\angnabla\conformalfactor\\
						\notag&+\lgensmoothfunction\cdot\left(\ACC,\rgeo^{-1}\right)\cdot\pfour\gfour\cdot\angnabla\conformalfactor\\
						\notag&+\lgensmoothfunction\cdot\ACC\cdot(\angnabla\chismall,\angnabla\modtorsion,\pfour^2\gfour)\\
						\notag&+\lgensmoothfunction\cdot\pfour\boxg\gfour+\lgensmoothfunction\cdot\boxg\gfour\cdot\left(\ACC,\rgeo^{-1}\right)\\
						\notag&+\lgensmoothfunction\cdot\ACC\cdot\ACC\cdot\left(\ACC,\rgeo^{-1}\right)
						,\end{align} 
					where we enjoy a remarkable cancellation of the term $4\angnabla^2\conformalfactor\cdot\hchi$.
					
					Let M be the first 11 terms on the RHS of (\ref{EQ:Transportmodmass1}). After rearranging the terms and performing cancellations, substituting $\Lunit\gtr\upchi$ by (\ref{EQ:Ltraceupchi}), and using (\ref{EQ:riccilunitlunit}), we derive:
					\begin{align}
						\label{EQ:trashtermofmodmass}M=&-2k_{\spherenormal\spherenormal}\Lunit\gtr\upchi-k_{\spherenormal\spherenormal}(\gtr\upchi)^2+\frac{1}{2}\Lunit\gtr\upchi\Chfour_{\uLunit}+\frac{1}{2}(\gtr\upchi)^2\Chfour_{\uLunit}\\
						\notag=&\left(2k_{\spherenormal\spherenormal}-\frac{1}{2}\Chfour_{\uLunit}\right)\left(\sgabs{\hchi}^2+k_{\spherenormal\spherenormal}\gtr\upchi+\Ricfour{\Lunit}{\Lunit}\right)+\frac{1}{4}(\gtr\upchi)^2\Chfour_{\Lunit}\\
						\notag=&\rgeo^{-2}\xi+\lgensmoothfunction\cdot\boxg\gfour\cdot\pfour\gfour+\lgensmoothfunction\cdot\left(\ACC,\rgeo^{-1}\right)\cdot\ACC\cdot\pfour\gfour
						,\end{align}
					where the term $\rgeo^{-2}\xi$ in (\ref{EQ:trashtermofmodmass}) comes from the term $\frac{1}{2}(\gtr\upchi)^2\Chfour_{\uLunit}$.
					Combining (\ref{EQ:Transportmodmass1}) and (\ref{EQ:trashtermofmodmass}), we obtain: 
					\begin{align}
						\label{EQ:Transportmodmass}\Lunit\modmass+\gtr\upchi\modmass=&\rgeo^{-2}\xi+\rgeo^{-1}\angdiv\xi\\
						\notag&+\lgensmoothfunction\cdot\left(\angnabla\pfour\gfour,\angnabla\chismall\right)\cdot\angnabla\conformalfactor\\
						\notag&+\lgensmoothfunction\cdot\left(\ACC,\rgeo^{-1}\right)\cdot\pfour\gfour\cdot\angnabla\conformalfactor\\
						\notag&+\lgensmoothfunction\cdot\ACC\cdot(\angnabla\chismall,\angnabla\modtorsion,\pfour^2\gfour)\\
						\notag&+\lgensmoothfunction\cdot\pfour\boxg\gfour+\lgensmoothfunction\cdot\boxg\gfour\cdot\left(\ACC,\rgeo^{-1}\right)\\
						\notag&+\lgensmoothfunction\cdot\ACC\cdot\ACC\cdot\left(\ACC,\rgeo^{-1}\right)
						.\end{align}
					We obtain the desired equation (\ref{EQ:Lunitmodmass}).
				\end{proof}
				\begin{proof}[Proof of the Hodge systems (\ref{EQ:angdivupzeta})-(\ref{EQ:angcurlhodgemass})]
					The desired Hodge system for $\upzeta$ is derived by invoking \eqref{EQ:ALuLA}-\eqref{EQ:antisym riem curv} and equation \eqref{EQ:mass} of $\mass$ in \eqref{EQ:angdivupzet}-\eqref{EQ:angcurlupzet}. Together with the definitions of $\modtorsion$ in \eqref{DE:defmodtorsion} and $\modmass$ in \eqref{DE:modmass}, this leads to the Hodge system \eqref{EQ:angdivmodtorsion}-\eqref{EQ:angcurlmodtorsion}. Here, note that the $\mass$-terms, and hence the $\angnabla\cphi$-terms as well, are cancelled in \eqref{EQ:angdivmodtorsion}. Finally, using the definition of $\hodgemass$ in \eqref{DE:hodgemass} and \eqref{EQ:angdivmodtorsion}-\eqref{EQ:angcurlmodtorsion}, we arrive at \eqref{EQ:angdivhodgemass}-\eqref{EQ:angcurlhodgemass}. 
				\end{proof}
				\begin{proof}[Proof of a decomposition of $\hodgemass$ and Hodge-transport systems related to this decomposition in (\ref{EQ:hodgemudecomposition})-(\ref{EQ:angcurl2})] Recall definition of $\hodgemass$ (\ref{DE:hodgemass}).
					We use the commutator formula (\ref{EQ:commutatorformulafortensor}) in the equations below:
					\begin{align}\label{EQ:angAhodgemass}
						\angnabla_A\left(\angD_{\Lunit}\hodgemass+\frac{1}{2}\gtr\upchi\hodgemass\right)_B=&\angnabla_A\angD_{\Lunit}\hodgemass_B-\angD_{\Lunit}\angnabla_A\hodgemass_B+\angD_{\Lunit}\angnabla_A\hodgemass_B+\frac{1}{2}\angnabla_A\gtr\upchi\hodgemass_B+\frac{1}{2}\gtr\upchi\angnabla_A\hodgemass_B\\
						\notag=&\left(\upchi_{AB}k_{C\spherenormal}-\upchi_{AC}k_{B\spherenormal}-\Riemfour{B}{C}{\Lunit}{A}\right)\hodgemass_C\\
						\notag&+\hchi_{AC}\angnabla_C\hodgemass_B+\frac{1}{2}\angnabla_A\gtr\upchi\hodgemass_B\\
						\notag&+\gtr\upchi\angnabla_A\hodgemass_B+\angD_{\Lunit}\angnabla_A\hodgemass_B
						.\end{align}
					Taking $\gsphere$-trace of (\ref{EQ:angAhodgemass}), using (\ref{EQ:Evolutionaverage}) and (\ref{EQ:Lunitmodmass}), we derive 
					\begin{align}
						\angdiv\left(\angD_{\Lunit}\hodgemass+\frac{1}{2}\gtr\upchi\hodgemass\right)=&\left(\gtr\upchi k_{B\spherenormal}-\upchi_{AB}k_{A\spherenormal}-\Riemfour{A}{B}{\Lunit}{A}\right)\hodgemass_B\\
						\notag&+\hchi_{AB}\angnabla_B\hodgemass_A+\frac{1}{2}\angnabla_A\gtr\upchi\hodgemass_A\\
						\notag&+\frac{1}{2}\left(\Lunit\left(\modmass-\average{\modmass}\right)+\gtr\upchi\left(\modmass-\average{\modmass}\right)\right)\\
						\notag=&\left(\gtr\upchi k_{B\spherenormal}-\upchi_{AB}k_{A\spherenormal}-\Riemfour{A}{B}{\Lunit}{A}\right)\hodgemass_B\\
						\notag&+\hchi_{AB}\angnabla_B\hodgemass_A+\frac{1}{2}\angnabla_A\gtr\upchi\hodgemass_A\\
						\notag&+\frac{1}{2}\left\{\left(\Lunit\modmass+\gtr\upchi\modmass\right)-\average{\Lunit\modmass+\gtr\upchi\modmass}-\left(\gtr\upchi-\average{\gtr\upchi}\right)\average{\modmass}\right\}\\
						\notag=&\left(\gtr\upchi k_{B\spherenormal}-\upchi_{AB}k_{A\spherenormal}-\Riemfour{A}{B}{\Lunit}{A}\right)\hodgemass_B\\
						\notag&+\hchi_{AB}\angnabla_B\hodgemass_A+\frac{1}{2}\angnabla_A\gtr\upchi\hodgemass_A\\
						\notag&+\mathfrak{J}^{(1)}-\average{\mathfrak{J}^{(1)}}+\mathfrak{J}^{(2)}-\average{\mathfrak{J}^{(2)}}
						-\frac{1}{2}\left(\gtr\upchi-\average{\gtr\upchi}\right)\average{\modmass}
						.\end{align}
					Moreover, we calculate
					\begin{align}
						\angcurl\left(\angD_{\Lunit}\hodgemass+\frac{1}{2}\gtr\upchi\hodgemass\right)&=\antisymmetic^{AB}\left(-\upchi_{AC}k_{B\spherenormal}-\Riemfour{B}{C}{\Lunit}{A}\right)\hodgemass_C\\
						\notag&+\antisymmetic^{AB}\hchi_{AC}\angnabla_C\hodgemass_B+\frac{1}{2}\antisymmetic^{AB}\angnabla_A\gtr\upchi\hodgemass_B
						.\end{align}
					We decompose $\hodgemass=\massone+\masstwo$, where $\massone$ and $\masstwo$ verify the following equations:
					\begin{subequations}	
						\begin{align}
							\angdiv\left(\angD_{\Lunit}\massone+\frac{1}{2}\gtr\upchi\massone\right)=&\mathfrak{J}_{(1)}-\average{\mathfrak{J}_{(1)}},\\
							\angdiv\left(\angD_{\Lunit}\masstwo+\frac{1}{2}\gtr\upchi\masstwo\right)=&\mathfrak{J}_{(2)}-\average{\mathfrak{J}_{(2)}}
							-\frac{1}{2}\left(\gtr\upchi-\average{\gtr\upchi}\right)\average{\modmass}\\
							\notag&+\left(\gtr\upchi k_{B\spherenormal}-\upchi_{AB}k_{A\spherenormal}-\Riemfour{A}{B}{\Lunit}{A}\right)\hodgemass_B\\
							\notag&+\hchi_{AB}\angnabla_B\hodgemass_A+\frac{1}{2}\angnabla_A\gtr\upchi\hodgemass_A\\
							\notag=&\mathfrak{J}^{(2)}-\average{\mathfrak{J}^{(2)}}+\left(\gtr\upchi-\average{\gtr\upchi}\right)\average{\modmass}\\
							\notag&+\lgensmoothfunction\cdot\left(\ACC,\rgeo^{-1}\right)\cdot\pfour\gfour\cdot\hodgemass\\
							\notag&+\lgensmoothfunction\cdot\left(\angnabla\pfour\gfour,\angnabla\chismall\right)\cdot\hodgemass+\lgensmoothfunction\cdot\hchi\cdot\angnabla\hodgemass
							.\end{align}
					\end{subequations}
					
					Since $\massone$ corresponds to the terms $\gtr\upchi\angnabla_A\hodgemass_B$ and $\angD_{\Lunit}\angnabla_A\hodgemass_B$ in (\ref{EQ:angAhodgemass}), by the definition of $\hodgemass$ (\ref{DE:hodgemass}), it is obvious that 
					\begin{align}
						\angcurl\left(\angD_{\Lunit}\massone+\frac{1}{2}\gtr\upchi\massone\right)=0
						.\end{align}
					For $\masstwo$, we have:
					\begin{align}
						\angcurl\left(\angD_{\Lunit}\masstwo+\frac{1}{2}\gtr\upchi\masstwo\right)=&\angcurl\left(\angD_{\Lunit}\hodgemass+\frac{1}{2}\gtr\upchi\hodgemass\right)\\
						\notag=&\antisymmetic^{AB}\left(-\upchi_{AC}k_{B\spherenormal}-\Riemfour{B}{C}{\Lunit}{A}\right)\hodgemass_C\\
						\notag&+\antisymmetic^{AB}\hchi_{AC}\angnabla_C\hodgemass_B+\frac{1}{2}\antisymmetic^{AB}\angnabla_A\gtr\upchi\hodgemass_B\\
						\notag=&\lgensmoothfunction\cdot\left(\ACC,\rgeo^{-1}\right)\cdot\pfour\gfour\cdot\hodgemass\\
						\notag&+\lgensmoothfunction\cdot\left(\angnabla\pfour\gfour,\angnabla\chismall\right)\cdot\hodgemass+\lgensmoothfunction\cdot\hchi\cdot\angnabla\hodgemass 
						.\end{align}
					We then have obtained the desired equations (\ref{EQ:hodgemudecomposition})-(\ref{EQ:angcurl2}).
				\end{proof}
				
				\section{Control of the Causal Geometry}\label{S:ControloftheGeometry}
				In this section, we prove the estimates for the geometric quantities which are listed in Proposition \ref{PR:mainproof}.
				\subsection{Restatement of Bootstrap Assumptions and Estimates for Quantities Constructed out of the Eikonal Equation}\label{SS:sectionRestatement}
				In this section, we restate the consequence of bootstrap assumptions (\ref{BA:bt1}), followed by the bootstrap assumptions for geometry. Then we state the main estimates for the Eikonal function quantities in Proposition \ref{PR:mainproof}. The estimates in Proposition \ref{PR:mainproof} are crucial for deriving our conformal energy estimates and hence for closing the whole bootstrap argument. We prove Proposition \ref{PR:mainproof} in Section \ref{PR:mainproof} via a bootstrap argument, where the bootstrap assumptions are listed in Section \ref{SS:BAGeo}. For the complete details of the proof, we refer readers to \cite[Section 10]{3DCompressibleEuler} and \cite[Section 5-6]{AGeoMetricApproach}.
				
				\subsubsection{The fixed number $p$}\label{SSS:choiceofp}
				In the rest of the article, let $p$ denote a fixed number satisfying 
				\begin{align}
					0<\delta_0<1-\frac{2}{p}<N-3,
				\end{align}
				where $\delta_0$ is defined in Section \ref{SS:ChoiceofParameters}.
				
				\subsubsection{Bootstrap assumptions for geometric quantities} After rescaling in Section \ref{SS:SelfRescaling}, we introduce several bootstrap assumptions for the geometric quantities. These assumptions will be recovered and improved by using the estimates in Proposition \ref{PR:mainproof}\label{SS:BAGeo}\\
				
				For the rescaled unknowns defined in Definition \ref{DE:rescaledquantities}, we reformulate the bootstrap assumptions (\ref{BA:bt1}) as below:
				
				\textbf{Estimates by using bootstrap assumptions of variables}
				\begin{align}
					\label{BA:bt1.1}&\twonorms{\pfour\vvariables}{t}{2}{x}{\infty}{\region}+\lambda^{\delta_0}\sqrt{\sum\limits_{\upnu>2}\upnu^{2\delta_0}\twonorms{\littlewood\pfour\vvariables}{t}{2}{x}{\infty}{\region}^2}\lesssim\lambda^{-1/2-4\varepsilon_0}
					.\end{align}
				
				We now make a few more bootstrap assumptions regarding the geometry (defined in Section \ref{S:causalgeometry}):
				
				\textbf{Bootstrap assumptions for the geometry}
				\begin{subequations}
					\begin{align}
						\label{bt2}\max\limits_{A,B=1,2}\norm{\rgeo^{-2}\gsphere\left(\deriasphere,\deribsphere\right)-\esphere\left(\deriasphere,\deribsphere\right)}_{L^\infty_x(\region)}&\lesssim\lambda^{-\varepsilon_0},\\
						\label{bt3}\max\limits_{A,B,C=1,2}\twonorms{\dericsphere\left(\rgeo^{-2}\gsphere\left(\deriasphere,\deribsphere\right)-\esphere\left(\deriasphere,\deribsphere\right)\right)}{t}{\infty}{\sangle}{p}{\coneu}&\lesssim\lambda^{-\varepsilon_0}
						.\end{align}
				\end{subequations}
				
				Also, for $\CC=(\chismall,\hchi,\upzeta)$ defined in Definition \ref{DE:conformalchangemetric}, we have:
				\begin{subequations}
					\begin{align}
						\label{bt4}\holdertwonorms{\CC}{t}{2}{\sangle}{0}{\delta_0}{\coneu}&\lesssim\lambda^{-1/2+2\varepsilon_0},\\
						\label{bt5}\onenorm{\rgeo\CC}{\sangle}{p}{\stu}&\leq 1
						.\end{align}
				\end{subequations}
				
				Moreover, suppose that
				\begin{subequations}
					\begin{align}
						\label{bt6}\norm{\nulllapse-1}_{L^\infty_\sangle(\stu)}\leq&\frac{1}{2}
						,\\
						\label{bt6.5}\norm{\gtr\spheresecondfund-\frac{2}{\rgeo}}_{L^3(\St)}\leq& 1.
					\end{align}
				\end{subequations}
				
				Finally, we assume that the following estimates hold in the interior region $\intregion$ (defined in Section \ref{SS:opticalfunction}):
				\begin{align}
					\label{bt7}\holderthreenorms{\chismall,\hchi}{t}{2}{u}{\infty}{\sangle}{0}{\delta_0}{\intregion}&\leq\lambda^{-1/2},&
					\twonorms{\upzeta}{t}{2}{x}{\infty}{\intregion}&\leq\lambda^{-1/2},& 
					\threenorms{\angnabla\conformalfactor}{u}{2}{t}{2}{\sangle}{\infty}{\intregion}&\leq 1
					.\end{align}

				\subsection{Preliminary Estimates for Controlling the Causal Geometry}
				In this section, we provide preliminary estimates for controlling the geometry. We begin with the trace estimates in Proposition \ref{PR:sobolevinequalities}, followed by a transport lemma in Lemma \ref{LE:transport}, and Hodge estimates on $\stu$ in Lemma \ref{LE:LemmaCZ}.
				
				\begin{proposition}[Hardy-Littlewood maximal function] If $f=f(t)$ is a scalar function defined on the interval $I$, we define the Hardy-Littlewood maximal function $\mathcal{M}(f)(t)$ as follows:
					\begin{align}
						\mathcal{M}(f)(t):=\sup\limits_{t^\prime\in I\bigcup(-\infty,t)}\frac{1}{\abs{t-t^\prime}}\int_{t^\prime}^{t}f(\tau)\diff\tau
						.\end{align}
					We will use the following well-known estimate, which is valid for $1<Q\leq\infty$:
					\begin{align}
						\label{ES:hardylittlewood}\norm{\mathcal{M}(f)}_{L^Q(I)}\lesssim\norm{f}_{L^Q(I)}
						.\end{align}
				\end{proposition}
				\subsubsection{Trace Estimates}\label{SSS:TraceEstimates}
				\begin{proposition}\cite[Proposition 10.2]{3DCompressibleEuler}, \cite[Lemma 5.5 - Lemma 5.7. Calculus and elliptic estimates]{Roughsolutionsofthe3-DcompressibleEulerequations}.\label{PR:sobolevinequalities} Under the bootstrap assumptions in Section \ref{SS:bootstrap}, we have the following estimates for the $\stu$-tangent tensor fields $\upxi$:\\
					
					\noindent	\textbf{Comparison of $\stu$-norms with different volume forms}: If $1\leq Q<\infty$, there holds
					\begin{align}
						\label{ES:comparisonofnorms}\onenorm{\upxi}{\gsphere}{Q}{\stu}\approx\onenorm{\rgeo^{\frac{2}{Q}}\upxi}{\sangle}{Q}{\stu}
						.\end{align}
					\textbf{Trace inequalities:}
					\begin{align}
						\label{ES:traceinequalities}\onenorm{\rgeo^{-\frac{1}{2}}\upxi}{\gsphere}{2}{\stu}+\onenorm{\upxi}{\gsphere}{4}{\stu}\lesssim\sobolevnorm{\upxi}{1}{\St}
						.\end{align}
					\textbf{Sobolev and Morrey-type inequalities:}
					\begin{subequations}
						\begin{align}
							\label{ES:sobolevone}\twonorms{\upxi}{u}{2}{\sangle}{2}{\St}&\lesssim\sobolevnorm{\upxi}{1}{\St},\\
							\label{ES:sobolevtwo}\twonorms{\rgeo^{\frac{1}{2}}\upxi}{\sangle}{2p}{t}{\infty}{\coneu}^2&\lesssim\left(\twonorms{\rgeo\angD_{\Lunit}\upxi}{\sangle}{p}{t}{2}{\coneu}+\twonorms{\upxi}{\sangle}{p}{t}{2}{\coneu}\right)\twonorms{\upxi}{\sangle}{\infty}{t}{2}{\coneu}
							.\end{align}
					\end{subequations}
					Furthermore, if $2< Q<\infty$, the following estimates are valid:
					\begin{subequations}
						\begin{align}
							\label{ES:sobolevthree}\onenorm{\upxi}{\sangle}{Q}{\stu}&\lesssim\onenorm{\rgeo\angnabla\upxi}{\sangle}{2}{\stu}^{1-\frac{2}{Q}}\onenorm{\upxi}{\sangle}{2}{\stu}^{\frac{2}{Q}}+\onenorm{\upxi}{\sangle}{2}{\stu},\\
							\label{ES:sobolevfour}\holdernorm{\upxi}{\sangle}{0}{1-\frac{2}{Q}}{\stu}&\lesssim\onenorm{\rgeo\angnabla\upxi}{\sangle}{Q}{\stu}+\onenorm{\upxi}{\sangle}{2}{\stu}
							.\end{align}
					\end{subequations}
					In addition, if $2\leq Q$, there holds
					\begin{align}
						\label{ES:sobolevfive}\twonorms{\rgeo^{\frac{1}{2}-\frac{1}{Q}}\upxi}{\gsphere}{2Q}{u}{\infty}{\St}^2\lesssim\left(\twonorms{\rgeo(\angD_\spherenormal,\angnabla)\upxi}{\sangle}{Q}{u}{2}{\St}+\twonorms{\upxi}{\sangle}{Q}{u}{2}{\St}\right)\twonorms{\upxi}{\sangle}{\infty}{u}{2}{\St}
						.\end{align}
					Finally, if $0<1-\frac{2}{Q}<\delta\leq N-3$, then for any scalar function $f$, we have:
					\begin{align}
						\label{ES:sobolevsix}\twonorms{\rgeo f}{u}{2}{\sangle}{Q}{\St}\lesssim&\sobolevnorm{f}{N-3}{\St}.
					\end{align}
				\end{proposition}
				
				\begin{lemma}\cite[Lemma 5.11. Transport lemma]{AGeoMetricApproach}.\label{LE:transport}
					Let $m$ be a constant, and let $\upxi$ and $\mathfrak{F}$ be $\stu$-tangent tensor fields, such that the following transport equation holds along the null cone portion $\coneu\subset\region$:
					\begin{align}\label{EQ:transportlemma2}
						\angD_{\Lunit}\upxi+m\gtr\upchi\upxi=\mathfrak{F}
						.\end{align}
					Then, the following identities hold for $\tstart:=max\{u,0\}$:
					\begin{align}
						\label{EQ:vmupxi}(\volume^m\upxi)(t,u,\sangle)=&\lim\limits_{\tau\downarrow\tstart}(\volume^m\upxi)(\tau,u,\sangle)+\int_{\tstart}^{t}(\volume^m\mathfrak{F})(\tau,u,\sangle)\diff\tau,\\
						\label{EQ:rtwomupxi}(\rgeo^{2m}\upxi)(t,u,\sangle)=&\lim\limits_{\tau\downarrow\tstart}(\rgeo^{2m}\upxi)(\tau,u,\sangle)\\
						\notag&+\int_{\tstart}^{t}(\rgeo^{2m}\mathfrak{F})(\tau,u,\sangle)+m\left(\rgeo^{2m}\left(\frac{2}{\rgeo}-\gtr\upchi\right)\upxi\right)(\tau,u,\sangle)\diff\tau
						.\end{align}
					Similarly, if $\upxi$, $\mathfrak{F}$ and $\mathfrak{G}$ are $\stu$-tangent tensor fields, such that the following transport equation holds:
					\begin{align}
						\label{EQ:transportlemma3}\angD_{\Lunit}\upxi+\frac{2m}{\rgeo}\upxi=\mathfrak{G}\cdot\upxi+\mathfrak{F}
						,\end{align}
					and if 
					\begin{align}
						\label{AS:mathfrakGC}\twonorms{\mathfrak{G}}{\sangle}{\infty}{t}{1}{\coneu}\leq C,
					\end{align}
					then, under the bootstrap assumptions, the following estimates hold:
					\begin{align}
						\label{ES:r2mupxi}\sgabs{\rgeo^{2m}\upxi}(t,u,\sangle)\lesssim\lim\limits_{\tau\downarrow\tstart}\sgabs{\rgeo^{2m}\upxi}(\tau,u,\sangle)+\int_{\tstart}^t\sgabs{\rgeo^{2m}\mathfrak{F}}(\tau,u,\sangle)\diff\tau
						,\end{align}
                         where the implicit constants in the estimates above depend on the constant C on (\ref{AS:mathfrakGC}).
				\end{lemma}
				\begin{lemma}\cite[Lemma 10.9. H\"older estimates along the integral curve of $\Lunit$]{3DCompressibleEuler}.
					Let $\mathfrak{F}$ be a scalar function on $\coneu$. Let $\varphi$ be the solution of the following transport equation with initial data $\mathring{\varphi}$ given on $S_{\tstart,u}$:
					\begin{subequations}
						\begin{align}
							\Lunit\varphi=&\mathfrak{F},\\
							\varphi(\tstart,u,\sangle)=&\mathring{\varphi}(\sangle).
						\end{align}
					\end{subequations}
					Then, $\varphi$ satisfies the following estimate:
					\begin{align}
						\holdernorm{\varphi}{\sangle}{0}{\delta_0}{\stu}\lesssim\holdernorm{\mathring{\varphi}}{\sangle}{0}{\delta_0}{S_{\tstart,u}}+\int_{\tstart}^{t}\holdernorm{\mathfrak{F}}{\sangle}{0}{\delta_0}{\sTauu}\diff\tau.
					\end{align}
				\end{lemma}
				\begin{lemma}\cite[Lemma 10.10. Calderon-Zygmund and Schauder-type Hodge estimates on $\stu$]{3DCompressibleEuler}.\label{LE:LemmaCZ}
					If $\upxi$ is an $\stu$-tangent one-form and $2\leq Q\leq p<\infty$, then we have
					\begin{align}
						\label{ES:hodgeone}\onenorm{\angnabla\upxi}{\gsphere}{Q}{\stu}+\onenorm{\rgeo^{-1}\upxi}{\gsphere}{Q}{\stu}\lesssim\onenorm{\angdiv\upxi}{\gsphere}{Q}{\stu}+\onenorm{\angcurl\upxi}{\gsphere}{Q}{\stu}.
					\end{align}
					Similarly, if $\upxi$ is an $\stu$-tangent type $\binom{0}{2}$ symmetric trace-free tensor field, there holds
					\begin{align}
						\label{ES:hodgetwo}\onenorm{\angnabla\upxi}{\gsphere}{Q}{\stu}+\onenorm{\rgeo^{-1}\upxi}{\gsphere}{Q}{\stu}\lesssim\onenorm{\angdiv\upxi}{\gsphere}{Q}{\stu}.
					\end{align}
					Moreover, let $\upxi$ be an $\stu$-tangent type $\binom{0}{2}$ symmetric trace-free tensor field, let $\mathfrak{F}_{(1)}$ be a scalar function and $\mathfrak{F}_{(2)}$ be a type $\binom{0}{2}$ symmetric trace-free tensor field, and let $\mathfrak{G}$ be an $\stu$-tangent one-form. Assume that
					\begin{align}
						\angdiv\upxi=\angnabla\mathfrak{F}_{(1)}+\angdiv\mathfrak{F}_{(2)}+\mathfrak{G}.
					\end{align}
					Let $2<Q<\infty$, and $Q^\prime$ be defined by $\frac{1}{2}+\frac{1}{Q}=\frac{1}{Q^\prime}$. Then, the following estimate holds:
					\begin{align}
						\label{ES:hodgethree}\onenorm{\upxi}{\gsphere}{Q}{\stu}\lesssim\sum\limits_{i=1,2}\onenorm{\mathfrak{F}_{(i)}}{\gsphere}{Q}{\stu}+\onenorm{\mathfrak{G}}{\gsphere}{Q^\prime}{\stu}.
					\end{align}
					In addition, if $2<Q<\infty$, then we have
					\begin{equation}
						\label{ES:hodgefive}\holdernorm{\upxi}{\sangle}{0}{1-\frac{2}{Q}}{\stu}\lesssim\sum\limits_{i=1,2}\holdernorm{\mathfrak{F}_{(i)}}{\sangle}{0}{1-\frac{2}{Q}}{\stu}+\onenorm{\rgeo\mathfrak{G}}{\sangle}{Q}{\stu}.
					\end{equation}
					Similarly, assume that $\upxi$, $\mathfrak{F}_{(1)}$, and $\mathfrak{F}_{(2)}$ are $\stu$-tangent one-forms, and $\mathfrak{G}_{(1)}$ and $\mathfrak{G}_{(2)}$ are scalar functions, such that the following Hodge system is verified:
					\begin{subequations}
						\begin{align}
							\label{EQ:angdivupxi}\angdiv\upxi&=\angdiv\mathfrak{F}_{(1)}+\mathfrak{G}_{(1)},\\
							\label{EQ:angcurlupxi}\angcurl\upxi&=\angcurl\mathfrak{F}_{(2)}+\mathfrak{G}_{(2)}.
						\end{align}
					\end{subequations}
					Then, $\upxi$ satisfies estimates (\ref{ES:hodgethree}) and (\ref{ES:hodgefive}) with $\mathfrak{G}=(\mathfrak{G}_{(1)},\mathfrak{G}_{(2)})$.\\

					Finally, let $\upxi$, $\mathfrak{F}=\left(\mathfrak{F}_{(1)},\mathfrak{F}_{(2)}\right)$, and $\mathfrak{G}$ be the $\stu$ tensor fields of the type from the previous paragraphs. Assume that $\mathfrak{F}$ is the $\stu$-projection of a space-time tensor field $\tilde{\mathfrak{F}}$, or is a contraction of a space-time tensor field $\tilde{\mathfrak{F}}$ against $\Lunit$, $\uLunit$, or $\spherenormal$. If $Q>2$, $1\leq c<\infty$, and $\delta^\prime>0$ is sufficiently small, then the following estimates hold:
					\begin{align}
						\label{ES:hodgesix}\onenorm{\upxi}{\sangle}{\infty}{\stu}\lesssim\norm{\upnu^{\delta^\prime}\littlewood\vec{\tilde{\mathfrak{F}}}}_{l^c_\upnu L_{\sangle}^\infty(\stu)}+\onenorm{\vec{\tilde{\mathfrak{F}}}}{\sangle}{\infty}{\stu}+\onenorm{\rgeo^{1-\frac{2}{Q}}\mathfrak{G}}{\gsphere}{Q}{\stu},
					\end{align}
					where $\vec{\tilde{\mathfrak{F}}}$ denotes the array of (scalar) Cartesian component function of $\tilde{\mathfrak{F}}$.
				\end{lemma}
				\begin{proposition}[Frequently used estimates]\label{frequentlyusedestimates} We have the following estimates, {which will be used silently in the rest of the article}.
					Recall that the following estimates hold for $t$ and $\rgeo$:
					\begin{subequations}
						\begin{align}
							\label{ES:boundnessoft}0&<t\lesssim\lambda^{1-8\varepsilon_0},\\
							\label{ES:boundnessofrgeo2}0&<\rgeo\lesssim\lambda^{1-8\varepsilon_0}.
						\end{align}
					\end{subequations}
					Also, for $\lgensmoothfunction$, we have
					\begin{align}\label{ES:lgensmooth}
						\holderthreenorms{\lgensmoothfunction}{t}{\infty}{u}{\infty}{\sangle}{0}{\delta_0}{\region}\lesssim 1.
					\end{align}
					Moreover, for a scalar functions $f$, and $1\leq Q<p\leq\infty$ , there holds
					\begin{align}\label{ES:LQnormoff}
						\onenorm{f}{\sangle}{Q}{\stu}\leq\onenorm{f}{\sangle}{p}{\stu}\abs{\stu}_{\esphere}^{\frac{1}{Q}-\frac{1}{p}}\lesssim\onenorm{f}{\sangle}{p}{\stu}.
					\end{align}
					Finally, we have the following estimates for the average of scalar functions among $\stu$ (defined in Definition \ref{DE:average}):
					\begin{align}
						\label{ES:averagef}|\average{f}|&\lesssim\onenorm{f}{\sangle}{Q}{\stu},&
						\onenorm{\average{f}}{\sangle}{Q}{\stu}&\lesssim\onenorm{f}{\sangle}{Q}{\stu}.
					\end{align}
				\end{proposition}
				\begin{proof}[Proof of Proposition \ref{frequentlyusedestimates}]
					(\ref{ES:boundnessoft}) and (\ref{ES:boundnessofrgeo2}) are derived from (\ref{ES:lengtht}), (\ref{ES:estimatesradial}) respectively.  (\ref{ES:LQnormoff}) and (\ref{ES:averagef}) follow directly by using the H\"older's inequality and the Minkowski's inequality.
					(\ref{ES:lgensmooth}) can be proved by using (\ref{EQ:LunitLunit}) and (\ref{EQ:derinormalLunit}), we refer readers to \cite[Lemma 10.6 - Corollary 10.8]{3DCompressibleEuler} for the detailed proofs. Compared to \cite{3DCompressibleEuler}, where $\gfour$ depends only on $\vvariables$, in this article, $\gfour$ depends on both $\vvariables$ and $\cvvariables$. Nevertheless, the same argument goes through, since $\cvvariables$ is well-behaved due to Theorem \ref{TH:LinearStrichartz}.
				\end{proof}
				Now, employing the above preliminary estimates, and the energy estimates established in Section \ref{S:EnergyEstimates} and Section \ref{S:connectioncoefficientsandpdes}, we derive the following estimates for $\vvariables,\cvvariables$ and $\boxg\gfour$:
				\begin{proposition}[Estimates for $\vvariables,\cvvariables$ and $\boxg\gfour$]
					\label{estimateoffluidvariables} The rescaled variables $\vvariables,\cvvariables,\gfour$ are  defined in Definition \ref{DE:rescaledquantities}. Under the bootstrap assumptions, for any $2\leq Q\leq p$, and
					$0<\delta_0<1-\frac{2}{p}<N-2$, the following estimates hold on $\region$:
					\begin{subequations}
						\begin{align}
							\label{fv1}\twonorms{\pfour\vvariables,\pfour\cvvariables}{u}{2}{\sangle}{p}{\St},\twonorms{\rgeo^{1/2}\left(\pfour\vvariables,\pfour\cvvariables\right)}{u}{\infty}{\sangle}{2p}{\St}&\lesssim\lambda^{-1/2},\\
							\label{fv2}\twonorms{\rgeo^{1-\frac{2}{Q}}\left(\pfour^2\vvariables,\pfour^2\cvvariables\right)}{u}{2}{\gsphere}{Q}{\St}&\lesssim\lambda^{-1/2},\\
							\label{fv3}\twonorms{\pfour\vvariables,\pfour\cvvariables}{t}{2}{\sangle}{\infty}{\coneu},\twonorms{\pfour\vvariables,\pfour\cvvariables}{t}{2}{\sangle}{p}{\coneu}&\lesssim\lambda^{-1/2-4\varepsilon_0},\\
							\label{fv5}\twonorms{\rgeo\left(\pfour\vvariables,\pfour\cvvariables\right)}{t}{2}{\sangle}{\infty}{\coneu}&\lesssim\lambda^{1/2-12\varepsilon_0},\\
							\label{fv6}\norm{(\angnabla,\angD_{\Lunit})(\pfour\vvariables,\pfour\cvvariables)}_{L^2(\coneu)},\twonorms{\rgeo^{1-\frac{2}{p}}\left\{(\angnabla,\angD_{\Lunit})\pfour\vvariables,\pfour^2\cvvariables\right\}}{t}{2}{\gsphere}{p}{\coneu}&\lesssim\lambda^{-1/2},\\
							\label{fv7}\norm{\pfour^3\cvvariables}_{L^2(\coneu)},\twonorms{\rgeo^{1-\frac{2}{p}}\pfour^3\cvvariables}{t}{2}{\gsphere}{p}{\coneu}&\lesssim\lambda^{-3/2},\\
							\label{fv8}\threenorms{\rgeo^{1/2}\left(\pfour\vvariables,\pfour\cvvariables\right)}{u}{2}{t}{\infty}{\sangle}{2p}{\region}&\lesssim\lambda^{-4\varepsilon_0}.
						\end{align}
					\end{subequations}
					Moreover, we derive the following estimates for $\boxg\gfour$:
					\begin{subequations}
						\begin{align}
							\label{fv9}\holdertwonorms{\boxg \gfour}{t}{1}{x}{0}{\delta_0}{\region},\twonorms{\boxg \gfour}{t}{1}{\sangle}{p}{\coneu}\lesssim&\lambda^{-1-8\varepsilon_0},\\
							\label{fv10}\twonorms{\rgeo\boxg\gfour}{t}{2}{\sangle}{p}{\coneu}\lesssim&\lambda^{-1/2-8\varepsilon_0},\\
							\label{fv11}\twonorms{\rgeo(\angnabla,\angD_{\Lunit})\boxg\gfour}{t}{1}{\sangle}{p}{\coneu}\lesssim&\lambda^{-1-4\varepsilon_0},\\
							\label{fv12}\threenorms{\rgeo\pfour\boxg\gfour}{t}{1}{u}{2}{\sangle}{p}{\region}\lesssim&\lambda^{-1/2-8\varepsilon_0}.
						\end{align}
					\end{subequations}
				\end{proposition}
				\begin{proof}[Proof of (\ref{fv2})]
					By using (\ref{ES:sobolevsix}), and rescaling the energy estimates (\ref{ES:energyestimates1}) and (\ref{ES:energyestimates2}), we obtain
					\begin{align}
						\twonorms{\rgeo^{1-\frac{2}{Q}}\left(\pfour^2\vvariables,\pfour^2\cvvariables\right)}{u}{2}{\gsphere}{Q}{\St}\approx&\twonorms{\rgeo\left(\pfour^2\vvariables,\pfour^2\cvvariables\right)}{u}{2}{\sangle}{Q}{\St}\\
						\notag\lesssim&\sobolevnorm{\pfour^2\vvariables,\pfour^2\cvvariables}{N-3}{\St}\\
						\notag\lesssim&\lambda^{-1/2}
						.\end{align}
				\end{proof}
				\begin{proof}[Proof of (\ref{fv1})]
					We begin with estimating the first norm in (\ref{fv1}). Employing (\ref{ES:sobolevone}), (\ref{ES:sobolevthree}), (\ref{ES:sobolevsix}), (\ref{fv2}), and the energy estimates (\ref{ES:energyestimates1}) and (\ref{ES:energyestimates2}), we deduce
					\begin{align}
						\twonorms{\pfour\vvariables,\pfour\cvvariables}{u}{2}{\sangle}{p}{\St}\lesssim&\twonorms{\rgeo\angnabla\left(\pfour\vvariables,\pfour\cvvariables\right)}{u}{2}{\sangle}{2}{\St}^{1-\frac{2}{p}}\twonorms{\pfour\vvariables,\pfour\cvvariables}{u}{2}{\sangle}{2}{\St}^{\frac{2}{p}}+\twonorms{\pfour\vvariables,\pfour\cvvariables}{u}{2}{\sangle}{2}{\St}\\
						\notag\lesssim&\sobolevnorm{\pfour\vvariables,\pfour\cvvariables}{1}{\St}\\
						\notag\lesssim&\lambda^{-1/2}
						.\end{align}
					For the estimate of the second norm in (\ref{fv1}), by using (\ref{ES:sobolevfive}), (\ref{ES:sobolevfour}), (\ref{ES:sobolevsix}), (\ref{ES:sobolevone}), (\ref{fv2}), and the energy estimates (\ref{ES:energyestimates1}) and (\ref{ES:energyestimates2}), we have
					\begin{align}
						&\twonorms{\rgeo^{1/2}\left(\pfour\vvariables,\pfour\cvvariables\right)}{\sangle}{2p}{u}{\infty}{\St}^2\\
                        \notag\lesssim&\twonorms{\left(\pfour\vvariables,\pfour\cvvariables\right)}{u}{2}{\sangle}{\infty}{\St}\cdot\left(\twonorms{\rgeo(\angD_\spherenormal,\angnabla)\left(\pfour\vvariables,\pfour\cvvariables\right)}{\sangle}{Q}{u}{2}{\St}+\twonorms{\left(\pfour\vvariables,\pfour\cvvariables\right)}{\sangle}{Q}{u}{2}{\St}\right)\\
						\notag\lesssim&\left(\twonorms{\rgeo\angnabla\left(\pfour\vvariables,\pfour\cvvariables\right)}{u}{2}{\sangle}{Q}{\St}+\twonorms{\left(\pfour\vvariables,\pfour\cvvariables\right)}{u}{2}{\sangle}{2}{\St}\right)\cdot\left(\twonorms{\rgeo(\angD_\spherenormal,\angnabla)\left(\pfour\vvariables,\pfour\cvvariables\right)}{\sangle}{Q}{u}{2}{\St}+\twonorms{\left(\pfour\vvariables,\pfour\cvvariables\right)}{\sangle}{Q}{u}{2}{\St}\right)\\
						\notag\lesssim&\sobolevnorm{\pfour\vvariables,\pfour\cvvariables}{N-2}{\St}^2\\
						\notag\lesssim&\lambda^{-1}
						.\end{align}
				\end{proof}
				\begin{proof}[Proof of (\ref{fv8})]
					Using the fact that $u\lesssim\lambda^{1-8\varepsilon_0}$ in $\region$ and (\ref{fv1}), we arrive at
					\begin{align}
						\threenorms{\rgeo^{1/2}\left(\pfour\vvariables,\pfour\cvvariables\right)}{u}{2}{t}{\infty}{\sangle}{2p}{\region}\lesssim&\lambda^{1/2-4\varepsilon_0}	\threenorms{\rgeo^{1/2}\left(\pfour\vvariables,\pfour\cvvariables\right)}{u}{\infty}{t}{\infty}{\sangle}{2p}{\region}
						\\
						\notag\lesssim&\lambda^{-4\varepsilon_0}
						.\end{align}
				\end{proof}
				
				\begin{proof}[Proof of (\ref{fv3})-(\ref{fv7})]
					
					Estimates (\ref{fv3})-(\ref{fv5}) follow directly from the rescaled bootstrap assumptions (\ref{BA:bt1})-(\ref{BA:bt3}).
					
					The first estimates of (\ref{fv6})-(\ref{fv7}) are direct results of energy estimates along null hypersurfaces (\ref{ES:coneenergy1})-(\ref{ES:coneenergy3}). 
					
					The proof of the second estimates of (\ref{fv6})-(\ref{fv7}) is similar to that in \cite[Lemma 5.5]{AGeoMetricApproach}. We use the Sobolev inequality (\ref{ES:sobolevthree}) and the following inequality in 
					\cite[Proposition 2.7]{Causalgeometryofrougheinsteincmcshspacetime}:
					\begin{align}
						\sum\limits_{l>1}\twonorms{l^{N-2}\anglittlewood_{l}(\angD_{\Lunit},\angnabla)f}{t}{2}{\gsphere}{2}{\coneu}^2\lesssim&\sum\limits_{\upnu>1}\upnu^{2(N-2)}\left(\mathcal{F}_{(wave)}[\littlewood f;\coneu]+\sobolevnorm{\littlewood f}{1}{\St}\right)\\
						\notag&+\mathcal{F}_{(wave)}[\littlewood f;\coneu]+\sobolevnorm{f}{1}{\St}^2,
					\end{align}
					where $\anglittlewood_{l}$ is the Littlewood-Paley projection operator on $\stu$. We obtain the results by applying the energy estimates (\ref{ES:energyestimates1})(\ref{ES:energyestimates2}) (along constant-time hypersurfaces) and (\ref{ES:coneenergy1})-(\ref{ES:coneenergy3}) (along null hypersurfaces).
				\end{proof}
				\begin{proof}[Proof of (\ref{fv9})]
					For the first estimate in (\ref{fv9}), we consider the following algebraic expression:
					\begin{align}\label{EQ:controlofboxg1.1}
						\boxg\cvvariables=(\pfour\vvariables,\pfour\cvvariables)\cdot\pfour\cvvariables+(\vvariables,\cvvariables)\cdot\pfour^2\cvvariables.
					\end{align} 
					Combining (\ref{EQ:EW5}) and (\ref{EQ:controlofboxg1.1}), we have: 
					\begin{align}\label{EQ:controlofboxg1.2}
						\boxg\gfour_{\alpha\beta}=\gensmoothfunction(\vvariables,\cvvariables)\cdot(\pfour\vvariables,\pfour\cvvariables)\cdot(\pfour\vvariables,\pfour\cvvariables)+\left(\gensmoothfunction(\vvariables,\cvvariables)+1\right)\cdot\pfour^2\cvvariables.
					\end{align}
					We consider the $L^1_tC_{x}^{0,\delta_0}$-norm of the first term on the RHS of (\ref{EQ:controlofboxg1.2}). By (\ref{BA:databa}), (\ref{BA:bt1}) and (\ref{BA:bt2.1}), we derive:
					\begin{align}
						\holdertwonorms{\gensmoothfunction(\vvariables,\cvvariables)\cdot(\pfour\vvariables,\pfour\cvvariables)\cdot(\pfour\vvariables,\pfour\cvvariables)}{t}{1}{x}{0}{\delta_0}{\region}\lesssim&\holdertwonorms{\pfour\cvvariables}{t}{2}{x}{0}{\delta_0}{\region}\holdertwonorms{\pfour\vvariables}{t}{2}{x}{0}{\delta_0}{\region}\lesssim\lambda^{-1-8\varepsilon_0}.
					\end{align}
					Then, for the second term, employing (\ref{BA:databa}), (\ref{ES:Trescale}) and (\ref{BA:bt3.1}), we obtain:
					\begin{align}
						\holdertwonorms{\left(\gensmoothfunction(\vvariables,\cvvariables)+1\right)\cdot\pfour^2\cvvariables}{t}{1}{x}{0}{\delta_0}{\region}\lesssim&\holdertwonorms{\pfour^2\cvvariables}{t}{2}{x}{0}{\delta_0}{\region}\Trescale^{1/2}\lesssim\lambda^{-1-12\varepsilon_0},
					\end{align}
					which then yields the first estimate in (\ref{fv9}). In a similar fashion, using \eqref{fv3}, we obtain the second estimates in \eqref{fv9}.
				\end{proof}
				\begin{proof}[Proof of (\ref{fv10})]
					By using (\ref{fv1}), (\ref{EQ:controlofboxg1.2}), and the bootstrap assumptions (\ref{BA:bt1})-(\ref{BA:bt3}), we have:
					\begin{align}
						\twonorms{\rgeo\boxg\gfour}{t}{2}{\sangle}{p}{\coneu}\lesssim&\lambda^{1/2-4\varepsilon_0}\twonorms{\rgeo^{1/2}(\pfour\vvariables,\pfour\cvvariables)}{t}{\infty}{\sangle}{p}{\coneu}\cdot\twonorms{\pfour\vvariables,\pfour\cvvariables}{t}{2}{\sangle}{\infty}{\coneu}+\lambda^{1-8\varepsilon_0}\twonorms{\pfour^2\cvvariables}{t}{2}{\sangle}{\infty}{\coneu}\\
						\notag\lesssim&\lambda^{-1/2-8\varepsilon_0}.
					\end{align}
				\end{proof}
				\begin{proof}[Proof of (\ref{fv11})]
					First, we take $(\angnabla,\angD_{\Lunit})$ of (\ref{EQ:controlofboxg1.2}) as below:
					\begin{align}\label{ES:fv12.0}
						(\angnabla,\angD_{\Lunit})\boxg\gfour=&\gensmoothfunction(\vvariables,\cvvariables)\cdot\left\{(\angnabla,\angD_{\Lunit})\pfour\vvariables,\pfour^2\cvvariables\right\}\cdot(\pfour\vvariables,\pfour\cvvariables)+\gensmoothfunction(\vvariables,\cvvariables)\cdot(\pfour\vvariables,\pfour\cvvariables)\cdot(\pfour\vvariables,\pfour\cvvariables)\cdot(\pfour\vvariables,\pfour\cvvariables)\\
						\notag&+\gensmoothfunction(\vvariables,\cvvariables)\cdot(\angnabla,\angD_{\Lunit})\pfour^2\cvvariables
					\end{align}
					Applying the bootstrap assumptions (\ref{BA:bt1}), the estimates (\ref{BA:bt2}), (\ref{fv1}), (\ref{fv6}), and (\ref{fv7}), we derive
					\begin{subequations}\label{ES:fv12.1}
						\begin{align}
							\twonorms{\rgeo\gensmoothfunction(\vvariables,\cvvariables)\cdot\left\{(\angnabla,\angD_{\Lunit})\pfour\vvariables,\pfour^2\cvvariables\right\}\cdot(\pfour\vvariables,\pfour\cvvariables)}{t}{1}{\sangle}{p}{\coneu}\lesssim&\twonorms{\pfour\vvariables,\pfour\cvvariables}{t}{2}{\sangle}{\infty}{\coneu}\twonorms{\rgeo\left\{(\angnabla,\angD_{\Lunit})\pfour\vvariables,\pfour^2\cvvariables\right\}}{t}{2}{\sangle}{p}{\coneu}\\
							\notag\lesssim&\lambda^{-1-4\varepsilon_0},\\
							\twonorms{\rgeo\gensmoothfunction(\vvariables,\cvvariables)\cdot(\pfour\vvariables,\pfour\cvvariables)\cdot(\pfour\vvariables,\pfour\cvvariables)\cdot(\pfour\vvariables,\pfour\cvvariables)}{t}{1}{\sangle}{p}{\coneu}\lesssim&\lambda^{1/2-4\varepsilon_0}\twonorms{\rgeo^{1/2}(\pfour\vvariables,\pfour\cvvariables)}{t}{\infty}{\sangle}{p}{\coneu}\twonorms{\pfour\vvariables,\pfour\cvvariables}{t}{2}{\sangle}{\infty}{\coneu}^2\\
							\notag\lesssim&\lambda^{-1-12\varepsilon_0},\\
							\twonorms{\rgeo\gensmoothfunction(\vvariables,\cvvariables)\cdot(\angnabla,\angD_{\Lunit})\pfour^2\cvvariables}{t}{1}{\sangle}{p}{\coneu}\lesssim&\lambda^{1/2-4\varepsilon_0}\twonorms{\rgeo(\angnabla,\angD_{\Lunit})\pfour^2\cvvariables}{t}{2}{\sangle}{p}{\coneu}\\
							\notag\lesssim&\lambda^{-1-4\varepsilon_0}.
						\end{align}
					\end{subequations}
					The desired estimate (\ref{fv11}) follows by combining (\ref{ES:fv12.0})-(\ref{ES:fv12.1}).
				\end{proof}
				\begin{proof}[Proof of (\ref{fv12})]
					Taking one derivative of (\ref{EQ:controlofboxg1.2}), we have:
					\begin{align}\label{ES:fv122}
						\pfour\boxg\gfour=&\gensmoothfunction(\vvariables,\cvvariables)\cdot(\pfour^2\vvariables,\pfour^2\cvvariables)\cdot(\pfour\vvariables,\pfour\cvvariables)+\left(\gensmoothfunction(\vvariables,\cvvariables)+1\right)\cdot\pfour^3\cvvariables\\
						\notag&+\gensmoothfunction(\vvariables,\cvvariables)\cdot(\pfour\vvariables,\pfour\cvvariables)\cdot(\pfour\vvariables,\pfour\cvvariables)\cdot(\pfour\vvariables,\pfour\cvvariables).
					\end{align}
					For the first two terms on the RHS of (\ref{ES:fv122}), by (\ref{ES:sobolevsix}), the bootstrap assumptions (\ref{BA:bt1}), estimate (\ref{BA:bt2}), and rescaling the energy estimates (\ref{ES:energyestimates1}) and (\ref{ES:energyestimates2}), we derive:
					\begin{subequations}
						\begin{align}\label{ES:fv121}
							\threenorms{\rgeo\gensmoothfunction(\vvariables,\cvvariables)\cdot(\pfour^2\vvariables,\pfour^2\cvvariables)\cdot(\pfour\vvariables,\pfour\cvvariables)}{t}{1}{u}{2}{\sangle}{p}{\region}\lesssim&\sobolevtwonorms{\pfour^2\vvariables,\pfour^2\cvvariables}{t}{2}{x}{N-3}{\region}\cdot\twonorms{\pfour\vvariables,\pfour\cvvariables}{t}{2}{x}{\infty}{\region}\lesssim\lambda^{-1/2-8\varepsilon_0},\\
							\threenorms{\rgeo\left(\gensmoothfunction(\vvariables,\cvvariables)+1\right)\cdot\pfour^3\cvvariables}{t}{1}{u}{2}{\sangle}{p}{\region}\lesssim&\sobolevtwonorms{\pfour^3\cvvariables}{t}{1}{x}{N-3}{\region}\lesssim\lambda^{-1/2-8\varepsilon_0}.
						\end{align}
					\end{subequations}
					For the last term on the RHS of (\ref{ES:fv122}), by utilizing the second estimate of (\ref{fv1}), the bootstrap assumptions (\ref{BA:bt1}), and estimate (\ref{BA:bt2}), we deduce
					\begin{align}\label{ES:fv123}
						&\threenorms{\rgeo\gensmoothfunction(\vvariables,\cvvariables)\cdot(\pfour\vvariables,\pfour\cvvariables)\cdot(\pfour\vvariables,\pfour\cvvariables)\cdot(\pfour\vvariables,\pfour\cvvariables)}{t}{1}{u}{2}{\sangle}{p}{\region}\\
						\notag\lesssim&\lambda^{1/2-4\varepsilon_0}\threenorms{\rgeo^{1/2}(\pfour\vvariables,\pfour\cvvariables)}{t}{\infty}{u}{\infty}{\sangle}{2p}{\region}^2\cdot\threenorms{\pfour\vvariables,\pfour\cvvariables}{t}{1}{u}{\infty}{\sangle}{\infty}{\region}\\
						\notag\lesssim&\lambda^{-1/2-12\varepsilon_0}.
					\end{align}
					Combining (\ref{ES:fv121}) and (\ref{ES:fv123}), we conclude the desired result (\ref{fv12}). 
				\end{proof}
				
				\subsection{Main Estimates for the Causal Geometry}
				\begin{proposition}[The main estimates for the eikonal function quantities] Under the bootstrap assumptions, for the geometric quantities defined in Definition \ref{DE:DEFSOFCONNECTIONCOEFFICIENTS}, Definition \ref{DE:conformalchangemetric} and Definition \ref{DE:mass}, we have the following estimates for $2<q\leq4$:\label{PR:mainproof}\\ 
					
					\noindent\textbf{Estimates for connection coefficients:}
					\begin{subequations}
						\begin{align}
							\label{mr1}\twonorms{\CC}{t}{2}{\sangle}{p}{\coneu},\twonorms{\rgeo\angD_{\Lunit}\CC}{t}{2}{\sangle}{p}{\coneu},\twonorms{\rgeo^{1/2}\CC}{t}{\infty}{\sangle}{p}{\coneu}&\lesssim\lambda^{-1/2},\\
							\label{mr3}\twonorms{\rgeo\CC}{t}{\infty}{\sangle}{p}{\coneu}&\lesssim\lambda^{-4\varepsilon_0}
							,\end{align}
					\end{subequations}
                    and
					\begin{subequations}
						\begin{align}
	\label{mr4}\rgeo\Restrace\reschi&\approx 1,&\twonorms{\rgeo^{1/2}\chismall}{u}{\infty}{\sangle}{2p}{\St}&\lesssim\lambda^{-1/2},\\
							\label{mr5}\norm{\rgeo^{1/2}\chismall}_{L^\infty(\region)}&\lesssim\lambda^{-1/2},&\threenorms{\rgeo^{3/2}\angnabla\chismall}{t}{\infty}{u}{\infty}{\sangle}{p}{\region}&\lesssim\lambda^{-1/2},\\
							\label{mr7}\twonorms{\rgeo\left(\angnabla\CC,\mass\right)}{t}{2}{\sangle}{p}{\coneu}&\lesssim\lambda^{-1/2},&\holdertwonorms{\CC}{t}{2}{\sangle}{0}{\delta_0}{\coneu}&\lesssim\lambda^{-1/2}
							.\end{align}
					\end{subequations}
					In addition, the null lapse $\nulllapse$ verifies the following estimates:
					\begin{align}
						\label{mr9}\twonorms{\frac{\nulllapse^{-1}-1}{\rgeo}}{t}{2}{x}{\infty}{\region},\threenorms{\frac{\nulllapse^{-1}-1}{\rgeo^{1/2}}}{t}{\infty}{u}{\infty}{\sangle}{2p}{\region},\twonorms{\rgeo(\angD_{\Lunit},\angnabla)\left(\frac{\nulllapse^{-1}-1}{\rgeo}\right)}{t}{2}{\sangle}{p}{\coneu}&\lesssim\lambda^{-1/2}
						.\end{align}
					Moreover, we obtain
					\begin{align}
						\label{mr10}\holdertwonorms{\lgensmoothfunction}{t}{\infty}{x}{0}{\delta_0}{\region}&\lesssim 1
						.\end{align}
					Furthermore, there hold
					\begin{subequations}
						\begin{align}
							\label{mr11}\holderthreenorms{\chismall,\hchi,\gtr\upchi-\frac{2}{\rgeo}}{t}{\frac{q}{2}}{u}{\infty}{\sangle}{0}{\delta_0}{\region}\lesssim&\lambda^{\frac{2}{q}-1-4\varepsilon_0\left(\frac{4}{q}-1\right)},\\
							\label{mr115}\twonorms{\upzeta}{t}{\frac{q}{2}}{x}{\infty}{\region}\lesssim&\lambda^{\frac{2}{q}-1-4\varepsilon_0\left(\frac{4}{q}-1\right)},\\
							\label{mr1152}\norm{\gtr\spheresecondfund-\frac{2}{\rgeo}}_{L^3(\St)}\lesssim& \lambda^{-4\varepsilon_0}
							.\end{align}
					\end{subequations}
					\textbf{Improved estimates in the interior region:}
					\begin{subequations}
						\begin{align}
							\label{mr12}\twonorms{\frac{\nulllapse^{-1}-1}{\rgeo}}{t}{2}{x}{\infty}{\intregion}&\lesssim\lambda^{-1/2-4\varepsilon_0},\\
							\label{mr13}\twonorms{\rgeo^{1/2}\CC}{\sangle}{2p}{t}{\infty}{\coneu}&\lesssim\lambda^{-1/2},& \text{if} \ \coneu&\subset\intregion,\\
							\label{mr14}\holderthreenorms{\chismall,\hchi,\gtr\upchi-\frac{2}{\rgeo}}{t}{2}{u}{\infty}{\sangle}{0}{\delta_0}{\intregion}&\lesssim\lambda^{-1/2-4\varepsilon_0},\\
							\label{mr145}\twonorms{\upzeta}{t}{2}{x}{\infty}{\intregion}&\lesssim\lambda^{-1/2-4\varepsilon_0}.
						\end{align}
					\end{subequations}
					\textbf{Estimates for the geometric angular coordinate components of $\gsphere$:}
					\begin{subequations}
						\begin{align}
							\label{mr15}\max\limits_{A,B=1,2}\norm{\rgeo^{-2}\gsphere\left(\deriasphere,\deribsphere\right)-\esphere\left(\deriasphere,\deribsphere\right)}_{L^\infty(\region)}&\lesssim\lambda^{-4\varepsilon_0},\\
							\label{mr16}\max\limits_{A,B,C=1,2}\twonorms{\dericsphere\left(\rgeo^{-2}\gsphere\left(\deriasphere,\deribsphere\right)-\esphere\left(\deriasphere,\deribsphere\right)\right)}{t}{\infty}{\sangle}{p}{\coneu}&\lesssim\lambda^{-4\varepsilon_0}
							.\end{align}
					\end{subequations}
					\textbf{Estimates for $\volume$ and $\nulllapse$:}
					\begin{subequations}
						\begin{align}
							\label{mr17}\volume&:=\frac{\sqrt{\mathrm{det}\gsphere}}{\sqrt{\mathrm{det}\esphere}}\approx\rgeo^2,\\
							\label{mr18}\norm{\nulllapse-1}_{L^\infty(\region)}&\lesssim\lambda^{-4\varepsilon_0}<\frac{1}{4}
							.\end{align}
					\end{subequations}
					Furthermore, we have
					\begin{align}
						\label{mr19}\threenorms{\rgeo^{1/2}\angnabla\cphi}{t}{\infty}{u}{\infty}{\sangle}{p}{\region}\lesssim\lambda^{-1/2},\  \twonorms{\angnabla\cphi}{t}{2}{\sangle}{p}{\coneu}\lesssim\lambda^{-1/2},\ 
						\twonorms{\rgeo\Lunit\angnabla\cphi}{t}{2}{\sangle}{p}{\coneu}&\lesssim\lambda^{-1/2}
						.\end{align}
					\textbf{Interior region estimates for $\conformalfactor$:}
					\begin{subequations}
						\begin{align}
							\label{mr21}\twonorms{\rgeo^{1/2}\Lunit\conformalfactor}{t}{\infty}{\sangle}{2p}{\coneu}\lesssim&\lambda^{-1/2-2\varepsilon_0},& \twonorms{\rgeo^{1/2}\angnabla\conformalfactor}{t}{\infty}{\sangle}{p}{\coneu}, \twonorms{\angnabla\conformalfactor}{t}{2}{\sangle}{p}{\coneu}&\lesssim\lambda^{-1/2},&  \text{if}\ \coneu&\subset\intregion,\\
							\label{mr22}\norm{\conformalfactor}_{L^\infty\left(\intregion\right)}&\lesssim\lambda^{-8\varepsilon_0},&\norm{\rgeo^{-1/2}\conformalfactor}_{L^\infty\left(\intregion\right)}&\lesssim\lambda^{-1/2-4\varepsilon_0}
							.\end{align}
					\end{subequations}
					\textbf{Interior region estimates for $\conformalfactor$, $\mass$, $\modtorsion$ and $\hodgemass$:}
					\begin{align}
						\label{mr24}\holderthreenorms{\angnabla\conformalfactor}{u}{2}{t}{2}{\sangle}{0}{\delta_0}{\intregion},\threenorms{\rgeo\modmass,\rgeo\angnabla\modtorsion}{u}{2}{t}{2}{\sangle}{p}{\intregion},\threenorms{\rgeo^{\frac{3}{2}}\modmass}{u}{2}{t}{\infty}{\sangle}{p}{\intregion}&\lesssim\lambda^{-4\varepsilon_0}
						.\end{align}
					In addition, we obtain
					\begin{align}
						\label{mr26}\threenorms{\rgeo\angnabla\hodgemass,\hodgemass}{t}{2}{u}{2}{\sangle}{p}{\intregion},\holderthreenorms{\hodgemass}{t}{2}{u}{2}{\sangle}{0}{\delta_0}{\intregion}&\lesssim\lambda^{-4\varepsilon_0}
						.\end{align}
					\textbf{Decomposition of $\angnabla\conformalfactor$ and corresponding estimates in the interior region}: In $\intregion$, we can decompose $\angnabla\conformalfactor$ as follows:
					\begin{align}\label{decompofangnablasigma}
						\angnabla\conformalfactor=-\upzeta+(\modtorsion-\hodgemass)+\massone+\masstwo
						,\end{align} 
					where the following asymptotic conditions near the cone-tip axis are satisfied:
					\begin{align}
						\label{mr265}\rgeo\massone(t,u,\sangle)=\zero(\rgeo), \ \rgeo\masstwo(t,u,\sangle)=\zero(\rgeo), \quad \text{as}\ t\downarrow u
						.\end{align}
					Moreover, we have
					\begin{subequations}
						\begin{align}
							\label{mr27}\threenorms{\modtorsion-\hodgemass}{t}{2}{u}{\infty}{\sangle}{\infty}{\intregion},\threenorms{\massone}{t}{2}{u}{\infty}{\sangle}{\infty}{\intregion}&\lesssim\lambda^{-1/2-4\varepsilon_0},\\
							\label{mr28}\threenorms{\masstwo}{u}{2}{t}{\infty}{\sangle}{\infty}{\intregion}&\lesssim\lambda^{-1/2-4\varepsilon_0}
							.\end{align}
					\end{subequations}
				\end{proposition}
				
				\subsection{Proof of Proposition \ref{PR:mainproof}}\label{SS:controlacoustic} With the help of the previous results, we are now ready to control the geometric quantities. Recall that $\tstart:=max\{u,0\}$. {\color{black}Before delving into the lengthy details, we first present a table outlining the steps and the underlying logic of the proof.
\newpage

					\begin{table}[h!]
						\centering
						\begin{tabular}{|c|p{5cm}|c|p{4cm}|}
							\hline
							\textbf{Order} & \textbf{Estimates} & \textbf{Depends on}  \\ \hline
							1 &Estimate in \eqref{mr17} & --   \\ \hline
							2 &Estimate in \eqref{mr18} & --   \\ \hline
							3 &Estimate of $\twonorms{\hchi}{t}{2}{\sangle}{p}{\coneu}$ in \eqref{mr1} & --  \\ \hline
							4 &Estimate of $\twonorms{\rgeo^{1/2}\hchi}{t}{\infty}{\sangle}{p}{\coneu}\lesssim\lambda^{-1/2}$ in \eqref{mr1} & --  \\ \hline
							5 &Estimate of $\twonorms{\rgeo\angD_{\Lunit}\hchi}{t}{2}{\sangle}{p}{\coneu}$ in \eqref{mr1} &  Estimates of \eqref{mr1} in Step 3,4\\ \hline
							6 &Estimate of $\twonorms{\upzeta}{t}{2}{\sangle}{p}{\coneu}$ in \eqref{mr1} & Estimate of \eqref{mr1} in Step 4 \\ \hline
							7 &Estimate of $\twonorms{\rgeo^{1/2}\upzeta}{t}{\infty}{\sangle}{p}{\coneu}$ in \eqref{mr1} &-- \\ \hline
							8 &Estimate of $\twonorms{\rgeo\angD_{\Lunit}\upzeta}{t}{2}{\sangle}{p}{\coneu}$ in \eqref{mr1}& Estimates of \eqref{mr1} in Step 6, 7 \\ \hline
							9 &Estimate $\twonorms{\rgeo^{1/2}(\hchi,\upzeta)}{\sangle}{2p}{t}{\infty}{\coneu}$ in \eqref{mr13}& Estimates in \eqref{mr1} \\ \hline
							10 &Estimate of $\rgeo\Restrace\reschi$ in \eqref{mr4}& -- \\ \hline
							11 &Estimate of $\norm{\rgeo^{1/2}\chismall}_{L^\infty(\region)}$ in \eqref{mr5}& -- \\ \hline
							12 &Estimates of $\twonorms{\rgeo^{1/2}\chismall}{t}{\infty}{\sangle}{p}{\coneu}$ in \eqref{mr1},  $\twonorms{\rgeo^{1/2}\chismall}{\sangle}{2p}{u}{\infty}{\St}$ in \eqref{mr4}, $\twonorms{\rgeo^{1/2}\chismall}{\sangle}{2p}{t}{\infty}{\coneu}$ in \eqref{mr13}& Estimate of \eqref{mr5} in step 11 \\ \hline
							13 &Estimate of $\norm{\gtr\spheresecondfund-\frac{2}{\rgeo}}_{L^3(\St)}$ in \eqref{mr1152} & Estimate of \eqref{mr4} in step 12 \\ \hline
							14 &Estimate of $\holdertwonorms{\chismall}{t}{2}{\sangle}{0}{\delta_0}{\coneu}$ in \eqref{mr7} &  \\ \hline
							15 &Estimate $\twonorms{\chismall}{t}{2}{\sangle}{p}{\coneu}$ in \eqref{mr1}& Estimate of \eqref{mr7} in step 14\\ \hline
							16 &Estimate of $\twonorms{\rgeo\angD_{\Lunit}\chismall}{t}{2}{\sangle}{p}{\coneu}$ in \eqref{mr1}& Estimates of $\twonorms{\rgeo^{1/2}\CC}{t}{\infty}{\sangle}{p}{\coneu}$ in \eqref{mr1}\\ \hline
							17 &Estimate $\twonorms{\rgeo\left(\angnabla\chismall,\angnabla\hchi\right)}{t}{2}{\sangle}{p}{\coneu}$ in \eqref{mr7}& \\ \hline
							18 &Estimate $\threenorms{\rgeo^{3/2}\angnabla\chismall}{t}{\infty}{u}{\infty}{\sangle}{p}{\region}$ in \eqref{mr5}& Estimate \eqref{mr7} for $\twonorms{\rgeo\angnabla\hchi}{t}{2}{\sangle}{p}{\coneu}$ in step 17\\ \hline
							19 &Estimate of $\holdertwonorms{\hchi}{t}{2}{\sangle}{0}{\delta_0}{\coneu}$ in \eqref{mr7}& Estimates \eqref{mr7} for $\rgeo\angnabla\hchi$ and \eqref{mr1} for $\hchi$\\ \hline
							20 &Estimate \eqref{mr14}&Proven result \eqref{mr10} (in Proposition \ref{frequentlyusedestimates})\\ \hline
							21 &Estimate \eqref{mr11}& \\ \hline
							22 &Estimates \eqref{mr15} and \eqref{mr16}& Estimates \eqref{mr1}, \eqref{mr7} and \eqref{mr11}\\ \hline
							23 &Estimate of $\twonorms{\frac{\nulllapse^{-1}-1}{\rgeo}}{t}{2}{x}{\infty}{\region}$ in \eqref{mr9}& \\ \hline
                                    \end{tabular}
					\end{table}
                    \newpage
                    \begin{table}[h!]
						\centering
						\begin{tabular}{|c|p{5cm}|c|p{4cm}|}
							\hline
							\textbf{Order} & \textbf{Estimates} & \textbf{Depends on}  \\ \hline
							24 &Estimate  $\twonorms{\rgeo(\angD_{\Lunit},\angnabla)\left(\frac{\nulllapse^{-1}-1}{\rgeo}\right)}{t}{2}{\sangle}{p}{\coneu}$ in \eqref{mr9}& Estimates  \eqref{mr1} and \eqref{mr18}\\ \hline
							25 &Estimate  $\threenorms{\frac{\nulllapse^{-1}-1}{\rgeo^{1/2}}}{t}{\infty}{u}{\infty}{\sangle}{2p}{\region}$ in \eqref{mr9}& Estimates of \eqref{mr9} in step 23 and 24\\ \hline   
							26 &Estimates in \eqref{mr19}& Proven estimates in \eqref{mr7}, \eqref{mr11} and \eqref{mr17}\\ \hline  
							27 &Estimates $\twonorms{\rgeo\angnabla\upzeta}{t}{2}{\sangle}{p}{\coneu},\holdertwonorms{\upzeta}{t}{2}{\sangle}{0}{\delta_0}{\coneu}$ in \eqref{mr7}&Estimates \eqref{mr1}, \eqref{mr3}, \eqref{mr11} and \eqref{mr19}\\ \hline   
							28 &Estimates of $\twonorms{\rgeo\mass}{t}{2}{\sangle}{p}{\coneu}$ in \eqref{mr7}&Estimates \eqref{mr1}, \eqref{mr3}, \eqref{mr7} for $\CC$ and \eqref{mr19}\\ \hline  
							29 &Estimates of $\twonorms{\upzeta}{t}{\frac{q}{2}}{x}{\infty}{\region}$ in \eqref{mr115} and $\twonorms{\upzeta}{t}{2}{x}{\infty}{\intregion}$ in \eqref{mr145}&Estimates \eqref{mr3}, \eqref{mr11} \eqref{mr14} and \eqref{mr19}\\ \hline  
							30 &Estimates of \eqref{mr21}-\eqref{mr22}&Estimates \eqref{mr7}\\ \hline 
							31 &Estimate  $\threenorms{\rgeo\left(\modmass,\angnabla\modtorsion\right)}{u}{2}{t}{2}{\sangle}{p}{\intregion}$ in \eqref{mr24}& Estimates \eqref{mr3}, \eqref{mr7}, \eqref{mr14}, \eqref{mr145} and \eqref{mr21} \\ \hline 
							32 &Estimate $\threenorms{\rgeo^{3/2}\modmass}{u}{2}{t}{\infty}{\sangle}{p}{\intregion}$ in \eqref{mr24}&  \\ \hline
							33 &Estimate $\holderthreenorms{\angnabla\conformalfactor}{u}{2}{t}{2}{\sangle}{0}{\delta_0}{\intregion}$ in \eqref{mr24}&Estimates \eqref{mr7}, \eqref{mr21}, \eqref{mr24} for $\threenorms{\rgeo\angnabla\modtorsion}{u}{2}{t}{2}{\sangle}{p}{\intregion}$  \\ \hline
							34 &Estimate \eqref{mr26}&Estimates in \eqref{mr24}  \\ \hline
							35 &Estimate $\threenorms{\modtorsion-\hodgemass}{t}{2}{u}{\infty}{\sangle}{\infty}{\intregion}$ in \eqref{mr27}&  \\ \hline
							36 &Estimate \eqref{mr265}&  \\ \hline
							37 &Estimate $\threenorms{\massone}{t}{2}{u}{\infty}{\sangle}{\infty}{\intregion}$ in \eqref{mr27}& Estimates \eqref{mr14} and \eqref{mr265}  \\ \hline
							38 &Estimate \eqref{mr28}& Estimates \eqref{mr7}, \eqref{mr14}, \eqref{mr24} and \eqref{mr26}  \\ \hline
						\end{tabular}
						\caption{Outline of proof for Proposition \ref{PR:mainproof}}
					\end{table}
					Following the order presented in the above table, we now prove the corresponding estimates in turn.}
				\begin{proof}[Proof of $\volume\approx\rgeo^2$ in (\ref{mr17})]
					Recall equation (\ref{EQ:Lunitcphi}), there holds
					\begin{align}\label{EQ:Lrvolume}
						\Lunit\left(\rgeo^{-2}\volume\right)=\left(\chismall-\Chfour_{\Lunit}\right)\left(\rgeo^{-2}\volume\right).
					\end{align}
					Integrating (\ref{EQ:Lrvolume}) along the integral curve of $\Lunit$, we deduce
					\begin{align}
						\left(\rgeo^{-2}\volume\right)(t,u,\sangle)=\lim\limits_{\tau\downarrow\tstart}\left(\rgeo^{-2}\volume\right)(\tau,u,\sangle)+\int_{\tstart}^t\left(\chismall-\Chfour_{\Lunit}\right)\left(\rgeo^{-2}\volume\right)\diff\tau,
					\end{align}
					where $\tstart:=\text{max}\{u,0\}$.
					Using the Gr\"onwall's inequality, the bootstrap assumptions for $\chismall$ and $\pfour\vvariables$ in Section \ref{SS:BAGeo}, the $L^2_tL_x^\infty$ bound for $\pfour\cvvariables$ (\ref{BA:bt2}), and the initial conditions for $\rgeo^{-2}\volume$ (\ref{ES:lapseminusone}), (\ref{ES:initialfoliation3}) and (\ref{ES:initialcontip3}), we obtain:
					\begin{align}
						\left(\rgeo^{-2}\volume\right)(t,u,\sangle)\lesssim\lim\limits_{\tau\downarrow\tstart}\left(\rgeo^{-2}\volume\right)\exp\left(\twonorms{\chismall-\Chfour_{\Lunit}}{t}{1}{\sangle}{\infty}{\coneu}\right)\lesssim 1.
					\end{align}
					We get $\left(\rgeo^{-2}\volume\right)(t)\gtrsim 1$ by applying a similar argument to: 
					\begin{align}
						\left(-\rgeo^{-2}\volume\right)(t,u,\sangle)=\lim\limits_{\tau\downarrow\tstart}\left(-\rgeo^{-2}\volume\right)+\int_{\tstart}^t\left(\chismall-\Chfour_{\Lunit}\right)\left(-\rgeo^{-2}\volume\right)\diff\tau.
					\end{align}
				\end{proof}
				\begin{proof}[Proof of $\norm{\nulllapse-1}_{L^\infty(\region)}\lesssim\lambda^{-4\varepsilon_0}<\frac{1}{4}$ in (\ref{mr18})]
					Recall equation (\ref{EQ:Lunitnulllapse}). Taking the initial conditions for the two separate cases $u<0$ and $u\geq 0$ into consideration, we derive that 
					\begin{subequations}
						\begin{align}
							\nulllapse&=1+\int_{\tstart}^{t}\nulllapse\cdot\lgensmoothfunction\cdot(\pfour\vvariables,\pfour\cvvariables)\diff\tau,& u&\geq 0,\\
							\nulllapse&=a+\int_{\tstart}^{t}\nulllapse\cdot\lgensmoothfunction\cdot(\pfour\vvariables,\pfour\cvvariables) \diff\tau,& u&< 0,
						\end{align}
					\end{subequations}
					where $a$ is defined in (\ref{DE:Defofa}).
					Using the Gr\"onwall's inequality, the bootstrap assumption for $\pfour\vvariables$ and the initial condition (\ref{ES:lapseminusone}) in Section \ref{SS:BAGeo}, the $L^2_tL_x^\infty$ bound for $\pfour\cvvariables$ (\ref{BA:bt2}), we conclude the desired estimate.
				\end{proof}
				\begin{proof}[Proof of $\twonorms{\hchi}{t}{2}{\sangle}{p}{\coneu}\lesssim\lambda^{-1/2}$ in (\ref{mr1})]
					Recalling equation (\ref{EQ:angDhchi}), using the transport identity (\ref{EQ:rtwomupxi}) with $m=1$, we get
					\begin{align}
						\label{r2hchi}\sgabs{\rgeo^2\hchi}(t,u,\sangle)&\leq\lim\limits_{\tau\downarrow\tstart}\sgabs{\rgeo^2\hchi}(\tau,u,\sangle)\\
						\notag&+\int_{\tstart}^{t}\sgabs{\rgeo^2\lgensmoothfunction\cdot\boxg\gfour}+\sgabs{\rgeo^2(\angnabla,\angD_{\Lunit})\upxi}+\sgabs{\rgeo^2\lgensmoothfunction\cdot\left(\ACC,\rgeo^{-1}\right)\cdot\pfour\gfour}\diff\tau.
					\end{align}
					Dividing (\ref{r2hchi}) by $\rgeo^2(t,u)$, and taking the $L^2_tL^p_\sangle$ norm, we have:
					\begin{align}
						\label{hchit2sanglep}\twonorms{\hchi}{t}{2}{\sangle}{p}{\coneu}\leq&\twonorms{\rgeo^{-2}\left(\lim\limits_{\tau\downarrow\tstart}\rgeo^2\hchi\right)}{t}{2}{\sangle}{p}{\coneu}+\twonorms{\rgeo^{-2}\int_{\tstart}^{t}\sgabs{\rgeo^2\left(\ACC,\rgeo^{-1}\right)\cdot\pfour\gfour}\diff\tau}{t}{2}{\sangle}{p}{\coneu}\\
						\notag&+\twonorms{\rgeo^{-2}\int_{\tstart}^{t}\sgabs{\rgeo^2\boxg\gfour}\diff\tau}{t}{2}{\sangle}{p}{\coneu}+\twonorms{\rgeo^{-2}\int_{\tstart}^{t}\sgabs{\rgeo^2(\angnabla,\angD_{\Lunit})\upxi}\diff\tau}{t}{2}{\sangle}{p}{\coneu}.
					\end{align}
					We now consider the initial conditions. When $u\geq0$, by the initial condition on the cone-tip (\ref{ES:initialcontipone}), there holds $\lim\limits_{\tau\downarrow\tstart}\sgabs{\rgeo^2\hchi}=0$. When $u<0$,  we have $\lim\limits_{\tau\downarrow\tstart}\sgabs{\rgeo^2\hchi}=w^2\sgabs{\hchi}(0,u,\sangle)$. Then, by (\ref{ES:initialfoliationtwo}) and (\ref{fv1}), we derive:
					\begin{align}
						\label{whchisanglep}\onenorm{w^{1/2}\hchi(0,u,\sangle)}{\sangle}{p}{S_w}\lesssim\onenorm{w^{1/2}\hat{\spheresecondfund}(0,u,\sangle)}{\sangle}{p}{S_w}+\onenorm{w^{1/2}\pfour\gfour(0,u,\sangle)}{\sangle}{p}{S_w}\lesssim\lambda^{-1/2}.
					\end{align}
					Then, we deduce
					\begin{align}
						\twonorms{\rgeo^{-2}w^2\hchi(0,u,\sangle)}{t}{2}{\sangle}{p}{\coneu}\lesssim&\left(\int_{0}^{\Trescale}\left(\rgeo^{-2}w^{3/2}\lambda^{-1/2}\right)^2\diff\tau\right)^{1/2}\\
						\notag=&\left\{\lambda^{-1}w^3\left(\left.-\frac{1}{(\tau+w)^3}\right\vert^{\Trescale}_0\right)\right\}^{1/2}\\
						\notag\lesssim&\lambda^{-1/2}.
					\end{align}
					Now we estimate the other terms in (\ref{hchit2sanglep}). By (\ref{fv9}), we have
					\begin{align}
						\twonorms{\rgeo^{-2}\int_{\tstart}^{t}\sgabs{\rgeo^2\boxg\gfour}\diff\tau}{t}{2}{\sangle}{p}{\coneu}&\lesssim\lambda^{-1/2-12\varepsilon_0}.
					\end{align}
					Using the estimate for the Hardy-Littlewood maximal function (\ref{ES:hardylittlewood}), and estimate (\ref{fv6}), we have
					\begin{align}
						\label{rD*upxit2sanglep}\twonorms{\rgeo^{-2}\int_{\tstart}^{t}\sgabs{\rgeo^2(\angnabla,\angD_{\Lunit})\upxi}\diff\tau}{t}{2}{\sangle}{p}{\coneu}\lesssim&\twonorms{\frac{1}{t-\tstart}\int_{\tstart}^{t}\sgabs{\rgeo(\angnabla,\angD_{\Lunit})\upxi}\diff\tau}{t}{2}{\sangle}{p}{\coneu}\\
						\notag\lesssim&\twonorms{\rgeo(\angnabla,\angD_{\Lunit})\upxi}{t}{2}{\sangle}{p}{\coneu}\\
						\notag\lesssim&\lambda^{-1/2}.
					\end{align}
					Similarly, by applying the Hardy-Littlewood maximal inequality  (\ref{ES:hardylittlewood}), (\ref{BA:bt2}), the bootstrap assumptions (\ref{BA:bt1}), (\ref{bt4}), and (\ref{bt5}), the following estimate holds
					\begin{align}
						\label{rzt2sanglep}&\twonorms{\rgeo^{-2}\int_{\tstart}^{t}\sgabs{\rgeo^2\left(\ACC,\rgeo^{-1}\right)\cdot\pfour\gfour}\diff\tau}{t}{2}{\sangle}{p}{\coneu}\\
						\notag&\lesssim\twonorms{\rgeo\left(\ACC,\rgeo^{-1}\right)\cdot(\pfour\vvariables,\pfour\cvvariables)}{t}{2}{\sangle}{p}{\coneu}\\
						\notag&\lesssim\twonorms{\rgeo\left(\CC,\rgeo^{-1}\right)}{t}{\infty}{\sangle}{p}{\coneu}\twonorms{\pfour\vvariables,\pfour\cvvariables}{t}{2}{\sangle}{\infty}{\coneu}+\lambda^{1/2-4\varepsilon_0}\twonorms{\rgeo^{1/2}(\pfour\vvariables,\pfour\cvvariables)}{t}{\infty}{\sangle}{p}{\coneu}\twonorms{\pfour\vvariables,\pfour\cvvariables}{t}{2}{\sangle}{\infty}{\coneu}\\
						\notag&\lesssim\lambda^{-1/2-4\varepsilon_0}.
					\end{align}
					Combining (\ref{r2hchi})-(\ref{rzt2sanglep}), we conclude the desired estimate.
				\end{proof}
				\begin{proof}[Proof of $\twonorms{\rgeo^{1/2}\hchi}{t}{\infty}{\sangle}{p}{\coneu}\lesssim\lambda^{-1/2}$ in (\ref{mr1})]
					Dividing equation (\ref{r2hchi}) by $\rgeo^{3/2}(t,u)$, and taking the $L_{t}^{\infty}L_{\sangle}^{p}$ norm, we obtain
					\begin{align}\label{ES101}
						\twonorms{\rgeo^{1/2}\hchi}{t}{\infty}{\sangle}{p}{\coneu}\leq&\twonorms{\rgeo^{-3/2}\left(\lim\limits_{\tau\downarrow\tstart}\rgeo^2\hchi\right)}{t}{\infty}{\sangle}{p}{\coneu}+\twonorms{\rgeo^{-3/2}\int_{\tstart}^{t}\sgabs{\rgeo^2\boxg\gfour}\diff\tau}{t}{\infty}{\sangle}{p}{\coneu}\\
						\notag&+\twonorms{\rgeo^{-3/2}\int_{\tstart}^{t}\sgabs{\rgeo^2(\angnabla,\angD_{\Lunit})\upxi}\diff\tau}{t}{\infty}{\sangle}{p}{\coneu}\\
						\notag&+\twonorms{\rgeo^{-3/2}\int_{\tstart}^{t}\sgabs{\rgeo^2\left(\ACC,\rgeo^{-1}\right)\cdot\pfour\gfour}\diff\tau}{t}{\infty}{\sangle}{p}{\coneu}.
					\end{align}
					We control the initial condition as in the previous proof. In particular, using (\ref{whchisanglep}) for $u<0$, we deduce
					\begin{align}
						\twonorms{\rgeo^{-3/2}\left(\lim\limits_{\tau\downarrow\tstart}\rgeo^2\hchi\right)}{t}{\infty}{\sangle}{p}{\coneu}\lesssim\lambda^{-1/2}.
					\end{align}
					By (\ref{fv9}), we get
					\begin{align}
						\twonorms{\rgeo^{-3/2}\int_{\tstart}^{t}\sgabs{\rgeo^2\boxg\gfour}\diff\tau}{t}{\infty}{\sangle}{p}{\coneu}\lesssim\lambda^{-1/2-12\varepsilon_0}.
					\end{align}
					Employing the Minkowski's integral inequality and H\"older's inequality, and (\ref{fv6}), we derive
					\begin{align}
						\label{rD*tinftysanglep}\twonorms{\rgeo^{-3/2}\int_{\tstart}^{t}\sgabs{\rgeo^2(\angnabla,\angD_{\Lunit})\upxi}\diff\tau}{t}{\infty}{\sangle}{p}{\coneu}&\leq\sup\limits_t\onenorm{\rgeo^{-1/2}(t,u)\int_{\tstart}^{t}\sgabs{\rgeo(\angnabla,\angD_{\Lunit})\upxi}\diff\tau}{\sangle}{p}{\stu}\\
						\notag&\lesssim\sup\limits_t\frac{(t-\tstart)^{1/2}}{\rgeo^{1/2}(t,u)}\twonorms{\rgeo(\angnabla,\angD_{\Lunit})\upxi}{t}{2}{\sangle}{p}{\coneu}\\
						\notag&\lesssim\lambda^{-1/2}.
					\end{align}
					By using the same argument as above and also (\ref{rzt2sanglep}), we have
					\begin{align}
						\label{rztinftysanglep}\twonorms{\rgeo^{-3/2}\int_{\tstart}^{t}\sgabs{\rgeo^2\left(\ACC,\rgeo^{-1}\right)\cdot\pfour\gfour}\diff\tau}{t}{\infty}{\sangle}{p}{\coneu}&\lesssim\sup\limits_t\frac{(t-\tstart)^{1/2}}{\rgeo^{1/2}(t,u)}\twonorms{\rgeo\left(\ACC,\rgeo^{-1}\right)\cdot(\pfour\vvariables,\pfour\cvvariables)}{t}{2}{\sangle}{p}{\coneu}\\
						\notag&\lesssim\lambda^{-1/2-4\varepsilon_0}.
					\end{align}
					Combining (\ref{ES101})-(\ref{rztinftysanglep}), we conclude the desired estimates.
				\end{proof}
				\begin{proof}[Proof of $\twonorms{\rgeo\angD_{\Lunit}\hchi}{t}{2}{\sangle}{p}{\coneu}\lesssim\lambda^{-1/2}$ in (\ref{mr1})]
					Consider equation (\ref{EQ:angDhchi}):
					\begin{align}\label{ES111}
						\rgeo\angD_{\Lunit}\hchi=&\left(-\rgeo\chismall+2\right)\hchi+\rgeo(\angnabla,\angD_{\Lunit})\xi+\rgeo\lgensmoothfunction\cdot\boxg\gfour+\rgeo\lgensmoothfunction\cdot\left(\ACC,\rgeo^{-1}\right)\cdot\pfour\gfour.
					\end{align}
					Employing the bootstrap assumption (\ref{bt4}), and the previously proved result (\ref{mr1}) for $\twonorms{\rgeo^{1/2}\hchi}{t}{\infty}{\sangle}{p}{\coneu}$, we obtain:
					\begin{align}
						\twonorms{\rgeo\chismall\cdot\hchi}{t}{2}{\sangle}{p}{\coneu}\lesssim\lambda^{1/2-4\varepsilon_0}\twonorms{\chismall}{t}{2}{\sangle}{\infty}{\coneu}\twonorms{\rgeo^{1/2}\hchi}{t}{\infty}{\sangle}{p}{\coneu}\lesssim\lambda^{-1/2-2\varepsilon_0}.
					\end{align}
					By (\ref{fv6}), (\ref{rzt2sanglep}), we arrive at
					\begin{align}\label{ES112}
						&\twonorms{\rgeo(\angnabla,\angD_{\Lunit})\xi}{t}{2}{\sangle}{p}{\coneu}, \twonorms{\rgeo\lgensmoothfunction\cdot\left(\ACC,\rgeo^{-1}\right)\cdot\pfour\gfour}{t}{2}{\sangle}{p}{\coneu}\lesssim\lambda^{-1/2}.
					\end{align}
					By (\ref{fv10}), it holds that
					\begin{align}\label{ES113}
						\twonorms{\rgeo\lgensmoothfunction\cdot\boxg\gfour}{t}{2}{\sangle}{p}{\coneu}\lesssim&\lambda^{-1/2-8\varepsilon_0}.
					\end{align}
					Combining (\ref{ES111})-(\ref{ES113}) with (\ref{mr1}) for $\twonorms{\hchi}{t}{2}{\sangle}{p}{\coneu}$, we conclude the desired estimates.
				\end{proof}
				\begin{proof}[Proof of $\twonorms{\upzeta}{t}{2}{\sangle}{p}{\coneu}\lesssim\lambda^{-1/2}$ in (\ref{mr1})]
					Considering equation (\ref{EQ:angDLunitupzeta}), using (\ref{EQ:rtwomupxi}) with $m=\frac{1}{2}$, we get
					\begin{align}
						\label{rupzeta}\sgabs{\rgeo\upzeta}\leq&\lim\limits_{\tau\downarrow\tstart}\sgabs{\rgeo\upzeta}(\tau,u,\sangle)\\
						\notag&+\int_{\tstart}^{t}\sgabs{\rgeo\boxg\gfour}+\sgabs{\rgeo(\angnabla,\angD_{\Lunit})\upxi}\\		\notag&+\sgabs{\rgeo\lgensmoothfunction\cdot\left(\ACC,\rgeo^{-1}\right)\cdot\pfour\gfour}+\sgabs{\rgeo\lgensmoothfunction\cdot\upzeta\cdot\hchi} \diff\tau.
					\end{align}
					Dividing (\ref{rupzeta}) by $\rgeo(t,u)$, and taking the $L^2_tL^p_\sangle$ norm, we derive
					\begin{align}
						\twonorms{\upzeta}{t}{2}{\sangle}{p}{\coneu}\leq&\twonorms{\rgeo^{-1}\left(\lim\limits_{\tau\downarrow\tstart}\rgeo\upzeta\right)}{t}{2}{\sangle}{p}{\coneu}+\twonorms{\rgeo^{-1}\int_{\tstart}^{t}\sgabs{\rgeo\boxg\gfour}\diff\tau}{t}{2}{\sangle}{p}{\coneu}\\
						\notag&+\twonorms{\rgeo^{-1}\int_{\tstart}^{t}\sgabs{\rgeo(\angnabla,\angD_{\Lunit})\upxi}\diff\tau}{t}{2}{\sangle}{p}{\coneu}+\twonorms{\rgeo^{-1}\int_{\tstart}^{t}\sgabs{\rgeo\lgensmoothfunction\cdot\left(\ACC,\rgeo^{-1}\right)\cdot\pfour\gfour}\diff\tau}{t}{2}{\sangle}{p}{\coneu}\\
						\notag&+\twonorms{\rgeo^{-1/2}\int_{\tstart}^{t}\sgabs{\rgeo\upzeta\cdot\hchi}\diff\tau}{t}{\infty}{\sangle}{p}{\coneu}..
					\end{align}
					We now consider the initial condition. When $u\geq 0$, by the initial condition on the cone-tip (\ref{ES:initialcontipone}), we have $\rgeo^{-1}\left(\lim\limits_{\tau\downarrow\tstart}\rgeo\upzeta\right)=0$. When $u<0$, we use (\ref{E:CONNECTIONCOEFFICIENT}), (\ref{ES:initialfoliationtwo}), the bootstrap assumptions (\ref{BA:bt1}), and  estimates (\ref{BA:bt2}) to deduce:
					\begin{align}
						\twonorms{\rgeo^{-1}\left(\lim\limits_{\tau\downarrow\tstart}\rgeo\upzeta\right)}{t}{2}{\sangle}{p}{\coneu}&\lesssim\twonorms{\rgeo^{-1}w\angnabla\ln\lapse}{t}{2}{\sangle}{p}{\coneu}+\twonorms{\rgeo^{-1}w(\pfour\vvariables,\pfour\cvvariables)}{t}{2}{\sangle}{p}{\coneu}\\
						\notag&\lesssim\lambda^{-1/2}\int_{0}^{\Trescale}\frac{w}{(\tau+w)^2}\diff\tau+\lambda^{-1/2-4\varepsilon_0}\\
						\notag&\lesssim\lambda^{-1/2}.
					\end{align}
					It follows from (\ref{fv9}) that
					\begin{align}
						\twonorms{\rgeo^{-1}\int_{\tstart}^{t}\sgabs{\rgeo\boxg\gfour}\diff\tau}{t}{2}{\sangle}{p}{\coneu}\lesssim\lambda^{-1/2-12\varepsilon_0}.
					\end{align}
					Using the same method as in (\ref{rD*upxit2sanglep}) and (\ref{rzt2sanglep}), we obtain
					\begin{align}
						\twonorms{\rgeo^{-1}\int_{\tstart}^{t}\sgabs{\rgeo(\angnabla,\angD_{\Lunit})\upxi}\diff\tau}{t}{2}{\sangle}{p}{\coneu}&\lesssim\lambda^{-1/2},\\
						\twonorms{\rgeo^{-1}\int_{\tstart}^{t}\sgabs{\rgeo\lgensmoothfunction\cdot\left(\ACC,\rgeo^{-1}\right)\cdot\pfour\gfour}\diff\tau}{t}{2}{\sangle}{p}{\coneu}&\lesssim\lambda^{-1/2-4\varepsilon_0}.
					\end{align}
					Using the estimate for the Hardy-Littlewood maximal function (\ref{ES:hardylittlewood}), the bootstrap assumption (\ref{bt4}), and the previously proven estimate for $\rgeo^{1/2}\hchi$ in (\ref{mr1}), we derive
					\begin{align}
						\label{rupzetahchit2sanglep}\twonorms{\rgeo^{-1}\int_{\tstart}^{t}\sgabs{\rgeo\upzeta\cdot\hchi}\diff\tau}{t}{2}{\sangle}{p}{\coneu}&\lesssim\twonorms{\rgeo\upzeta\cdot\hchi}{t}{2}{\sangle}{p}{\coneu}\\
						\notag&\lesssim\lambda^{1/2-4\varepsilon_0}\twonorms{\rgeo^{1/2}\hchi}{t}{\infty}{\sangle}{p}{\coneu}\twonorms{\upzeta}{t}{2}{\sangle}{\infty}{\coneu}\\
						\notag&\lesssim\lambda^{-1/2-2\varepsilon_0}.
					\end{align}
					Combining (\ref{rupzeta})-(\ref{rupzetahchit2sanglep}), we conclude the desired estimate.
				\end{proof}
				\begin{proof}[Proof of $\twonorms{\rgeo^{1/2}\upzeta}{t}{\infty}{\sangle}{p}{\coneu}\lesssim\lambda^{-1/2}$ in (\ref{mr1})]
					Dividing the equation (\ref{rupzeta}) by $\rgeo^{1/2}(t,u)$, and taking the $L^\infty_tL^p_\sangle$ norm, we deduce
					\begin{align}\label{ES131}
						\twonorms{\rgeo^{1/2}\upzeta}{t}{\infty}{\sangle}{p}{\coneu}\leq&\twonorms{\rgeo^{-1/2}\left(\lim\limits_{\tau\downarrow\tstart}\rgeo\upzeta\right)}{t}{\infty}{\sangle}{p}{\coneu}\\
						\notag&+\twonorms{\rgeo^{-1/2}\int_{\tstart}^{t}\sgabs{\rgeo\boxg\gfour}\diff\tau}{t}{\infty}{\sangle}{p}{\coneu}+\twonorms{\rgeo^{-1/2}\int_{\tstart}^{t}\sgabs{\rgeo(\angnabla,\angD_{\Lunit})\upxi}\diff\tau}{t}{\infty}{\sangle}{p}{\coneu}\\
						\notag&+\twonorms{\rgeo^{-1/2}\int_{\tstart}^{t}\sgabs{\rgeo\lgensmoothfunction\cdot\left(\ACC,\rgeo^{-1}\right)\cdot\pfour\gfour}\diff\tau}{t}{\infty}{\sangle}{p}{\coneu}\\
						\notag&+\twonorms{\rgeo^{-1/2}\int_{\tstart}^{t}\sgabs{\rgeo\upzeta\cdot\hchi}\diff\tau}{t}{\infty}{\sangle}{p}{\coneu}.
					\end{align}
					We first consider the initial condition. When $u\geq 0$, by the initial condition on the cone-tip (\ref{ES:initialcontipone}), there holds 
					\begin{align}
						\rgeo^{-1/2}\left(\lim\limits_{\tau\downarrow\tstart}\rgeo\upzeta\right)=0.
					\end{align} When $u<0$, we use  (\ref{E:CONNECTIONCOEFFICIENT}), (\ref{ES:initialfoliationtwo}), and the estimate (\ref{fv1}) to deduce:
					\begin{align}
						\twonorms{\rgeo^{-1/2}\left(\lim\limits_{\tau\downarrow\tstart}\rgeo\upzeta\right)}{t}{\infty}{\sangle}{p}{\coneu}&\lesssim\twonorms{w^{1/2}\angnabla\ln\lapse}{t}{\infty}{\sangle}{p}{\coneu}+\twonorms{\rgeo^{1/2}(\pfour\vvariables,\pfour\cvvariables)}{t}{\infty}{\sangle}{p}{\coneu}\lesssim\lambda^{-1/2}.
					\end{align}
					By (\ref{fv9}), we have:
					\begin{align}
						\twonorms{\rgeo^{-1/2}\int_{\tstart}^{t}\sgabs{\rgeo\boxg\gfour}\diff\tau}{t}{\infty}{\sangle}{p}{\coneu}\lesssim\lambda^{-1/2-12\varepsilon_0}.
					\end{align}
					Using the same method as in (\ref{rD*tinftysanglep}) and (\ref{rztinftysanglep}), we can prove
					\begin{align}
						\twonorms{\rgeo^{-1/2}\int_{\tstart}^{t}\sgabs{\rgeo(\angnabla,\angD_{\Lunit})\upxi}\diff\tau}{t}{\infty}{\sangle}{p}{\coneu}&\lesssim\lambda^{-1/2},\\
						\label{upzetapfourvvariable}\twonorms{\rgeo^{-1/2}\int_{\tstart}^{t}\sgabs{\rgeo\lgensmoothfunction\cdot\left(\ACC,\rgeo^{-1}\right)\cdot\pfour\gfour}\diff\tau}{t}{\infty}{\sangle}{p}{\coneu}&\lesssim\lambda^{-1/2-4\varepsilon_0}.
					\end{align}
					With the help of (\ref{rupzetahchit2sanglep}), we obtain
					\begin{align}\label{ES132}
						\twonorms{\rgeo^{-1/2}\int_{\tstart}^{t}\sgabs{\rgeo\upzeta\cdot\hchi}\diff\tau}{t}{\infty}{\sangle}{p}{\coneu}&\lesssim\sup\limits_t\rgeo^{-1/2}\int_{\tstart}^{t}\onenorm{\rgeo\upzeta\cdot\hchi}{\sangle}{p}{\sTauu} \diff\tau\\
						\notag&\lesssim\sup\limits_t\frac{(t-\tstart)^{1/2}}{(t-u)^{1/2}}\twonorms{\rgeo\upzeta\cdot\hchi}{t}{2}{\sangle}{p}{\coneu}\\
						\notag&\lesssim\lambda^{-1/2-2\varepsilon_0}.
					\end{align}
					Combining (\ref{ES131})-(\ref{ES132}), we conclude the desired estimate.
				\end{proof}
				\begin{proof}[Proof of $\twonorms{\rgeo\angD_{\Lunit}\upzeta}{t}{2}{\sangle}{p}{\coneu}\lesssim\lambda^{-1/2}$ in (\ref{mr1})]
					We consider equation (\ref{EQ:angDLunitupzeta}):
					\begin{align}\label{ES141}
						\rgeo\angD_{\Lunit}\upzeta=&\upzeta+\rgeo(\angnabla,\angD_{\Lunit})\xi+\rgeo\lgensmoothfunction\cdot\boxg\gfour+\rgeo\lgensmoothfunction\cdot\left(\ACC,\rgeo^{-1}\right)\cdot\pfour\gfour+\rgeo\lgensmoothfunction\cdot\upzeta\cdot\ACC.
					\end{align}
					Using the bootstrap assumptions (\ref{BA:bt1}), (\ref{bt4}), and the proven result (\ref{mr1}) for $\twonorms{\rgeo^{1/2}\upzeta}{t}{\infty}{\sangle}{p}{\coneu}$, there holds
					\begin{align}
						\label{rzupzeta}\twonorms{\rgeo\lgensmoothfunction\cdot\ACC\cdot\upzeta}{t}{2}{\sangle}{p}{\coneu}&\lesssim\lambda^{1/2-4\varepsilon_0}\twonorms{\ACC}{t}{2}{\sangle}{\infty}{\coneu}\twonorms{\rgeo^{1/2}\upzeta}{t}{\infty}{\sangle}{p}{\coneu}\lesssim\lambda^{-1/2-2\varepsilon_0}.
					\end{align}
					By (\ref{fv6}) (\ref{fv10}), and (\ref{rzt2sanglep}), we derive
					\begin{align}\label{ES142}
						\twonorms{\rgeo(\angnabla,\angD_{\Lunit})\xi}{t}{2}{\sangle}{p}{\coneu}, \twonorms{\rgeo\lgensmoothfunction\cdot\boxg\gfour}{t}{2}{\sangle}{p}{\coneu},\twonorms{\rgeo\lgensmoothfunction\cdot\left(\ACC,\rgeo^{-1}\right)\cdot\pfour\gfour}{t}{2}{\sangle}{p}{\coneu}\lesssim\lambda^{-1/2}.
					\end{align}
					Combining (\ref{ES141})-(\ref{ES142}) with (\ref{mr1}) for $\twonorms{\upzeta}{t}{2}{\sangle}{p}{\coneu}$, we conclude the desired estimates.
				\end{proof}
				\begin{proof}[Proof of $\twonorms{\rgeo^{1/2}(\hchi,\upzeta)}{\sangle}{2p}{t}{\infty}{\coneu}\lesssim\lambda^{-1/2}$ when $\coneu\subset\intregion$ in (\ref{mr13})]
					Utilizing the Sobolev inequality (\ref{ES:sobolevtwo}), the previously proven estimate (\ref{mr1}), and bootstrap assumptions (\ref{bt7}), we conclude
					\begin{align}\label{ES151}
						\twonorms{\rgeo^{\frac{1}{2}}(\hchi,\upzeta)}{\sangle}{2p}{t}{\infty}{\coneu}^2&\lesssim\left(\twonorms{\rgeo\angD_{\Lunit}(\hchi,\upzeta)}{\sangle}{p}{t}{2}{\coneu}+\twonorms{(\hchi,\upzeta)}{\sangle}{p}{t}{2}{\coneu}\right)\twonorms{(\hchi,\upzeta)}{\sangle}{\infty}{t}{2}{\coneu}\lesssim\lambda^{-1}.
					\end{align}
					In deriving (\ref{ES151}), we also used the Minkowski's inequality for integrals to switch the order of the $L_t$ and $L_\sangle$ norms. 
				\end{proof}
				\begin{proof}[Proof of $\rgeo\Restrace\reschi\approx 1$ in (\ref{mr4})]
					It is sufficient to show $\sgabs{\rgeo\chismall}\lesssim\lambda^{-4\varepsilon_0}$. We plug (\ref{EQ:Lz}) into (\ref{EQ:transportlemma3}), where $\mathfrak{G}=(\chismall,\hchi,\pfour\gfour)$. Then, dividing the estimate (\ref{ES:r2mupxi}) by $\rgeo(t,u)$, we obtain
					\begin{align}
						\label{rgeoz}\sgabs{\rgeo\chismall}\lesssim&\rgeo^{-1}\lim\limits_{\tau\downarrow\tstart}\sgabs{\rgeo^2\chismall}\\
						\notag&+\rgeo^{-1}\int_{\tstart}^{t}\sgabs{\rgeo^2\lgensmoothfunction\cdot\boxg\gfour+\rgeo^2\lgensmoothfunction\cdot(\hchi,\pfour\gfour)\cdot(\hchi,\pfour\gfour)+\rgeo\lgensmoothfunction\cdot\pfour\gfour}\diff\tau.
					\end{align}
					We now consider the initial conditions. When $u\geq 0$, by the initial condition on the cone-tip (\ref{ES:initialcontipone}), it holds that $\rgeo^{-1}\lim\limits_{\tau\downarrow\tstart}\sgabs{\rgeo^2\chismall}\leq\lim\limits_{\tau\downarrow\tstart}\sgabs{\rgeo\chismall}=0$. When $u<0$, by (\ref{EQ:zinitial}) and (\ref{ES:lapseminusone}), we have:
					\begin{align}
						\rgeo^{-1}\lim\limits_{\tau\downarrow\tstart}\sgabs{\rgeo^2\chismall}=\frac{1}{\rgeo}\frac{2w(1-\lapse)}{\lapse}\leq\frac{2(1-\lapse)}{\lapse}\lesssim\lambda^{-4\varepsilon_0}.
					\end{align} 
					By (\ref{fv9}), we deduce
					\begin{align}
						\rgeo^{-1}\int_{\tstart}^{t}\abs{\rgeo^2\lgensmoothfunction\cdot\boxg\gfour}\diff\tau\lesssim\lambda^{-16\varepsilon_0}.
					\end{align}
					By the bootstrap assumptions (\ref{BA:bt1}), (\ref{bt4}) and estimate (\ref{BA:bt2}), we obtain
					\begin{align}
						\rgeo^{-1}\int_{\tstart}^{t}\sgabs{\rgeo^{2}\lgensmoothfunction\cdot(\hchi,\pfour\gfour)\cdot(\hchi,\pfour\gfour)} \diff\tau\lesssim\lambda^{1-8\varepsilon_0}\twonorms{\hchi,\pfour\vvariables,\pfour\cvvariables}{t}{2}{\sangle}{\infty}{\coneu}\twonorms{\hchi,\pfour\vvariables,\pfour\cvvariables}{t}{2}{\sangle}{\infty}{\coneu}\lesssim\lambda^{-4\varepsilon_0}.
					\end{align}
					\begin{align}\label{ES161}
						\rgeo^{-1}\int_{\tstart}^{t}\sgabs{\rgeo\lgensmoothfunction\cdot\pfour\gfour} \diff\tau\lesssim\twonorms{(\pfour\vvariables,\pfour\cvvariables)}{t}{2}{\sangle}{\infty}{\coneu}\lambda^{1/2-4\varepsilon_0}\lesssim\lambda^{-8\varepsilon_0}.
					\end{align}
					Combining (\ref{rgeoz})-(\ref{ES161}), we conclude the desired estimate. The estimate $\twonorms{\rgeo\chismall}{t}{\infty}{\sangle}{p}{\coneu}\lesssim\lambda^{-4\varepsilon_0}$ follows from the proof of $\rgeo\Restrace\reschi\approx 1$.
				\end{proof}
				\begin{proof}[Proof of $\norm{\rgeo^{1/2}\chismall}_{L^\infty(\region)}\lesssim\lambda^{-1/2}$ in (\ref{mr5})]
					Dividing (\ref{rgeoz}) by $\rgeo^{1/2}$, we get
					\begin{align}\label{ES171}
						\sgabs{\rgeo^{1/2}\chismall}\lesssim&\rgeo^{-3/2}\lim\limits_{\tau\downarrow\tstart}\sgabs{\rgeo^2\chismall}\\
						\notag&+\rgeo^{-3/2}\int_{\tstart}^{t}\sgabs{\rgeo^2\lgensmoothfunction\cdot\boxg\gfour+\rgeo^2\lgensmoothfunction\cdot\left(\hchi,\pfour\gfour\right)\cdot(\hchi,\pfour\gfour)+\rgeo\lgensmoothfunction\cdot\pfour\gfour}\diff\tau.
					\end{align}
					We consider the initial conditions. When $u\geq 0$, by the initial condition on the cone-tip (\ref{ES:initialcontipone}),  we have
					\begin{align}
						\rgeo^{-3/2}\lim\limits_{\tau\downarrow\tstart}\sgabs{\rgeo^2\chismall}\leq\rgeo^{-1/2}\lim\limits_{\tau\downarrow\tstart}\sgabs{\rgeo\chismall}\lesssim\lim\limits_{\tau\downarrow\tstart}\rgeo^{1/2}=0.
					\end{align} 
					When $u<0$, by (\ref{EQ:zinitial}) and (\ref{ES:lapseminusone}), there holds
					\begin{align}
						\rgeo^{-3/2}\lim\limits_{\tau\downarrow\tstart}\sgabs{\rgeo^2\chismall}=\frac{1}{\rgeo^{3/2}}\frac{2w(1-\lapse)}{\lapse}\leq w^{-1/2}(1-\lapse)\lesssim\lambda^{-1/2}.
					\end{align} 
					By (\ref{fv9}), we obtain
					\begin{align}
						\norm{\rgeo^{-\frac{3}{2}}\int_{\tstart}^t\sgabs{\rgeo^2\boxg\gfour}\diff\tau}_{L^\infty(\region)}&\lesssim\lambda^{-1/2-12\varepsilon_0}.
					\end{align}
					By the bootstrap assumptions (\ref{BA:bt1}), (\ref{bt4}) and estimate (\ref{BA:bt2}), we derive
					\begin{align}
						\rgeo^{-3/2}\int_{\tstart}^{t}\sgabs{\rgeo^{2}\lgensmoothfunction\cdot(\hchi,\pfour\gfour)\cdot(\hchi,\pfour\gfour)} \diff\tau\lesssim&\lambda^{1/2-4\varepsilon_0}\twonorms{\hchi,\pfour\vvariables,\pfour\cvvariables}{t}{2}{\sangle}{\infty}{\coneu}\twonorms{\hchi,\pfour\vvariables,\pfour\cvvariables}{t}{2}{\sangle}{\infty}{\coneu}\\
						\notag\lesssim&\lambda^{-1/2}.
					\end{align}
					\begin{align}\label{ES172}
						\rgeo^{-3/2}\int_{\tstart}^{t}\sgabs{\rgeo\lgensmoothfunction\cdot\pfour\gfour} \diff\tau\lesssim&\rgeo^{-1/2}\int_{\tstart}^{t}\onenorm{\lgensmoothfunction\cdot(\pfour\vvariables,\pfour\cvvariables)}{x}{\infty}{\St}\diff\tau\\
						\notag\lesssim&\frac{(t-\tstart)^{1/2}}{(t-u)^{1/2}}\twonorms{\pfour\vvariables,\pfour\cvvariables}{t}{2}{x}{\infty}{\region}\\
						\notag
						\lesssim&\lambda^{-1/2-4\varepsilon_0}.
					\end{align}
					Combining (\ref{ES171})-(\ref{ES172}), we conclude the desired estimate.
				\end{proof}
				{Following from the previously proven estimate $\norm{\rgeo^{1/2}\chismall}_{L^\infty(\region)}\lesssim\lambda^{-1/2}$, we can then show the estimates, $\twonorms{\rgeo^{1/2}\chismall}{t}{\infty}{\sangle}{p}{\coneu}\lesssim\lambda^{-1/2}$ in (\ref{mr1}),  $\twonorms{\rgeo^{1/2}\chismall}{\sangle}{2p}{u}{\infty}{\St}\lesssim\lambda^{-1/2}$ in (\ref{mr4}), and $\twonorms{\rgeo^{1/2}\chismall}{\sangle}{2p}{t}{\infty}{\coneu}\lesssim\lambda^{-1/2}$ in (\ref{mr13}).}
				\begin{proof}[Proof of $\norm{\gtr\spheresecondfund-\frac{2}{\rgeo}}_{L^3(\St)}\lesssim \lambda^{-4\varepsilon_0}$ in (\ref{mr1152})]
					By the definition of $\spheresecondfund$ in (\ref{E:CONNECTIONCOEFFICIENT}), it holds that
					\begin{align}
						\gtr\spheresecondfund-\frac{2}{\rgeo}=\chismall+\pfour\gfour.
					\end{align}
					Then, the desired estimate (\ref{mr1152}) follows by using H\"older inequality, $\twonorms{\rgeo^{1/2}\chismall}{\sangle}{2p}{u}{\infty}{\St}\lesssim\lambda^{-1/2}$ in (\ref{mr4}), and the second estimate of (\ref{fv1}).
				\end{proof}
				\begin{proof}[Proof of $\holdertwonorms{\chismall}{t}{2}{\sangle}{0}{\delta_0}{\coneu}\lesssim\lambda^{-1/2}$ in (\ref{mr7})]
					Plugging equation (\ref{EQ:Lz}) into (\ref{EQ:transportlemma2}), and dividing (\ref{EQ:rtwomupxi}) by $\rgeo^2(t,u)$ (with $m=1$), we derive:
					\begin{align}\label{ES182}
						\sgabs{\chismall}\leq&\rgeo^{-2}\lim\limits_{\tau\downarrow\tstart}\sgabs{\rgeo^2\chismall}+\rgeo^{-2}\int_{\tstart}^{t}\sgabs{\rgeo^2\lgensmoothfunction\cdot\boxg\gfour\diff\tau}\\
						\notag&+\rgeo^{-2}\int_{\tstart}^{t}\left|\rgeo^2\lgensmoothfunction\cdot\ACC\cdot\ACC+\rgeo\lgensmoothfunction\cdot\pfour\gfour \right|_{\gsphere}\diff\tau
						.\end{align}
					We now consider the initial condition. When $u\geq 0$, by (\ref{ES:initialcontipone}), there holds $\lim\limits_{\tau\downarrow\tstart}\left(\rgeo^2\chismall\right)=0$.
					When $u<0$, by (\ref{EQ:zinitial}), we have:
					\begin{align}
						\holdertwonorms{\rgeo^{-2}\lim\limits_{\tau\downarrow\tstart}\left(\rgeo^2\chismall\right)}{t}{2}{\sangle}{0}{\delta_0}{\coneu}^2\lesssim\int_{0}^{\Trescale}\left(\frac{w^{3/2}}{\rgeo^2}\lambda^{-1/2}\right)^2\diff\tau\leq\lambda^{-1}
						.\end{align}
					Using the Hardy-Littlewood maximal inequality (\ref{ES:hardylittlewood}), the bootstrap assumptions (\ref{BA:bt1}), and estimate (\ref{BA:bt2.1}), we get
					\begin{align}
						\holdertwonorms{\rgeo^{-2}\int_{\tstart}^{t}\left|\rgeo\lgensmoothfunction\cdot\pfour\gfour \right|_{\gsphere}\diff\tau}{t}{2}{\sangle}{0}{\delta_0}{\coneu}\lesssim\norm{\mathcal{M}\left(\holdernorm{\pfour\vvariables,\pfour\cvvariables}{\sangle}{0}{\delta_0}{\stu}\right)}_{L^2_t}\lesssim \holdertwonorms{\pfour\vvariables,\pfour\cvvariables}{t}{2}{\sangle}{0}{\delta_0}{\coneu}\lesssim\lambda^{-1/2-4\varepsilon_0}
						.\end{align}
					By the bootstrap assumptions (\ref{BA:bt1}), (\ref{bt4}), and estimates (\ref{fv9}), (\ref{BA:bt2}), we have:
					\begin{align}
						\holdertwonorms{\rgeo^{-2}\int_{\tstart}^{t}\sgabs{\rgeo^2\lgensmoothfunction\cdot\boxg\gfour}\diff\tau}{t}{2}{\sangle}{0}{\delta_0}{\coneu}&\lesssim\lambda^{-1/2-12\varepsilon_0},\\
						\label{ES183}\holdertwonorms{\rgeo^{-2}\int_{\tstart}^{t}\sgabs{\rgeo^2\lgensmoothfunction\cdot\ACC\cdot\ACC}\diff\tau}{t}{2}{\sangle}{0}{\delta_0}{\coneu}&\lesssim\lambda^{-1/2}
						.\end{align}
					Combining (\ref{ES182})-(\ref{ES183}), we conclude the desired result.
				\end{proof}
				The estimate $\twonorms{\chismall}{t}{2}{\sangle}{p}{\coneu}\lesssim\lambda^{-1/2}$ in (\ref{mr1}) follows directly from $\twonorms{\chismall}{t}{2}{\sangle}{\infty}{\coneu}\lesssim\lambda^{-1/2}$.
				\begin{proof}[Proof of $\twonorms{\rgeo\angD_{\Lunit}\chismall}{t}{2}{\sangle}{p}{\coneu}\lesssim\lambda^{-1/2}$ in (\ref{mr1})]
					Using equation (\ref{EQ:Lz}), we obtain
					\begin{align}\label{ES191}
						\rgeo\angD_{\Lunit}\chismall=&-2\chismall+\rgeo\lgensmoothfunction\cdot\boxg\gfour+\rgeo\lgensmoothfunction\cdot\ACC\cdot\ACC+\lgensmoothfunction\cdot\pfour\gfour
						.\end{align}
					By (\ref{fv10}), we deduce
					\begin{align}
						\twonorms{\rgeo\boxg\gfour}{t}{2}{\sangle}{p}{\coneu}&\lesssim^{-1/2-8\varepsilon_0}
						.\end{align}
					By estimate (\ref{fv3}), we can show
					\begin{align}
						\twonorms{\lgensmoothfunction\cdot\pfour\gfour}{t}{2}{\sangle}{p}{\coneu}\lesssim\lambda^{-1/2-4\varepsilon_0}
						.\end{align}
					By the bootstrap assumptions (\ref{BA:bt1}), (\ref{bt4}), and the previously proven results (\ref{BA:bt2}), (\ref{fv1}), (\ref{mr1}) for $\twonorms{\rgeo^{1/2}\CC}{t}{\infty}{\sangle}{p}{\coneu}$, we derive
					\begin{align}
						\label{zzhchihchi}\twonorms{\rgeo\lgensmoothfunction\cdot\ACC\cdot\ACC}{t}{2}{\sangle}{p}{\coneu}&\lesssim\lambda^{1/2-4\varepsilon_0}\twonorms{\ACC}{t}{2}{\sangle}{\infty}{\coneu}\twonorms{\rgeo^{1/2}\ACC}{t}{\infty}{\sangle}{p}{\coneu}\\
						\notag&\lesssim\lambda^{-1/2-2\varepsilon_0}
						.\end{align}
					Combining (\ref{ES191})-(\ref{zzhchihchi}), we conclude the desired estimate.
				\end{proof}
				\begin{remark}
					From (\ref{fv3}), (\ref{mr1}), and (\ref{zzhchihchi}), it follows that 
					\begin{align}
						\label{rmk1}\twonorms{\rgeo\lgensmoothfunction\cdot\left(\ACC,\rgeo^{-1}\right)\cdot\ACC}{t}{2}{\sangle}{p}{\coneu}\lesssim\lambda^{-1/2}
						.\end{align}
				\end{remark}
				\begin{proof}[Proof of $\twonorms{\rgeo\left(\angnabla\chismall,\angnabla\hchi\right)}{t}{2}{\sangle}{p}{\coneu}\lesssim\lambda^{-1/2}$ in (\ref{mr7})]
					First, we bound $\twonorms{\rgeo\angnabla\chismall}{t}{2}{\sangle}{p}{\coneu}$. {Plugging (\ref{EQ:Lunitangnablaz}) into (\ref{EQ:transportlemma3}), and letting $\mathfrak{G}=\lgensmoothfunction\cdot\ACC$ in (\ref{EQ:transportlemma3}), then by Lemma \ref{LE:transport}, we arrive at \eqref{ES:r2mupxi}. Dividing (\ref{ES:r2mupxi}) by $\rgeo^2(t,u)$, we derive}
					\begin{align}\label{ES201}
						\sgabs{\rgeo\angnabla\chismall}\lesssim&\rgeo^{-2}\lim\limits_{\tau\downarrow\tstart}\sgabs{\rgeo^3\angnabla\chismall}+\rgeo^{-2}\int_{\tstart}^{t}\sgabs{\rgeo^3\lgensmoothfunction\cdot\angnabla\boxg\gfour}\diff\tau\\
						\notag&+\rgeo^{-2}\int_{\tstart}^{t}\sgabs{\rgeo^3\lgensmoothfunction\cdot\boxg\gfour\cdot\left(\ACC,\rgeo^{-1}\right)}\diff\tau\\
						\notag&+\rgeo^{-2}\int_{\tstart}^{t}\sgabs{\rgeo^3\lgensmoothfunction\cdot\angnabla\pfour\gfour\cdot\left(\ACC,\rgeo^{-1}\right)+\rgeo^3\lgensmoothfunction\cdot\angnabla\hchi\cdot\ACC} \diff\tau\\
						\notag&+\rgeo^{-2}\int_{\tstart}^{t}\sgabs{\rgeo^3\lgensmoothfunction\cdot\left(\ACC,\rgeo^{-1}\right)\cdot\left(\ACC,\rgeo^{-1}\right)\cdot\pfour\vvariables }\diff\tau
						.\end{align}
					Consider the initial conditions. When $u\geq0$, we use the estimate for the Hardy-Littlewood maximal function (\ref{ES:hardylittlewood}) , and (\ref{ES:initialcontipone}), to deduce
					$\twonorms{\rgeo^{-2}\lim\limits_{\tau\downarrow\tstart}(\rgeo^3\angnabla\chismall)}{t}{2}{\sangle}{p}{\coneu}=0$. When $u<0$, we use the initial condition (\ref{ES:in7}) to deduce:
					\begin{align}
						\twonorms{\rgeo^{-2}\lim\limits_{\tau\downarrow\tstart}(\rgeo^3\angnabla\chismall)}{t}{2}{\sangle}{p}{\coneu}\lesssim\lambda^{-1/2}\twonorms{\rgeo^{-2}w^{3/2}}{t}{2}{\sangle}{p}{\coneu}\lesssim\lambda^{-1/2}
						.\end{align}
					It follows from (\ref{fv11}) that
					\begin{align}
						\twonorms{\rgeo^{-2}\int_{\tstart}^t\sgabs{\rgeo^3\lgensmoothfunction\cdot\angnabla\boxg\gfour}\diff\tau}{t}{2}{\sangle}{p}{\coneu}&\lesssim\lambda^{-1/2-8\varepsilon_0}
						.\end{align}
					Via using (\ref{fv9}), (\ref{fv10}), (\ref{BA:bt1}), (\ref{BA:bt2}), (\ref{bt4}), we derive
					\begin{align}
						\twonorms{\rgeo^{-2}\int_{\tstart}^t\sgabs{\rgeo^3\lgensmoothfunction\cdot\boxg\gfour\cdot\left(\ACC,\rgeo^{-1}\right)}\diff\tau}{t}{2}{\sangle}{p}{\coneu}&\lesssim\lambda^{-1/2-10\varepsilon_0}
						.\end{align}
					Appealing to the estimate for the Hardy-Littlewood maximal function (\ref{ES:hardylittlewood}), the bootstrap assumptions (\ref{BA:bt1}), (\ref{bt4}), and estimate (\ref{BA:bt2}), (\ref{fv6}), we deduce
					\begin{align}
						&\twonorms{\rgeo^{-2}\int_{\tstart}^{t}\rgeo^3\lgensmoothfunction\cdot\angnabla\pfour\gfour\cdot\left(\ACC,\rgeo^{-1}\right)\diff\tau}{t}{2}{\sangle}{p}{\coneu}\\
						\notag&\lesssim\norm{\int_{\tstart}^{t}\onenorm{\rgeo(\angnabla\pfour\vvariables,\angnabla\pfour\cvvariables)}{\sangle}{p}{\sTauu}\onenorm{\ACC}{\sangle}{\infty}{\sTauu}\diff\tau}_{L^2_t}+\norm{\rgeo^{-1}\int_{\tstart}^{t}\onenorm{\rgeo(\angnabla\pfour\vvariables,\angnabla\pfour\cvvariables)}{\sangle}{p}{\sTauu}\diff\tau}_{L^2_t}\\
						\notag&\lesssim\lambda^{1/2-4\varepsilon_0}\twonorms{\rgeo(\angnabla\pfour\vvariables,\angnabla\pfour\cvvariables)}{t}{2}{\sangle}{p}{\coneu}\twonorms{\ACC}{t}{2}{\sangle}{\infty}{\coneu}+\twonorms{\rgeo(\angnabla\pfour\vvariables,\angnabla\pfour\cvvariables)}{t}{2}{\sangle}{p}{\coneu}\\
						\notag&\lesssim\lambda^{-1/2}
						.\end{align}
					By the estimate for the Hardy-Littlewood maximal function (\ref{ES:hardylittlewood}), the bootstrap assumptions (\ref{BA:bt1}), (\ref{bt4}), (\ref{bt5}), and estimate (\ref{BA:bt2}), (\ref{fv5}), there holds
					\begin{align}
						&\twonorms{\rgeo^{-2}\int_{\tstart}^{t}\rgeo^3\lgensmoothfunction\cdot\left(\ACC,\rgeo^{-1}\right)\cdot\left(\ACC,\rgeo^{-1}\right)\cdot\pfour\gfour\diff\tau}{t}{2}{\sangle}{p}{\coneu}\\
						\notag\lesssim&\left\|\int_{\tstart}^{t}\onenorm{\rgeo\left(\ACC,\rgeo^{-1}\right)}{\sangle}{p}{\sTauu}\cdot\onenorm{\ACC}{\sangle}{\infty}{\sTauu}\cdot\onenorm{\pfour\vvariables,\pfour\cvvariables}{\sangle}{\infty}{\sTauu}\diff\tau\right\|_{L^2_t}\\
						\notag&+\norm{\rgeo^{-1}\int_{\tstart}^{t}\onenorm{\pfour\vvariables,\pfour\cvvariables}{\sangle}{\infty}{\sTauu}\diff\tau}_{L^2_t}\\
						\notag\lesssim&\lambda^{1/2-4\varepsilon_0}\twonorms{\ACC}{t}{2}{\sangle}{\infty}{\coneu}\twonorms{\pfour\vvariables,\pfour\cvvariables}{t}{2}{\sangle}{\infty}{\coneu}+\twonorms{\pfour\vvariables,\pfour\cvvariables}{t}{2}{\sangle}{\infty}{\coneu}\\
						\notag\lesssim&\lambda^{-1/2-4\varepsilon_0}
						.\end{align}
					Employing bootstrap assumption (\ref{bt4}), we get
					\begin{align}\label{ES202}
						\twonorms{\rgeo^{-2}\int_{\tstart}^{t}\rgeo^3\lgensmoothfunction\cdot\angnabla\hchi\cdot\ACC \diff\tau}{t}{2}{\sangle}{p}{\coneu}&\lesssim\norm{\twonorms{\rgeo\angnabla\hchi}{t}{2}{\sangle}{p}{\coneu}\twonorms{\ACC}{t}{2}{\sangle}{\infty}{\coneu}}_{L_t^2}\\
						\notag&\lesssim\lambda^{-2\varepsilon_0}\twonorms{\rgeo\angnabla\hchi}{t}{2}{\sangle}{p}{\coneu}
						.\end{align}
					Combining (\ref{ES201})-(\ref{ES202}), we arrive at
					\begin{align}
						\label{rangnablazlessrangnablahchi}\twonorms{\rgeo\angnabla\chismall}{t}{2}{\sangle}{p}{\coneu}\lesssim\lambda^{-1/2}+\lambda^{-2\varepsilon_0}\twonorms{\rgeo\angnabla\hchi}{t}{2}{\sangle}{p}{\coneu}
						.\end{align}
					Now we consider $\twonorms{\rgeo\angnabla\hchi}{t}{2}{\sangle}{p}{\coneu}$. Plugging equation (\ref{EQ:angdvhchi}) into the Hodge estimate (\ref{ES:hodgetwo}), we derive
					\begin{align}
						\twonorms{\rgeo\angnabla\hchi}{t}{2}{\sangle}{p}{\coneu}\lesssim&\twonorms{\rgeo\angnabla\chismall}{t}{2}{\sangle}{p}{\coneu}\\
						\notag&+\twonorms{\rgeo(\angnabla\pfour\vvariables,\angnabla\pfour\cvvariables)}{t}{2}{\sangle}{p}{\coneu}+\twonorms{\rgeo\lgensmoothfunction\cdot\left(\ACC,\rgeo^{-1}\right)\cdot(\pfour\vvariables,\pfour\cvvariables)}{t}{2}{\sangle}{p}{\coneu}
						.\end{align}
					Then, by utilizing (\ref{fv6}) and (\ref{rmk1}), we obtain
					\begin{align}
						\label{rangnablahchilessrangnablaz}\twonorms{\rgeo\angnabla\hchi}{t}{2}{\sangle}{p}{\coneu}\lesssim\lambda^{-1/2}+\twonorms{\rgeo\angnabla\chismall}{t}{2}{\sangle}{p}{\coneu}
						.\end{align}
					Combining (\ref{rangnablazlessrangnablahchi}) and (\ref{rangnablahchilessrangnablaz}), we prove the desired estimates $\twonorms{\rgeo\left(\angnabla\chismall,\angnabla\hchi\right)}{t}{2}{\sangle}{p}{\coneu}\lesssim\lambda^{-1/2}$.
				\end{proof}
				\begin{proof}[Proof of $\threenorms{\rgeo^{3/2}\angnabla\chismall}{t}{\infty}{u}{\infty}{\sangle}{p}{\region}\lesssim\lambda^{-1/2}$ in (\ref{mr5})]
					{Applying Lemma \ref{LE:transport} first to equation (\ref{EQ:Lunitangnablaz}),  then dividing (\ref{ES:r2mupxi}) by $\rgeo^{3/2}(t,u)$, we derive:}
					\begin{align}
						&\sgabs{\rgeo^{3/2}\angnabla\chismall}\\
						\notag\lesssim&\rgeo^{-3/2}\lim\limits_{\tau\downarrow\tstart}\sgabs{\rgeo^3\angnabla\chismall}+\rgeo^{-3/2}\int_{\tstart}^{t}\sgabs{\rgeo^3\lgensmoothfunction\cdot\angnabla\boxg\gfour}\diff\tau\\
						\notag&+\rgeo^{-3/2}\int_{\tstart}^{t}\sgabs{\rgeo^3\lgensmoothfunction\cdot\boxg\gfour\cdot\left(\ACC,\rgeo^{-1}\right)}\diff\tau\\
						\notag&+\rgeo^{-3/2}\int_{\tstart}^{t}\sgabs{\rgeo^3\lgensmoothfunction\cdot\angnabla\pfour\gfour\cdot\left(\ACC,\rgeo^{-1}\right)+\rgeo^3\lgensmoothfunction\cdot\angnabla\hchi\cdot\ACC} \diff\tau\\
						\notag&+\rgeo^{-3/2}\int_{\tstart}^{t}\sgabs{\rgeo^3\lgensmoothfunction\cdot\left(\ACC,\rgeo^{-1}\right)\cdot\left(\ACC,\rgeo^{-1}\right)\cdot\pfour\gfour} \diff\tau
						.\end{align}
					We now consider the initial conditions. When $u\geq0$, we use (\ref{ES:initialcontipone}) to deduce
					$\rgeo^{-3/2}\lim\limits_{\tau\downarrow\tstart}\left(\rgeo^3\angnabla\chismall\right)=0$. When $u<0$, we use the initial condition (\ref{ES:in7}) to deduce $\twonorms{\rgeo^{-3/2}\lim\limits_{\tau\downarrow\tstart}\sgabs{\rgeo^3\angnabla\chismall}}{u}{\infty}{\sangle}{p}{\Szero}\lesssim\lambda^{-1/2}$.\\
					By (\ref{fv11}), we have:
					\begin{align}
						\threenorms{\rgeo^{-3/2}\int_{\tstart}^{t}\sgabs{\rgeo^3\lgensmoothfunction\cdot\angnabla\boxg\gfour}\diff\tau}{t}{\infty}{u}{\infty}{\sangle}{p}{\region}\lesssim\lambda^{-1/2-8\varepsilon_0}
						.\end{align}
					Combining (\ref{fv9}), (\ref{fv10}), (\ref{BA:bt1}), (\ref{BA:bt2}), (\ref{bt4}), we obtain
					\begin{align}
						\threenorms{\rgeo^{-3/2}\int_{\tstart}^{t}\sgabs{\rgeo^3\lgensmoothfunction\cdot\boxg\gfour\cdot\left(\ACC,\rgeo^{-1}\right)}\diff\tau}{t}{\infty}{u}{\infty}{\sangle}{p}{\region}\lesssim\lambda^{-1/2-10\varepsilon_0}
						.\end{align}
					And the following estimate results from the bootstrap assumptions (\ref{BA:bt1}), (\ref{bt4}), and estimates (\ref{fv6}), (\ref{BA:bt2}):
					\begin{align}
						&\onenorm{\rgeo^{-3/2}\int_{\tstart}^{t}\sgabs{\rgeo^3\lgensmoothfunction\cdot\angnabla\pfour\gfour\cdot\left(\ACC,\rgeo^{-1}\right)}\diff\tau}{\sangle}{p}{\stu}\\
						\notag&\lesssim\lambda^{1/2-4\varepsilon_0}\twonorms{\rgeo(\angnabla\pfour\vvariables,\angnabla\pfour\cvvariables)}{t}{2}{\sangle}{p}{\stu}\twonorms{\ACC}{t}{2}{\sangle}{\infty}{\stu}+\frac{(t-\tstart)^{1/2}}{(t-u)^{1/2}}\twonorms{\rgeo(\angnabla\pfour\vvariables,\angnabla\pfour\cvvariables)}{t}{2}{\sangle}{p}{\coneu}\\
						\notag&\lesssim\lambda^{-1/2}
						,\end{align}
					Employing the bootstrap assumptions (\ref{BA:bt1}), (\ref{bt4}), and estimates (\ref{fv1}), (\ref{BA:bt2}), we deduce
					\begin{align}
						&\onenorm{\rgeo^{-3/2}\int_{\tstart}^{t}\sgabs{\rgeo^3\lgensmoothfunction\cdot\left(\ACC,\rgeo^{-1}\right)\cdot\left(\ACC,\rgeo^{-1}\right)\cdot\pfour\gfour} \diff\tau}{\sangle}{p}{\stu}\\
						\notag\lesssim&\lambda^{1/2-4\varepsilon_0}\int_{\tstart}^{t}\onenorm{\rgeo\left(\ACC,\rgeo^{-1}\right)}{\sangle}{p}{\sTauu}\onenorm{\ACC}{\sangle}{\infty}{\sTauu}\onenorm{\pfour\vvariables,\pfour\cvvariables}{\sangle}{\infty}{\sTauu}\diff\tau\\
						\notag&+\rgeo^{-1/2}\int_{\tstart}^{t}\onenorm{\pfour\vvariables,\pfour\cvvariables}{\sangle}{\infty}{\sTauu} \diff\tau\\
						\notag\lesssim&\lambda^{1/2-4\varepsilon_0}\twonorms{\ACC}{t}{2}{\sangle}{\infty}{\coneu}\twonorms{\pfour\vvariables,\pfour\cvvariables}{t}{2}{\sangle}{\infty}{\coneu}+\frac{(t-\tstart)^{1/2}}{(t-u)^{1/2}}\twonorms{\pfour\vvariables,\pfour\cvvariables}{t}{2}{\sangle}{\infty}{\coneu}\\
						\notag\lesssim&\lambda^{-1/2-4\varepsilon_0}
						.\end{align}
					Utilizing the bootstrap assumptions (\ref{BA:bt1}), (\ref{bt4}), estimate (\ref{BA:bt2}), and the previously proven result (\ref{mr7}) for $\twonorms{\rgeo\angnabla\hchi}{t}{2}{\sangle}{p}{\coneu}$, we get
					\begin{align}
						\onenorm{\rgeo^{-3/2}\int_{\tstart}^{t}\sgabs{\rgeo^3\lgensmoothfunction\cdot\angnabla\hchi\cdot\ACC} \diff\tau}{\sangle}{p}{\sTauu}\lesssim\lambda^{1/2-4\varepsilon_0}\twonorms{\rgeo\angnabla\hchi}{t}{2}{\sangle}{p}{\coneu}\twonorms{\ACC}{t}{2}{\sangle}{\infty}{\coneu}\lesssim\lambda^{-1/2-2\varepsilon_0}
						.\end{align}
				\end{proof}
				\begin{proof}[Proof of $\holdertwonorms{\hchi}{t}{2}{\sangle}{0}{\delta_0}{\coneu}\lesssim\lambda^{-1/2}$ in (\ref{mr7})]
					Using the Sobolev inequality (\ref{ES:sobolevfour}) with $Q=p$, and the previously proven estimates (\ref{mr7}) for $\rgeo\angnabla\hchi$ and (\ref{mr1}) for $\hchi$, we arrive at
					\begin{align}
						\holdertwonorms{\hchi}{t}{2}{\sangle}{0}{\delta_0}{\coneu}\lesssim\twonorms{\rgeo\angnabla\hchi}{t}{2}{\sangle}{p}{\coneu}+\twonorms{\hchi}{t}{2}{\sangle}{2}{\coneu}\lesssim\lambda^{-1/2}
						.\end{align}
				\end{proof}
				\begin{proof}[Proof of $\holderthreenorms{\chismall,\hchi, \gtr\upchi-\frac{2}{\rgeo}}{t}{2}{u}{\infty}{\sangle}{0}{\delta_0}{\intregion}\lesssim\lambda^{-1/2-4\varepsilon_0}$ in (\ref{mr14})]
					We first bound $\chismall$. Using equation (\ref{EQ:Lz}), and dividing $\Lunit(\rgeo^2\chismall)$ by $\rgeo^2(t,u)$, we obtain
					\begin{align}\label{ES211}
						\sgabs{\chismall}\lesssim&\rgeo^{-2}\lim\limits_{\tau\downarrow\tstart}\sgabs{\rgeo^2\chismall}+\rgeo^{-2}\int_{\tstart}^{t}\sgabs{\rgeo^2\lgensmoothfunction\cdot\boxg\gfour}\diff\tau\\
						\notag&+\rgeo^{-2}\int_{\tstart}^{t}\left|\rgeo^2\lgensmoothfunction\cdot\rgeo^{-1}\cdot\pfour\gfour\right.\\
						\notag&+\left.\rgeo^2(\chismall,\hchi)\cdot(\chismall,\hchi,\pfour\gfour) \right|\diff\tau
						.\end{align}
					Applying the initial condition (\ref{ES:initialcontipone}), and taking the $L^\infty_x$-norm of (\ref{ES211}), we deduce
					\begin{align}\label{ES:211}
						\holdertwonorms{\chismall}{u}{\infty}{\sangle}{0}{\delta_0}{\Stint}\lesssim&\holderthreenorms{\boxg\gfour}{t}{1}{u}{\infty}{\sangle}{0}{\delta_0}{\intregion}\\
						\notag&+\holderthreenorms{\chismall,\hchi,\pfour\vvariables,\pfour\cvvariables}{t}{2}{u}{\infty}{\sangle}{0}{\delta_0}{\intregion}^2+\mathcal{M}\left(\holdertwonorms{\pfour\vvariables,\pfour\cvvariables}{u}{\infty}{\sangle}{0}{\delta_0}{\Stint}\right)
						.\end{align}
					Via using the Hardy-Littlewood maximal inequality (\ref{ES:hardylittlewood}), the bootstrap assumptions (\ref{BA:bt1}), and estimate (\ref{BA:bt2.1}), we have
					\begin{align}
						\norm{\mathcal{M}\left(\holdertwonorms{\pfour\vvariables,\pfour\cvvariables}{u}{\infty}{\sangle}{0}{\delta_0}{\Stint}\right)}_{L^2_t}\lesssim \holdertwonorms{\pfour\vvariables,\pfour\cvvariables}{t}{2}{x}{0}{\delta_0}{\intregion}\lesssim\lambda^{-1/2-4\varepsilon_0}
						.\end{align}
					Appealing to the bootstrap assumption (\ref{BA:bt1}), (\ref{bt7}) and estimate (\ref{BA:bt2.1}), (\ref{fv9}), we get
					\begin{subequations}
						\begin{align}
							\holdertwonorms{\boxg\gfour}{t}{1}{x}{0}{\delta_0}{\intregion}&\lesssim\lambda^{-1-8\varepsilon_0},\\
							\label{ES212}\holderthreenorms{\chismall,\hchi,\pfour\vvariables,\pfour\cvvariables}{t}{2}{u}{\infty}{\sangle}{0}{\delta_0}{\intregion}^2&\lesssim\lambda^{-1}
							.\end{align}
					\end{subequations}
					Taking the $L^2_t$-norm of (\ref{ES:211}), combining (\ref{ES211})-(\ref{ES212}), we conclude the desired result for $\chismall$. 
					
					Now, by using (\ref{DE:tracechismall}), and $\gtr\upchi-\frac{2}{\rgeo}=\chismall-\Chfour_{\Lunit}$, we obtain that
					\begin{align}\label{ES213}
						\holderthreenorms{\gtr\upchi-\frac{2}{\rgeo}}{t}{2}{u}{\infty}{\sangle}{0}{\delta_0}{\intregion}=\holderthreenorms{\chismall}{t}{2}{u}{\infty}{\sangle}{0}{\delta_0}{\intregion}+\holderthreenorms{\lgensmoothfunction\cdot\pfour\gfour}{t}{2}{u}{\infty}{\sangle}{0}{\delta_0}{\intregion}
						.\end{align}
					By the previously proven (\ref{mr10}) (in Proposition \ref{frequentlyusedestimates}), the bootstrap assumption (\ref{BA:bt1}), (\ref{bt7}), estimate (\ref{BA:bt2.1}), there holds
					\begin{align}\label{ES214}
						\holderthreenorms{\lgensmoothfunction\cdot\pfour\gfour}{t}{2}{u}{\infty}{\sangle}{0}{\delta_0}{\intregion}\lesssim\holderthreenorms{\lgensmoothfunction}{t}{\infty}{u}{\infty}{\sangle}{0}{\delta_0}{\intregion}\holderthreenorms{\pfour\vvariables,\pfour\cvvariables}{t}{2}{u}{\infty}{\sangle}{0}{\delta_0}{\intregion}\lesssim\lambda^{-1/2-4\varepsilon_0}
						.\end{align}
					Combining (\ref{ES213})-(\ref{ES214}) with the result for $\chismall$, we conclude the proof for $\gtr\upchi-\frac{2}{\rgeo}$.\\
					
					We now prove $\holderthreenorms{\hchi}{t}{2}{u}{\infty}{\sangle}{0}{\delta_0}{\intregion}\lesssim\lambda^{-1/2-4\varepsilon_0}$.
					Plugging equation (\ref{EQ:angdvhchi}) into the Hodge estimate (\ref{ES:hodgefive}) with $Q=p$, we obtain
					\begin{align}\label{ES215}
						\holderthreenorms{\hchi}{t}{2}{u}{\infty}{\sangle}{0}{\delta_0}{\intregion}\lesssim&\holderthreenorms{\chismall,\pfour\gfour}{t}{2}{u}{\infty}{\sangle}{0}{\delta_0}{\intregion}+\threenorms{\rgeo\lgensmoothfunction\cdot(\ACC,\rgeo^{-1})\cdot\pfour\gfour}{t}{2}{u}{\infty}{\sangle}{p}{\intregion}
						.\end{align}
					By the bootstrap assumptions (\ref{BA:bt1}), (\ref{bt5}), (\ref{bt2}), and estimates (\ref{BA:bt2.1}), (\ref{fv1}), we have:
					\begin{align}
						\label{hchit2uinftysangledelta0mint}&\threenorms{\rgeo\lgensmoothfunction\cdot\left(\ACC,\rgeo^{-1}\right)\cdot\pfour\gfour}{t}{2}{u}{\infty}{\sangle}{p}{\intregion}\\
						\notag\lesssim&\norm{\left(\lambda^{1/2-4\varepsilon_0}\onenorm{\rgeo^{1/2}(\pfour\vvariables,\pfour\cvvariables)}{\sangle}{p}{\stu},\onenorm{\rgeo\left(\CC,\rgeo^{-1}\right)}{\sangle}{p}{\stu}\right)\onenorm{\pfour\vvariables,\pfour\cvvariables}{\sangle}{\infty}{\stu}}_{L^2_tL^\infty_u}\\
						\notag\lesssim&\threenorms{\pfour\vvariables,\pfour\cvvariables}{t}{2}{u}{\infty}{\sangle}{\infty}{\region}\\
						\notag\lesssim&\lambda^{-1/2-4\varepsilon_0}
						.\end{align}
					Combining (\ref{ES215})-(\ref{hchit2uinftysangledelta0mint}) with the previously proven result for $\chismall$, we conclude the desired estimate.
				\end{proof}
				
				The proof of $\holderthreenorms{\chismall,\hchi,\gtr\upchi-\frac{2}{\rgeo}}{t}{\frac{q}{2}}{u}{\infty}{\sangle}{0}{\delta_0}{\region}\lesssim\lambda^{\frac{2}{q}-1-4\varepsilon_0\left(\frac{4}{q}-1\right)}$ in (\ref{mr11}) follows a similar fashion as the previous proof of $\holderthreenorms{\chismall,\hchi, \gtr\upchi-\frac{2}{\rgeo}}{t}{2}{u}{\infty}{\sangle}{0}{\delta_0}{\intregion}\lesssim\lambda^{-1/2-4\varepsilon_0}$, provided the following bound on initial hypersurface $\Szero$:
				
				When $u<0$, for the first term on the RHS of (\ref{ES211}), by (\ref{EQ:zinitial}), we have:
				\begin{align}
					\holderthreenorms{\rgeo^{-2}\lim\limits_{\tau\downarrow\tstart}\left(\rgeo^2\chismall\right)}{t}{\frac{q}{2}}{u}{\infty}{\sangle}{0}{\delta_0}{\region}&\lesssim\left(\int_{0}^{\Trescale}\left(\rgeo^{-1/2}\lambda^{-1/2}\right)^{\frac{q}{2}}\diff\tau\right)^{\frac{2}{q}}\\
					\notag&\lesssim\left(\lambda^{(1-8\varepsilon_0)(1-\frac{q}{4})}\lambda^{-\frac{q}{4}}\right)^{\frac{2}{q}}=\lambda^{\frac{2}{q}-1-4\varepsilon_0\left(\frac{4}{q}-1\right)}
					.\end{align}
				
				\begin{proof}[Proof of (\ref{mr15}) and (\ref{mr16})]
					We first prove (\ref{mr15}). Integrate equation (\ref{EQ:dtrminus2gsphere}) along the integral curves of $\Lunit$ as follows:
					\begin{align}
						&\left\{\rgeo^{-2}\gsphere\left(\deriasphere,\deribsphere\right)-\esphere\left(\deriasphere,\deribsphere\right)\right\}\\
						\notag=&\lim\limits_{t\downarrow\tstart}	\left\{\rgeo^{-2}\gsphere\left(\deriasphere,\deribsphere\right)-\esphere\left(\deriasphere,\deribsphere\right)\right\}\\
						\notag&+\int_{\tstart}^{t}\left(\chismall-\Chfour_{\Lunit}\right)\left\{\rgeo^{-2}\gsphere\left(\deriasphere,\deribsphere\right)-\esphere\left(\deriasphere,\deribsphere\right)\right\}\diff\tau\\
						\notag&+\int_{\tstart}^{t}\left(\chismall-\Chfour_{\Lunit}\right)\esphere\left(\deriasphere,\deribsphere\right)+\frac{2}{\rgeo^2}\hchi\left(\deriasphere,\deribsphere\right)\diff\tau
						.\end{align}
					Applying the initial conditions (\ref{ES:initialfoliation3}) and (\ref{ES:initialcontip3}), noticing that $\hchi\left(\deriasphere,\deribsphere\right)=\rgeo^2\hchi_{AB}$, and using the Gr\"onwall's inequality, we derive
					\begin{align}
						&\sgabs{\rgeo^{-2}\gsphere\left(\deriasphere,\deribsphere\right)-\esphere\left(\deriasphere,\deribsphere\right)}\\
						\notag&\lesssim\left(\lambda^{-4\varepsilon_0}+\twonorms{\chismall,\pfour\gfour,\hchi}{t}{1}{\sangle}{\infty}{\coneu}\right)\exp\left(\twonorms{\chismall,\pfour\gfour}{t}{1}{\sangle}{\infty}{\coneu}\right)
						.\end{align}
					Employing the bootstrap assumptions (\ref{BA:bt1}), estimates (\ref{BA:bt2}), H\"older inequality, and the previously proven result (\ref{mr11}) with $q=4$, it follows that
					\begin{align}
						\norm{\rgeo^{-2}\gsphere\left(\deriasphere,\deribsphere\right)-\esphere\left(\deriasphere,\deribsphere\right)}_{L^\infty(\region)}\lesssim\lambda^{-4\varepsilon_0}
						.\end{align}
					We now prove (\ref{mr16}), we first apply $\dericsphere$ to equation (\ref{EQ:dtrminus2gsphere}). Note that it follows from (\ref{DE:tracechismall}) that $\chismall-\Chfour_{\Lunit}=\gtr\upchi-\frac{2}{\rgeo}$. Since $\Lunit$ and $\dericsphere$ commute, there holds
					\begin{align}
						&\tderivative\left\{\dericsphere\left[\rgeo^{-2}\gsphere\left(\deriasphere,\deribsphere\right)-\esphere\left(\deriasphere,\deribsphere\right)\right]\right\}=\dericsphere\gtr\upchi\rgeo^{-2}\gsphere\left(\deriasphere,\deribsphere\right)\\
						\notag&+\left(\gtr\upchi-\frac{2}{\rgeo}\right)\dericsphere\left\{\rgeo^{-2}\gsphere\left(\deriasphere,\deribsphere\right)-\esphere\left(\deriasphere,\deribsphere\right)\right\}\\
						\notag&+\left(\gtr\upchi-\frac{2}{\rgeo}\right)\dericsphere\esphere\left(\deriasphere,\deribsphere\right)+2\rgeo^{-2}\dericsphere\hchi\left(\deriasphere,\deribsphere\right)
						.\end{align}
					Integrating the above equation along the integral curves of $\Lunit$, taking the $L^p_\sangle$ norm, and using the initial condition (\ref{ES:initialfoliation4}), (\ref{ES:initialcontip4}), and then applying the Gr\"onwall's inequality, we have:
					\begin{align}
						&\onenorm{\dericsphere\left\{\rgeo^{-2}\gsphere\left(\deriasphere,\deribsphere\right)-\esphere\left(\deriasphere,\deribsphere\right)\right\}}{\sangle}{p}{\stu}\\
						\notag\lesssim&\left(\lambda^{-4\varepsilon_0}+\twonorms{\rgeo\angnabla\gtr\upchi}{t}{1}{\sangle}{p}{\coneu}+\twonorms{\left(\gtr\upchi-\frac{2}{\rgeo}\right)\dericsphere\esphere}{t}{1}{\sangle}{p}{\coneu}\right.\\
						\notag&\left.+\twonorms{\rgeo\angnabla\hchi}{t}{1}{\sangle}{p}{\coneu}+\twonorms{\Chfour\cdot\hchi}{t}{1}{\sangle}{p}{\coneu}\right)\cdot\exp\left(\twonorms{\gtr\upchi-\frac{2}{\rgeo}}{t}{1}{\sangle}{\infty}{\coneu}\right)
						.\end{align}
					By using the bootstrap assumptions (\ref{BA:bt1}), estimates (\ref{BA:bt2}), and the previously proven results (\ref{mr1}), (\ref{mr7}) and (\ref{mr11}) with $q=4$, we conclude
					\begin{align}
						\onenorm{\dericsphere\left\{\rgeo^{-2}\gsphere\left(\deriasphere,\deribsphere\right)-\esphere\left(\deriasphere,\deribsphere\right)\right\}}{\sangle}{p}{\stu}\lesssim\lambda^{-4\varepsilon_0}
						.\end{align}
				\end{proof}
				\begin{proof}[Proof of $\twonorms{\frac{\nulllapse^{-1}-1}{\rgeo}}{t}{2}{x}{\infty}{\region}\lesssim\lambda^{-1/2}$ in (\ref{mr9})]
					We first bound $\twonorms{\frac{\nulllapse^{-1}-1}{\rgeo}}{t}{2}{x}{\infty}{\intregion}$. Integrating equation (\ref{EQ:Lunitnulllapse}) along along the integral curves of $\Lunit$ emanating from the cone-tip, and using the initial condition (\ref{ES:initialcontipone}), we have:
					\begin{align}\label{ES235}
						\nulllapse^{-1}-1=-\int_{u}^{t}\left(\nulllapse^{-1}-1\right)\lgensmoothfunction\cdot\pfour\gfour \diff\tau-\int_{u}^{t}\lgensmoothfunction\cdot\pfour\gfour \diff\tau
						.\end{align}
					Using the Gr\"onwall's inequality, we obtain
					\begin{align}\label{ES234}
						\sgabs{\frac{\nulllapse^{-1}-1}{\rgeo}}\lesssim\mathcal{M}\left(\onenorm{\pfour\vvariables,\pfour\cvvariables}{\sangle}{\infty}{\stu}\right)\cdot\exp\left(\int_{u}^{t}\lgensmoothfunction\cdot\pfour\gfour\diff\tau\right)
						.\end{align}
					Hence, by the bootstrap assumption (\ref{BA:bt1}), and estimate (\ref{BA:bt2}), there holds
					\begin{align}\label{ES230}
						\twonorms{\frac{\nulllapse^{-1}-1}{\rgeo}}{t}{2}{x}{\infty}{\intregion}\lesssim\twonorms{\pfour\vvariables,\pfour\cvvariables}{t}{2}{x}{\infty}{\region}\lesssim\lambda^{-1/2-4\varepsilon_0}
						.\end{align}
					Now we consider the case when $u<0$. Integration of (\ref{EQ:Lunitnulllapse}) along along the integral curves of $\Lunit$ emanating from $\Sigma_{0}$ yields that
					\begin{align}
						\nulllapse^{-1}-1=\lapse^{-1}-1-\int_{0}^{t}\left(\nulllapse^{-1}-1\right)\lgensmoothfunction\cdot\pfour\gfour \diff\tau-\int_{u}^{t}\lgensmoothfunction\cdot\pfour\gfour \diff\tau.
					\end{align}
					 Applying the Gr\"onwall's inequality, we derive
					\begin{align}\label{ES231}
						\sgabs{\frac{\nulllapse^{-1}-1}{\rgeo}}\leq\sgabs{\frac{\nulllapse^{-1}-\lapse^{-1}}{\rgeo}}+\sgabs{\frac{\lapse^{-1}-1}{\rgeo}}\lesssim\left\{\mathcal{M}\left(\onenorm{\pfour\vvariables,\pfour\cvvariables}{\sangle}{\infty}{\stu}\right)+\sgabs{\frac{\lapse^{-1}-1}{\rgeo}}\right\}\cdot\exp\left(\int_{u}^{t}\lgensmoothfunction\cdot\pfour\gfour\diff\tau\right)
						.\end{align}
					Using the initial condition (\ref{ES:lapseminusone}), we get
					\begin{align}\label{ES232}
						\twonorms{\frac{\lapse^{-1}-1}{\rgeo}}{t}{2}{x}{\infty}{\region}\lesssim\onenorm{\frac{\lapse^{-1}-1}{w^{1/2}}}{x}{\infty}{\region}\left(\int_{0}^{\Trescale}\frac{w}{(\tau+w)^2}\diff\tau\right)^{1/2}\lesssim\lambda^{-1/2}
						.\end{align}
					Together with (\ref{ES231}), these estimates then lead to
					\begin{align}\label{ES233}
						\twonorms{\frac{\nulllapse^{-1}-1}{\rgeo}}{t}{2}{x}{\infty}{\extregion}\lesssim\lambda^{-1/2}
						.\end{align}
					Combining (\ref{ES230}) and (\ref{ES233}), we conclude the desired estimate.
				\end{proof}
				\begin{proof}[Proof of $\twonorms{\rgeo(\angD_{\Lunit},\angnabla)\left(\frac{\nulllapse^{-1}-1}{\rgeo}\right)}{t}{2}{\sangle}{p}{\coneu}\lesssim\lambda^{-1/2}$ in (\ref{mr9})] We first prove $\twonorms{\rgeo\angD_{\Lunit}\left(\frac{\nulllapse^{-1}-1}{\rgeo}\right)}{t}{2}{\sangle}{p}{\coneu}\lesssim\lambda^{-1/2}$.
					Noting that by (\ref{EQ:Lunitnulllapse}) and the fact that $\Lunit(\rgeo)=1$, the following identity is valid.
					\begin{align}
						\rgeo\angD_{\Lunit}\left(\frac{\nulllapse^{-1}-1}{\rgeo}\right)=-\nulllapse^{-1}\lgensmoothfunction\cdot\pfour\gfour-\frac{\nulllapse^{-1}-1}{\rgeo}
						.\end{align}
					Employing the bootstrap assumption (\ref{BA:bt1}), estimate (\ref{BA:bt2}), the previously proven result (\ref{mr18}), and the result (\ref{ES234}) for $\frac{\nulllapse^{-1}-1}{\rgeo}$ in the last proof, we deduce
					\begin{align}
						\twonorms{\rgeo\angD_{\Lunit}\left(\frac{\nulllapse^{-1}-1}{\rgeo}\right)}{t}{2}{\sangle}{p}{\coneu}\lesssim\twonorms{\pfour\vvariables,\pfour\cvvariables}{t}{2}{\sangle}{p}{\coneu}+\twonorms{\frac{\nulllapse^{-1}-1}{\rgeo}}{t}{2}{\sangle}{p}{\coneu}\lesssim\lambda^{-1/2}
						.\end{align}
					Now we bound $\rgeo\angnabla\left(\frac{\nulllapse^{-1}-1}{\rgeo}\right)$. Recall (\ref{E:CONNECTIONCOEFFICIENT}), $\upzeta=\angnabla\ln\nulllapse+\lgensmoothfunction\cdot\pfour\gfour$. Therefore, we derive
					\begin{align}
						\rgeo\angnabla\left(\frac{\nulllapse^{-1}-1}{\rgeo}\right)=-\nulllapse^{-1}\angnabla\ln\nulllapse=\nulllapse^{-1}\left(-\upzeta+\lgensmoothfunction\cdot\pfour\gfour\right)
						.\end{align}
					By estimates (\ref{fv3}), the proven results (\ref{mr1}) and (\ref{mr18}), we conclude
					\begin{align}
						\twonorms{\rgeo\angnabla\left(\frac{\nulllapse^{-1}-1}{\rgeo}\right)}{t}{2}{\sangle}{p}{\coneu}\lesssim\lambda^{-1/2}
						.\end{align}
				\end{proof}
				\begin{proof}[Proof of $\threenorms{\frac{\nulllapse^{-1}-1}{\rgeo^{1/2}}}{t}{\infty}{u}{\infty}{\sangle}{2p}{\region}\lesssim
					\lambda^{-1/2}$ in (\ref{mr9})]
					By the Sobolev inequality (\ref{ES:sobolevtwo}), we have:
					\begin{align}
						&\threenorms{\frac{\nulllapse^{-1}-1}{\rgeo^{1/2}}}{t}{\infty}{u}{\infty}{\sangle}{2p}{\region}^2\lesssim\left(\threenorms{\rgeo\angD_{\Lunit}\left(\frac{\nulllapse^{-1}-1}{\rgeo}\right)}{t}{2}{u}{\infty}{\sangle}{p}{\region}+\threenorms{\frac{\nulllapse^{-1}-1}{\rgeo}}{t}{2}{u}{\infty}{\sangle}{p}{\region}\right)\threenorms{\frac{\nulllapse^{-1}-1}{\rgeo}}{t}{2}{u}{\infty}{\sangle}{\infty}{\region}
						.\end{align}
					Appealing to the proven first and third estimates of (\ref{mr9}), we conclude the desired result.
				\end{proof}
				\begin{proof}[Proof of (\ref{mr19})]
					We first prove $\twonorms{\angnabla\cphi}{t}{2}{\sangle}{p}{\coneu}\lesssim\lambda^{-1/2}$. Plugging equation (\ref{EQ:Lunitangnablacphi}) into equation (\ref{EQ:transportlemma2}) and using estimate (\ref{EQ:rtwomupxi}) with $m=\frac{1}{2}$, by (\ref{mr17}), we have:
					\begin{align}
						\sgabs{\rgeo\angnabla\cphi}=&\lim\limits_{t\downarrow\tstart}\sgabs{\rgeo\angnabla\cphi}+\int_{\tstart}^{t}\sgabs{\rgeo\left(\frac{1}{2}\left(\chismall-\Chfour_{\Lunit}\right)+\lgensmoothfunction\cdot\hchi\right)\angnabla\cphi} \diff\tau\\
						\notag&+\int_{\tstart}^{t}\sgabs{\rgeo\angnabla\left(\chismall-\Chfour_{\Lunit}\right)}\diff\tau
						.\end{align}
					Applying the Gr\"onwall's inequality, and using the proven results (\ref{mr11}) with $q=4$, we derive
					\begin{align}
						\label{rcphi}\sgabs{\rgeo\angnabla\cphi}\lesssim&\left(\lim\limits_{t\downarrow\tstart}\sgabs{\rgeo\angnabla\cphi}+\int_{\tstart}^{t}\sgabs{\rgeo\angnabla\left(\chismall-\Chfour_{\Lunit}\right)}\diff\tau\right)\\
						\notag&\cdot\exp\left(\twonorms{\gtr\upchi-\frac{2}{\rgeo},\hchi}{t}{1}{\sangle}{\infty}{\coneu}\right)\\
						\notag\lesssim&\lim\limits_{t\downarrow\tstart}\sgabs{\rgeo\angnabla\cphi}+\int_{\tstart}^{t}\sgabs{\rgeo\angnabla\left(\chismall-\Chfour_{\Lunit}\right)}\diff\tau
						.\end{align}
					We now consider the initial conditions. When $u\geq 0$, by initial condition (\ref{ES:initialcontip3}) and (\ref{ES:initialcontip4}),  we have $\lim\limits_{t\downarrow\tstart}\sgabs{\rgeo\angnabla\cphi}=0$. When $u<0$, by (\ref{ES:in5}), $\onenorm{\lim\limits_{t\downarrow0}\rgeo^{1/2}\angnabla\cphi}{\sangle}{p}{S_w}\lesssim\lambda^{-1/2}$.
					Now dividing both sides of (\ref{rcphi}) by $\rgeo(t,u)$, then taking the $L^2_tL^p_\sangle$-norm and using the Hardy-Littlewood maximal inequality (\ref{ES:hardylittlewood}), we deduce:
					\begin{align}
						\twonorms{\angnabla\cphi}{t}{2}{\sangle}{p}{\coneu}\lesssim\lambda^{-1/2}\left(\int_{0}^{\Trescale}\frac{w}{(\tau+w)^2}\diff\tau\right)^{1/2}+\twonorms{\rgeo\left(\angnabla\chismall,\angnabla\pfour\gfour\right)}{t}{2}{\sangle}{p}{\coneu}
						.\end{align}
					By estimate (\ref{fv6}) and proven result (\ref{mr7}), we conclude the desired estimate. 
					
					Now we bound $\rgeo^{1/2}\angnabla\cphi$. Dividing both sides of (\ref{rcphi}) by $\rgeo^{1/2}$, taking the $L_\sangle^p$ norm, by estimate (\ref{fv6}) and the proven result (\ref{mr7}), we arrive at
					\begin{align}
						\onenorm{\rgeo^{1/2}\angnabla\cphi}{\sangle}{p}{\stu}\lesssim\lambda^{-1/2}+\frac{\left(t-\tstart\right)^{1/2}}{(t-u)^{1/2}}\twonorms{\rgeo\left(\angnabla\chismall,\angnabla\pfour\gfour\right)}{t}{2}{\sangle}{p}{\coneu}\lesssim\lambda^{-1/2}
						.\end{align}
					
					Now we prove for $\rgeo\Lunit\angnabla\cphi$. Using equation (\ref{EQ:Lunitangnablacphi}), we derive:
					\begin{align}
						&\sgabs{\rgeo\Lunit\angnabla\cphi}\lesssim\sgabs{\rgeo\angnabla\left(\chismall-\Chfour_{\Lunit}\right)+\rgeo\left(\chismall,\hchi,\pfour\gfour\right)\cdot\angnabla\cphi+\angnabla\cphi}
						.\end{align}
					Employing the bootstrap assumptions (\ref{BA:bt1}), estimates (\ref{BA:bt2}), (\ref{fv6}), and the proven results (\ref{mr11}) with $q=4$, (\ref{mr7}), we obtain
					\begin{align}
						\twonorms{\rgeo\Lunit\angnabla\cphi}{t}{2}{\sangle}{p}{\coneu}\lesssim&\twonorms{\rgeo\left(\angnabla\chismall,\angnabla\pfour\gfour\right)}{t}{2}{\sangle}{p}{\coneu}+\twonorms{\angnabla\cphi}{t}{2}{\sangle}{p}{\coneu}\\
						\notag&+\lambda^{1/2-4\varepsilon_0}\twonorms{\chismall,\hchi,\pfour\gfour}{t}{2}{\sangle}{\infty}{\coneu}\twonorms{\rgeo^{1/2}\angnabla\cphi}{t}{\infty}{\sangle}{p}{\coneu}\\
			\notag\lesssim&\lambda^{-1/2}
						.\end{align}
				\end{proof}
				\begin{proof}[Proof of $\twonorms{\rgeo\angnabla\upzeta}{t}{2}{\sangle}{p}{\coneu},\holdertwonorms{\upzeta}{t}{2}{\sangle}{0}{\delta_0}{\coneu}\lesssim\lambda^{-1/2}$ in (\ref{mr7})]
					Plugging equation (\ref{EQ:angdivupzeta}) and (\ref{EQ:angcurlupzeta}) into the Hodge estimate (\ref{ES:hodgeone}) with $Q=p$, we derive
					\begin{align}
						\twonorms{\rgeo\angnabla\upzeta}{t}{2}{\sangle}{p}{\coneu}\lesssim&	\twonorms{\rgeo\lgensmoothfunction\cdot\boxg\gfour}{t}{2}{\sangle}{p}{\coneu}+	\twonorms{\rgeo\angnabla\pfour\gfour}{t}{2}{\sangle}{p}{\coneu}\\
						\notag&+	\twonorms{\rgeo\lgensmoothfunction\cdot(\ACC,\angnabla\cphi,\rgeo^{-1})\cdot\ACC}{t}{2}{\sangle}{p}{\coneu}
						.\end{align}
					Using estimates (\ref{fv6}), (\ref{fv10}), the proven results (\ref{mr1}), (\ref{mr3}), (\ref{mr11}), (\ref{mr19}), (\ref{rmk1}), we deduce
					\begin{align}\label{ES241}
						\twonorms{\rgeo\angnabla\upzeta}{t}{2}{\sangle}{p}{\coneu}\lesssim\lambda^{-1/2}+\lambda^{-4\varepsilon_0}\twonorms{\upzeta}{t}{2}{\sangle}{\infty}{\coneu}
						.\end{align}
					Applying the Sobolev inequality (\ref{ES:sobolevfour}) with $Q=p$, and the already proven result (\ref{mr1}), we get
					\begin{align}\label{ES242}
						\holdertwonorms{\upzeta}{t}{2}{\sangle}{0}{\delta_0}{\coneu}\lesssim\twonorms{\rgeo\angnabla\upzeta}{t}{2}{\sangle}{p}{\coneu}+\twonorms{\upzeta}{t}{2}{\sangle}{2}{\coneu}\lesssim\lambda^{-1/2}+\twonorms{\rgeo\angnabla\upzeta}{t}{2}{\sangle}{p}{\coneu}
						.\end{align}
					Combining (\ref{ES241}) and (\ref{ES242}), we conclude the desired estimates.
				\end{proof}
				\begin{proof}[Proof of $\twonorms{\rgeo\mass}{t}{2}{\sangle}{p}{\coneu}\lesssim\lambda^{-1/2}$ in (\ref{mr7})]
					Using equation (\ref{EQ:mass}), by the estimates (\ref{fv6}), (\ref{fv10}), and the proven results (\ref{mr1}), (\ref{mr3}), (\ref{mr7}), (\ref{mr19}), (\ref{rmk1}), we have
					\begin{align}
						\twonorms{\rgeo\mass}{t}{2}{\sangle}{p}{\coneu}\lesssim&\twonorms{\rgeo\lgensmoothfunction\cdot\boxg\gfour}{t}{2}{\sangle}{p}{\coneu}+\twonorms{\rgeo\angnabla\pfour\gfour}{t}{2}{\sangle}{p}{\coneu}
						\notag&+\twonorms{\rgeo\lgensmoothfunction\cdot(\ACC,\angnabla\cphi,\rgeo^{-1})\cdot\ACC}{t}{2}{\sangle}{p}{\coneu}\\
						\notag\lesssim&\lambda^{-1/2}
						.\end{align}
				\end{proof}
				\begin{proof}[Proof of $\twonorms{\upzeta}{t}{2}{x}{\infty}{\intregion}\lesssim\lambda^{-1/2-4\varepsilon_0}$ in (\ref{mr145}) and $\twonorms{\upzeta}{t}{\frac{q}{2}}{x}{\infty}{\region}\lesssim\lambda^{\frac{2}{q}-1-4\varepsilon_0\left(\frac{4}{q}-1\right)}$ in (\ref{mr115})]
					Plugging equations (\ref{EQ:angdivupzeta}) and (\ref{EQ:angcurlupzeta}) into the Hodge estimate (\ref{ES:hodgesix}) with $Q=p$, $c=2$, $\mathfrak{F}=\lgensmoothfunction\cdot\pfour\gfour$, we obtain
					\begin{align}
						\label{upzetat2xinfty}\twonorms{\upzeta}{t}{2}{x}{\infty}{\intregion}\lesssim&\norm{\upnu^{\delta^\prime}\littlewood\left(\lgensmoothfunction\cdot\pfour\gfour\right)}_{L^2_tL^\infty_ul^2_\upnu L^\infty_\sangle(\intregion)}+\threenorms{\lgensmoothfunction\cdot\pfour\gfour}{t}{2}{u}{\infty}{\sangle}{\infty}{\intregion}\\
						\notag&+\threenorms{\rgeo\lgensmoothfunction\cdot\boxg\gfour}{t}{2}{u}{\infty}{\sangle}{p}{\intregion}+\threenorms{\rgeo\lgensmoothfunction\cdot(\hchi,\upzeta)\cdot(\hchi,\upzeta)}{t}{2}{u}{\infty}{\sangle}{p}{\coneu}\\
						\notag&+\threenorms{\rgeo\lgensmoothfunction\cdot\angnabla\cphi\cdot(\pfour\gfour,\upzeta)}{t}{2}{u}{\infty}{\sangle}{p}{\coneu}\\
						\notag&+\threenorms{\rgeo\lgensmoothfunction\cdot(\ACC,\rgeo^{-1})\cdot\pfour\gfour}{t}{2}{u}{\infty}{\sangle}{p}{\coneu}
						.\end{align}
					Via employing the bootstrap assumptions (\ref{BA:bt1}), estimates (\ref{BA:bt2.1}), (\ref{fv10}), the already proven results (\ref{mr3}), (\ref{mr14}), (\ref{mr19}), and (\ref{rzt2sanglep}), we arrive at
					\begin{align}\label{ES251}
						\twonorms{\upzeta}{t}{2}{x}{\infty}{\intregion}\lesssim\lambda^{-1/2-4\varepsilon_0}+\lambda^{-4\varepsilon_0}\twonorms{\upzeta}{t}{2}{x}{\infty}{\intregion}
						.\end{align}
					By absorbing the second term on the right of (\ref{ES251}) into the left, we conclude the desired result $\twonorms{\upzeta}{t}{2}{x}{\infty}{\intregion}\lesssim\lambda^{-1/2-3\varepsilon_0}$.
                    
                    To bound $\twonorms{\upzeta}{t}{\frac{q}{2}}{x}{\infty}{\region}$, we apply the H\"older's inequality to (\ref{upzetat2xinfty}) and {use the similar argument given above (with (\ref{mr14}) replaced by (\ref{mr11})) to deduce:}
					\begin{align}\label{ES252}
						\twonorms{\upzeta}{t}{\frac{q}{2}}{x}{\infty}{\region}\lesssim\lambda^{\frac{2}{q}-1-4\varepsilon_0\left(\frac{4}{q}-\frac{1}{4}\right)}+\lambda^{-4\varepsilon_0}\twonorms{\upzeta}{t}{\frac{q}{2}}{x}{\infty}{\region}
						.\end{align}
					By absorbing the second term on the right of (\ref{ES252}) into the left, we conclude the desired result.
				\end{proof}
				\begin{proof}[Proof of (\ref{mr21})-(\ref{mr22})]
					By (\ref{DE:Lunitsigma}) and the initial condition for $\conformalfactor$ (\ref{DE:initialsigma}), it holds that
					\begin{align}
						\conformalfactor=\frac{1}{2}\int_{u}^{t}\Chfour_{\Lunit}\diff\tau
						.\end{align}
					Using the bootstrap assumptions (\ref{BA:bt1}) and estimate (\ref{BA:bt2}), we find that:
					\begin{align}
						\norm{\conformalfactor}_{L^\infty(\intregion)}&\lesssim\lambda^{-8\varepsilon_0},\\
						\norm{\rgeo^{-1/2}\conformalfactor}_{L^\infty(\intregion)}&\lesssim\sup\limits_{t,u}\frac{(t-u)^{1/2}}{\rgeo^{1/2}}\lambda^{-1/2-4\varepsilon_0}\lesssim\lambda^{-1/2-4\varepsilon_0}
						.\end{align}
					
					To prove $\twonorms{\angnabla\conformalfactor}{t}{2}{\sangle}{p}{\coneu}\lesssim\lambda^{-1/2}$, we plug equation (\ref{EQ:Lunitangnablasigma}) into (\ref{EQ:transportlemma2}) and use (\ref{EQ:rtwomupxi}) with $m=\frac{1}{2}$, by the initial condition (\ref{ES:initialcontipone}), we deduce
					\begin{align}
						\sgabs{\rgeo\angnabla\conformalfactor}&\lesssim\int_{u}^{t}\rgeo\left(\sgabs{\angnabla\pfour\gfour}+\sgabs{\ACC\cdot\angnabla\conformalfactor}\right)\diff\tau
						.\end{align}
					Using the Gr\"onwall's inequality, the bootstrap assumption (\ref{BA:bt1}), estimate (\ref{BA:bt2}) and the proven results (\ref{mr7}), we derive
					\begin{align}
						\label{rconformalfactor}\sgabs{\rgeo\angnabla\conformalfactor}\lesssim\int_{u}^{t}\rgeo\sgabs{\angnabla\pfour\gfour}\diff\tau\cdot\exp\left(\twonorms{\ACC}{t}{1}{\sangle}{\infty}{\coneu}\right)\lesssim\int_{u}^{t}\rgeo\sgabs{\angnabla\pfour\gfour}\diff\tau
						.\end{align}
					Dividing both sides of (\ref{rconformalfactor}) by $\rgeo(t,u)$, taking the $L^2_tL^p_\sangle$-norm, using the estimate for Hardy-Littlewood maximal function (\ref{ES:hardylittlewood}) and the estimate (\ref{fv6}), we obtain the desired bound:
					\begin{align}
						\twonorms{\angnabla\conformalfactor}{t}{2}{\sangle}{p}{\coneu}\lesssim\twonorms{\rgeo(\angnabla\pfour\vvariables,\angnabla\pfour\cvvariables)}{t}{2}{\sangle}{p}{\coneu}\lesssim\lambda^{-1/2}
						.\end{align}

					To show $\twonorms{\rgeo^{1/2}\angnabla\conformalfactor}{\sangle}{p}{t}{\infty}{\coneu}\lesssim\lambda^{-1/2}$, first dividing estimate (\ref{rconformalfactor}) by $\rgeo^{1/2}(t,u)$, then taking the $L^p_\sangle$-norm, we have:
					\begin{align}
						\onenorm{\rgeo^{1/2}\angnabla\conformalfactor}{\sangle}{p}{\stu}\lesssim\frac{(t-u)^{1/2}}{\rgeo^{1/2}}\twonorms{\rgeo(\angnabla\pfour\vvariables,\angnabla\pfour\cvvariables)}{t}{2}{\sangle}{p}{\coneu}\lesssim\lambda^{-1/2}
						.\end{align}\\
					
					We now proceed to $\twonorms{\rgeo^{1/2}\Lunit\conformalfactor}{t}{\infty}{\sangle}{2p}{\coneu}\lesssim\lambda^{-1/2-2\varepsilon_0}$. Recall $\Lunit\conformalfactor=\frac{1}{2}\Chfour_{\Lunit}$. Applying the Sobolev equality (\ref{ES:sobolevtwo}), the estimate (\ref{fv3}) and (\ref{fv6}), we get
					\begin{align}
						\twonorms{\rgeo^{1/2}\Lunit\conformalfactor}{t}{\infty}{\sangle}{2p}{\coneu}^2\lesssim\left(\twonorms{\rgeo(\angD_{\Lunit}\pfour\vvariables,\angD_{\Lunit}\pfour\cvvariables)}{\sangle}{p}{t}{2}{\coneu}+\twonorms{\pfour\vvariables,\pfour\cvvariables}{\sangle}{p}{t}{2}{\coneu}\right)\twonorms{\pfour\vvariables,\pfour\cvvariables}{\sangle}{\infty}{t}{2}{\coneu}\lesssim\lambda^{-1-4\varepsilon_0}
						.\end{align}
				\end{proof}
				\begin{proof}[Proof of $\threenorms{\rgeo\left(\modmass,\angnabla\modtorsion\right)}{u}{2}{t}{2}{\sangle}{p}{\intregion}\lesssim\lambda^{-4\varepsilon_0}$ in (\ref{mr24})]
					In this paragraph, we assume $\coneu\subset\intregion$. We will implicitly use the fact that $0\leq u\leq t\leq\Trescale\lesssim\lambda^{1-8\varepsilon_0}$. We start by deriving a preliminary estimate for $\threenorms{\rgeo\angnabla\modtorsion}{u}{2}{t}{2}{\sangle}{p}{\intregion}$. By substituting equations (\ref{EQ:angdivmodtorsion}) and (\ref{EQ:angcurlmodtorsion}) into the Hodge estimate (\ref{ES:hodgeone}), we obtain
					\begin{align}
						\threenorms{\rgeo\angnabla\modtorsion}{u}{2}{t}{2}{\sangle}{p}{\intregion}\lesssim&\threenorms{\rgeo\lgensmoothfunction\cdot(\angnabla\pfour\vvariables,\angnabla\pfour\cvvariables)}{u}{2}{t}{2}{\sangle}{p}{\intregion}+\threenorms{\rgeo\lgensmoothfunction\cdot\boxg\gfour}{u}{2}{t}{2}{\sangle}{p}{\intregion}\\
						\notag&+\threenorms{\rgeo\lgensmoothfunction\cdot\left(\ACC,\rgeo^{-1}\right)\cdot\ACC}{u}{2}{t}{2}{\sangle}{p}{\intregion}+\threenorms{\rgeo\modmass}{u}{2}{t}{2}{\sangle}{p}{\intregion}
						.\end{align}
					By estimates (\ref{fv6}), (\ref{fv10}), (\ref{rmk1}), there holds
					\begin{align}
						\label{ES260}\threenorms{\rgeo\angnabla\modtorsion}{u}{2}{t}{2}{\sangle}{p}{\intregion}&\lesssim\lambda^{-4\varepsilon_0}+\threenorms{\rgeo\modmass}{u}{2}{t}{2}{\sangle}{p}{\intregion}
						.\end{align}
					
					We now bound $\threenorms{\rgeo\modmass}{u}{2}{t}{2}{\sangle}{p}{\intregion}$. Plugging (\ref{EQ:Lunitmodmass}) into equation (\ref{EQ:transportlemma2}) and using (\ref{EQ:rtwomupxi}) with $m=1$, invoking the initial condition on cone-tip (\ref{ES:initialcontipone}), we derive
					\begin{align}
						\sgabs{\rgeo^2\modmass}\lesssim\int_{u}^{t}\rgeo^2\left(\sgabs{\mathfrak{J}_{(1)}}+\sgabs{\mathfrak{J}_{(2)}}\right)\diff\tau+\int_{u}^{t}\rgeo^2\sgabs{\ACC}\cdot\sgabs{\modmass} \diff\tau
						.\end{align}
					Using the Gr\"onwall's inequality, estimates (\ref{fv3}) and the already proven estimate (\ref{mr14}), we have:
					\begin{align}
						\label{r2modmass}\sgabs{\rgeo^2\modmass}&\lesssim\int_{u}^{t}\rgeo^2\left(\sgabs{\mathfrak{J}_{(1)}}+\sgabs{\mathfrak{J}_{(2)}}\right)\diff\tau\exp\left(\twonorms{\ACC}{t}{1}{\sangle}{\infty}{\coneu}\right)\\
						\notag&\lesssim\int_{u}^{t}\rgeo^2\left(\sgabs{\mathfrak{J}_{(1)}}+\sgabs{\mathfrak{J}_{(2)}}\right)\diff\tau
						.\end{align}
					We divide estimate (\ref{r2modmass}) by $\rgeo(t,u)$, and then take the $L^2_uL^2_tL^p_\sangle$ norm. We now estimate terms in $\mathfrak{J}_{(1)}$ and $\mathfrak{J}_{(2)}$ defined in (\ref{DE:J1J2}). First, for $\mathfrak{J}_{(1)}$, by the Hardy-Littlewood maximal inequality  (\ref{ES:hardylittlewood}), estimates (\ref{fv1}), (\ref{fv6}), we conclude:
					\begin{subequations}
						\begin{align}
							\label{ES261}\threenorms{\rgeo^{-1}\int_u^t\rgeo\angnabla\pfour\gfour \diff\tau}{u}{2}{t}{2}{\sangle}{p}{\intregion}&\lesssim\threenorms{\rgeo(\angnabla\pfour\vvariables,\angnabla\pfour\cvvariables)}{u}{2}{t}{2}{\sangle}{p}{\intregion}\lesssim\lambda^{-4\varepsilon_0},\\
							\label{ES262}\threenorms{\rgeo^{-1}\int_{u}^{t}\pfour\gfour \diff\tau}{u}{2}{t}{2}{\sangle}{p}{\intregion}&\lesssim\threenorms{\pfour\vvariables,\pfour\cvvariables}{u}{2}{t}{2}{\sangle}{p}{\intregion}\lesssim\lambda^{-4\varepsilon_0}
							.\end{align}
					\end{subequations}
					Next, we estimate terms in $\mathfrak{J}_{(2)}$. For the first term on the RHS of (\ref{DE:J2}), by estimate (\ref{fv12}), there holds
					\begin{align}
						\label{ES263}\threenorms{\rgeo^{-1}\int_{u}^{t}\sgabs{\rgeo^2\pfour\boxg\gfour}\diff\tau}{u}{2}{t}{2}{\sangle}{p}{\intregion}&\lesssim\norm{\threenorms{\rgeo\pfour\boxg\gfour}{t}{1}{u}{2}{\sangle}{p}{\intregion}}_{L_t^2}\lesssim\lambda^{-12\varepsilon_0}.
					\end{align}
					For the second term on the RHS of (\ref{DE:J2}), by (\ref{fv9}), (\ref{mr3}), we deduce
					\begin{align}
						&\threenorms{\rgeo^{-1}\int_{\tstart}^t\sgabs{\rgeo^2\boxg\gfour\cdot\left(\ACC,\rgeo^{-1}\right)}\diff\tau}{u}{2}{t}{2}{\sangle}{p}{\intregion}\\
						\notag\lesssim&\lambda^{1/2-4\varepsilon_0}\norm{\threenorms{\boxg\gfour}{t}{1}{u}{\infty}{\sangle}{\infty}{\intregion}\cdot\threenorms{\rgeo(\ACC,\rgeo^{-1})}{t}{\infty}{u}{\infty}{\sangle}{p}{\intregion}}_{L^2_t}\\
						\notag\lesssim&\lambda^{-16\varepsilon_0}
						.\end{align}
					For the third term on the right of (\ref{DE:J2}),	by the proven estimates (\ref{mr7}), (\ref{fv2}), we obtain
					\begin{subequations}
						\begin{align}
							\threenorms{\rgeo^{-1}\int_{u}^{t}\rgeo^2\sgabs{(\angnabla\modtorsion,\angnabla\chismall)\cdot\ACC}\diff\tau}{u}{2}{t}{2}{\sangle}{p}{\intregion}\lesssim&\norm{\twonorms{\rgeo(\angnabla\modtorsion,\angnabla\chismall)}{t}{2}{\sangle}{p}{\coneu}\twonorms{\ACC}{t}{2}{\sangle}{\infty}{\coneu}}_{L^2_uL^2_t}\\
							\notag\lesssim&\lambda^{-8\varepsilon_0}+\lambda^{-4\varepsilon_0}\threenorms{\rgeo\angnabla\modtorsion}{u}{2}{t}{2}{\sangle}{p}{\intregion},\\
							\threenorms{\rgeo^{-1}\int_{u}^{t}\rgeo^2\sgabs{\pfour^2\gfour\cdot\ACC}\diff\tau}{u}{2}{t}{2}{\sangle}{p}{\intregion}\lesssim&\lambda^{1/2-4\varepsilon_0}\norm{\threenorms{\rgeo\pfour^2\gfour}{t}{\infty}{u}{2}{\sangle}{p}{\coneu}\threenorms{\ACC}{t}{2}{u}{\infty}{\sangle}{\infty}{\coneu}}_{L^2_t}\\
							\notag\lesssim&\lambda^{-8\varepsilon_0}
							.\end{align}
					\end{subequations}
					For the fourth and the fifth term on the right of (\ref{DE:J2}), by the bootstrap assumptions (\ref{bt7}), and the estimates (\ref{fv6}), (\ref{mr7}), (\ref{mr21}) and (\ref{rmk1}), we derive:
					\begin{align}\label{ES2635}
						&\threenorms{\rgeo^{-1}\int_{u}^{t}\sgabs{\rgeo^2\angnabla\conformalfactor\cdot\left(\angnabla\pfour\gfour,\angnabla\chismall,\left(\ACC,\rgeo^{-1}\right)\cdot\pfour\gfour\right)} \diff\tau}{u}{2}{t}{2}{\sangle}{p}{\intregion}\\
						\notag&\lesssim
						\norm{\threenorms{\angnabla\conformalfactor}{t}{2}{u}{2}{\sangle}{\infty}{\coneu}\threenorms{\rgeo\left(\angnabla\pfour\gfour,\angnabla\chismall,\left(\ACC,\rgeo^{-1}\right)\cdot\pfour\gfour\right)}{t}{2}{u}{\infty}{\sangle}{p}{\coneu}}_{L^2_t}\\
						&\notag\lesssim\lambda^{-4\varepsilon_0}
						.\end{align}
					For the last term on the right of (\ref{DE:J2}), by the bootstrap assumptions (\ref{BA:bt1}), estimate (\ref{BA:bt2}), already proven results (\ref{mr14}), (\ref{mr145}), and the estimate (\ref{rmk1}), we have:
					\begin{align}\label{ES264}
						&\threenorms{\rgeo^{-1}\int_{u}^{t}\rgeo^2\sgabs{\left(\ACC,\rgeo^{-1}\right)\cdot\ACC\cdot\ACC}\diff\tau}{u}{2}{t}{2}{\sangle}{p}{\intregion}\\
						\notag&\lesssim\norm{\twonorms{\rgeo\left(\ACC,\rgeo^{-1}\right)\cdot\ACC}{t}{2}{\sangle}{p}{\coneu}\twonorms{\ACC}{t}{2}{\sangle}{\infty}{\coneu}}_{L^2_uL^2_t}\\
						\notag&\lesssim\lambda^{-8\varepsilon_0}
						.\end{align}
					Note that $\mathfrak{J}_{(1)}$ is bounded by (\ref{ES261})-(\ref{ES262}), and $\mathfrak{J}_{(2)}$ is bounded by (\ref{ES263})-(\ref{ES264}). Therefore, we conclude 
					\begin{align}\label{ES265}
						\threenorms{\rgeo\modmass}{u}{2}{t}{2}{\sangle}{p}{\intregion}\lesssim\lambda^{-4\varepsilon_0}+\lambda^{-4\varepsilon_0}\threenorms{\rgeo\angnabla\modtorsion}{u}{2}{t}{2}{\sangle}{p}{\intregion}
						.\end{align}
					Combining estimate (\ref{ES265}) with (\ref{ES260}), and soaking the second term on the right of (\ref{ES260}) to the left of (\ref{ES265}), we conclude the desired results (\ref{mr24}).
				\end{proof}
				\begin{proof}[Proof of $\threenorms{\rgeo^{3/2}\modmass}{u}{2}{t}{\infty}{\sangle}{p}{\intregion}\lesssim\lambda^{-4\varepsilon_0}$ in (\ref{mr24})]
					We divide (\ref{r2modmass}) by $\rgeo^{1/2}(t,u)$, and then take the $L^2_uL_t^\infty L^p_\sangle$ norms. The estimates of the terms in $\mathfrak{J}_{(1)}$ and $\mathfrak{J}_{(2)}$, defined in (\ref{DE:J1J2}), follow a similar pattern to the previous proof of $\threenorms{\rgeo\angnabla\modmass}{u}{2}{t}{2}{\sangle}{p}{\intregion}\lesssim\lambda^{-4\varepsilon_0}$ in (\ref{mr24}) except for the term in $\mathfrak{J}_{(1)}$, for which we leverage on estimates (\ref{fv6}), (\ref{fv1}) to deduce:
					\begin{subequations}
						\begin{align}
							\label{ES271}\threenorms{\rgeo^{-1/2}\int_u^t\rgeo\angnabla\pfour\gfour \diff\tau}{u}{2}{t}{\infty}{\sangle}{p}{\intregion}&\lesssim\threenorms{\rgeo\angnabla\pfour\gfour}{u}{2}{t}{2}{\sangle}{p}{\intregion}\lesssim\lambda^{-4\varepsilon_0},\\
							\label{ES272}\threenorms{\rgeo^{-1/2}\int_{u}^{t}\pfour\gfour \diff\tau}{u}{2}{t}{\infty}{\sangle}{p}{\intregion}&\lesssim\threenorms{\pfour\gfour}{u}{2}{t}{2}{\sangle}{p}{\intregion}\lesssim\lambda^{-4\varepsilon_0}
							.\end{align}
					\end{subequations}
				\end{proof}
				\begin{proof}[Proof of $\holderthreenorms{\angnabla\conformalfactor}{u}{2}{t}{2}{\sangle}{0}{\delta_0}{\intregion}\lesssim\lambda^{-4\varepsilon_0}$ in (\ref{mr24})]
					Plugging (\ref{DE:defmodtorsion}) into the Sobolev equality (\ref{ES:sobolevfour}) with $Q=p$, we obtain
					\begin{align}
						\holderthreenorms{\angnabla\conformalfactor}{u}{2}{t}{2}{\sangle}{0}{\delta_0}{\intregion}\lesssim\threenorms{\rgeo(\angnabla\modtorsion,\angnabla\upzeta)}{u}{2}{t}{2}{\sangle}{p}{\intregion}+\threenorms{\angnabla\conformalfactor}{u}{2}{t}{2}{\sangle}{2}{\intregion}
						.\end{align}
					Together with the proven estimates (\ref{mr7}), (\ref{mr21}), (\ref{mr24}) for $\threenorms{\rgeo\angnabla\modtorsion}{u}{2}{t}{2}{\sangle}{p}{\intregion}$, we conclude the desired estimate.
				\end{proof}
				\begin{proof}[Proof of (\ref{mr26})]
					Plugging (\ref{DE:hodgemass}) into the Hodge estimate (\ref{ES:hodgeone}), using (\ref{ES:averagef}) and the proven estimate (\ref{mr24}), we get
					\begin{align}
						\threenorms{\rgeo\angnabla\hodgemass,\hodgemass}{t}{2}{u}{2}{\sangle}{p}{\intregion}\lesssim\threenorms{\rgeo\modmass}{t}{2}{u}{2}{\sangle}{p}{\intregion}+\threenorms{\rgeo\umodmass}{t}{2}{u}{2}{\sangle}{p}{\intregion}\lesssim\lambda^{-4\varepsilon_0}
						.\end{align}
					By the Sobolev inequality (\ref{ES:sobolevfour}) with $Q=p$, we deduce
					\begin{align}
						\holderthreenorms{\hodgemass}{t}{2}{u}{2}{\sangle}{0}{\delta_0}{\intregion}\lesssim\threenorms{\rgeo\angnabla\hodgemass}{t}{2}{u}{2}{\sangle}{p}{\intregion}+\threenorms{\hodgemass}{t}{2}{u}{2}{\sangle}{2}{\intregion}\lesssim\lambda^{-4\varepsilon_0}
						.\end{align}
				\end{proof}
				\begin{proof}[Proof of $\threenorms{\modtorsion-\hodgemass}{t}{2}{u}{\infty}{\sangle}{\infty}{\intregion}\lesssim\lambda^{-1/2-4\varepsilon_0}$ in (\ref{mr27})]
					Plugging equations (\ref{EQ:angdivhodgemass}) and (\ref{EQ:angcurlhodgemass}) into the Hodge estimate (\ref{ES:hodgesix}) with $Q=p$, $c=2$ and $\delta^\prime\leq\delta_0$, we have:
					\begin{align}
						\threenorms{\modtorsion-\hodgemass}{t}{2}{u}{\infty}{\sangle}{\infty}{\intregion}\lesssim&\norm{\upnu^{\delta^\prime}\littlewood\pfour\gfour}_{L^2_tL^\infty_ul^2_\upnu L_\sangle^\infty\left(\intregion\right)}+\threenorms{\pfour\gfour}{t}{2}{u}{\infty}{\sangle}{\infty}{\intregion}\\
						\notag&+\threenorms{\rgeo\boxg\gfour}{t}{2}{u}{\infty}{\sangle}{p}{\intregion}+\threenorms{\rgeo(\upzeta,\hchi)\cdot(\upzeta,\hchi))}{t}{2}{u}{\infty}{\sangle}{p}{\intregion}\\
						\notag&+\threenorms{\rgeo(\ACC,\rgeo^{-1})\cdot\pfour\gfour}{t}{2}{u}{\infty}{\sangle}{p}{\intregion}
						.\end{align}
					By the bootstrap assumptions (\ref{BA:bt1}), estimates (\ref{BA:bt2})(\ref{fv10})(\ref{rzt2sanglep}), we conclude the desired bound:
					\begin{align}
						\threenorms{\modtorsion-\hodgemass}{t}{2}{u}{\infty}{\sangle}{\infty}{\intregion}&\lesssim\lambda^{-1/2-4\varepsilon_0}
						.\end{align}
				\end{proof}
				\begin{proof}[Proof of (\ref{mr265})]
					By (\ref{EQ:initialhodgemass}), it suffices to show $\lim\limits_{t\downarrow u}\rgeo\hodgemass=\zero(\rgeo)$. Plugging the definition (\ref{DE:hodgemass}) into the Hodge estimate (\ref{ES:hodgesix}) with $\mathfrak{F}=0$ and $Q=p$, and together using the initial condition (\ref{ES:initialcontipone}), we obtain
					\begin{align}
						\onenorm{\lim\limits_{t\downarrow u}\rgeo\hodgemass}{\sangle}{\infty}{\stu}\lesssim\onenorm{\lim\limits_{t\downarrow u}\rgeo^2\modmass}{\sangle}{p}{\stu}=\zero(\rgeo)
						.\end{align}
				\end{proof}
				\begin{proof}[Proof of $\threenorms{\massone}{t}{2}{u}{\infty}{\sangle}{\infty}{\intregion}\lesssim\lambda^{-1/2-4\varepsilon_0}$ in (\ref{mr27})]
					Let 
					\begin{align}\label{ES280}
						\mathfrak{H}:=\angD_{\Lunit}\massone+\frac{1}{2}\gtr\upchi\massone.
					\end{align} Plugging (\ref{ES280}) into (\ref{EQ:transportlemma2}), then applying  (\ref{EQ:rtwomupxi}) with $m=\frac{1}{2}$ and initial condition (\ref{mr265}) for $\rgeo\massone$, we get
					\begin{align}
						\sgabs{\rgeo\massone}=\int_{u}^{t}\sgabs{\rgeo\mathfrak{H}+\rgeo(\chismall,\pfour\gfour)\cdot\massone} \diff\tau
						.\end{align}
					Employing the Gr\"onwall's inequality, the bootstrap assumption (\ref{BA:bt1}), estimate (\ref{BA:bt2}) and the proven estimate (\ref{mr14}), we derive
					\begin{align}\label{ES281}
						\sgabs{\rgeo\massone}\lesssim\int_{u}^{t}\sgabs{\rgeo\mathfrak{H}}\diff\tau\exp\left(\twonorms{\chismall,\pfour\gfour}{t}{1}{\sangle}{\infty}{\coneu}\right)\lesssim\int_{u}^{t}\sgabs{\rgeo\mathfrak{H}}\diff\tau
						.\end{align}
					Dividing both sides of (\ref{ES281}) by $\rgeo$, and then by substituting equations (\ref{EQ:angdiv1}) and (\ref{EQ:angcurl1}) into the Hodge estimate (\ref{ES:hodgesix}) with $\mathfrak{F}=\rgeo^{-1}\upxi$, $Q=p$, $c=2$ and $\delta^\prime\leq\delta_0$, we deduce
					\begin{align}
						\threenorms{\massone}{t}{2}{u}{\infty}{\sangle}{\infty}{\intregion}\lesssim&\left\|\rgeo^{-1}\int_{u}^{t}\rgeo\left(\rgeo^{-1}\norm{\upnu^{\delta^\prime}\littlewood\pfour\gfour}_{l^2_\upnu L^\infty_\sangle(\sTauu)}+\onenorm{\rgeo^{-1}\pfour\gfour}{\sangle}{\infty}{\sTauu}\right.\right.\\
						\notag&\left.\left.+\onenorm{\rgeo^{-1}\pfour\gfour}{\sangle}{p}{\sTauu}\right)\diff\tau\right\|_{L^2_tL^{\infty}_u}
						.\end{align}
					Appealing to the Hardy-Littlewood maximal inequality (\ref{ES:hardylittlewood}), the bootstrap assumptions (\ref{BA:bt1}), estimates (\ref{BA:bt2}), (\ref{fv1}), we conclude the desired estimate.
				\end{proof}
				\begin{proof}[Proof of $\threenorms{\masstwo}{u}{2}{t}{\infty}{\sangle}{\infty}{\intregion}\lesssim\lambda^{-1/2-4\varepsilon_0}$ in (\ref{mr28})]
					Following the same argument as in the proof of (\ref{mr27}), it holds that
					\begin{align}\label{ES291}
						\sgabs{\rgeo\masstwo}\lesssim\int_{u}^{t}\sgabs{\rgeo\mathfrak{H}}\diff\tau
						,\end{align}
					where $\mathfrak{H}:=\angD_{\Lunit}\masstwo+\frac{1}{2}\gtr\upchi\masstwo$. Now divide both sides of (\ref{ES291}) by $\rgeo$, and take the $L^2_uL^\infty_tL^\infty_\sangle$-norm. Notice that: 
					\begin{align}\label{ES292}
						\threenorms{\masstwo}{u}{2}{t}{\infty}{\sangle}{\infty}{\intregion}\lesssim\threenorms{\rgeo^{-1}\int_{u}^{t}\sgabs{\rgeo\mathfrak{H}}\diff\tau}{u}{2}{t}{\infty}{\sangle}{\infty}{\intregion}\lesssim\threenorms{\mathfrak{H}}{u}{2}{t}{1}{\sangle}{\infty}{\intregion}
						.\end{align}
					Applying the Hodge estimate (\ref{ES:hodgesix})  to equations (\ref{EQ:angdiv2}) and (\ref{EQ:angcurl2}) with $\mathfrak{F}=0$ and $Q=p$, we have:
					\begin{align}
						\threenorms{\mathfrak{H}}{u}{2}{t}{1}{\sangle}{\infty}{\intregion}\lesssim&\threenorms{\rgeo\mathfrak{J}_{(2)}}{u}{2}{t}{1}{\sangle}{p}{\intregion}+\threenorms{\rgeo\hchi\cdot\angnabla\hodgemass}{u}{2}{t}{1}{\sangle}{p}{\intregion}\\
						\notag&+\threenorms{\rgeo\left(\angnabla\pfour\gfour,\angnabla\chismall\right)\cdot\hodgemass}{u}{2}{t}{1}{\sangle}{p}{\intregion}\\
						\notag&+\threenorms{\rgeo\left(\ACC,\rgeo^{-1}\right)\cdot\pfour\gfour\cdot\hodgemass}{u}{2}{t}{1}{\sangle}{p}{\intregion}\\
						\notag&+\threenorms{\rgeo\left(\gtr\upchi-\average{\gtr\upchi}\right)\modmass}{u}{2}{t}{1}{\sangle}{p}{\intregion}
						.\end{align}
					For the terms in $\mathfrak{J_{(2)}}$, the proof is similar to that for $\threenorms{\rgeo\angnabla\modmass}{u}{2}{t}{2}{\sangle}{p}{\intregion}$ (specifically, the arguments in (\ref{ES263})-(\ref{ES264})). Adopting this approach, we derive
					\begin{align} \label{12.216}
						\threenorms{\rgeo\mathfrak{J}_{(2)}}{u}{2}{t}{1}{\sangle}{p}{\intregion}\lesssim\lambda^{-1/2-4\varepsilon_0}.
					\end{align}
					In particular, in the estimate of \eqref{12.216}, we note that, by using $\holderthreenorms{\angnabla\conformalfactor}{u}{2}{t}{2}{\sangle}{0}{\delta_0}{\intregion}\lesssim\lambda^{-4\varepsilon_0}$ in (\ref{mr24}), the analogous estimate of the corresponding term of $\mathfrak{J}_{(2)}$ in (\ref{ES2635}) is improved as follows:
					\begin{align}
						&\threenorms{\rgeo\angnabla\conformalfactor\cdot\left(\angnabla\pfour\gfour,\angnabla\chismall,\left(\ACC,\rgeo^{-1}\right)\cdot\pfour\gfour\right)} {u}{2}{t}{1}{\sangle}{p}{\intregion}\\
						\notag&\lesssim
						\threenorms{\angnabla\conformalfactor}{t}{2}{u}{2}{\sangle}{\infty}{\intregion}\threenorms{\rgeo\left(\angnabla\pfour\gfour,\angnabla\chismall,\left(\ACC,\rgeo^{-1}\right)\cdot\pfour\gfour\right)}{t}{2}{u}{\infty}{\sangle}{p}{\intregion}\\
						&\notag\lesssim\lambda^{-1/2-4\varepsilon_0}
						.\end{align}
					Using (\ref{mr26}) and (\ref{mr14}), we deduce
					\begin{align}
						\threenorms{\rgeo\hchi\cdot\angnabla\hodgemass}{u}{2}{t}{1}{\sangle}{p}{\intregion}\lesssim\threenorms{\rgeo\angnabla\hodgemass}{u}{2}{t}{2}{\sangle}{p}{\intregion}\threenorms{\hchi}{u}{\infty}{t}{2}{\sangle}{\infty}{\intregion}\lesssim\lambda^{-1/2-8\varepsilon_0}
						.\end{align}
					Combining (\ref{fv6}), (\ref{mr7}), (\ref{rmk1}) and (\ref{mr26}), we have
					\begin{align}
						&\threenorms{\rgeo\hodgemass\cdot\left(\angnabla\pfour\gfour,\angnabla\chismall,\left(\ACC,\rgeo^{-1}\right)\cdot\pfour\gfour\right)}{u}{2}{t}{1}{\sangle}{p}{\intregion}\\
						\notag\lesssim&\threenorms{\hodgemass}{u}{2}{t}{2}{\sangle}{\infty}{\intregion}\cdot\threenorms{\rgeo\left(\angnabla\pfour\gfour,\angnabla\chismall,\left(\ACC,\rgeo^{-1}\right)\cdot\pfour\gfour\right)}{u}{\infty}{t}{2}{\sangle}{p}{\intregion}\\
						\notag\lesssim&\lambda^{-1/2-4\varepsilon_0}
						.\end{align}
					Note that $\average{\frac{1}{\rgeo}}=\frac{1}{\rgeo}$ (see Definition \ref{DE:average} for the definition of $\average{\frac{1}{\rgeo}}$). Employing (\ref{BA:bt1}), (\ref{BA:bt2}), (\ref{mr14}), and (\ref{mr24}), we obtain
					\begin{align}\label{ES293}
						\threenorms{\rgeo\left(\gtr\upchi-\average{\gtr\upchi}\right)\modmass}{u}{2}{t}{1}{\sangle}{p}{\intregion}&\lesssim\threenorms{\pfour\gfour,\chismall}{u}{\infty}{t}{2}{\sangle}{\infty}{\intregion}\threenorms{\rgeo\modmass}{u}{2}{t}{2}{\sangle}{p}{\intregion}\\
						\notag&\lesssim\lambda^{-1/2-8\varepsilon_0}
						.\end{align}
					Combining (\ref{ES292})-(\ref{ES293}), we conclude the desired estimate (\ref{mr28}).
				\end{proof}

				\appendix
				\section*{Appendix}
				\addcontentsline{toc}{section}{Appendix}
				\renewcommand{\thesubsection}{\Alph{subsection}}
				\subsection{Notations}\label{AP:notations}
				In this appendix, we compile the notations that are frequently used throughout the article.
				\begin{center}

					\begin{tabular}{|c|c|} 
						\hline
						Symbol & Reference  \\ 
						\hline\hline
						$\lesssim,\p,\pfour,\antisymmetic_{\alpha\beta\gamma\delta},\antisymmetic_{ijk},\dive,\curl,\St$ & Section \ref{SS:Notations} \\
						\hline
						$\boxg,\resboxg$ & Definition \ref{DE:waveoperators} \\
						\hline
						$\vec{v},\vec{Y},F,W,\vec{U},\speedone,\speedtwo$ & Section \ref{SS:elasticwaveeq} \\
						\hline
						$G$ & Section \ref{SS:ModifiedElasticWaves} \\
						\hline
						$\vvardiv,\vvarcurl,\gfour_{E},\hfour_{E}$ & Section \ref{S:GeometricFormulation} \\
						\hline
						$\gfour,\gfour^{-1},\hfour,\hfour^{-1}$ & Section \ref{S:PDEs} \\
						\hline
						$\Chfour$ & Definition \ref{DE:Christoffelsymbols} \\
						\hline
						$\vvariables,\cvvariables$ & Definition \ref{DE:vvariablecvvariables} \\
						\hline
						$\littlewood$ & Section \ref{SS:definitionlittlewood} \\
						\hline
						$N,q,\varepsilon_0,\delta_0,\delta,\delta_1$ & Section \ref{SS:ChoiceofParameters} \\
						\hline
						$\Domain,D,\Int$& Section \ref{SS:Data}\\
						\hline
						$\Tstar$& Section \ref{SS:bootstrap}\\
						\hline
						$Q_{\mu\nu},\Jenarg{\mathbf{X}}{\alpha} ,\deform{\mathbf{X}},\Dfour,\Timelike,\mathbb{E}[\varphi]$ & Definition \ref{DE:energymomentum} \\	
						\hline
						$P_N$ & Section \ref{SSS:proofofEnergy}\\
						\hline
						$\lambda,I_k,\varphi$ & Section \ref{SS:partitioning}\\
						\hline
					\end{tabular}
					\begin{tabular}{|c|c|} 
						\hline
						Symbol & Reference  \\ 
						\hline\hline
						
						$\Trescale,(\dive\vvardiv)_{(\lambda)}$ & Section \ref{SS:SelfRescaling}\\
						\hline
						$B_R(p),B_\rho(p,\gt)$ & Theorem \ref{TH:Spatiallylocalizeddecay}\\
						\hline
						$u,w_*,\tip,\lapse,w,\coneu,\intregion,\extregion,\stu$ & Section \ref{SS:opticalfunction}\\
						\hline
						$\sangle^A,\esphere$ & Definition \ref{DE:norms}\\
						\hline
						$\rgeo,\Lgeo,\nulllapse,\spherenormal,\Lunit,\uLunit,\gsphere,e_A,\angnabla$ & Definition \ref{DE:Nullframe}\\
						\hline
						$\Conenormal$& Definition \ref{DE:Conenormal}\\
						\hline
						$\sphereproject,\gtr$& Definition \ref{D:stutensor fields}\\
						\hline
						$\Energyconformal[\varphi]$& Section \ref{SS:SectionSetupconformal}\\
						\hline
						$\breve{Q}_{\mu\nu},\mathcal{F}_{(1)}[\varphi;\coneu],\mathcal{F}_{(2)}[\varphi;\coneu]$& Section \ref{energyalongnullhypersurfaces}\\
						\hline
						$\angD,\Riemfour{\gamma}{\alpha}{\delta}{\beta},\Ricfour{\alpha}{\beta},[\cdot,\cdot]$& Section \ref{SSS:levi-civita2}\\
						\hline		
						$k,\chi,\spheresecondfund,\uchi,\upzeta,\uzeta,\Lie$& Definition \ref{DE:DEFSOFCONNECTIONCOEFFICIENTS}\\
						\hline	
						$\conformalfactor,\rescaledgfour,\congsphere,\reschi,\resuchi,\Restrace,\chismall,\CC,\ACC$& Definition \ref{DE:conformalchangemetric}\\
						\hline	
						$\average{f},\diff\gvol,\diff\tvol,\diff\spherevol,\volume,\conformalvol$& Definition \ref{DE:average}\\
						\hline	
						$\mass,\modmass,\hodgemass,\modtorsion$& Definition \ref{DE:mass}\\
						\hline		
						$\massone,\masstwo$& \eqref{EQ:hodgemudecomposition}\\
						\hline	
					\end{tabular}
				\end{center}
				\bibliography{Low-Regularity_Local_Well-Posedness_for_the_Elastic_Wave_System.bbl}
				\bibliographystyle{plain}
			\end{document}